\documentclass[a4paper,english,11pt,twoside,dvips]{article}
\usepackage{citehack}
\usepackage{amsmath}
\usepackage{amsfonts}
\usepackage{amssymb}
\usepackage{longtable}
\usepackage{amsthm}
\usepackage{euscript}
\usepackage{array}

\usepackage{graphicx}

\newcommand{\REM}[1]{\relax}
\makeatletter%

\newcommand{\begthrm}[2]{\begin{trivlist}\it\item[\hspace{\labelsep}{\bf #1~#2.}]}
\newcommand{\enthrm}{\end{trivlist}}

\renewcommand{\thepage}{\roman{page}}

\newcommand{\OrdinaryPageStyle}{
\setcounter{page}{1}%
\renewcommand{\thepage}{\arabic{page}}
\renewcommand{\@evenhead}{\underline{\hbox to\textwidth{\it\thepage\hfil{\mychap}\hfil}}}%
\renewcommand{\@oddfoot}{\relax}%
\renewcommand{\@oddhead}{\underline{\hbox to\textwidth{\it\hfil{\mypar}\hfil\thepage}}}%
\renewcommand{\@evenfoot}{\relax}%
}

\def\mypar{\relax}

\renewcommand{\@makefntext}[1]{{\normalsize\parindent=1em\noindent\hbox to 1.8em{\hss
$^{\@thefnmark}$}#1}}

\renewcommand{\part}[1]{\ifodd\arabic{page}\newpage\hbox{~~~~}\thispagestyle{empty}\fi%
\setcounter{section}{0}\newpage%
\vspace*{-5ex}\begin{center}\@startsection{part}{1}%
{\parindent}{3.5ex plus 1ex minus .2ex}%
{2.3ex plus .2ex}{\large\bf Chapter~\Roman{part}.~#1}{}\end{center}\nopagebreak%
\vspace{20pt}\addcontentsline{toc}{partm}{Chapter~\Roman{part}.
#1}\gdef\mychap{Chapter~\Roman{part}. #1}\thispagestyle{plain}\gdef\mypar{\mychap}}

\newcommand{\parti}[1]{\ifodd\arabic{page}\newpage\thispagestyle{empty}\hbox{~~~~}\fi%
\setcounter{section}{0}\newpage%
\newpage\vspace*{-5ex}\begin{center}\@startsection{part}{1}%
{\parindent}{3.5ex plus 1ex minus .2ex}%
{2.3ex plus .2ex}{\large\bf #1}{}\end{center}\nopagebreak%
\vspace{20pt}\addcontentsline{toc}{partm}{#1}\addtocounter{part}{-1}\gdef\mychap{#1}
\gdef\mypar{\mychap}\thispagestyle{empty}}

\newcommand{\partii}[1]{\ifodd\arabic{page}\newpage\thispagestyle{empty}\hbox{~~~~}\fi%
\setcounter{section}{0}\newpage%
\newpage\vspace*{-5ex}\begin{center}\@startsection{part}{1}%
{\parindent}{3.5ex plus 1ex minus .2ex}%
{2.3ex plus .2ex}{\large\bf #1}{}\end{center}\nopagebreak%
\vspace{20pt}\addtocounter{part}{-1}\gdef\mychap{#1}\thispagestyle{plain}\gdef\mypar{\mychap}}

\renewcommand{\thesection}{\arabic{part}.\arabic{section}}

\renewcommand{\section}{\vskip15mm\@startsection{section}{2}%
{\parindent}{3.5ex plus 1ex minus .2ex}%
{2.3ex plus .2ex}{\large\bf}}

\renewcommand{\l@section}{\@dottedtocline{2}{1.5em}{2.3em}}

\newcommand{\l@partm}[2]{\@dottedtocline{1}{3em}{0em}{%
#1}{#2}}
\renewcommand{\l@part}[2]{\relax}
\renewcommand{\tableofcontents}{\partii{Contents}\@starttoc{toc}}
\def\thebibliography#1{\parti{References}\thispagestyle{plain}\list
 {\@biblabel{\arabic{enumiv}}}{\settowidth\labelwidth{\@biblabel{#1}}%
 \leftmargin\labelwidth\advance\leftmargin\labelsep
 \usecounter{enumiv}\let\p@enumiv\@empty
 \def\theenumiv{\arabic{enumiv}}}%
 \def\newblock{\hskip .11em plus.33em minus.07em}%
 \sloppy\clubpenalty4000\widowpenalty4000\sfcode`\.=\@m\relax}

\renewcommand{\p@equation}{\relax}

\renewcommand{\l@subsection}{\@dottedtocline{3}{1.5em}{3.3em}}

\expandafter\def\csname @tocrmarg\endcsname{4em}

\thm@notefont{\rm}

\newtheorem{theorem}{Theorem}[section]
\newtheorem{prop}{Proposition}[section]

\newtheorem{lem}{Lemma}[section]
\newtheorem{ex}{Example}[section]
\newtheorem{tab}{Table}[section]
\newtheorem{corol}{Corollary}[section]

\def\zr{\ltimes}

\def\ZR{\rightthreetimes}

\def\spa{\mathop\text{{\rm span}}\nolimits}
\def\spar{\mathop\text{\rm span}_\mathbb{R}\nolimits}
\def\spac{\mathop\text{\rm span}_\mathbb{C}\nolimits}
\def\dimr{\mathop\text{\rm dim}_\mathbb{R}\nolimits}
\def\dimc{\mathop\text{\rm dim}_\mathbb{C}\nolimits}
\def\Hom{\mathop\text{\rm Hom}\nolimits}
\def\Ric{\mathop\text{\rm Ric}\nolimits}

\def\pr{\mathop\text{\rm pr}\nolimits}
\def\grad{\mathop\text{\rm grad}\nolimits}
\def\tr{\mathop\text{\rm tr}\nolimits}
\def\Aut{\mathop\text{\rm Aut}\nolimits}
\def\End{\mathop\text{\rm End}\nolimits}

\def\Real{\mathbb{R}}
\def\Co{\mathbb{C}}
\def\a{\mathfrak{a}}
\def\g{\mathfrak{g}}
\def\h{\mathfrak{h}}
\def\so{\mathfrak{so}}
\def\sod{\mathfrak{sod}}
\def\SOd{SOD}
\def\Ga{\Gamma}

\def\gl{\mathfrak{gl}}
\def\su{\mathfrak{su}}
\def\spin{\mathfrak{spin}}
\def\u{\mathfrak{u}}
\def\un{\mathfrak{u}}
\def\hol{\mathfrak{hol}}

\def\f{\mathfrak{f}}

\def\m{\mathfrak{m}}
\def\z{\mathfrak{z}}
\def\R{\mathcal{R}}
\def\P{\mathcal{P}}

\def\Conf{\mathop\text{\rm Conf}\nolimits}
\def\Simil{\mathop\text{\rm Sim}\nolimits}
\def\Isom{\mathop\text{\rm Isom}\nolimits}
\def\ad{\mathop\text{\rm ad}\nolimits}
\def\A{\mathcal {A}}
\def\B{\mathcal {B}}

\def\N{\mathcal {N}}
\def\K{\mathcal {K}}
\def\H{\mathcal {H}}
\def\F{\mathcal {F}}
\def\C{\mathcal {C}}

\def\id{\mathop\text{\rm id}\nolimits}

\def\im{\mathop\text{\rm Im}\nolimits}
\def\re{\mathop\text{\rm Re}\nolimits}
\def\lie{\mathop\text{\rm LA}\nolimits}
\def\pq{p_1\wedge q_1+p_2\wedge q_2}
\def\qp{p_1\wedge q_2-p_2\wedge q_1}

\def\i{{\hat{i}}}
\def\j{{\hat{j}}}
\def\k{{\hat{k}}}
\def\l{{\hat{l}}}
\def\ii{{\hat{\hat{i}}}}
\def\jj{{\hat{\hat{j}}}}
\def\kk{{\hat{\hat{k}}}}
\def\ll{{\hat{\hat{l}}}}
\def\iv{{\tilde{i}}}
\def\jv{{\tilde{j}}}
\def\kv{{\tilde{k}}}
\def\lv{{\tilde{l}}}
\def\iiv{{\tilde{\tilde{i}}}}
\def\jjv{{\tilde{\tilde{j}}}}
\def\kkv{{\tilde{\tilde{k}}}}
\def\llv{{\tilde{\tilde{l}}}}
\def\ib{{\breve{i}}}
\def\jb{{\breve{j}}}
\def\iib{{\breve{\breve{i}}}}
\def\jjb{{\breve{\breve{j}}}}
\def\p{\partial}

\newcommand{\blem}{\begin{lem}}
\newcommand{\elem}{\end{lem}}
\newcommand{\beq}{\begin{equation}}
\newcommand{\eeq}{\end{equation}}

\def\arcth{\mathrm{artanh\,}}
\def\sh{\mathrm{sinh\,}}
\def\ch{\mathrm{cosh\,}}

\def\th{\mathrm{tanh\,}}

\def\bR{{\mathbb R}}
\def\na{\nabla}
\def\n{\nabla}
\def\lam{\lambda}
\def\t{\tilde}
\def\pp{{\hat p_2}}

\def\yy{\hat y}
\def\pyy{\p_{\hat y}}
\def\py{\p_{y}}
\def\J{\t J}
\def\s{\sigma}
\newcommand\pK{ pseudo-K\"ahlerian }
\newcommand\pR{ pseudo-Riemannian }

\title{Holonomy groups and special geometric structures  of pseudo-K\"ahlerian manifolds of index 2}

\author{Anton S. Galaev \\ Humboldt-Universit\"at zu Berlin,
 Institut f\"ur Mathematik,\\ Rudower Chaussee 25,
 12489 Berlin}

\begin{document}


\thispagestyle{empty}
\begin{center}

$\quad$
\bigskip
{\Large\bfseries
Holonomy groups and  special  geometric  structures    of pseudo-K\"ahlerian  manifolds of index 2
}

\bigskip
\bigskip

Dissertation

zur Erlangung des akademischen Grades

doctor rerum naturalium

im Fach Mathematik

(endg\"ultige Fassung)

\end{center}

\bigskip
\bigskip
\bigskip
\bigskip
\bigskip
\noindent
eingereicht an der

\noindent
Mathematisch-Naturwissenschaftlichen Fakult\"at II

\noindent
der Humboldt-Universit\"at zu Berlin

\bigskip
\noindent
von Dipl. Math. Anton Galaev

\noindent
geb. am 02.02.1981 in Saratov

\bigskip
\bigskip
\bigskip
\bigskip
\noindent
Der Pr\"asident der Humboldt-Universit\"at zu Berlin

\noindent
Prof. Dr. Christoph Markschies

\bigskip
\noindent
Der Dekan der Mathematisch-Naturwissenschaftlichen Fakult\"at II

\noindent
Prof. Dr. Wolfgang Coy

\bigskip
\noindent Erste Gutachterin: Frau Prof. Dr. Helga Baum (Humboldt-Universit\"at zu Berlin)

\noindent Zweiter Gutachter: Herr Prof. Dr. Dmitri Alekseevsky (University of Edinburg)

\noindent Dritter Gutachter: Herr Prof. Dr. Vicente Cort\'es (Universit\"at Hamburg)

\bigskip
\noindent Tag der Verteidigung: 8. Dezember 2006

\newpage

{\small

 \centerline{\bf Abstract}

\bigskip

The holonomy group of a pseudo-Riemannian manifold, t.m. the group
of  parallel displacements along all loops at a fixed point,
is an invariant of the Levi-Civita connection and gives much
information about the special geometric structure of the manifold.
The problem to classify all possible holonomy groups is still
open. In the case of simply connected manifolds this problem
reduces to the classification of possible holonomy algebras. The
classification of {\em Riemannian} holonomy algebras is a
classical result. The classification of {\em Lorentzian}
holonomy algebras was achieved only recently.

In this thesis we classify holonomy algebras of \pK manifolds of
index 2 and consider some examples and applications.

Chapter I contains an introduction to the theory of holonomy
groups. After fixing the necessary notation we give an outline of
the methods and results which leads to the classification of the
holonomy algebras for Riemannian and Lorentzian manifolds.
Furthermore, for pseudo-Riemannian manifolds, we reduce the
classification problem to the classification of
weakly-irreducible not irreducible holonomy algebras.

In Chapter II we classify weakly-irreducible not irreducible
subalgebras of $\su(1,n+1)$ ($n\geq 0$).

Chapter III deals with the classification of weakly-irreducible
not irreducible holonomy algebras of pseudo-K\"ahlerian and
special pseudo-K\"ahlerian manifolds. First we classify
weakly-irreducible not irreducible {\em Berger} subalgebras of
$\u(1,n+1)$. These algebras are generated by the images of their
algebraic curvature tensors and are candidates for the holonomy
algebras. We describe the spaces of curvature tensors for all of
these algebras. After that, we construct a \pK metric of index 2
for each of these Berger algebras, i.e. we show that all Berger
subalgebras of $\u(1,n+1)$ are in fact holonomy algebras.

In Chapter IV we consider some examples and applications. In the
first part we describe examples of 4-dimensional Lie groups with
left-invariant \pK metrics and determine their holonomy algebras.
In the second part we use our classification of holonomy algebras
to give a new proof for the classification of simply connected \pK
symmetric spaces of index 2 with weakly-irreducible not
irreducible holonomy algebras. Finally we consider time-like cones
over Lorentzian Sasaki manifolds. These cones are also \pK
manifolds of index 2. We describe the local DeRham-Wu
decomposition of the cone in terms of the initial Lorentzian
Sasaki manifold and we describe all possible weakly-irreducible
not irreducible holonomy algebras of
such cones.

\newpage

\centerline{\bf Zusammenfassung}

\bigskip

Die Holonomiegruppe einer pseudo-Riemannschen Mannigfaltigkeit,
d.h. die Gruppe der Parallelverschiebungen entlang aller in einem
fixierten Punkt geschlossenen Kurven, ist eine Invariante des
Levi-Civita Zusammenhangs. Sie gibt viele Informationen \"uber die
spezielle geometrische Struktur der Mannigfaltigkeit. Das Problem,
alle m\"oglichen Holonomiegruppen zu klassifizieren ist nach wie
vor offen. Im Falle einfach-zusammenh\"angender Mannigfaltigkeiten
reduziert sich dieses Probelm auf die Klassifikation der
m\"oglichen Holonomiealgebren. Die Klassifikation der {\em
Riemannschen} Holonomiealgebren ist ein klassisches Resultat. Vor
kurzem
wurde die Klassifikation der {\em Lorentzschen} Holonomiealgebren abgeschlossen.

In dieser Dissertation klassifizieren wir die Holonomiealgebren
von pseudo-K\"ahlerschen Mannigfaltigkeiten vom Index 2 und
betrachten einige Beispiele und Anwendungen.

Kapitel I enth\"alt eine Einf\"uhrung in die Theorie der
Holonomiegruppen. Nach der Definition der n\"otigen Begriffe
beschreiben wir die Methode und die Ergebnisse der Klassifikation
der Holonomiealgebren von Riemannschen und Lorentzschen
Mannigfaltigkeiten. Insbesondere reduzieren wir das
Klassifikation-Problem auf den Fall von schwach-irreduziblen nicht
irreduziblen Holonomiealgebren.

Im Kapitel II klassifizieren wir alle schwach-irreduziblen nicht
irreduziblen Unteralgebren von $\su(1,n+1)$ ($n\geq 0$).

Das Kapitel III ist der Klassifikation der schwach-irreduziblen
nicht irreduziblen Holonomiealgebren von pseudo-K\"ahlerschen und
speziellen pseudo-K\"ahlerschen Mannigfaltigkeiten gewidmet. Dazu
klassifizieren wir zuerst schwach-irreduzible nicht irreduzible
{\em Berger}-Unteralgebren von  $\u(1,n+1)$. Solche Algebren
werden durch die Bilder ihrer algebraischen Kr\"ummungs\-tensoren
erzeugt und sind Kandidaten f\"ur die Holonomiealgebren. Danach
konstruieren wir f\"ur jede dieser Berger-Algebren eine
pseudo-K\"ahlersche Metrik vom Index 2, d.h. wir zeigen, dass alle
Berger-Unteralgebren von $\u(1,n+1)$ tats\"achlich
Holonomiealgebren von K\"ahler-Metriken vom Index 2 sind.

Im Kapitel IV betrachten wir einige Beispiele und Anwendungen. Im
1.Teil beschreiben wir Beispiele 4-dimensionaler Lie-Gruppen mit
links-invarianter pseudo-K\"ahlerscher Metrik vom Index 2 und
bestimmen ihre Holonomiealgebren. Im 2. Teil benutzen wir unsere
Klassifikation der Holonomiealgebren, um einen neuen Beweis f\"ur
die Klassifikation der einfach-zusammenh\"angenden symmetrischen
pseudo-K\"ahlersch Mannigfaltigkeiten vom Index 2 mit
schwach-irreduzibler nicht irreduzibler Holonomiealgebra
anzugeben. Zum Abschluss betrachten wir den zeitartigen Kegel
\"uber Lorentz-Sasaki Mannigfaltigkeiten. Diese Kegel sind
ebenfalls pseudo-K\"ahlersche Mannigfaltigkeiten vom Index 2. Wir
beschreiben die lokale DeRham-Wu-Zerlegung f\"ur den Kegel in
Abh\"angigkeit von den Eigenschaften der Kegelbasis und bestimmen
die m\"oglichen schwach-irreduziblen nicht-irredu\-zib\-len
Holonomiealgebren f\"ur solche Kegel.
}




\tableofcontents

\parti{Introduction}

\OrdinaryPageStyle

The Levi-Civita connection of a pseudo-Riemannian manifold $(M,g)$ defines the {\it parallel displacement}, i.e. for any
piecewise smooth curve $\gamma:[a,b]\subset\Real\to M$ we have an isomorphism $\tau_\gamma:T_{\gamma(a)}M\to
T_{\gamma(b)}M$ preserving the metric $g$. {\it The holonomy group} of a \pR manifold $(M,g)$ of signature $(r,s)$ at a
point $x\in M$ is the Lie subgroup of the pseudo-orthogonal Lie group $O(T_xM,g_x)\simeq O(r,s)$ that consists of parallel
displacements along all piecewise smooth loops at the point $x$. The corresponding subalgebra of $\so(T_xM,g_x)\simeq
\so(r,s)$ is called {\it the holonomy algebra} of the manifold $(M,g)$ at the point $x$. The holonomy group is an
invariant of the Levi-Civita connection of a pseudo-Riemannian manifold. In particular,  it allows to find all parallel
geometric objects on the manifold (e.g. tensor fields or distributions, see Theorems \ref{FP1}, \ref{FP2}). Knowing the
holonomy group of a pseudo-Riemannian manifold we can say whether the manifold is flat, orientable, pseudo-K\"ahlerian,
\pK and Ricci-flat or locally  decomposable (locally a product of \pR manifolds). The holonomy groups allow also to find
parallel spinors on \pR spin manifolds \cite{Wang, HelgaInes, Le4}. Using the cone constructions, one can  find Killing
spinors on \pR spin manifolds \cite{Baer,Helga,Kat99,Boh03}. Similarly, the holonomy algebra gives information about
locally defined parallel objects and it gives global information if the manifold is simply connected (then the holonomy
group is connected). The holonomy algebra gives also information about the curvature tensor of the manifold. Thus, \pR
manifolds with different holonomy groups have different kinds of geometries and to know all possible geometries we need a
classification of holonomy groups. Note that this problem depends  on the topology of the manifold and for simplicity we
restrict ourself only to connected holonomy groups (equivalently, holonomy algebras).

In {\bf Chapter I} we give  an introduction to the theory of holonomy groups.
We provide definitions,  general facts, examples and ideas of some proofs that illustrate the methods of the holonomy.
We recall the classifications of the  holonomy algebras for Riemannian and Lorentzian manifolds.

The classification of connected holonomy groups of Riemannian manifolds is a classical result.
First in 1952 A. Borel and A. Lichnerowicz proved that
a Riemannian manifold is locally a product of  Riemannian manifolds with irreducible holonomy groups,
see \cite{Bo-Li}.
In 1955 M. Berger gave a list of  possible connected irreducible
holonomy groups of Riemannian manifolds, see \cite{Ber}.
Later, in 1987 R. Bryant constructed metrics for the  exceptional
groups of this list, see \cite{Bryant}. See also Section \ref{sec1.2}.

 In the case of pseudo-Riemannian manifolds appears the situation that the holonomy group
preserves a degenerate vector subspace of the tangent space. In this situation the
Borel-Lichnerowicz theorem does not work. A subgroup $G\subset SO(p,q)$
is called {\it weakly-irreducible} if it does not  preserve any non-degenerate
proper subspace of $\Real^{r,s}$.
The Wu theorem states that {\it
a pseudo-Riemannian manifold is locally a product of pseudo-Riemannian manifolds
with weakly-irreducible holonomy groups},
see \cite{Wu}.
This reduces the problem of classification of the holonomy groups of pseudo-Riemannian manifolds
to the weakly-irreducible case.
If a holonomy group is irreducible, then it is weakly-irreducible. In \cite{Ber}
M. Berger gave also  a list  of possible connected  irreducible holonomy groups for
pseudo-Riemannian manifolds. This list was refined and completed by several people, see \cite{Br3,CMS,MS99}.

Thus the first problem is to classify weakly-irreducible not irreducible subgroups of $SO(p,q)$.
This was completely done only for connected groups in the Lorentzian case, i.e. for the signature $(1,n+1)$,
in 1993 by L. Berard Bergery and A. Ikemakhen, who divided connected weakly-irreducible not irreducible
subgroups of $SO(1,n+1)$ into 4 types, see \cite{B-I}. In \cite{Gal2} a more geometric proof of this result was given, see
Section \ref{sec1.3}. To each  weakly-irreducible not irreducible subgroup of  $G\subset SO(1,n+1)$
can be associated a subgroup of $SO(n)$, which is called the orthogonal part of $G$. Just recently T. Leistner showed that
the orthogonal part of a weakly-irreducible not irreducible holonomy group of a Lorentzian manifold must be the holonomy
group of a Riemannian manifold, see \cite{Le1,Le2,Le3}. In \cite{Gal5} metrics for all possible connected holonomy  groups
of Lorentzian manifold were constructed. This completes
the classification of connected holonomy groups for Lorentzian manifolds. See also Section \ref{sec1.3}.
In 1998 A. Ikemakhen classified connected weakly-irreducible
subgroups of $SO(2,N)$ that preserve an isotropic plane and satisfy an additional condition, see \cite{Ik22}.
There are  partial results for signature $(n,n)$, see \cite{B-Inn}.

In this thesis  we study connected holonomy groups of pseudo-K\"ahlerian manifolds of signature  $(2,2n+2)$, i.e.
holonomy groups contained in $U(1,n+1)\subset SO(2,2n+2)$. From the Wu theorem it follows that
any such group is a product of irreducible holonomy groups of K\"ahlerian manifolds and of the
weakly-irreducible  holonomy group of a  pseudo-K\"ahlerian manifold of signature  $(2,2k+2)$.

First of all,  we classify connected holonomy groups of pseudo-K\"ahlerian manifolds of index~2.

Let $\Real^{2,2n+2}$ be a $2n+4$-dimensional real  vector space endowed with a
complex structure $J$ and with a $J$-invariant metric $\eta$ of
signature $(2,2n+2)$ ($n\geq 0$).
In {\bf Chapter II}  we classify (up to conjugacy) all  connected subgroups of $SU(1,n+1)$ that
act weakly-irreducibly  and not irreducibly on $\Real^{2,2n+2}$, that is equivalent to the
classification of the corresponding subalgebras of $\su(1,n+1)$. Any such subgroup preserves a 2-dimensional isotropic $J$-invariant subspace of  $\Real^{2,2n+2}$.
We use a generalization of the method from \cite{Gal2} (see Section \ref{sec1.3}).

As the first case, we consider all subalgebras of $\su(1,1)$ that preserve a 2-dimensional isotropic $J$-invariant subspace of  $\Real^{2,2}$
and show which of these subalgebras are weakly-irreducible.

Then we consider the case $n\geq 1$. We denote by $\Co^{1,n+1}$ the $n+2$-dimensional complex vector space given by $(\Real^{2,2n+2},J,\eta)$.
Let $g$ be the  pseudo-Hermitian metric on $\Co^{1,n+1}$ of signature
$(1,n+1)$ corresponding to $\eta$.
If a subgroup $G\subset U(1,n+1)$ acts weakly-irreducibly on $\Real^{2,2n+2}$, then
$G$ acts  weakly-irreducibly on $\Co^{1,n+1}$, i.e. does not preserve any proper $g$-non-degenerate complex vector subspace.

We consider the boundary $\p\mathbf{H}^{n+1}_\mathbb{C}$ of
the complex hyperbolic space $\mathbf{H}^{n+1}_\mathbb{C}$ and
identify $\p\mathbf{H}^{n+1}_\mathbb{C}$ with the $2n+1$-dimensional
sphere $S^{2n+1}$. We fix  a complex isotropic line $l\subset \Co^{1,n+1}$ and
 denote by  $U(1,n+1)_l\subset U(1,n+1)$ the connected Lie  subgroup that
preserves the line $l$. Any connected subgroup $G\subset U(1,n+1)$ that acts on $\Co^{1,n+1}$
weakly-irreducibly and not irreducibly is conjugated to a subgroup of $U(1,n+1)_l$.

We identify the set $\p\mathbf{H}^{n+1}_\mathbb{C}\backslash\{l\}=S^{2n+1}\backslash\{point\}$ with
the Heisenberg space $\H_n=\Co^n\oplus\Real$. Any element $f\in U(1,n+1)_l$ induces
a transformation $\Gamma(f)$ of $\H_n$, moreover, $\Gamma(f)\in\Simil \H_n$,
where $\Simil \H_n$ is the group of the  Heisenberg similarity transformations of $\H_n$.
We show that $\Gamma:U(1,n+1)_l\to\Simil\H_n$ is a surjective Lie group homomorphism with
the kernel $\mathbb{T}$, where $\mathbb{T}$ is  the 1-dimensional subgroup generated by the complex structure
$J\in U(1,n+1)_l$. In particular,  $\mathbb{T}$ is the center of $U(1,n+1)_l$.
Let $SU(1,n+1)_{l}=U(1,n+1)_{l}\cap SU(1,n+1)$.
Then  $U(1,n+1)_{l}=SU(1,n+1)_l\cdot\mathbb{T}$
and the restriction
$$\Gamma|_{SU(1,n+1)_{l}}:SU(1,n+1)_{l}\to\Simil\H_n$$
is a Lie group isomorphism.

We consider the projection
$\pi:\Simil{\H_n}\to\Simil \Co^n$,
where $\Simil \Co^n$ is the group of similarity transformations of $\Co^n$.
The homomorphism  $\,\pi\,$ is surjective and its kernel is 1-dimensional.

We prove that {\it
if a subgroup $G\subset U(1,n+1)_{l}$ acts weakly-irreducibly on $\Co^{1,n+1}$,
then}
\begin{itemize}
\item[(1)] {\it the subgroup $\pi(\Gamma(G))\subset \Simil \Co^n$
does not preserve any proper
complex affine subspace of $\Co^n$;}
\item[(2)] {\it if $\pi(\Gamma(G))\subset \Simil \Co^n$ preserves a proper
non-complex affine subspace $L\subset\Co^n$, then the minimal
complex affine subspace of $\Co^n$ containing $L$ is $\Co^n$.}
\end{itemize}
This is the key statement for our classification.

Since we are interested in  connected Lie groups, it is enough to classify the corresponding Lie algebras. The
classification is done in the following way:\\ $\bullet$ First we describe non-complex vector subspaces $L\subset\Co^n$
with $\spa_{\Co}L=\Co^n$ (it is enough to consider only vector subspaces, since we do the classification  up to
conjugacy). Any such non-complex vector subspace has the form $L=\Co^m\oplus\Real^{n-m}$, where $0\leq m\leq n$. Here we
have 3 types of subspaces: 1) $m=0$ ($L$ is a real form of $\Co^n$); 2) $0< m<n$; 3) $m=n$ ($L=\Co^n$).\\ $\bullet$ We
describe the Lie algebras $\f$  of the connected  Lie subgroups $F\subset\Simil \Co^n$ preserving $L$. Without loss of
generality,  we can assume that each Lie group $F$ does not preserve any proper affine subspace of  $L$. This means that
$F$ acts irreducibly on $L$. By a theorem of D.V. Alekseevsky \cite{Al2,A-V-S}, $F$ acts transitively on $L$. In our
recent paper \cite{Gal2} we  divided transitive similarity transformation groups of Euclidean spaces into 4 types. Here we
unify two of the types. The group $F$ is contained in $(\Real^+\times SO(L)\times SO(L^{\bot_\eta}))\ZR L$, where
$\Real^+$ is the group of real dilations of $\Co^n$ about the origin and $L$ is the group of all translations in $\Co^n$
by vectors of $L$. In general situation we know only the projection of $F$ on $\Simil L=(\Real^+\times SO(L))\ZR L$, but
in our case the projection of $F$ on $SO(L)\times SO(L^{\bot_\eta})$ is also contained in $U(n)$ and we know the full
information about $F$. On this step we obtain $9$ types of Lie algebras.\\ $\bullet$  Then we describe subalgebras
$\a\subset\lie (\Simil\H_n)$ with $\pi(\a)=\f$. For each $\f$ we have 2 possibilities: $\a=\f+\ker\pi$ or
$\a=\{x+\zeta(x)|x\in\f\}$, where $\zeta:\f\to \ker\pi$ is a linear map. Using the  isomorphism
$(\Gamma|_{\su(1,n+1)_{l}})^{-1}$ we obtain a list of subalgebras $\g\subset \su(1,n+1)_{l}$. This gives us  12 types of
Lie algebras.\\ $\bullet$ Finally we check which of the obtained subalgebras of $\su(1,n+1)_l\subset\so(2,2n+2)$ are
weakly-irreducible. It turns out that some of the types contain Lie algebras that are not weakly-irreducible. Giving new
definitions to these types we obtain 11 types of weakly-irreducible Lie algebras. Unifying some of the types we obtain 7
types of weakly-irreducible subalgebras of $\su(1,n+1)_l\subset\so(2,2n+2)$.

\vskip 0.3cm

The result can be stated as  follows.

Let $n=0$. The Lie algebra $\su(1,1)_l$ is 2-dimensional nilpotent, we have $\su(1,1)_l=\Real\zr\Real$ and $[(a,0),(0,c)]=(0,2ac)$.
There are 2 weakly-irreducible subalgebras of $\su(1,1)_l$: $\{(0,c)|c\in\Real\}$ and the whole $\su(1,1)_l$.

Let $n>0$. For the Lie algebra $\su(1,n+1)_l$ we have the Iwasawa decomposition
$$\su(1,n+1)_l=(\Real\oplus\u(n))\zr\lie\H_n,$$ where $\lie\H_n=\Co^n\zr\Real$ is the Lie algebra of the Lie group $\H_n$
of the Heisenberg translations of the Heisenberg space $\H_n$.

Let $0\leq m\leq n$ be an integer. Consider the decomposition $\Co^n=\Co^m\oplus\Co^{n-m}$.
Let $\h\subset\u(m)\oplus\sod(n-m)$ be a subalgebra, here
$\sod(n-m)=\left.\left\{\left(\begin{smallmatrix}B&0\\0&B\end{smallmatrix}\right)\right|B\in\so(n-m)\right\}
\subset\su(n-m)$.
The Lie algebras of one of the types of weakly-irreducible subalgebras of $\su(1,n+1)_l\subset\so(2,2n+2)$ have the form
$$\g^{m,\h,\A^1}=(\Real\oplus\h)\zr((\Co^m\oplus\Real^{n-m})\zr\Real),$$
where $\Real^{n-m}\subset\Co^{n-m}$ is a real form.
The other types of weakly-irreducible subalgebras of $\su(1,n+1)_l\subset\so(2,2n+2)$ can be obtained from this one using some twisting
and they have the following forms:
\begin{description}

\item[$\g^{m,\h,\varphi}$=]$\{\varphi(A)+A|A\in\h\}\zr((\Co^m\oplus\Real^{n-m})\zr\Real) $,\\ where
  $\varphi:\h\to\Real$ is a linear map with $\varphi|_{\h'}=0$;

\item[$\g^{n,\h,\psi,k,l}$=]$\{A+\psi(A)|A\in\h\}\zr((\Co^k\oplus\Real^{n-l})\zr\Real) $,\\
 where $k$ and $l$ are integers such that $0< k\leq l\leq n$,
we have the decomposition $\Co^n=\Co^k\oplus\Co^{l-k}\oplus\Co^{n-l}$,
$\h\subset\u(k)$ is a subalgebra with $\dim\z(\h)\geq n+l-2k$  and
$\psi:\h\to\Co^{l-k}\oplus i\Real^{n-l}$ is a surjective linear map with  $\psi|_{\h'}=0$;

\item[$\g^{m,\h,\psi,k,l,r}$=]$\{A+\psi(A)|A\in\h\}\zr((\Co^k\oplus\Real^{m-l}\oplus\Real^{r-m})\zr\Real) $,\\
where $k$, $l$, $r$ and $m$ are integers such that $0< k\leq l\leq m\leq r\leq n$ and $m<n$,
we have the decomposition $\Co^n=\Co^k\oplus\Co^{l-k}\oplus\Co^{m-l}\oplus\Co^{r-m}\oplus\Co^{n-r}$,
$\h\subset\u(k)\oplus\sod(r-m)$ is a subalgebra with $\dim\z(\h)\geq n+m+l-2k-r$  and
$\psi:\h\to\Co^{l-k}\oplus \Real^{n-r}\oplus i\Real^{m-l}$ is a surjective linear map with  $\psi|_{\h'}=0$;

\item[$\g^{0,\h,\psi,k}$=]$ \{A+\psi(A)|A\in\h\}\zr(\Real^{k}\zr\Real) $,\\
 where $0<k<n$, we have the decomposition $\Co^n=\Co^k\oplus\Co^{n-k}$,
$\h\subset\sod(k)$ is a subalgebra such that $\dim\z(\h)\geq n-k$, $\psi:\h\to \Real^{n-k}$ is a surjective linear map with $\psi|_{\h'}=0$;


\item[$\g^{0,\h,\zeta}$=]
$\{A+\zeta(A)|A\in \h\}\zr\Real^{n}$,\\ where  $\h\subset\sod(n)$ is a subalgebra  with $\z(\h)\neq\{0\}$,
$\zeta:\h\to\Real\subset\lie\H_n$ is a non-zero linear map with $\zeta|_{\z(\h)}\neq 0$;

\item[$\g^{0,\h,\psi,k,\zeta}$=]$\{A+\psi(A)+\zeta(A)|A\in\h\}\zr\Real^{k}$,\\
where $1 \leq k<n$, we have the decomposition $\Co^n=\Co^{k}\oplus\Co^{n-k}$, $\h\subset\sod(k)$ is a subalgebra with
$\dim\z(\h)\geq n-k$, $\psi:\h\to \Real^{n-k}$ is a surjective linear map with $\psi|_{\h'}=0$,
$\zeta:\h\to\Real\subset\lie\H_n$ is a non-zero linear map with $\zeta|_{\h'}=0$.
\end{description}

Note that the last two types of  weakly-irreducible subalgebras $\g\subset\su(1,n+1)_l\subset\so(2,2n+2)$ were not
considered by A. Ikemakhen in \cite{Ik22}.

\vskip0.3cm

For each $\f\subset\lie(\Simil\H_n)$ as above and for each $\g\subset\su(1,n+1)_l$ with
$\pi(\Gamma(\g))=\f$ we  consider the  Lie algebras
$\g^J=\g\oplus \Real J$ and $\g^\xi=\{x+\xi(x)|x\in\g\}$, where $\xi:\g\to\Real$ is a non-zero linear map.
As we claimed above, any weakly-irreducible subalgebra of $\u(1,n+1)_l\subset\so(2,2n+2)$ is of the form $\g$, $\g^J$ or $\g^\xi$.
These subalgebras are candidates for the weakly-irreducible subalgebras of  $\u(1,n+1)_l\subset\so(2,2n+2)$.
We associate with each of these subalgebras  an integer $0\leq m\leq n$.
If $m>0$, then the  subalgebras of the form $\g$, $\g^J$ and $\g^\xi\subset\u(1,n+1)_l$  are weakly-irreducible.
We have  inclusions $\u(m)\subset\u(n)\subset\u(1,n+1)_l$ and  projection maps $\pr_{\u(m)}:\u(1,n+1)_l\to\u(m)$, $\pr_{\u(n)}:\u(1,n+1)_l\to\u(n)$.

\vskip0.3cm

In {\bf Chapter III} we classify weakly-irreducible Berger subalgebras of $\u(1,n+1)_l$. These subalgebras are candidates for the holonomy algebras.
Then we show that all these Berger algebras are holonomy algebras.

More precisely, for any subalgebra $\g\subset\u(1,n+1)_l$ consider the space $\R(\g)$ of curvature tensors of type $\g$,
$$\R(\g)=\left\{R\in\text{Hom}(\Real^{2,2n+2}\wedge \Real^{2,2n+2},\g)\,\left|\,\begin{array}{c}
R(u\wedge v)w+R(v\wedge w)u+R(w\wedge u)v=0\\ \text{ for all } u,v,w\in \Real^{2,2n+2}\end{array}\,\right\}\right..$$
Denote by $L(\R(\g))$ the vector subspace of $\g$ spanned by $R(u\wedge v)$ for all $R\in\R(\g),$ $u,v\in \Real^{2,2n+2},$
$$L(\R(\g))=\spa\{R(u\wedge v)|R\in\R(\g),\,u,v\in \Real^{2,2n+2}\}.$$
A subalgebra $\g\subset\u(1,n+1)_l$ is called {\it a Berger algebra} if $L(\R(\g))=\g$.
From the Ambrose-Singer theorem \cite{Am-Si} it follows that {\it if $\g\subset\u(1,n+1)_{l}$ is the holonomy algebra
of a pseudo-K\"ahlerian manifold, then $\g$ is a Berger algebra} (here we identify the tangent space to the manifold at some point with $\Real^{2,2n+2}$).

First we consider all subalgebras of $\u(1,1)_l$ and show which of these subalgebras are weakly-irreducible Berger subalgebras.

Then we consider the case $n\geq 1$. For any integer  $0\leq m\leq n$ and subalgebra $\un\subset\u(m)\oplus\sod(m+1,...,n)$
we consider a subalgebra $\g^{m,\un}\subset\u(1,n+1)_l$ and describe the space $\R(\g^{m,\un})$.
The Lie algebras of the form  $\g^{m,\un}$  contain all candidates for the weakly-irreducible subalgebras of $\u(1,n+1)_l$.
For any  subalgebra $\g\subset\g^{m,\un}$ the space $\R(\g)$ can be found from the following condition
$$R\in\R(\g)\text{ if and only if } R\in\R(\g^{m,\un}) \text{ and } R(\Real^{2,2n+2}\wedge \Real^{2,2n+2})\subset\g.$$
Using this, we easily find all weakly-irreducible not irreducible Berger subalgebras of $\u(1,n+1)_l$.

\vskip0.3cm

As the last step of the classification, we construct  metrics on $\Real^{2n+4}$ that realize all Berger algebras obtained above as holonomy algebras.
The coefficients of the metrics are polynomial functions, hence  the corresponding Levi-Civita connections
are analytic and in each case the holonomy algebra at the point $0\in\Real^{2n+4}$ is generated by the operators
$$R(X,Y)_0,\nabla_{Z_1} R(X,Y)_0,\nabla_{Z_2}\nabla_{Z_1} R(X,Y)_0,...$$
where   $X$, $Y$, $Z_1$, $Z_2$,... are vectors at the point $0$.
We explicitly compute for each metric the components of the curvature tensor and its derivatives. Then using the
induction, we find the holonomy algebra for each of the metrics.

\vskip 0.3cm

Thus we obtain the classification of weakly-irreducible not irreducible holonomy algebras contained in $\u(1,n+1)$.
The result can be stated as  follows.

Let $n=0$. The Lie algebra $\u(1,1)_l$ is a 3-dimensional nilpotent real Lie algebra, we have $\u(1,1)_l=\Co\zr\Real$ and
$[(a+ib,0),(0,c)]=(0,2ac)$, $a$, $b$, $c\in\Real$. There are three weakly-irreducible holonomy algebras contained in
$\u(1,1)_l$:

$$ \hol_{n=0}^1=\u(1,1)_l,\quad
\hol^{\gamma_1,\gamma_2}_{n=0}=\Real(\gamma_1+i\gamma_2)\zr\Real\,(\gamma_1,\gamma_2\in\Real),\quad \hol_{n=0}^2=\Co.$$

Let $n>0$. For the Lie algebra $\u(1,n+1)_l$ we have the Iwasawa decomposition
$$\u(1,n+1)_l=(\Co\oplus\u(n))\zr\lie\H_n.$$

Consider the following  type of weakly-irreducible holonomy algebras contained in $\u(1,n+1)_l\subset\so(2,2n+2)$:
$$\hol^{m,\un,\A^1,\t\A^2}=(\Real\oplus\Real(i+J_{n-m})\oplus\un)\zr((\Co^m\oplus\Real^{n-m})\zr\Real),$$ where
$0\leq m\leq n$ is an integer, we have the decompositions $\Co^n=\Co^m\oplus\Co^{n-m}$, $\Co=\Real\oplus i\Real$,
  $\un\subset\u(m)$ is a subalgebra,   $\Real^{n-m}\subset\Co^{n-m}$ is a real form and $J_{n-m}\subset\u(n-m)\subset\u(n)\subset\u(1,n+1)_l$
is the complex structure on $\Co^{n-m}$.

The other types of weakly-irreducible holonomy algebras contained in $\u(1,n+1)_l\subset\so(2,2n+2)$ can be obtained from this  type using some twisting.
 Let $\varphi,\phi:\un\to\Real$ be  linear maps with $\varphi|_{\un'}=\phi|_{\un'}=0$.
The other types have the following forms:

\begin{description}

\item[$\hol^{m,\un,\A^1,\phi}$=]$(\Real\oplus\{\phi(A)(i+J_{n-m})+A|A\in\un\})\zr((\Co^m\oplus\Real^{n-m})\zr\Real),$
\item[$\hol^{m,\un,\varphi,\phi}$=]$\{\varphi(A)+\phi(A)(i+J_{n-m})+A|A\in\un\}\zr((\Co^m\oplus\Real^{n-m})\zr\Real),$
\item[$\hol^{m,\un,\varphi,\t\A^2}$=]$(\Real(i+J_{n-m})\oplus\{\varphi(A)+A|A\in\un\})\zr((\Co^m\oplus\Real^{n-m})\zr\Real),$
\item[$\hol^{m,\un,\lambda}$=]$(\Real(1+\lambda(i+J_{n-m}))\oplus\un)\zr((\Co^m\oplus\Real^{n-m})\zr\Real),$ where $\lambda\in\Real$,
\item[$\hol^{n,\un,\psi,k,l}$=]$\{A+\psi(A)|A\in\un\}\zr((\Co^k\oplus\Real^{n-l})\zr\Real) $,\\
where $k$ and $l$ are integers such that $0< k\leq l\leq n$,
we have the decomposition $\Co^n=\Co^k\oplus\Co^{l-k}\oplus\Co^{n-l}$,
$\un\subset\u(k)$ is a subalgebra with $\dim\z(\un)\geq n+l-2k$  and
$\psi:\un\to\Co^{l-k}\oplus i\Real^{n-l}$ is a surjective linear map with  $\psi|_{\un'}=0$,
\item[$\hol^{m,\un,\psi,k,l,r}$=]$\{A+\psi(A)|A\in\un\}\zr((\Co^k\oplus\Real^{m-l}\oplus\Real^{r-m})\zr\Real) $,\\
where $k$, $l$, $r$ and $m$ are integers such that $0< k\leq l\leq m\leq r\leq n$ and $m<n$,
we have the decomposition $\Co^n=\Co^k\oplus\Co^{l-k}\oplus\Co^{m-l}\oplus\Co^{r-m}\oplus\Co^{n-r}$,
$\un\subset\u(k)$ is a subalgebra with $\dim\z(\un)\geq n+m+l-2k-r$  and
$\psi:\un\to\Co^{l-k}\oplus i\Real^{m-l}\oplus\Real^{n-r}$ is a surjective linear map with  $\psi|_{\un'}=0$.
\end{description}

\vskip0.3cm
As a corollary, we get the classification of  weakly-irreducible not irreducible holonomy algebras contained in $\su(1,n+1)$
(i.e. of the  holonomy algebras of special pseudo-K\"ahlerian manifolds). For $n=0$ these algebras are exhausted by $\{(0,c)|c\in\Real\}$ and $\su(1,1)_l$.
For $n>0$ these algebras are the following:

\begin{description}
\item[] $\hol^{m,\un,\A^1,\phi}$, $\hol^{m,\un,\varphi,\phi}$ with $\phi(A)=-\frac{1}{n-m+2}\tr_\Co A$;
\item[] $\hol^{n,\un,\psi,k,l}$, $\hol^{m,\un,\psi,k,l,r}$ with $\un\subset\su(k)$.
\end{description}

\vskip0.3cm

The above result  together with the Wu theorem and the classification  of irreducible holonomy algebras of M. Berger
gives us the classification of holonomy algebras (or equivalently, of connected holonomy groups) for
pseudo-K\"ahlerian manifolds of signature  $(2,2n+2)$.

Remark that we do not have any additional condition on the $\u(m)$-projection of a holonomy algebra, while in the Lorentzian case
an analogous subalgebra $\h\subset\so(n)$ associated to a holonomy algebra must be the holonomy algebra of a Riemannian manifold.
This shows the principal difference between our case and the case of Lorentzian manifolds.

\vskip0.2cm

In {\bf Chapter IV} we consider some examples and applications.

In Section \ref{secLie} we  give examples of 4-dimensional Lie groups with left-invariant \pK metrics that have holonomy
algebras $\hol_{n=0}^2$ and $\hol^{\gamma_1=0,\gamma_2=1}_{n=0}$.

\vskip0.1cm

In \cite{Ber57} M.~Berger obtained a classification of simply connected semi-simple pseudo-Riemannian symmetric spaces.
The holonomy algebra of any such manifold is either irreducible or it is weakly-irreducible and it preserves two
complementary isotropic subspaces (in the last case the manifold has signature $(n,n)$). In particular, there are three
such spaces with weakly-irreducible not irreducible holonomy algebras  of signature $(2,2)$. In \cite{KathOldrich} I.~Kath
and M.~Olbrich  obtained a classification of simply connected pseudo-Riemannian symmetric spaces of index 2 that are not
semi-simple. In particular, the holonomy algebras of these spaces are weakly-irreducible and not irreducible. In Section
\ref{secsym} we use our classification of holonomy algebras to give a new proof for the classification of simply connected
\pK symmetric spaces of index 2 with weakly-irreducible not irreducible holonomy algebras.
 We use the fact that any simply connected \pK symmetric spaces of index 2 is uniquely defined by a pair
 $(\hol,R)$, where $\hol\subset\u(1,n+1)$ is a  holonomy algebra and
$R\in\R_(\hol)$ is a curvature tensor that is annihilated by the natural representation of  the Lie algebra $\hol$ in the
vector space $\R_(\hol)$ and such that  $R(\Real^{2,2n+2}\wedge\Real^{2,2n+2})=\hol$. We find all such pairs and  check
which of them define isometric simply connected symmetric spaces. We show that all possible weakly-irreducible not
irreducible holonomy algebras of locally symmetric \pK manifolds of index 2 are exhausted by the following Lie algebras:
$\hol^2_{n=0}$, $\hol^{\gamma_1=0,\gamma_2=0}_{n=0}$, $\hol^{m=n-1=0,\{0\},\varphi=0,\phi=0}$ and $\hol^{m,\{\Real
J_m\},\varphi=0,\phi=2}$.

\vskip0.1cm

Finally in Section \ref{secSasaki} we consider  Lorentzian  manifold $(M,g)$ such that the time-like cone $(\t
M^-=\Real^+\times M,\t g^-=-dr^2+r^2g)$ over $(M,g)$  is a \pK manifold of index 2. Such Lorentzian  manifold is called
{\it a Lorentzian Sasaki manifold}.

First we describe the local DeRham-Wu decomposition of the cone in terms of the initial Lorentzian Sasaki manifold.
We show that {\it if the cone $(\t{M}^-,\t{g}^-)$ is locally decomposable, then $M$ locally has the form
$$((a,b)\times N_1\times N_2,ds^2+ \ch^2(s) g_1+\sh^2(s)g_2),\quad (a,b)\subset\bR^+,$$
where $(N_1,g_1)$ is a Lorentzian Sasaki manifold and $(N_2,g_2)$ is a Riemannian Sasaki manifold
(i.e. the space-like cone $(\t N_2^+=\Real^+\times N_2,\t g^+=dr^2+r^2g_2))$ over $(N_2,g_2)$  is a K\"ahlerian manifold).}

Then we study Lorentzian Sasaki manifold $(M,g)$ such that the holonomy algebra $\hol(\t M^-)$ of its time-like cone $(\t{M}^-,\t{g}^-)$ (with
a \pK structure $\J$)  preserves a 2-dimensional $\J$-invariant isotropic subspace.
We show that {\it in this situation  $\hol(\t M^-)$ annihilates this subspace,
i.e. there exist locally on $\t M^-$ two isotropic parallel vector fields $p_1$ and $p_2=\t J p_1$.
Furthermore, $(M,g)$  locally has the form
 $$((a,b)\times N,ds^2+e^{-2s}g_N),\quad (a,b)\subset\Real,$$
 where $(N,g_N)$ is a Lorentzian manifold with a parallel isotropic vector field.
There exist local coordinates $x,y,\hat x,\hat y,x_1,...,x_{2n}$ on $\t M^-$ such that
$\hat x,\hat y,x_1,...,x_{2n}$ are coordinates on $N$ and
$$\t g^-=2dxdy+y^2g_N=2dxdy+y^2(2d\hat xd\hat y+g_1),$$
where $g_1$ is a  $\hat y$-family   of K\"ahlerian metrics on the integral manifolds corresponding
to the coordinates $x_1,...,x_{2n}$. In these coordinates, $p_1=\p_x$ and $p_2=\hat y \p_x-\frac{1}{y}\p_{\hat x}$.

Moreover, if the local holonomy algebra $\hol(\t M^-)$ of $(\t M^-,\t g^-)$ is weakly-irreducible,
then the holonomy algebra $\hol(N)$ of $(N,g_N)$ is weakly-irreducible and
\begin{itemize}\item[1)] if $\hol(N)$ is of type $\g^{2,\u},$ $\u\subset\u(n)$ (see Section \ref{sec1.3} below), then $\hol(\t M^-)$ is
of type $\hol^{n,\u,\varphi=0,\phi=0}$;
\item[2)] if $\hol(N)$ is of type $\g^{4,\u,k,\psi},$  then $\hol(\t M^-)$ is
of type $\hol^{n,\u,\psi,k,l}$.\end{itemize}     }

\parti{Acknowledgments}


I would like to express sincere gratitude to my advisor Helga Baum for all-round help, guidance, support and encouragement.\\

I wish to acknowledge  Dmitri Vladimirovich Alekseevsky  for introducing me to the theory of holonomy groups and many useful suggestions.\\

I am deeply grateful to Mark Volfovich Losik for all he taught me at Saratov University and for his constant encouragement since 1998. \\

I thank Ines Kath and Martin Olbrich, who helped me to find many mathematical mistakes in this thesis.\\

I thank Lionel Berard Bergery,  Vicente Cort\'es, Thomas Leistner and Thomas Neukirchner for useful discussions.\\

This work was supported by DFG International Research Training group ``Arithmetic and Geometry'' (GRK 870) at
Humboldt-Universit\"at zu Berlin. I am thankful to the speaker of GRK 870  J\"urg Kramer, secretary Marion Thomma
and all fellows of GRK 870 for support, good atmosphere and opportunity to work.


\part{Holonomy groups of pseudo-Riemannian manifolds}

This chapter contains an introduction to the theory of  holonomy groups of pseudo-Riemannian manifolds.
In Section \ref{sec1.1} we provide definitions,
general facts, examples and ideas of some proofs that illustrate the methods of the holonomy.
In Section \ref{sec1.2} we recall the classification of  connected holonomy groups for Riemannian manifolds and
of connected irreducible holonomy groups for pseudo-Riemannian manifolds.
In Section \ref{sec1.3} we explain the classification of
connected holonomy groups for Lorentzian manifolds.

\section{Definitions and facts}\label{sec1.1}

All definitions and statements of this section can be found in
\cite{K-N} or in \cite{Be}. We assume that all manifolds are connected.

{\bf Definition.} {\it A pseudo-Riemannian manifold of signature $(r,s)$
is a differential manifold  $M$ equipped with a smooth field $g$ of
symmetric non-degenerate bilinear forms of signature $(r,s)$ at each point.
We assume that $r$ is the number of minuses and we call it the index of $(M,g)$.
If $r=0$,  then $(M,g)$ is a Riemannian manifold;
if $r=1$,  then $(M,g)$ is a Lorentzian manifold.}

On any pseudo-Riemannian manifold  $(M,g)$  exists the Levi-Civita connection
$\nabla$ which is defined by two conditions: $g$ is parallel $(\nabla  g=0)$  and
the torsion is zero $(Tor=0)$. The Levi-Civita connection
$\nabla$  is explicitly  given by the Koszul formulae
\begin{equation}\label{Koszul} \begin{array}{rl} 2g(\n_XY,Z)=&Xg(Y,Z)+Yg(X,Z)-Zg(X,Y)\\
&+g([X,Y],Z)+g([Z,X],Y)+g(X,[Z,Y]),\end{array}
\end{equation}
where $X$, $Y$ and $Z$ are vector fields on $M$.

It is known that for any  smooth curve $\gamma: [a,b]\subset\Real\to M$ and any vector $X_0\in T_{\gamma(a)}M$ there exists a unique
vector field $X$ defined along the curve $\gamma$ and satisfying the differential equation
$\nabla_{\dot\gamma(s)}X=0$ with the initial condition $X_{\gamma(a)}=X_0$.
Consequently, for any  smooth curve
$\gamma: [a,b]\subset\Real\to M$ we obtain the isomorphism
$\tau_{\gamma}:T_{\gamma(a)}M\to T_{\gamma(b)}M$ defined by $\tau_{\gamma}:X_0\mapsto X_{\gamma(b)}$.
The isomorphism $\tau_{\gamma}$ is called {\it the parallel displacement along the curve} $\gamma$.
The parallel displacement can be defined in the obvious way also for piecewise smooth curves.
Note that for the curve  $\gamma_{-}: [a,b]\to M$
($\gamma_{-}(t)=\gamma(a+b-t)$)
the isomorphism $\tau_{\gamma_{-}}:T_{\gamma(b)}M\to T_{\gamma(a)}M$
is the inverse to $\tau_{\gamma}$; for any piecewise smooth curve
$\gamma': [b,c]\to M$ such that $\gamma'(b)=\gamma(b)$
we have $\tau_{\gamma'*\gamma}=\tau_{\gamma'}\circ\tau_{\gamma}$,
where $\gamma'*\gamma$ is the curve given by $\gamma'*\gamma(t)=\gamma(t)$ if $a \leq t\leq b$ and
$\gamma'*\gamma(t)=\gamma'(t)$ if $b \leq t\leq c$.
For the constant curve $\gamma(t)\equiv x\in M$ we have $\tau_{\gamma}=\id_{T_xM}$.

 Let $x\in M$. We  denote by $Hol_x$ the set of parallel displacements along all
piecewise smooth loops at the point $x\in M$. Let  $Hol^0_x$ be
the set of parallel displacements along all
piecewise smooth null-homotopic loops at the point $x\in M$.
From the properties of the parallel displacements it follows that $Hol_x$ and $Hol^0_x$ are groups.
Obviously, $Hol^0_x\subset Hol_x$ is a subgroup. If the manifold
$M$ is simply connected, then $Hol^0_x=Hol_x$.

{\bf Definition.} {\it The group $Hol_x$ is called the holonomy group of the manifold
$(M,g)$ at the point $x$. The group $Hol^0_x$ is called the restricted holonomy group of the manifold $(M,g)$ at the point $x$.}

Let $O(T_xM,g_x)$ be the Lie group of isomorphisms of the vector space $T_xM$ that preserve the form $g_x$
and $\so(T_xM,g_x)$ the corresponding Lie algebra. Since $g$ is parallel, we see that $Hol_x\subset O(T_xM,g_x)$.

\begin{theorem}
 The group $Hol_x$ is a Lie subgroup of the Lie group $O(T_xM,g_x)$.
The group $Hol^0_x$ is the connected identity component of the Lie group $Hol_x$.
\end{theorem}

By $\hol_x$ we denote the Lie algebra of the Lie group $Hol_x$ (and of $Hol^0_x$).

{\bf Definition.} {\it The Lie subalgebra $\hol_x\subset\so(T_xM,g_x)$ is called the holonomy algebra
of the manifold $(M,g)$ at the point $x$.}

Remark that by the  holonomy group (resp. holonomy algebra) we understand not just the Lie group $Hol_x$  (resp. Lie algebra $\hol_x$),
but the Lie group $Hol_x$ with the representation $Hol_x\hookrightarrow O(T_xM,g_x)$ (resp. the  Lie algebra $\hol_x$ with the representation
$\hol\hookrightarrow\so(T_xM,g_x)$). This representation is called {\it the holonomy representation}.

Let $\gamma$ be a piecewise smooth curve in $M$ beginning at the point $x$ and ending
at a point $y\in M$. Then we have
$Hol_y=\tau_{\gamma}\circ Hol_x\circ\tau_{\gamma}^{-1}$,
$Hol^0_y=\tau_{\gamma}\circ Hol^0_x\circ\tau_{\gamma}^{-1}$
and $\hol_y=\tau_{\gamma}\circ\hol_x\circ\tau_{\gamma}^{-1}$.
This means that the holonomy groups at all points of the manifold
are isomorphic and we can speak about the holonomy group  $Hol$ of the manifold $(M,g)$
 (or about the restricted holonomy group $Hol^0$ of $(M,g)$ or about the holonomy algebra $\hol$  of $(M,g)$).
The holonomy algebra $\hol$ and the restricted
holonomy group $Hol^0$ uniquely define each other.
If the manifold $M$ is simply connected, then the holonomy algebra $\hol$
uniquely defines the holonomy group $Hol$.

Denote by $\Real^{r,s}$ the $(r+s)$-dimensional pseudo-Euclidean space
endowed with a non-degenerate symmetric bilinear form  $\eta$ of signature $(r,s)$.
For a pseudo-Riemannian manifold $(M,g)$ of signature $(r,s)$
we identify the tangent space  $(T_xM,g_x)$ at a point $x\in M$ with the space
$\Real^{r,s}$. Then the holonomy algebra $\hol_x$ can be identified with a subalgebra
of the pseudo-orthogonal Lie algebra $\so(r,s)$ and the holonomy group $Hol_x$ can be identified
with a Lie subgroup of the pseudo-orthogonal Lie group $O(r,s)$. We may write $\hol_x\subset\so(r,s)$ and $Hol_x\subset O(r,s)$.
As we claimed above, these representations are defined up to conjugacy.

Let $A$ be a tensor field of type $(p,q)$ on the manifold  $(M,g)$.
Recall that $A$ is called parallel if $\nabla A=0$.
This is equivalent to the condition that for any piecewise smooth curve
 $\gamma:[a,b]\to M$ holds
$\tau_{\gamma q}^pA_{\gamma(a)}=A_{\gamma(b)}$, where
$$\tau_{\gamma q}^p=
\underbrace{\tau_{\gamma}\otimes\cdots\otimes\tau_{\gamma}}_{p
\text{ times}}\otimes \underbrace{
(\tau_{\gamma}^*)^{-1}\otimes\cdots\otimes(\tau_{\gamma}^*)^{-1}}
_{q \text{ times}}:T_{\gamma(a)q}^pM\to T_{\gamma(b)q}^pM$$
is the tensorial extension of the isomorphism
$\tau_{\gamma}:T_{\gamma(a)}M\to T_{\gamma(b)}M$ to the isomorphism of the spaces of tensors of type $(p,q)$.
A distribution $E\subset TM$ is called parallel if for any vector field $X$ with values in
$E$ and any vector field $Y$ on $M$ the vector field
 $\nabla_YX$ is also with values in  $E$. This is equivalent to the condition that
for any piecewise smooth curve $\gamma:[a,b]\to M$ holds
$\tau_{\gamma}E_{\gamma(a)}=E_{\gamma(b)}$.

The importance of the holonomy groups shows the following theorem.

\begin{theorem}\label{FP1} ({\rm Fundamental principle 1.})
 For a pseudo-Riemannian manifold $(M,g)$ the following conditions are equivalent:
\begin{description}
\item[1)] There exists  a parallel tensor field $A$ of type $(p,q)$ on  $(M,g)$.

\item[2)] For some $x\in M$ there exists  a tensor $A_x$ of type $(p,q)$ on $T_xM$
that is invariant under the tensorial extension of the holonomy  representation
of the holonomy group $Hol_x\subset O(T_xM,g_x)$.
\end{description}\end{theorem}

 {\it The idea of the proof} is the following.
For a given parallel tensor field $A$ take the value $A_x$ at a chosen point $x\in M$. Since $A$  is invariant under the
parallel displacements, we see that the tensor $A_x$ is invariant under the parallel displacements along the loops at the
point $x$, i.e.
 under the tensorial extension of the holonomy representation.
Conversely, for a given tensor $A_x$ define the tensor field $A$ on $M$ such that at any point  $y\in M$ holds $A_y=\tau_{\gamma q}^pA_x$,
 where $\gamma$ is any curve beginning at  $x$ and ending at $y$. From the condition 2) it follows that $A_y$ does not depend on the curve $\gamma$. $\Box$

 The local analog of the previous theorem is the following

\begin{theorem}\label{FP1'} ({\rm Fundamental principle 1'.})
 For a pseudo-Riemannian manifold $(M,g)$ the following conditions are equivalent:
\begin{description}

\item[1)] There exists  a parallel tensor field $A$ of type $(p,q)$ on an open neighbourhood of a point $x\in M$.

\item[2)] There exists   a tensor $A_x$ of type $(p,q)$ on $T_xM$
that is invariant under  the tensorial extension of the holonomy representation
of the restricted holonomy group $Hol^0_x\subset O(T_xM,g_x)$.

\item[3)] There exists   a tensor $A_x$ of type $(p,q)$ on $T_xM$
 that is annihilated by the tensorial extension of the holonomy representation
of the holonomy algebra $\hol_x\subset \so(T_xM,g_x)$.

\end{description}\end{theorem}

Analog results hold for distributions.

\begin{theorem}\label{FP2} ({\rm Fundamental principle 2.})
 For a pseudo-Riemannian manifold $(M,g)$ the following conditions are equivalent:
\begin{description}

\item[1)] There exists  a parallel distribution $E$ of rang $p$ on  $(M,g)$.

\item[2)] For some $x\in M$ there exists  a vector subspace $E_x\subset T_xM$ of dimension $p$
that is invariant under the tensorial extension of the holonomy representation of the holonomy group $Hol_x\subset O(T_xM,g_x)$.
\end{description}\end{theorem}

\begin{theorem}\label{FP2'} ({\rm Fundamental principle 2'.})
 For a pseudo-Riemannian manifold $(M,g)$ the following conditions are equivalent:
\begin{description}

\item[1)] There exists  a parallel distribution $E$ of rang $p$ on an open neighbourhood of a point $x\in M$.

\item[2)] There exists a vector subspace $E_x\subset T_xM$  of dimension $p$
that is invariant under  the tensorial extension of the  holonomy representation
of the restricted holonomy group $Hol^0_x\subset O(T_xM,g_x)$.

\item[3)] There exists a vector subspace $E_x\subset T_xM$  of dimension $p$
that is invariant under to the tensorial extension of the holonomy representation
of the holonomy algebra $\hol_x\subset O(T_xM,g_x)$.

\end{description}\end{theorem}

Note that parallel distributions are involutive. Indeed, since
the torsion of the Levi-Civita connection is zero,
for any vector fields $X$ and $Y$ with values in  a parallel distribution $E$
 we have $[X,Y]=\nabla_XY-\nabla_YX\in E$.

Thus if we know the holonomy group of a pseudo-Riemannian manifold, then
the geometric problem to find  parallel tensor fields or parallel distributions
on the manifold can be reduced to the algebraic problem to find  invariant
subspaces for some representations of the holonomy group.

Consider several examples of using of the above theorems.

\begin{ex} A pseudo-Riemannian manifold $(M,g)$ of signature $(r,s)$
is orientable if and only if $Hol\subset SO(r,s)$.
\end{ex}
{\it Proof.} Recall that the manifold  $M$ is orientable if and only if
there exists  a nowhere vanishing differential form of degree $n=\dim M$ on $M$.
Furthermore, $SO(r,s)=\{A\in O(r,s)|\det
A=1\}$. Let $x\in M$. Suppose that $Hol_x\subset
SO(T_xM,g_x)$. Let $e_1,...,e_n$ be a basis of the vector space $T_xM$.
Consider the differential form $w_x=e_1^*\wedge\cdots\wedge e_n^*$. For $A\in
SO(T_xM,g_x)$ we have $A\cdot w_x=(Ae_1)^*\wedge\cdots\wedge (Ae_n)^*=(\det A)w_x=w_x$. From Theorem \ref{FP1} it follows that  $w_x$ defines
a parallel nowhere vanishing $n$-form $w$ on $M$.
Conversely, suppose that  $M$
is orientable, then there exists  a nowhere vanishing differential $n$-form $w$ on $M$.
Consider the differential form $\bar w=fw$, where $f$ is a function on $M$.
The condition $\nabla\bar w=0$ is equivalent to the condition $(Xf)w+f\nabla_Xw=0$ for all $X\in TM$. Since $w$ is nowhere vanishing, we have
$Xf=-f\frac{\nabla_Xw}{w}$. This system of differential equations has a unique global solution satisfying  $f(x)=1$. From Theorem \ref{FP1} it follows that the group
 $Hol_x$ preserves the form $\bar w_x$. Thus for any $A\in Hol_x$ we have  $\det A=1$, i.e. $Hol_x\subset SO(T_xM,g_x)$. $\Box$

Recall that a pseudo-Riemannian manifold is called flat if it admits  parallel local
fields of frames. From Theorem \ref{FP1'} we obtain

\begin{ex}\label{ex2} A pseudo-Riemannian manifold   $(M,g)$ is flat
if and only if $Hol^0=\{\id\}$.
\end{ex}

A pseudo-Riemannian manifold $(M,g)$ is called {\it pseudo-K\"ahlerian}
if  there exists  a parallel smooth  field of endomorphisms $J$ of the tangent bundle of $M$
that satisfies  $J^2=-\id$ and $g(JX,Y)+g(X,JY)=0$ for all vector fields $X$ and $Y$ on $M$.

\begin{ex}\label{exIn3} A pseudo-Riemannian manifold  $(M,g)$ of signature $(2r,2s)$
is pseudo-K\"ahlerian if and only if  $Hol\subset U(r,s)$.
\end{ex}
{\it Proof.} The inclusion $Hol_x\subset U(r,s)$ is equivalent to the existence of an endomorphism $J_x$ of the vector space $T_xM$ that satisfies
$J_x^2=-\id_{T_xM}$ and $g_x(J_xX,Y)+g(X,J_xY)=0$ for all $X,Y\in
T_xM$. Now the statement follows from Theorem \ref{FP1}. $\Box$

A pseudo-Riemannian manifold $(M,g)$ is called {\it special pseudo-K\"ahlerian}
if it is pseudo-K\"ahlerian and Ricci-flat. The following statement is well known.

\begin{ex}\label{exIn4} A pseudo-K\"ahlerian  manifold  $(M,g)$ of signature $(2r,2s)$
is special pseudo-K\"ahlerian if and only if  $\hol\subset \su(r,s)$.
\end{ex}

Let $(M,g)$ be a pseudo-Riemannian manifold.
The curvature tensor  $R$ of $(M,g)$ has the following properties:
\begin{description}
\item[a)] $R$ is a tensor field of type $(1,3)$, i.e. for any vector fields
 $X$, $Y$ and $Z$ on $M$,  $R(X,Y)Z$ is a vector field;

\item[b)] $R$ is skew-symmetric with respect to the first two variables, i.e.
$R(X,Y)=-R(Y,X)$;

\item[c)] $R(X,Y)Z+R(Y,Z)X+R(Z,X)Y=0$ (the Bianchi identity).

\end{description}

The next theorem gives relation between the holonomy algebra and  the curvature tensor of the manifold. It was proved by W.~Ambrose and I.~M.~Singer in 1953 in
\cite{Am-Si}.

\begin{theorem}\label{Am-Si} Let $(M,g)$ be a pseudo-Riemannian manifold and $x\in M$. The holonomy algebra $\hol_x$ is spanned by the following endomorphisms
$$\tau_{\gamma}^{-1}\circ R(\tau_{\gamma}U,\tau_{\gamma}V)\circ\tau_{\gamma}:
T_xM\to T_xM,$$
where $U,V\in T_xM$ and $\gamma$ is a piecewise smooth curve beginning at the point $x$.
\end{theorem}

Using Example \ref{ex2} and Theorem \ref{Am-Si}, we obtain
\begin{ex} A pseudo-Riemannian manifold $(M,g)$ is flat if and only if $R=0$.
\end{ex}

{\bf Definition.} {\it Let $\g\subset \so(r,s)$ be a subalgebra.
The space of curvature tensors $\R(\g)$ of type $\g$ is defined as follows
$$\R(\g)=\left\{R\in\Hom(\Real^{r,s}\wedge \Real^{r,s},\g)\,\left|\begin{array}{c}R(u\wedge v)w+R(v\wedge w)u+R(w\wedge u)v=0\\
\text{ }\text{\it for all  } u,v,w\in \Real^{r,s}\end{array}\right\}\right..$$
Denote by $L(\R(\g))$ the vector subspace of $\g$ spanned by the elements
 $R(u\wedge v)$ for all $R\in\R(\g)$ and $u,v\in \Real^{r,s}$.}

Let $\g\subset\so(r,s)$.
It is known that any $R\in\R(\g)$ satisfies
\begin{equation}
 \eta (R(u\wedge v)z,w)=\eta
(R(z\wedge w)u,v)\, \label{property}
\end{equation}
for all $u,v,z,w\in \Real^{r,s}$.

The next proposition follows from Theorem \ref{Am-Si}.

\begin{prop}\label{curv10} If a subalgebra $\g\subset\so(r,s)$
is the holonomy algebra of a pseudo-Riemannian manifold, then $L(\R(\g))=\g$.\end{prop}

{\bf Definition.} {\it A subalgebra $\g\subset \so(r,s)$ is called a Berger algebra if  $L(\R(\g))=\g$.}

The condition $L(\R(\g))=\g$ is quite strong and by Proposition \ref{curv10},  the Berger algebras can be
considered as candidates for the holonomy algebras. We use the notion of Berger algebras, as such irreducible
subalgebras of $\so(n)$ were classified by M.~Berger (see Theorem \ref{Ber} below).

Theorem \ref{Am-Si} in general does not allow to find the holonomy algebra, since it involves all parallel displacements.
The following theorem can be used to find the holonomy algebra of an analytic pseudo-Riemannian manifold.

\begin{theorem} \label{analythol} If a pseudo-Riemannian manifold $(M,g)$ is analytic,
then the holonomy algebra $\hol_x$ is generated by the following operators
$$R(X,Y)_x,\nabla_{Z_1} R(X,Y)_x,
\nabla_{Z_2}\nabla_{Z_1} R(X,Y)_x,...\in\so(T_xM,g_x),$$
where  $X$, $Y$, $Z_1$, $Z_2,...\in T_xM$.
\end{theorem}

We say that a subspace $U\subset \Real^{r,s}$ is non-degenerate if the restriction of the pseudo-Euclidean metric  $\eta$ to $U$ is non-degenerate.

{\bf Definition.} {\it  A Lie subgroup $G\subset O(r,s)$
(or a subalgebra $\g\subset \so(r,s)$)
is called irreducible if it does not preserve any proper vector subspace of  $
\Real^{r,s}$; $G$ (or $\g$) is called weakly-irreducible if it does not preserve any proper non-degenerate vector subspace of  $\Real^{r,s}$.}

It is clear that $\g\subset \so(r,s)$ is irreducible (resp. weakly-irreducible)
if and only if the corresponding connected Lie subgroup  $G\subset SO(r,s)$ is
irreducible (resp. weakly-irreducible).
If a subgroup $G\subset O(r,s)$ is irreducible, then it is weakly-irreducible.
The converse holds only for positively and negatively definite metrics $\eta$.

Let $(M,g)$ and $(N,h)$ be pseudo-Riemannian manifolds. Let $x\in M$,
$y\in N$  and $G_x$, $H_y$ be the corresponding holonomy groups.
The product $M\times N$  is a pseudo-Riemannian manifold with the metric $g+h$.
Denote by $F_{(x,y)}$ the holonomy group of the manifold  $M\times N$ at the point $(x,y)$.

\begin{theorem}\label{H} We have $F_{(x,y)}=G_x\times H_y$. \end{theorem}

{\bf Definition.} {\it A pseudo-Riemannian manifold is called locally indecomposable
if its restricted holonomy group is weakly-irreducible.}

Theorem \ref{H} has the following inversion.

\begin{theorem}\label{Wu}
 Let $(M,g)$ be a pseudo-Riemannian manifold and $x\in M$.
Then there exists a decomposition of $T_xM$ into an orthogonal direct sum of non-degenerate vector subspaces
$$T_xM=E_0\oplus E_1\oplus\cdots\oplus E_t$$ such that $Hol^0_x$ acts trivially on $E_0$,
$Hol^0_x(E_i)\subset E_i$ ($i=1,...,t$),
$Hol^0_x$ acts weakly-irreducibly on $E_i$ ($i=1,...,t$) and
$$Hol^0_x=\{\id\}\times H_1\times\cdots\times H_t,$$
where $H_i=Hol_x^0|_{E_i}\subset SO(E_i)$ ($i=1,...,t$).

Furthermore, if $(M,g)$ is simply connected and complete, then
exist a flat pseudo-Riemannian submanifold $N_0\subset M$
and locally indecomposable pseudo-Riemannian submanifolds
$N_1,...,N_t\subset M$ such that $T_xN_i=E_i$ ($i=0,...,t$) and $(M,g)$ is isometric
to the product $$(N_0\times N_1\times\cdots\times N_t,g|_{N_0}+g|_{N_1}+\cdots+g|_{N_t}).$$
In general, such isometry of $(M,g)$  exists locally.
\end{theorem}

The local version of Theorem \ref{Wu} for Riemannian manifolds was proved in 1952 by A.~Borel and A.~Lichnerowicz in \cite{Bo-Li}.
The global version  for Riemannian manifolds was proved in the same year by G.~DeRham in \cite{Rha}.
For pseudo-Riemannian manifolds Theorem \ref{Wu}  was proved by H.~Wu in 1967 in \cite{Wu}.

The following proposition is important for us.

\begin{prop}\label{utv1} The Lie algebras of the Lie groups  $H_i$ from Theorem \ref{Wu}
are Berger algebras. \end{prop}

We explain shortly  the ideas of the proofs of the local part of Theorem  \ref{Wu} and of proposition \ref{utv1}.
Suppose that $Hol^0_x$ preserves a non-degenerate proper vector subspace  $E_x\subset T_xM$.
The vector subspace $E_x^{\bot}\subset T_xM$ is also non-degenerate and preserved.
 We obtain the orthogonal direct sum $T_xM=E_x\oplus E_x^\bot$ of $Hol_x^0$-invariant vector subspaces.
From Theorem \ref{FP2'} it follows that in a neighbourhood of the point $x$ exist parallel distributions
$E$ and $E^\bot$. As it was remarked, these distributions are involutive. Let $M_1$ and $M_2$
be the corresponding maximal connected  integral  submanifolds of $M$ through the point $x$. Since the distributions  $E$ and $E^\bot$ are parallel, the submanifolds
 $M_1$ and $M_2$ are totally geodesic. Let  $V\subset T_xM$ be an open neighbourhood of the point $0$ such that the restriction of the exponential map to
$V$ is a diffeomorphism, then $\exp(V)=\exp(E\cap V)\times(E^\bot\cap V)$, where
$\exp(E\cap V)\subset M_1$ and $\exp(E^\bot\cap V)\subset M_2$ are open submanifolds.

Consider now  the subalgebras $\g_1=\{\xi\in\hol_x| \xi(E^\bot)=\{0\}\}
\subset\hol_x$ and $\g_2=\{\xi\in\hol_x| \xi(E)=\{0\}\}
\subset\hol_x$. Obviously $\g_1$ and $\g_2$ are ideals and
$\g_1\cap\g_2=\{0\}$. Let $R\in\R(\hol_x),$
$X_1,X_2,X_3\in E_x$ and $Y_1,Y_2\in E_x^\bot$.
Since $R(X_1\wedge X_2)Y_1+R(X_2\wedge Y_1)X_1+R(Y_1\wedge X_1)X_2=0$,
we have
$R(X_1\wedge X_2)Y_1=0$, i.e. $R(X_1\wedge X_2)\in\g_1$. Similarly we can show that
 $R(Y_1\wedge Y_2)\in\g_2$.
Using \eqref{property}, we get
$g(R(X_1\wedge Y_1)X_2,X_3)=g(R(X_2\wedge X_3)X_1,Y_1)=0$.
Therefore, $R(X_1\wedge Y_1)=0$. Hence,
$R(E\wedge E^\bot)=0$, $R(E\wedge E)\subset \g_1$ and
$R(E^\bot\wedge E^\bot)\subset \g_2$. Clearly, $R|_{E_x\wedge
E_x}\in\R(\g_1)$ and $R|_{E_x^\bot\wedge E_x^\bot}\in\R(\g_2)$.
All this yields that $\R(\hol_x)=\R(\g_1)\oplus\R(\g_2)$.
Hence,
$\hol_x=L(\R(\hol_x))=L(\R(\g_1))\oplus
L(\R(\g_2))\subset\g_1\oplus\g_2$. Since
$\g_1\oplus\g_2\subset\hol_x$, we have $\hol_x=\g_1\oplus\g_2$,
$L(\R(\g_1))=\g_1$ and $L(\R(\g_2))=\g_2$. $\Box$

All the above shows that  to get a classification of connected holonomy groups
(equivalently of holonomy algebra) for pseudo-Riemannian manifolds of signature $(r,s)$ one must solve the following problems
\begin{description}\item[Problem 1)] classify weakly-irreducible
subalgebras $\g\subset\so(r,s)$;
\item[Problem 2)] check which Lie algebras of  Problem 1) are Berger algebras;
\item[Problem 3)] find a pseudo-Riemannian manifold with the holonomy algebra $\g$
for each weakly-irreducible Berger algebra $\g\subset\so(r,s)$.
\end{description}

In the next two sections we will recall the solution of these problems for Riemannian and Lorentzian manifolds.\newpage


\section{Connected holonomy groups of Riemannian manifolds and
connected irreducible holonomy groups of pseudo-Riemannian manifolds}\label{sec1.2}
\def\mypar{\thesection~Connected holonomy groups of Riemannian manifolds and \ldots}

Consider now Riemannian manifolds. Weakly-irreducible subalgebras $\g\subset\so(n)$
are irreducible and Problem 1) is equivalent to the problem of classification
of irreducible representations of compact Lie algebras, this problem was solved by E.~Cartan in  \cite{Car1,Car2}.

Recall that for any locally symmetric Riemannian manifold
there exists a simply connected Riemannian symmetric space with the same restricted holonomy group (see Section \ref{secsym} for definitions).
 Simply connected Riemannian symmetric spaces were classified
by E.~Cartan  (\cite{Car6,Be,Helgason}). If the holonomy group of such  space
is irreducible, then it coincides with the isotropy representation.
Thus connected irreducible holonomy groups of locally symmetric Riemannian spaces are known.

In  1955 M.\,Berger obtained a list of possible connected irreducible holonomy groups of Riemannian manifolds, \cite{Ber}.

\begin{theorem}\label{Ber}
Let $G\subset SO(n)$ be a connected irreducible Lie subgroup
such that its Lie algebra $\g\subset\so(n)$ satisfies  $L(\R(\g))=\g$,
then either $G$ is the holonomy group of a locally symmetric Riemannian manifold, or
$G$  is one of the following groups:\\ $SO(n);$\\ $U(m)$, $SU(m)$, $n=2m$;\\
$Sp(m)$, $Sp(m)\cdot Sp(1)$, $n=4m$;\\ $Spin(7)$, $n=8$;\\
$G_2$, $n=7$.
\end{theorem}

The initial  list of M.~Berger contained also the Lie group
$Spin(9)\subset SO(16)$. In \cite{Al}  D.~V.~Alekseevsky showed that
Riemannian manifolds with the holonomy group $Spin(9)$  are locally symmetric.
The list of Theorem \ref{Ber} coincides with the list of connected Lie groups
$G\subset SO(n)$ acting transitively on the sphere $S^{n-1}\subset \Real^n$
(if we exclude from the last list  the Lie groups  $Spin(9)$ and $Sp(m)\cdot T$, where $T$ is the circle).
Using this, in  1962 J.~Simons gave in \cite{Sim} a geometric proof of
Theorem \ref{Ber}.

To prove  Theorem  \ref{Ber},  M.~Berger
used the classification of irreducible representations of  compact Lie algebras.
Each such representation can be obtained using tensor
products of the fundamental representations.
Berger  showed that most of these representations can not appear as holonomy representations: if the representation contains more then one
tensorial factor, then $\R(\g)=\{0\}$.
Using complicated computations, Berger showed that
for the fundamental representations that are not in Theorem \ref{Ber}
from the Bianchi identity
it follows  $\nabla R=0$ or  $R=0$.

Thus   Problem 2) for Riemannian manifolds is solved.

Examples of Riemannian manifolds with the holonomy groups
$U(\frac{n}{2})$, $SU(\frac{n}{2})$, $Sp(\frac{n}{4})$ and $Sp(\frac{n}{4})\cdot Sp(1)$
were constructed by E.~Calabi, S.~T.~Yau and D.~V.~Alekseevsky.
In 1987 in \cite{Bryant} R.~Bryant constructed examples of Riemannian manifolds with the holonomy groups $Spin(7)$ and $G_2$.
Many constructions are given in the book of D.~Joyce \cite{Jo}.
This finishes the classification of connected holonomy groups of Riemannian manifolds.

Consider some special geometric structures on Riemannian manifolds with the
holonomy groups from Theorem  \ref{Ber}.
\begin{description}
\item[$SO(n)$:] This is the holonomy group of Riemannian manifolds of "general position".
There are no geometric structures induced by the holonomy group on such manifolds.

\item[$U(m)$] ($n=2m$): Riemannian manifolds with this holonomy group are
 K\"ahlerian, they admit parallel Hermitian structures (Example \ref{exIn3}).

\item[$SU(m)$] ($n=2m$): Riemannian manifolds with this holonomy group are called
special K\"ahlerian or Calabi-Yau manifolds, they are K\"ahlerian
and Ricci-flat (Example \ref{exIn4}).

\item[$Sp(m)$] ($n=4m$): Riemannian manifolds with this holonomy group are called
hyper-K\"ahlerian, they are Ricci-flat, each such manifold admits a parallel quaternionic structure,
i.e. parallel Hermitian structures $I$, $J$ and $K$ such that $IJ=-JI=K$.

\item[$Sp(m)\cdot Sp(1)$] ($n=4m$): Riemannian manifolds with this holonomy group are called quaternionic-K\"ahlerian, each such manifold admits a
parallel subbundle of the bundle of the endomorphisms of the tangent spaces
that  locally is generated by a quaternionic structure.

\item[$Spin(7)$] ($n=8$), $G_2$ ($n=7$): Riemannian manifolds with these holonomy groups
are Ricci-flat. On a manifold with the holonomy group $Spin(7)$
exists a parallel 4-form, on a manifold with the holonomy group $G_2$
exists a parallel 3-form.
\end{description}

The next theorem gives a classification of possible connected irreducible holonomy groups
of pseudo-Riemannian manifolds. The  first version of this theorem was obtained by M.~Berger in 1955 in \cite{Ber}.
Later, the initial list of M.~Berger  was refined and completed by R.~Bryant, Q.-S.~Chi, S.~Merkulov and L.~Schwachh\"ofer, see \cite{Br3,CMS,MS99}.

\begin{theorem}\label{Ber1}
Let $G\subset SO(r,s)$ is a connected irreducible Lie subgroup
such that its Lie algebra $\g\subset\so(r,s)$ satisfies  $L(\R(\g))=\g$,
then either $G$ is the holonomy group of a locally symmetric pseudo-Riemannian manifold, or $G$  is one of the following groups:
\\ $SO(r,s);$\\ $U(p,q)$, $SU(p,q)$, $r=2p$, $s=2q$;\\
$Sp(p,q)$, $Sp(p,q)\cdot Sp(1)$, $r=4p$, $s=4q$;\\
$SO(r,\mathbb{C})$, $s=r$;\\
$Sp(p)\cdot SL(2,\Real)$, $r=s=2p$;\\
$Sp(p,\mathbb{C})\cdot SL(2,\Co)$, $r=s=4p$;\\
$Spin(7)$, $r=0$, $s=8$;\\
$Spin(4,3)$, $r=s=4$;\\
$Spin(7)^{\mathbb{C}}$, $r=s=8$;\\
$G_2$, $r=0$, $s=7$;\\
$G_{2(2)}^*$, $r=4$, $s=3$;\\
$G_2^\mathbb{C}$, $r=s=7$.
\end{theorem}

In 1957 in \cite{Ber57}
M.~Berger got a classification of simply connected pseudo-Riemannian symmetric
spaces with irreducible holonomy groups, thus the holonomy groups of these spaces are known.


Thus in the pseudo-Riemannian case we are left with the problem of classification
of connected weakly-irreducible not irreducible holonomy groups.

\section{Connected holonomy groups of Lorentzian manifolds}\label{sec1.3}

In this section we explain the classification of the holonomy algebras of Lorentzian manifolds.
We assume that the dimensions of the Lorentzian manifolds are $n+2$, where $n\geq 0$.
From Berger's classification of irreducible holonomy algebras of pseudo-Riemannian manifolds  it follows that {\it the
only irreducible holonomy algebra of Lorentzian
manifolds is  $\so(1,n+1)$}, see \cite{Disc-Ol} and \cite{Bo-Ze} for direct proofs of this fact.

Consider weakly-irreducible not irreducible holonomy algebras. First set some notation.

Let $(\Real^{1,n+1},\eta)$ be the Minkowski space of dimension $n+2$,
where  $\eta$ is a metric on $\Real^{n+2}$ of
signature $(1,n+1)$. We fix a basis
$p,e_1,...,e_n,q$ of $\Real^{1,n+1}$ such that the Gram matrix
of $\eta$ has the form $\left(\scriptsize{\begin{array}{ccc} 0 & 0 & 1\\ 0 & E_n & 0 \\ 1 & 0 & 0 \\ \end{array}}\right)$,
where $E_n$ is the $n$-dimensional identity matrix. We will denote by $E=\Real^n\subset\Real^{1,n+1}$ the Euclidean subspace spanned by the vectors $e_1,...,e_n$.

Denote  by $\so(1,n+1)_{\Real p}$ the subalgebra of $\so(1,n+1)$ that preserves the isotropic line $\Real p$. The Lie
algebra $\so(1,n+1)_{\Real p}$ can be identified with the following matrix algebra $$\so(1,n+1)_{\Real p}=\left\{\left.
\left (\scriptsize{\begin{array}{ccc} a &X^t & 0\\ 0 & A &-X \\ 0 & 0 & -a \\
\end{array}}\right)\right|\, a\in \Real,\, X\in \Real^n,\,A \in \so(n) \right\} .$$
We identify the above matrix  with the triple $(a,A,X)$.
Define the following subalgebras of $\so(1,n+1)_{\Real p}$:
$$\A=\{(a,0,0)|a\in \Real\},\quad \K=\{(0,A,0)|A\in \mathfrak{so}(n)\},\quad \mathcal N=\{(0,0,X)|X\in \Real^n\}.$$
We see that $\A$ commutes with $\K$, and $\N$ is a commutative ideal. We also see that $$[(a,A,0),(0,0,X)]=(0,0,aX+AX).$$
We obtain the decomposition $$\so(1,n+1)_{\Real p}=(\mathcal A\oplus\mathcal K)\zr\mathcal N=(\Real\oplus\so(n))\zr\Real^n.$$

If a weakly-irreducible subalgebra $\g\subset \so(1,n+1)$ preserves a degenerate proper subspace $U\subset \Real^{1,n+1}$, then it preserves the
isotropic line $U\cap U^\bot$, and $\g$ is conjugated to a weakly-irreducible subalgebra of $\so(1,n+1)_{\Real p}$.

Let $\h\subset\so(n)$ be a subalgebra. Recall that $\h$ is a compact Lie
algebra and we have the decomposition $\h=\h'\oplus\z(\h)$, where $\h'$ is the commutant of $\h$ and $\z(\h)$ is the center of $\h$,
see for example \cite{V-O}.

The first step towards the classification of weakly-irreducible not irreducible holonomy algebras of Lorentzian manifolds was done by L.~Berard~Bergery and A.~Ikemakhen who
proved in 1993 in \cite{B-I} the following theorem.

\begin{theorem}\label{B-I}
A subalgebra $\g\subset \so(1,n+1)_{\Real p}$ is  weakly-irreducible if and only if
$\g$ belongs to one of the following types
\begin{description}
\item[type 1.] $\g^{1,\h}=(\Real\oplus\h)\zr\Real^n=\left\{\left. \left (\scriptsize{\begin{array}{ccc}
a &X^t & 0\\ 0 & A &-X \\ 0 & 0 & -a \\
\end{array}}\right)\right|\, a\in \Real,\,\,A \in \h,\, X\in \Real^n
 \right\}$, where $\h\subset\so(n)$ is a subalgebra;

\item[type 2.] $\g^{2,\h}=\h\zr\Real^n=\left\{\left. \left (\scriptsize{\begin{array}{ccc}
0 &X^t & 0\\ 0 & A &-X \\ 0 & 0 & 0 \\
\end{array}}\right)\right|\,\,A \in \h,\,  X\in \Real^n  \right\}$;

\item[type 3.] $\g^{3,\h,\varphi}=\{(\varphi(A),A,0)|A\in\h\}\zr\Real^n=\left\{\left. \left (\scriptsize{\begin{array}{ccc}
\varphi(A) &X^t & 0\\ 0 & A &-X \\ 0 & 0 & -\varphi(A) \\
\end{array}}\right)\right|\,\,A \in \h,\,  X\in \Real^n \right\}$,
where $\h\subset\so(n)$ is a subalgebra with $\z(\h)\neq\{0\}$, and  $\varphi :\h\to\Real$ is a non-zero linear map with $\varphi|_{\h'}=0$;

\item[type 4.] $\g^{4,\h,m,\psi}=\{(0,A,X+\psi(A))|A\in\h,X\in\Real^m\}\\=\left\{\left. \left (\scriptsize{\begin{array}{cccc}
0 &X^t&\psi(A)^t & 0\\ 0 & A&0 &-X \\ 0 & 0 & 0 &-\psi(A) \\ 0&0&0&0\\ \end{array}}\right)\right|\,\,A \in \h,\,  X\in
\Real^{m} \right\}$, where  $0<m<n$ is an integer, $\h\subset\so(m)$ is a subalgebra with $\dim\z(\h)\geq n-m$, and
$\psi:\h\to \Real^{n-m}$ is a surjective linear map with $\psi|_{\h'}=0$.
\end{description}\end{theorem}

{\bf Definition.} {\it
The subalgebra $\h\subset\so(n)$ associated to a
weakly-irreducible subalgebra $\g\subset \so(1,n+1)_{\Real p}$ in Theorem \ref{B-I}  is called  the orthogonal part of $\g$.}

\vskip0.2cm

The proof of Theorem \ref{B-I} given by L.~Berard~Bergery and A.~Ikemakhen  was purely algebraic.
We describe now another more geometric proof   of this theorem from \cite{Gal2},
since in Chapter II we will generalize this idea in order to classify weakly-irreducible not irreducible subalgebras of $\su(1,n+1)\subset\so(2,2n+2)$.

Let $(E,\eta)$ be an Euclidean space. A map $f:E\to E$
is called a {\it similarity transformation} of $E$ if there exists
a $\lambda >0$ such that $\|f(x_1)-f(x_2)\|=\lambda\|x_1-x_2\|$ for all
$x_1,x_2\in E$, where $\|x\|^2=\eta(x,x)$. If $\lambda=1$, then $f$
is called an {\it isometry}. Denote by $\Simil E$ and $\Isom E$ the groups of all similarity transformations and isometries of $E$, respectively.
A subgroup $G\subset \Simil E$ is called {\it irreducible}
if it preserves no proper affine subspace of $E$.

Let $p,e_1,...,e_n,q$ be a basis of the vector space $\Real^{1,n+1}$ as
above.
Consider the basis $e_0,e_1,...,e_n,e_{n+1}$ of $\Real^{1,n+1}$,
where  $e_0=\frac{\sqrt{2}}{2}(p-q)$ and $e_{n+1}=\frac{\sqrt{2}}{2}(p+q)$.
With respect to this basis the Gram matrix of $\eta$ has the form
 $\left(\scriptsize{
\begin{array}{cc}
-1 & 0 \\
0 & E_{n+1}\\
 \end{array}} \right),$ where $E_{n+1}$ is the $n+1$-dimensional
 identity matrix.

The vector model of the $n+1$-dimensional {\it hyperbolic space} is
defined  in the following way:
$$H^{n+1}=\{x\in \Real^{1,n+1}|\ \eta(x,x)=-1,\ x_0>0\},$$
where $x_0$ is the first coordinate with respect to the basis $e_0,...,e_{n+1}$.
Recall that $H^{n+1}$ is an $n+1$-dimensional Riemannian submanifold of $\Real^{1,n+1}$.
The tangent space at a point $x\in H^{n+1}$ is identified with the
vector subspace $(x)^{\bot_\eta}\subset \Real^{1,n+1}$ and the restriction of the form
$\eta$ to this subspace is positively definite.

Let $O'(1,n+1)\subset O(1,n+1)$ be the subgroup preserving the space $H^{n+1}$. Then for any $f\in O'(1,n+1)$, the
restriction $f|_{H^{n+1}}$ is an isometry of $H^{n+1}$ and any isometry of $H^{n+1}$ can be obtained in this way. Hence we
have the isomorphism $$O'(1,n+1)\simeq\Isom H^{n+1},$$ where $\Isom H^{n+1}$ is the group of all isometries of $H^{n+1}$.

Consider the {\it light-cone} of $\Real^{1,n+1}$,
$$C=\{x\in \Real^{1,n+1}|\ \eta(x,x)=0\}.$$
The subset of the $n+1$-dimensional projective space $\mathbb{P}\Real^{1,n+1}$ that
consists of all {\it isotropic lines} $l\subset C$ is called the
{\it boundary of the hyperbolic space} $H^{n+1}$ and is
denoted by $\partial H^{n+1}$.

We identify $\partial H^{n+1}$ with the $n$-dimensional unit sphere $S^n$ in
the following way.
Consider the vector subspace $E_1=E\oplus\Real e_{n+1}$.
Each isotropic line intersects the hyperplane $e_0+E_1$ at a unique
point. The intersection $(e_0+E_1)\cap C$ is the set $$\{x\in \Real^{1,n+1} | x_0=1,\,
x_1^2+\cdots+ x_{n+1}^2=1\},$$ which is the $n$-dimensional sphere $S^n$.
This gives us the identification $\partial H^{n+1}\simeq S^n$.

Denote by $\Conf S^n$ the group of all conformal transformations
of $S^n$. Any transformation $f\in SO(1,n+1)$ takes isotropic lines to
isotropic lines. Moreover, under the above identification, we have
$f|_{\partial H^{n+1}}\in \Conf \partial H^{n+1}$ and any transformation from
$\Conf \partial H^{n+1}$ can be obtained in this way. Hence we have the
isomorphism $$SO(1,n+1)\simeq\Conf\partial H^{n+1}.$$

Denote by $SO(1,n+1)_{\Real p}$ the subgroup of $SO(1,n+1)$ that preserves the isotropic line $\Real p$.  Suppose that
$f\in SO(1,n+1)_{\Real p}$. The corresponding element $f\in\Conf S^n$ preserves the point $p_0=\Real p\cap (e_0+ E_1)$.
Clearly, $p_0=\sqrt 2 p$. Denote by $s_0$ the stereographic projection $s_0: S^n\backslash\{p_0\}\to e_0+E$. Since
$f\in\Conf S^n$, we see that the element $\Gamma(f)$ defined by $\Gamma(f)=s_0\circ f \circ s^{-1}_0:E\to E$ (here we
identify $e_0+E$ with $E$) is a similarity transformation of the Euclidean space $E$. Conversely, any similarity
transformation of $E$ can be obtained in this way. Thus we have the isomorphism $$\Gamma:SO(1,n+1)_{\Real p}\to \Simil
E.$$

Then  we prove that {\it a connected  subgroup $G\subset SO(V)_{\Real p}$
acts weakly-irreducibly on $\Real^{1,n+1}$ if and
only if $\Gamma(G)\subset\Simil E$ acts transitively on the Euclidean space
$E$.} From this and Theorem \ref{Al} below we obtain

\begin{theorem} \label{wir} Let $G$ be a connected  subgroup
of $SO(1,n+1)_{\Real p}$. Then $G$ acts weakly-irreducibly on $\Real^{1,n+1}$ if and
only if $\Gamma(G)$ acts transitively on the Euclidean space  $E=\partial H^{n+1}\backslash\{\Real p\}$.\end{theorem}

Using results of \cite{Al2} and \cite{A-V-S}, we divide transitively acting connected subgroups of $\Simil E$ into 4 types
(see Theorem \ref{trans} below) and show that the Lie algebras of the corresponding subgroups of $SO(1,n+1)_{\Real p}$
have the same types introduced by L.~Berard~Bergery and A.~Ikemakhen. This finishes the proof of Theorem \ref{B-I}. $\Box$

In order to get a classification of weakly-irreducible Berger subalgebras of
$\so(1,n+1)_{\Real p}$, in \cite{Gal1} we study the spaces $\R(\g)$ for the Lie algebras of Theorem \ref{B-I} in terms of their orthogonal parts $\h\subset\so(n)$.

\vskip0.2cm

{\bf Definition.} {\it
The vector  space
$$\P(\h)=\left\{P\in \Hom (\Real^n,\h)\,\left|\,\begin{array}{c}
\eta(P(u)v,w)+\eta(P(v)w,u)+\eta(P(w)u,v)=0\\ \text{ for all }u,v,w\in \Real^n\end{array}\right\}\right.$$
is called the space of  weak-curvature tensors of type $\h$.
A subalgebra $\h\subset\so(n)$ is called a weak-Berger algebra if $L(\P(\h))=\h$, where
$$L(\P(\h))=\spa\{P(u)|P\in\P(\h),\,u\in \Real^n\}$$ is the vector subspace of $\h$ spanned by
$P(u)$ for all $P\in\P(\h)$ and $u\in \Real^n$. }

\vskip0.2cm

In \cite{Gal1} we proved the following theorem.

\begin{theorem} A weakly-irreducible subalgebra $\g\subset\so(1,n+1)_{\Real p}$ is a Berger algebra if and only if
its orthogonal part $\h\subset\so(n)$ is a weak-Berger algebra.\end{theorem}

Thus we need a classification of weak-Berger algebras.
In \cite{Le1,Le2,Le3,Le4} T.~Leistner proved the following theorem.

\begin{theorem} \label{Le}
A subalgebra $\h\subset\so(n)$ is a weak-Berger algebra if and only if $\h$ is a Berger algebra.\end{theorem}

First note that any weak-Berger algebra is a direct sum of irreducible weak-Berger
algebras. To prove Theorem \ref{Le}, T.~Leistner used the classification of irreducible
representations of compact Lie algebras and complicated computations.
Recall that from the classification of Riemannian holonomy algebras it follows that {\it a subalgebra $\h\subset\so(n)$ is a Berger algebra if and only if
$\h$ is the holonomy algebra of a Riemannian manifold.}

Thus {\it a subalgebra $\g\subset\so(1,n+1)$ is a  weakly-irreducible not irreducible Berger algebra if and only if $\g$ is conjugated to one
of the subalgebras $\g^{1,\h},\g^{2,\h},\g^{3,\h,\varphi},\g^{4,\h,m,\psi}\subset\so(1,n+1)_{\Real p}$, where $\h\subset\so(n)$ is the
holonomy algebra of a Riemannian manifold.}

The last step to complete the classification of Lorentzian holonomy algebras
is to realize all weakly-irreducible Berger subalgebra of $\so(1,n+1)_{\Real p}$
as the holonomy algebras of Lorentzian manifolds.  This was completely done in \cite{Gal5}.
First consider two examples.

\vskip0.3cm

{\bf Example \cite{B-I}.} In 1993 L. Berard Bergery and A. Ikemakhen realized
the weakly-irreducible Berger subalgebra of $\so(1,n+1)_{\Real p}$ of type 1 and 2
as the holonomy algebras of Lorentzian manifolds. They constructed the following metrics.
Let $\h\subset\so(n)$ be the holonomy algebra of a Riemannian manifold.
Let $x^{0},x^{1},...,x^{n},x^{n+1}$ be the standard coordinates on $\Real^{n+2}$,
$h$ be a metric on $\Real^n$ with the holonomy algebra $\h$, and $f(x^{0},...,x^{n+1})$ be  a function with
$\frac{\partial f}{\partial x^1}\neq 0,$...$,\frac{\partial f}{\partial x^n}\neq 0$.
If $\frac{\partial f}{\partial x^0}\neq 0$, then the holonomy algebra of the metric
$$g=2dx^0dx^{n+1}+h+f\cdot(dx^{n+1})^2$$ is $\g^{1,\h}$. If  $\frac{\partial f}{\partial x^0}=0$, then the holonomy algebra of the metric $g$ is $\g^{2,\h}$.

\vskip0.3cm

The idea of the general constructions came from the following example of A. Ikemakhen.

\vskip0.3cm

{\bf Example \cite{Ik}.}
Let $x^{0},x^{1},...,x^{5},x^{6}$ be the standard coordinates on $\Real^{7}$. Consider the following metric
$$g=2dx^0dx^{6}+\sum^{5}_{i=1}(dx^i)^2+2\sum^{5}_{i=1}u^i dx^{i} dx^{6},$$
where $$\begin{array}{ll}u^1=-(x^3)^2-4(x^4)^2-(x^5)^2,& u^2=u^4=0,\\ u^3=-2\sqrt{3}x^2x^3-2x^4x^5,& u^5=2\sqrt{3}x^2x^5+2x^3x^4. \end{array}$$

The holonomy algebra  of this metric at the point $0$ is $\g^{2,\rho(\so(3))}\subset\so(1,6)$, where
$\rho:\so(3)\to\so(5)$ is the representation given by the highest irreducible
component of the representation $\otimes^2\id:\so(3)\to\otimes^2\so(3)$.
The image $\rho(\so(3))\subset\so(5)$ is spanned by the matrices
$$A_1=\left(\begin{smallmatrix}0&0&-1&0&0\\0&0&\sqrt{3}&0&0\\1&-\sqrt{3}&0&0&0\\0&0&0&0&-1\\0&0&0&1&0\end{smallmatrix}\right),\quad
A_2=\left(\begin{smallmatrix}0&0&0&-4&0\\0&0&0&0&0\\0&0&0&0&-2\\4&0&0&0&0\\0&0&2&0&0\end{smallmatrix}\right),\quad
A_3=\left(\begin{smallmatrix}0&0&0&0&-1\\0&0&0&0&-\sqrt{3}\\0&0&0&-1&-1\\0&0&1&0&0\\1&\sqrt{3}&1&0&0\end{smallmatrix}\right).$$

We have $\pr_{\so(n)}\left(R\left(\frac{\partial}{\partial x^{3}},\frac{\partial}{\partial x^{6}}\right)_0\right)=A_1$,
$\pr_{\so(n)}\left(R\left(\frac{\partial}{\partial x^{4}},\frac{\partial}{\partial x^{6}}\right)_0\right)=A_2$,
$\pr_{\so(n)}\left(R\left(\frac{\partial}{\partial x^{5}},\frac{\partial}{\partial x^{6}}\right)_0\right)=A_3$
and $\pr_{\so(n)}\left(R\left(\frac{\partial}{\partial x^{1}},\frac{\partial}{\partial x^{6}}\right)_0\right)=
\pr_{\so(n)}\left(R\left(\frac{\partial}{\partial x^{2}},\frac{\partial}{\partial x^{6}}\right)_0\right)=0$.

Note the following.
Let $P\in\Hom(\Real^n,\h)$ be a linear map defined by: $P(e_1)=P(e_2)=0$, $P(e_3)=A_1$, $P(e_4)=A_2$ and  $P(e_5)=A_3$.
Then $P\in\P(\h)$, $P(\Real^n)=\h$ and
$\pr_{\so(n)}\left(R\left(\frac{\partial}{\partial x^{i}},\frac{\partial}{\partial x^{6}}\right)_0\right)=P(e_i)$ for all $1\leq i\leq 5$.

\vskip0.3cm

Now we explain the general construction.
Let $\h\subset\so(n)$ be the holonomy algebra of a Riemannian manifold.
Theorem \ref{Wu} states that we have an orthogonal
decomposition \begin{equation}\label{LM0A}\Real^{n}=\Real^{n_1}\oplus\cdots\oplus\Real^{n_s}\oplus\Real^{n_{s+1}}\end{equation} and the corresponding
decomposition of $\h$ into the direct sum of ideals
\begin{equation}\label{LM0B}\h=\h_1\oplus\cdots\oplus\h_s\oplus\{0\}\end{equation} such that $\h$ annihilates $\Real^{n_{s+1}}$,
$\h_i(\Real^{n_j})=0$ for $i\neq j$, and $\h_i\subset\so(n_i)$ is an
irreducible subalgebra for $1\leq i\leq s$. Moreover, the Lie algebras $\h_i$ are the holonomy algebras of Riemannian manifolds.
In \cite{Gal1} it was proved that \begin{equation}\label{LM0C}\P(\h)=\P(\h_1)\oplus\cdots\oplus\P(\h_s).\end{equation}
We choose the basis $e_1,...,e_n$ of $\Real^n$  compatible  with the above decomposition of $\Real^n$.

Let $n_0=n_1+\cdots+n_s=n-n_{s+1}$. Obviously, $\h\subset\so(n_0)$ and $\h$ does not annihilate any proper subspace of $\Real^{n_0}$.
Note that in the case of the Lie algebra $\g^{4,\h,m,\psi}$ we have $0<n_0\leq m$.

We will always assume that the indices $b,c,d,f$ run from $0$ to $n+1$,
the indices $i,j,k,l$ run from $1$ to $n$, the indices $\i,\j,\k,\l$ run from $1$ to $n_0$,
the indices $\ii,\jj,\kk,\ll$ run from $n_0+1$ to $n$,  and the indices $\alpha,\beta,\gamma$
run from $1$ to $N$. In case of the Lie algebra  $\g^{4,\h,m,\psi}$ we will also  assume  that
the indices $\iv,\jv,\kv,\lv$ run from $n_0+1$ to $m$ and the indices $\iiv,\jjv,\kkv,\llv$ run from $m+1$ to $n$.
We will use the Einstein rule for sums.

Let $(P_\alpha)_{\alpha=1}^N$ be linearly independent elements of $\P(\h)$ such that the subset
$\{P_\alpha(u)|1\leq\alpha\leq N,\,u\in \Real^{n_0}\}\subset\h$ generates the Lie algebra $\h$. For example, it can be any basis
of the vector space $\P(\h)$.
We have $P_\alpha|_{\Real^{n_{s+1}}}=0$ and $P_\alpha$ can be considered as linear maps $P_\alpha:\Real^{n_0}\to\h\subset\so(n_0)$.
For each $P_\alpha$ define the numbers $P_{\alpha\j\i}^\k$ such that $P_\alpha(e_\i)e_\j=P_{\alpha\j\i}^\k e_\k.$
Since $P_\alpha\in\P(\h)$, we have
\begin{equation}\label{LM1} P_{\alpha\k\i}^\j=-P_{\alpha\j\i}^\k\quad \text{and}\quad
 P_{\alpha\j\i}^\k+P_{\alpha\k\j}^\i+P_{\alpha\i\k}^\j=0.\end{equation}

Define the following numbers
\begin{equation}\label{LM2A} a_{\alpha\j\i}^\k=\frac{1}{3\cdot(\alpha-1)!}\left(P_{\alpha\j\i}^\k+P_{\alpha\i\j}^\k\right).\end{equation}
Then \begin{equation}\label{LM3} a_{\alpha\j\i}^\k=a_{\alpha\i\j}^\k.\end{equation}
From \eqref{LM1} it follows that
\begin{equation}\label{LM4} P_{\alpha\j\i}^\k=(\alpha-1)!\left(a_{\alpha\j\i}^\k-a_{\alpha\k\i}^\j\right)\quad \text{and} \quad
 a_{\alpha\j\i}^\k+a_{\alpha\k\j}^\i+a_{\alpha\i\k}^\j=0.\end{equation}

Let $x^0,...,x^{n+1}$ be the standard coordinates on $\Real^{n+2}$. Consider the metric
\begin{equation}\label{LM6} g=2dx^0dx^{n+1}+\sum^{n}_{i=1}(dx^i)^2+2\sum^{n_0}_{\i=1}u^\i dx^\i dx^{n+1}+f\cdot(dx^{n+1})^2,\end{equation}
where \begin{equation}\label{LM6A}u^\i=a_{\alpha\j\k}^\i x^\j x^\k(x^{n+1})^{\alpha-1}\end{equation}
and $f$ is a function which depends on the type of the holonomy algebra that we wish to obtain.

For the Lie algebra $\g^{3,\h,\varphi}$ (if it exists)  define the numbers
\begin{equation}\label{LM7}\varphi_{\alpha\i}=\frac{1}{(\alpha-1)!}\varphi(P_{\alpha}(e_\i)).\end{equation}

For the Lie algebra $\g^{4,\h,m,\psi}$ (if it exists)  define the numbers $\psi_{\alpha\i\iiv}$ such that

\begin{equation}\label{LM8}\frac{1}{(\alpha-1)!}\psi(P_{\alpha}(e_\i))=\sum^{n}_{\iiv=m+1}\psi_{\alpha\i\iiv}e_\iiv.\end{equation}

If we choose $f$ such that $f(0)=0$, then $g_0=\eta$ and we can identify the tangent space to $\Real^{n+2}$ at $0$ with
the vector space $\Real^{1,n+1}$.

\begin{theorem}\label{LMth1}
The holonomy algebra $\hol_0$ of the metric $g$ at the point $0\in \Real^{n+2}$
depends on the function $f$  in the following way
$$
\begin{array}{|cc|l|}\hline f& &\hol_0 \\\hline
(x^0)^2+\sum_{\ii=n_0+1}^{n}(x^\ii)^2& &\g^{1,\h}\\
\sum_{\ii=n_0+1}^{n}(x^\ii)^2& &\g^{2,\h}\\
2x^0\varphi_{\alpha\i}x^\i(x^{n+1})^{\alpha-1}+\sum_{\ii=n_0+1}^{n}(x^\ii)^2& &\g^{3,\h,\varphi} (\text{if } \z(\h)\neq\{0\}) \\
2\psi_{\alpha\i\iiv}x^\i x^\iiv(x^{n+1})^{\alpha-1}+\sum_{\iv=n_0+1}^{m}(x^\iv)^2& & \g^{4,\h,m,\psi}(\text{if } \dim\z(\h)\geq n-m)\\\hline\end{array}$$
\end{theorem}

In order to prove  Theorem \ref{LMth1}, we computed in \cite{Gal5} the components of all covariant derivatives of the
curvature tensors of the metrics \eqref{LM6} with the above chosen $f$, then we used Theorem \ref{analythol}.

As the corollary we obtain the classification of the weakly-irreducible not irreducible Lorentzian holonomy algebras.

\begin{theorem}\label{LMth2}
A weakly-irreducible not irreducible subalgebra  $\g\subset\so(1,n+1)$ is the holonomy algebra of a Lorentzian manifold if and only if $\g$ is conjugated to one
of the subalgebras $\g^{1,\h},\g^{2,\h},\g^{3,\h,\varphi},\g^{4,\h,m,\psi}\subset\so(1,n+1)_{\Real p}$, where $\h\subset\so(n)$ is the
holonomy algebra of a Riemannian manifold.
\end{theorem}

From Theorem \ref{LMth2}, Wu's theorem and Berger's list it follows that
{\it the holonomy algebra $\hol\subset\so(1,N+1)$ of any Lorentzian manifold of dimension $N+2$
has the form $\hol=\g\oplus\h_1\oplus\cdots\oplus\h_r$,
where either $\g=\so(1,n+1)$ or $\g$ is a Lie algebra from Theorem \ref{LMth2}, and $\h_i\subset\so(n_i)$ are the irreducible holonomy
algebras of Riemannian manifolds ($N=n+n_1+\cdots+n_r$).}

\vskip 0.3cm

Now we compare our method of constructions of the metric \eqref{LM6} with the example of A. Ikemakhen.
Let us construct the metric for $\g^{2,\rho(\so(3))}\subset\so(1,6)$ by
our method. Take $P\in\P(\rho(\so(3)))$ defined as $P(e_1)=P(e_2)=0$, $P(e_3)=A_1$, $P(e_4)=A_2$ and $P(e_5)=A_3$.
By our constructions,
$$g=2dx^0dx^{6}+\sum^{5}_{i=1}(dx^i)^2+2\sum^{5}_{i=1}u^i dx^{i} dx^{6},$$
where $$\begin{array}{ll}u^1=-\frac{2}{3}((x^3)^2+4(x^4)^2+(x^5)^2),& u^2=\frac{2\sqrt{3}}{3}((x^3)^2-(x^5)^2),\\
u^3=\frac{2}{3}(x^1x^3-\sqrt{3}x^2x^3-3x^4x^5-(x^5)^2),& u^4=\frac{8}{3}x^1x^4,\\
u^5=\frac{2}{3}(x^1x^5+\sqrt{3}x^2x^5+3x^3x^4+x^3x^5).& \end{array}$$
We still have $\pr_{\so(n)}\left(R\left(\frac{\partial}{\partial x^{3}},\frac{\partial}{\partial x^{6}}\right)_0\right)=A_1$,
$\pr_{\so(n)}\left(R\left(\frac{\partial}{\partial x^{4}},\frac{\partial}{\partial x^{6}}\right)_0\right)=A_2$, \\
$\pr_{\so(n)}\left(R\left(\frac{\partial}{\partial x^{5}},\frac{\partial}{\partial x^{6}}\right)_0\right)=A_3$
and $\pr_{\so(n)}\left(R\left(\frac{\partial}{\partial x^{1}},\frac{\partial}{\partial x^{6}}\right)_0\right)=
\pr_{\so(n)}\left(R\left(\frac{\partial}{\partial x^{2}},\frac{\partial}{\partial x^{6}}\right)_0\right)=0$.

The reason why we obtain another metric is the following. The idea of our constructions is to find the constants $a_{\alpha\j\i}^\k$ such that
$$\begin{array}{l}
\pr_{\so(n)}\left(R\left(\frac{\partial}{\partial x^{\i}},\frac{\partial}{\partial x^{n+1}}\right)_0\right)=P_1(e_\i),\\
...,\\
\pr_{\so(n)}\left(\nabla^{N-1}  R\left(\frac{\partial}{\partial x^{\i}},\frac{\partial}{\partial x^{n+1}};
\frac{\partial}{\partial x^{n+1}};\cdots ;\frac{\partial}{\partial x^{n+1}}\right)_0\right)=P_N(e_\i).
\end{array}$$
These conditions give us the system of equations
$$\left\{\begin{array}{rl}(\alpha-1)!\left(a_{\alpha\j\i}^\k-a_{\alpha\k\i}^\j\right)=&P_{\alpha\j\i}^\k,\\
a_{\alpha\j\i}^\k-a_{\alpha\i\j}^\k=&0.\end{array}\right.$$
One of the solutions of this system is given by \eqref{LM2A}, but this system can have other solutions.
We use the solution given by \eqref{LM2A}, taking another solution of the above system, we can obtain the metric constructed by
A. Ikemakhen.

Thus the choice of the functions $u^\i$ given by \eqref{LM6A} guarantees us that the orthogonal part of the holonomy algebra
$\hol_0$ coincides with the given Riemannian  holonomy algebra $\h\subset\so(n)$ (the other values of
$\pr_{\so(n)}(\nabla^{r}  R)$ does not give us anything new). This also guarantees us the inclusion $\Real^{n_0}\subset\hol_0$.
The reason why we choose the function $f$ as in Theorem \ref{LMth1} can be easily understood from
the following formulas

$$\pr_{\Real}\left(R\left(\frac{\partial}{\partial x^{0}},\frac{\partial}{\partial x^{n+1}}\right)_0\right)=
\frac{1}{2}\frac{\partial^2 f}{(\partial x^{0})^2} \text{ (we use this for } \g^{1,\h}),$$
$$\pr_{\Real}\left(\nabla^{\alpha-1}  R\left(\frac{\partial}{\partial x^{\i}},\frac{\partial}{\partial x^{n+1}};
\frac{\partial}{\partial x^{n+1}};\cdots ;\frac{\partial}{\partial x^{n+1}}\right)_0\right)=
\frac{1}{2}\frac{\partial^{\alpha+1} f}{\partial x^{0}\partial x^{\i}(\partial x^{n+1})^{\alpha-1}} $$
$$\text{ (we use this for } \g^{3,\h,\varphi}),$$
$$\pr_{\Real^n}\left(\nabla^{\alpha-1}  R\left(\frac{\partial}{\partial x^{\ii}},\frac{\partial}{\partial x^{n+1}};
\frac{\partial}{\partial x^{n+1}};\cdots ;\frac{\partial}{\partial x^{n+1}}\right)_0\right)=
\frac{1}{2}\sum_{\jj=n_0+1}^{n}
\frac{\partial^{\alpha+1} f}{\partial x^{\ii}\partial x^{\jj}(\partial x^{n+1})^{\alpha-1}}e_\jj
$$
$$\text{ (for $\alpha=0$ we use this for all algebras, for $\alpha\geq 0$ we use this for }
\g^{4,\h,m,\psi}).$$

\vskip0.3cm

As application let us construct  metrics for the Lie algebras $\g^{2,\g_2}\subset\so(1,8)$ and $\g^{2,\spin(7)}\subset\so(1,9)$.

\vskip0.3cm

{\bf Example (Metric with the holonomy algebra $\g^{2,\g_2}\subset\so(1,8)$).}
Consider the Lie subalgebra $\g_2\subset\so(7)$. The vector subspace $\g_2\subset\so(7)$ is spanned by the following matrices
(\cite{HelgaInes}):
$$\begin{array}{llll}
A_1=E_{12}-E_{34},&A_2=E_{12}-E_{56},&A_3=E_{13}+E_{24},&A_4=E_{13}-E_{67},\\
A_5=E_{14}-E_{23},&A_6=E_{14}-E_{57},&A_7=E_{15}+E_{26},&A_8=E_{15}+E_{47},\\
A_9=E_{16}-E_{25},&A_{10}=E_{16}+E_{37},&A_{11}=E_{17}-E_{36},&A_{12}=E_{17}-E_{45},\\
A_{13}=E_{27}-E_{35},&A_{14}=E_{27}+E_{46},& &
\end{array}$$ where $E_{ij}\in\so(7)$ ($i<j$) is the
skew-symmetric matrix such that $(E_{ij})_{ij}=1$, $(E_{ij})_{ji}=-1$ and $(E_{ij})_{kl}=0$ for other $k$ and $l$.

Consider the linear map $P\in\Hom(\Real^7,\g_2)$ defined as
$$\begin{array}{llll}P(e_1)=A_6,&P(e_2)=A_4+A_5,&P(e_3)=A_1+A_7,&P(e_4)=A_1,\\P(e_5)=A_4,&P(e_6)=-A_5+A_6,&P(e_7)=A_7.&\end{array}$$
It can be checked that $P\in\P(\g_2)$. Moreover,  the elements $A_1,A_4,A_5,A_6,A_7\in \g_2$ generate the Lie algebra
$\g_2$.

The holonomy algebra of the metric
$$g=2dx^0dx^{8}+\sum^{7}_{i=1}(dx^i)^2+2\sum^{7}_{i=1}u^i dx^{i} dx^{8},$$
where
$$\begin{array}{l}
u^1=\frac{2}{3}(2x^2x^3+x^1x^4+2x^2x^4+2x^3x^5+x^5x^7),\\
u^2=\frac{2}{3}(-x^1x^3-x^2x^3-x^1x^4+2x^3x^6+x^6x^7),\\
u^3=\frac{2}{3}(-x^1x^2+(x^2)^2-x^3x^4-(x^4)^2-x^1x^5-x^2x^6),\\
u^4=\frac{2}{3}(-(x^1)^2-x^1x^2+(x^3)^2+x^3x^4),\\
u^5=\frac{2}{3}(-x^1x^3-2x^1x^7-x^6x^7),\\
u^6=\frac{2}{3}(-x^2x^3-2x^2x^7-x^5x^7),\\
u^7=\frac{2}{3}(x^1x^5+x^2x^6+2x^5x^6),\end{array}$$
at the point $0\in\Real^9$ is $\g^{2,\g_2}\subset\so(1,8)$.

\vskip0.3cm

{\bf Example (Metric with the holonomy algebra $\g^{2,\spin(7)}\subset\so(1,9)$).}
Consider the Lie subalgebra $\spin(7)\subset\so(8)$.
The vector subspace $\spin(7)\subset\so(8)$ is spanned by the following matrices
(\cite{HelgaInes}):
$$\begin{array}{llll}
A_{1}=E_{12}+E_{34},&A_{2}=E_{13}-E_{24},&A_{3}=E_{14}+E_{23},&A_{4}=E_{56}+E_{78},\\
A_{5}=-E_{57}+E_{68},&A_{6}=E_{58}+E_{67},&A_{7}=-E_{15}+E_{26},&A_{8}=E_{12}+E_{56},\\
A_{9}=E_{16}+E_{25},&A_{10}=E_{37}-E_{48},&A_{11}=E_{38}+E_{47},&A_{12}=E_{17}+E_{28},\\
A_{13}=E_{18}-E_{27},&A_{14}=E_{35}+E_{46},&A_{15}=E_{36}-E_{45},&A_{16}=E_{18}+E_{36},\\
A_{17}=E_{17}+E_{35},&A_{18}=E_{26}-E_{48},&A_{19}=E_{25}+E_{38},&A_{20}=E_{23}+E_{67},\\
A_{21}=E_{24}+E_{57}.
\end{array}$$
Consider the linear map $P\in\Hom(\Real^8,\spin(7))$ defined as
$$\begin{array}{llll}
P(e_1)=0,&P(e_2)=-A_{14},&P(e_3)=0,&P(e_4)=A_{21},\\P(e_5)=A_{20},&P(e_6)=A_{21}-A_{18},&P(e_7)=A_{15}-A_{16},&P(e_7)=A_{14}-A_{17}.
\end{array}$$
It can be checked that $P\in\P(\spin(7))$. Moreover,  the elements $A_{14},A_{15}-A_{16},A_{17},A_{18},A_{20},A_{21}\in  \spin(7)$
generate the Lie algebra
$\spin(7)$.

The holonomy algebra of the metric
$$g=2dx^0dx^{9}+\sum^{8}_{i=1}(dx^i)^2+2\sum^{8}_{i=1}u^i dx^{i} dx^{9},$$
where
$$\begin{array}{ll}
u^1=-\frac{4}{3}x^7x^8,&
u^2=\frac{2}{3}((x^4)^2+x^3x^5+x^4x^6-(x^6)^2),\\
u^3=-\frac{4}{3}x^2x^5,&
u^4=\frac{2}{3}(-x^2x^4-2x^2x^6-x^5x^7+2x^6x^8),\\
u^5=\frac{2}{3}(x^2x^3+2x^4x^7+x^6x^7),&
u^6=\frac{2}{3}(x^2x^4+x^2x^6+x^5x^7-x^4x^8),\\
u^7=\frac{2}{3}(-x^4x^5-2x^5x^6+x^1x^8),&
u^8=\frac{2}{3}(-x^4x^6+x^1x^7),
\end{array}$$
at the point $0\in\Real^9$ is $\g^{2,\spin(7)}\subset\so(1,9)$.


\part{Weakly-irreducible not irreducible subalgebras  of $\su(1,n+1)$}

In this chapter we classify weakly-irreducible not irreducible subalgebras  of $\su(1,n+1)$. The result is stated in Section \ref{w-ir}.
In the other sections we give preliminaries  and the proof of the main theorem.

\section{Classification of weakly-irreducible not irreducible subalgebras  of $\su(1,n+1)$}\label{w-ir}

Let $\Real^{2,2n+2}$ be a $2n+4$-dimensional real  vector space endowed with a complex structure $J\in\Aut \Real^{2,2n+2},
J^2=-\id$ and with a $J$-invariant metric $\eta$ of signature $(2,2n+2)$, i.e. $\eta(Jx,Jy)=\eta(x,y)$ for all $x,y\in\Real^{2,2n+2}$.
We fix a basis $p_1$, $p_2$, $e_1$,...,$e_n$,  $f_1$,...,$f_n$, $q_1$, $q_2$ of the vector space $\Real^{2,2n+2}$
such that the Gram matrix of the metric $\eta$  and the complex structure $J$ have the forms
 $${\scriptsize \left (\begin{array}{ccccc}
0 &0 &0 & 1 & 0\\
0 &0 &0 & 0& 1\\
0 &0 & E_{2n} & 0 & 0 \\
1 &0 &0 &0 & 0 \\
0 &1 &0 &0 & 0 \\
\end{array}\right)} \text{ \rm{ and }}
{\scriptsize \left (\begin{array}{cccccc}
0&-1 &0 &0 & 0 & 0\\
1& 0 &0 &0 & 0 & 0\\
0 &0 &0& -E_n & 0 & 0 \\
0 &0 &E_n & 0 & 0& 0 \\
0 &0 &0 &0 & 0&-1 \\
0 &0 &0 &0 &1 &0 \\
\end{array}\right),}\text{ \rm{ respectively}}.$$



We denote by $\u(1,n+1)_{<p_1,p_2>}$ the subalgebra of $\u(1,n+1)$ that preserves the $J$-invariant 2-dimensional
isotropic  subspace  $\Real p_1\oplus\Real p_2\subset\Real^{2,2n+2}$.
The Lie algebra $\u(1,n+1)_{<p_1,p_2>}$ can be identified with the following matrix algebra
$${\scriptsize  \u(1,n+1)_{<p_1,p_2>}=\left\{\left. \left (\begin{array}{cccccc}
a_1&-a_2 &-z^t_1 & -z^t_2 &0 &-c\\
a_2&a_1 &z^t_2 & -z^t_1 &c &0\\
0 &0&B&-C&z_1&-z_2\\
0 &0&C&B&z_2&z_1\\
0&0&0&0&-a_1&-a_2\\
0&0&0&0&a_2&-a_1\\
\end{array}\right)\right|\,
\begin{array}{c}
a_1,a_2,c\in \Real,\\ z_1, z_2\in \Real^n,\\
\bigl(\begin{smallmatrix}B&-C\\C&B\end{smallmatrix}\bigr)\in \u(n) \end{array} \right\}.}$$
Recall that $$\u(n)=\left.\left\{\bigl(\begin{smallmatrix}B&-C\\C&B\end{smallmatrix}\bigr)\right|B\in\so(n),C\in\gl(n),C^t=C\right\}$$
and $$\su(n)=\left.\left\{\bigl(\begin{smallmatrix}B&-C\\C&B\end{smallmatrix}\bigr)\in\u(n)\right|\tr C=0\right\}.$$

We identify the above element of $\u(1,n+1)_{<p_1,p_2>}$  with the $7$-tuple
$(a_1,a_2,B,C,z_1,z_2,c)$. Define the following vector subspaces of
$\u(1,n+1)_{<p_1,p_2>}$:
$$\begin{array}{lcrlcr}
\A^1&=&\{(a_1,0,0,0,0,0,0)|a_1\in\Real\},&\A^2&=&\{(0,a_2,0,0,0,0,0)|a_2\in\Real\},\\
\N^1&=&\{(0,0,0,0,z_1,0,0)|z_1\in \Real^n\},&\N^2&=&\{(0,0,0,0,0,z_2,0)|z_2\in \Real^n\} \end{array}$$
and $$\C=\{(0,0,0,0,0,0,c)|c\in\Real\}.$$ We consider $\u(n)$ as a subalgebra of $\u(1,n+1)_{<p_1,p_2>}$.

We see that $\C$ is a  commutative ideal, which commutes with $\A^2$, $\N^1$, $\N^2$ and $\u(n)$,
and $\A^1\oplus\A^2$  is a commutative subalgebras, which commutes with $\u(n)$.

Furthermore, for $a_1,a_2,c\in \Real$, $z_1,z_2,w_1,w_2\in \Real^n$ and $\bigl(\begin{smallmatrix}B&-C\\C&B\end{smallmatrix}\bigr)\in \u(n)$
we obtain
$$\begin{array}{rcl}
[(a_1,0,0,0,0,0,0),(0,0,0,0,z_1,z_2,c)]&=&(0,0,0,0,a_1z_1,a_1z_2,2a_1 c),\\
{[(0,a_2,0,0,0,0,0),(0,0,0,0,z_1,z_2,0)]}&=&(0,0,0,0,a_2z_2,-a_2z_1,0),\\
{[(0,0,B,C,0,0,0),(0,0,0,0,z_1,z_2,0)]}&=&(0,0,0,0,Bz_1-Cz_2,Cz_1+Bz_2,0),\\
{[(0,0,0,0,z_1,z_2,0),(0,0,0,0,w_1,w_2,0)]}&=&(0,0,0,0,0,0,2(-z_1w_2^t+z_2w_1^t)).
\end{array}$$
Hence  we obtain the decomposition\footnote{We denote by $\ZR$ and $\zr$
the semi-direct product and semi-direct sum for Lie groups and Lie algebras, respectively. If a Lie algebra is decomposed
into a direct sum of vector subspaces, then we use $+$.}
$$\u(1,n+1)_{<p_1,p_2>}=(\A^1\oplus\A^2\oplus\u(n))\zr(\N^1+\N^2+\C).$$

Denote by $\su(1,n+1)_{<p_1,p_2>}$ the subalgebra of $\su(1,n+1)$ that preserves the subspace $\Real p_1\oplus\Real p_2\subset\Real^{2,2n+2}$.
Then
$$\su(1,n+1)_{<p_1,p_2>}=\{(a_1,a_2,B,C,z_1,z_2,c)\in\u(1,n+1)_{<p_1,p_2>}|2a_2+\tr_\Real C=0\}$$
and $$\u(1,n+1)_{<p_1,p_2>}=\su(1,n+1)_{<p_1,p_2>}\oplus\Real J.$$
Therefore we obtain the decomposition
$$\su(1,n+1)_{<p_1,p_2>}=(\A^1\oplus\su(n)\oplus\Real I_0)\zr(\N^1+\N^2+\C),$$
where  {\scriptsize$$I_0=  \left (\begin{array}{cccccc}
0&\frac{n}{n+2} &0 &0 & 0 & 0\\
-\frac{n}{n+2}& 0 &0 &0 & 0 & 0\\
0 &0 &0& -\frac{2}{n+2}E_n & 0 & 0 \\
0 &0 &\frac{2}{n+2}E_n & 0 & 0& 0 \\
0 &0 &0 &0 & 0&\frac{n}{n+2} \\
0 &0 &0 &0 &-\frac{n}{n+2} &0 \\
\end{array}\right).$$}

Note that {\scriptsize$$\u(1,1)_{<p_1,p_2>}=\left\{\left. \left (\begin{array}{cccc}
a_1&-a_2  &0 &-c\\
a_2&a_1  &c &0\\
0&0&-a_1&-a_2\\
0&0&a_2&-a_1\\
\end{array}\right)\right|\, a_1,a_2,c\in \Real \right\}=(\A^1\oplus \A^2)\zr\C$$}
and $\su(1,1)_{<p_1,p_2>}=\A^1\zr\C$.

If a weakly-irreducible subalgebra $\g\subset \u(1,n+1)$ preserves a degenerate proper subspace $W\subset \Real^{2,2n+2}$, then $\g$
preserves the $J$-invariant 2-dimensional isotropic  subspace  $W_1\subset\Real^{2,2n+2}$,
where $W_1=(W\cap JW)\cap(W\cap JW)^\bot$ if $W\cap JW\neq\{0\}$ and $W_1=(W\oplus JW)\cap(W\oplus JW)^\bot$ if $W\cap JW=\{0\}$.
Therefore $\g$ is conjugated to a weakly-irreducible subalgebra of $\u(1,n+1)_{<\Real p_1,\Real p_2>}$.

Let $E=\spa\{e_1,...,e_n,f_1,...,f_n\}$, $E^1=\spa\{e_1,...,e_n\}$ and  $E^2=\spa\{f_1,...,f_n\}$.
For any integers $k$ and $l$ with $1\leq k\leq l\leq n$ we consider the following subspaces:
$$E^1_{k,...,l}=\spa\{e_k,...,e_l\}\subset E^1,E^2_{k,...,l}=\spa\{f_k,...,f_l\}\subset E^2, E_{k,...,l}=E^1_{k,...,l}\oplus E^2_{k,...,l}\subset E,$$
$$\N^1_{k,...,l}=\{(0,0,0,0,z_1,0,0)|z_1\in E^1_{k,...,l}\}\subset \N^1$$  and
$$\N^2_{k,...,l}=\{(0,0,0,0,0,z_2,0)|z_2\in E^2_{k,...,l}\}\subset \N^2.$$
Clearly, $E^1=E^1_{1,...,n}$, $E^2=E^2_{1,...,n}$, $E=E_{1,...,n}$, $\N^1=\N^1_{1,...,n}$ and  $\N^2=\N^2_{1,...,n}$.

We denote by $\u(e_k,...,e_l)$ the subalgebra of $\u(n)$ that preserves the  vector subspace $E_{k,...,l}\subset E$
and annihilates the orthogonal complement to this subspace. We denote $\u(1,...,l)$ just by $\u(l)$.
Furthermore, let $J_{k,...,l}$ be the element of $\u(1,n+1)_{<\Real p_1,\Real p_2>}$
defined by  $J_{k,...,l}|_{E_{k,...,l}}=J|_{E_{k,...,l}}$ and $J_{k,...,l}|_{E_{k,...,l}^\bot}=0$.
We denote $J_{1,...,l}$ just by $J_l$.
Consider the following Lie algebra $${\scriptsize
\sod(k,...,l)=\left\{\left. \left (\begin{array}{cccccc}
0&0&0&0&0&0\\
0&B&0&0&0&0\\
0&0&0&0&0&0\\
0&0&0&0&0&0\\
0&0&0&0&B&0\\
0&0&0&0&0&0\\
\end{array}\right)\right|B\in\so(l-k+1)\right\}\subset\su(k,...,l),}$$ where the matrices of the
operators are written with respect to the decomposition $$E=E^1_{1,...,k-1}\oplus E^1_{k,...,l}\oplus E^1_{l+1,...,n}\oplus
E^2_{1,...,k-1}\oplus E^2_{k,...,l}\oplus E^2_{l+1,...,n}.$$
The subalgebra  $\sod(k,...,l)\subset\u(n)$ annihilates the orthogonal complement to $E^1_{k,...,l}\oplus E^2_{k,...,l}$ and
acts diagonally on $E^1_{k,...,l}\oplus E^2_{k,...,l}$.

For any $0\leq m\leq n$ define the following vector space
$$\t\A^2=\{(0,a_2,0,0,0,0,0)+a_2J_{m+1,...,n}|a_2\in\Real\}.$$
Clearly, if $m=n$, then $\t\A^2=\A^2$; if $m=0$, then $\t\A^2=\Real J$.

Let $\h\subset\u(n)$ be a subalgebra. Recall that $\h$ is a compact Lie
algebra and we have the decomposition $\h=\h'\oplus\z(\h)$, where $\h'$ is the commutant of $\h$ and $\z(\h)$ is the center of $\h$, see for example \cite{V-O}.

\begin{theorem}\label{Tw-ir}
{\bf 1)} A subalgebra $\g\subset\su(1,1)\subset\so(2,2)$ is weakly-irreducible and not irreducible if and only if
$\g$ is conjugated to the subalgebra $\C\subset\su(1,1)_{<p_1,p_2>}$ or to $\su(1,1)_{<p_1,p_2>}$.

{\bf 2)} Let $n\geq 1$. Then a subalgebra $\g\subset\su(1,n+1)\subset\so(2,2n+2)$ is weakly-irreducible and  not irreducible
if and only if $\g$ is conjugated to one of the following
subalgebras of $\su(1,n+1)_{<p_1,p_2>}$:

\begin{description}

\item[$\g^{m,\h,\A^1}$=]$(\A^1\oplus\{(0,-\frac{1}{2}\tr C,B,C,0,0,0)|\bigl(\begin{smallmatrix}B&-C\\C&B\end{smallmatrix}\bigr)\in \h\})\zr(\N^1+\N^2_{1,...,m}+\C)$,\\
 where $0\leq m\leq n$ and  $\h\subset\su(m)\oplus(J_m-\frac{m}{n+2}J_n)\oplus\sod(m+1,...,n)$ is a subalgebra;

\item[$\g^{m,\h,\varphi}$=]$\{(\varphi(B,C),-\frac{1}{2}\tr C,B,C,0,0,0)|\bigl(\begin{smallmatrix}B&-C\\C&B\end{smallmatrix}\bigr)\in \h\}\zr(\N^1+\N^2_{1,...,m}+\C)$,\\
 where $0\leq m\leq n$,  $\h\subset\su(m)\oplus(J_m-\frac{m}{n+2}J_n)\oplus\sod(m+1,...,n)$ is a subalgebra and
 $\varphi:\h\to\Real$ is a linear map with $\varphi|_{\h'}=0$;

\item[$\begin{array}{rl}\g^{n,\h,\psi,k,l}=
&\{(0,-\frac{1}{2}\tr C,B,C,\psi_1(B,C),\psi_2(B,C)+\psi_3(B,C),0)|\bigl(\begin{smallmatrix}B&-C\\C&B\end{smallmatrix}\bigr)\in \h\}\\
&\zr(\N^1_{1,...,k}+\N^2_{1,...,k}+\N^1_{l+1,...,n}+\C),\end{array}$] $\qquad$\\      where
$k$ and $l$ are integers such that $0< k\leq l\leq n$, $\h\subset\su(k)\oplus\Real(J_k-\frac{k}{n+2}J_n)$ is a subalgebra with $\dim\z(\h)\geq n+l-2k$,
$\psi:\h\to E^1_{k+1,...,l}+E^2_{k+1,...,l}+E^2_{l+1,...,n}$ is a surjective linear map with  $\psi|_{\h'}=0$,
$\psi_1=\pr_{E^1_{k+1,...,l}}\circ\psi$, $\psi_2=\pr_{E^2_{k+1,...,l}}\circ\psi$ and $\psi_3=\pr_{E^2_{l+1,...,n}}\circ\psi$;

\item[$\begin{array}{rl}\g^{m,\h,\psi,k,l,r}=
&\{(0,-\frac{1}{2}\tr C,B,C,\psi_1(B,C)+\psi_4(B,C),\psi_2(B,C)+\psi_3(B,C),0)|\bigl(\begin{smallmatrix}B&-C\\C&B\end{smallmatrix}\bigr)\in \h\}\\
&\zr(\N^1_{1,...,k}+\N^2_{1,...,k}+\N^1_{l+1,...,m}+\N^1_{m+1,...,r}+\C),\end{array}$] $\qquad$\\ where
$k$, $l$, $r$ and $m$ are integers such that $0< k\leq l\leq m\leq r\leq n$ and $m<n$,
$\h\subset\su(k)\oplus\Real(J_k-\frac{k}{n+2}J_n)\oplus\sod(m+1,...,r)$ is a subalgebra with $\dim\z(\h)\geq n+m+l-2k-r$,
$\psi:\h\to E^1_{k+1,...,l}+E^2_{k+1,...,l}+E^2_{l+1,...,m}+E^1_{r+1,...,n}$ is a surjective linear map with $\psi|_{\h'}=0$,
$\psi_1=\pr_{E^1_{k+1,...,l}}\circ\psi$, $\psi_2=\pr_{E^2_{k+1,...,l}}\circ\psi$,  $\psi_3=\pr_{E^2_{l+1,...,m}}\circ\psi$ and $\psi_4=\pr_{E^1_{r+1,...,n}}\circ\psi$;

\item[$\g^{0,\h,\psi,k}$=]$\{(0,0,B,0,\psi(B),0,0)|\bigl(\begin{smallmatrix}B&0\\0&B\end{smallmatrix}\bigr)\in \h\}\zr(\N^1_{1,...,k}+\C),$\\ where
$0<k<n$, $\h\subset\sod(1,...,k)$ is a subalgebra such that $\dim\z(\h)\geq n-k$, $\psi:\h\to E^1_{k+1,...,n}$ is a surjective linear map with $\psi|_{\h'}=0$;


\item[$\g^{0,\h,\zeta}$=]
$\{(0,0,B,0,0,0,\zeta(B))|\bigl(\begin{smallmatrix}B&0\\0&B\end{smallmatrix}\bigr)\in \h\} \zr\N^1$,\\ where
$\h\subset\sod(1,...,n)$ is a subalgebra with $\z(\h)\neq\{0\}$ and  $\zeta:\h\to\Real$ is a non-zero linear map with
$\zeta|_{\z(\h)}\neq 0$;

\item[$\g^{0,\h,\psi,k,\zeta}$=]$
\{(0,0,B,0,\psi(B),0,\zeta(B))|\bigl(\begin{smallmatrix}B&0\\0&B\end{smallmatrix}\bigr)\in\h,\}\zr\N^1_{1,...,k}$,\\ where
 $1 \leq k<n$, $\h\subset\sod(1,...,k)$ is a subalgebra with $\dim\z(\h)\geq n-k$,
$\psi:\h\to E^1_{k+1,...,n}$ is a surjective linear map with $\psi|_{\h'}=0$ and $\zeta:\h\to\Real$ is a non-zero linear
map with $\zeta|_{\h'}=0$.
\end{description}
\end{theorem}

{\bf Remark.}
In \cite{Ik22} A.~Ikemakhen  classified weakly-irreducible subalgebras of $\so(2,N+2)_{<p_1,p_2>}$ that contain the ideal
$\C$, i.e.  the last two types of Lie algebras from the above theorem were not considered in \cite{Ik22}.


\section{Transitive similarity transformation groups of Euclidian spaces}


Let $(E,\eta)$ be an Euclidean space, $\|x\|^2=\eta(x,x)$. A map $f:E\to E$
is called a {\it similarity transformation} of $E$ if there exists
a $\lambda >0$ such that $\|f(x_1)-f(x_2)\|=\lambda\|x_1-x_2\|$ for all
$x_1,x_2\in E$. If $\lambda=1$, then $f$
is called an {\it isometry}.
Denote by $\Simil E$ and $\Isom E$ the  groups of all similarity
transformations and isometries of $E$, respectively.
A subgroup $G\subset \Simil E$ such that $G\not\subset\Isom E$ is called
{\it essential}. A subgroup $G\subset \Simil E$ is called {\it irreducible}
if it does not preserve any proper affine subspace of $E$.

The following theorem is due to D.V. Alekseevsky (see \cite{Al2} or \cite{A-V-S}).
\begin{theorem}\label{Al}
A subgroup $G\subset\Simil  E$ acts irreducibly on $E$ if
and only if it acts transitively.
\end{theorem}

We denote by $E$ the group of all translations in $E$ and by
$A^1=\Real^+$ the identity component of the group of all dilations of $E$ about
the origin.
For the connected  identity component of the  Lie group   $\Isom E$ we have the decomposition
$$\Isom^0 E=SO(E)\ZR E,$$ where  $E$ is a normal subgroup of $\Isom^0 E$.
For $\Simil^0 E$ we obtain $$\Simil^0 E=A^1\ZR\Isom^0 E=(A^1\times SO(E))\ZR E,$$
where $E$ is a normal subgroup of $\Simil^0 E$ and $A^1$ commutes with $SO(E)$.

In \cite{Gal2} we  deduced from results of \cite{Al2} and \cite{A-V-S}
the following theorem.

\begin{theorem} \label{trans}
Let $G\subset\Simil E$ be a
transitively acting connected subgroup. Then $G$ belongs to one of the following types
\begin{description}
\item[type 1.] $G=(A\times H)\ZR E$, where $H\subset SO(E)$ is a connected Lie subgroup;
\item[type 2.] $G=H\ZR E$;
\item[type 3.] $G=(A^\Phi \times H)\ZR E$,
where $\Phi: A\to SO(E)$ is a non-trivial homomorphism and
$$A^\Phi=\{a\cdot\Phi(a)|a\in A\}\subset A\times SO(E)$$
 is a group of screw
dilations of $E$ that commutes with $H$;
\item[type 4.] $G=(H\times U^{\Psi})\ZR W,$ where $H\subset SO(W)$ is a connected Lie subgroup, $E=W\oplus U$ is an orthogonal decomposition,
$\Psi:U\to SO(W)$ is a homomorphism with $\ker \text{\rm d\,}\Psi=\{0\}$, and
$$U^{\Psi}=\{\Psi(u)\cdot u|u\in U\}\subset SO(W)\times U$$
is a group of screw isometries of $E$ that commutes with $H$.
\end{description}
\end{theorem}


To the above decomposition of the Lie group $\Simil E$ corresponds the following decomposition of its
Lie algebra $\lie(\Simil E)$:  $$\lie(\Simil E)=(\A^1\oplus\so(E))\zr E,$$
where $\A^1=\Real$ is the Lie algebra of the Lie group $A^1$ and $E$ is the Lie algebra of the Lie group $E$.
We see that $\A^1$ commutes with $\so(E)$, and $E$ is a commutative ideal in $\lie(\Simil E)$.

Let $\B$ be a Lie algebra. We denote by $\B'$ the commutant of $\B$ and by $\z(\B)$ the center of $\B$.
If $\B\subset\so(n)$, then $\B$ is  a compact Lie algebra and we have $\B=\B'\oplus \z(\B).$

The Lie algebras of the Lie groups from the above theorem have the following forms (see \cite{Gal2}):
\begin{description}
\item[type 1.] $\g^{1,\B}=(\mathcal A^1\oplus\mathcal B)\zr E$,
where $\B\subset\so(E)$ is a subalgebra;

\item[type 2.] $\g^{2,\B}=\mathcal B \zr E$;

\item[type 3.] $\g^{3,\B,\varphi}=(\mathcal B'\oplus
 \{(\varphi(x)+x)|x\in \mathcal \z(\B)\})\zr E$, where $\varphi
 :\z(\B)\to\A^1$ is a non-zero linear map;

\item[type 4.] $\g^{4,k,\B,\psi}=(\mathcal B'\oplus
\{(\psi(x)+x)|x\in \z(B)\})\zr W$, where we have a non-trivial   orthogonal
decomposition $E=W\oplus U$ such that $\B\subset\so(W)$, and
$\psi:\z(\B)\to  U$ is a surjective linear map; $k=\dim W$.
\end{description}

It is convenient for us to extend the maps $\varphi$ and $\psi$ to $\B$ and to unify the Lie algebras of type 2 and 3 assuming
that $\varphi=0$ for the Lie algebras of type 2. We obtain the following

\begin{theorem} \label{transalg}
Let $G\subset\Simil E$ be a
transitively acting connected subgroup. Then the Lie algebra of  $G$ belongs to one of the following types
\begin{description}
\item[type $\A^1$.] $\g^{\B,\A^1}=(\mathcal A^1\oplus\mathcal B)\zr E$,
where $\B\subset\so(E)$ is a subalgebra;

\item[type $\varphi$.] $\g^{\B,\varphi}= \{\varphi(x)+x|x\in \B\}\zr E$, where $\varphi
 :\B\to\A^1$ is a linear map with $\varphi|_{\B'}=0$;

\item[type $\psi$.] $\g^{\B,\psi,k}=\{\psi(x)+x|x\in \B\})\zr W$, where we have a non-trivial   orthogonal
decomposition $E=W\oplus U$ such that $\B\subset\so(W)$, and
$\psi:\B\to  U$ is a surjective linear map with $\psi|_{\B'}=0$; $k=\dim W$.
\end{description}
\end{theorem}

For each Lie algebra $\g$ from Theorem \ref{transalg}  we call $\B\subset\so(E)$ {\it the orthogonal part of} $\g$.

\section{Similarity transformations of the  Heisenberg spaces}

In this section we explain notation  from \cite{Gold}, which we will use later.

Let $\tilde E$ be a complex vector space of
dimension $n$ endowed with a Hermitian metric $g$.

By definition, {\it the Heisenberg space associated to $n$}
is the direct sum $\H_n=\tilde E\oplus \Real$.
The line $\Real$ is called the {\it vertical axis}.

We consider $\H_n$ also as a group with respect to the operation
$$(z,u)\cdot(w,v)=(z+w,u+v+2\im g(z,w)),$$
where $z,w\in\tilde E$ and $u,v\in \Real$.
The group $\H_n$ is nilpotent and $\Real\subset\H_n$ is a normal subgroup.

We consider the action of the group  $\H_n$ on itself by left translations.
These transformations are called  {\it the  Heisenberg translations}.

The unitary group $U(\tilde E)$ acts on $\H_n$ by
$$A:(z,u)\mapsto (Az,u),$$
where $A\in U(\tilde E)$. These transformations are called
{\it the Heisenberg rotations about the vertical axis.}

The group $\Co^*$ of non-zero complex numbers acts on $\H_n$
by $$\lambda:(z,u)\mapsto (\lambda z,|\lambda|^2u),$$
where $\lambda\in\Co^*$. These transformations are called
{\it the complex Heisenberg dilations about the origin}.

The intersection of the group of the  Heisenberg rotations about the vertical axis
and the group of the complex Heisenberg dilations about the origin is
the group of scalar multiplications by unit complex numbers $\mathbb{T}=\Co^*\cap U(\tilde E)$.
The groups $\Co^*$ and $U(\tilde E)$
generate the group $\Real^+\times U(\tilde E),$ where $\Real^+$ is the group of
{\it the real Heisenberg dilations about the origin}, i.e. with $\lambda\in \Real^+$.

All the above transformations generate {\it the Heisenberg similarity transformation group}
$$\Simil\H_n=(\Real^+ \times U(\tilde E))\ZR \H_n,$$
where the subgroup $\H_n\subset\Simil\H_n$ is normal.

By definition, {\it the similarity transformation group of the Hermitian  space} $\tilde E$
is $$\Simil\tilde E=(\Real^+ \times U(\tilde E))\ZR \tilde E,$$
where $\tilde E\subset\Simil\tilde E$ is a normal subgroup
that consists of translations in $\tilde E$.

We have the natural projection
$$\pi:\Simil{\H_n}\to\Simil\tilde E.$$
The  kernel of $\ \pi\ $ is 1-dimensional
and  consists of the Heisenberg translations
$(z,u)\mapsto (z,u+c)$, $c\in\Real$.

To the above decomposition of the Lie group $\Simil\H_n$ corresponds the following decomposition of
its Lie algebra
 $$\lie(\Simil \H_n)=(\A^1\oplus\u(n))\zr\lie(\H_n),$$
where $\A^1=\Real$ is the Lie algebra of the Lie group of the  real Heisenberg dilations about the origin.

Let us denote by $\C$ the Lie algebra of the group of the Heisenberg translations $(0,u):(w,v)\mapsto (w,u+v)$.
Let $\tilde E\subset\lie(\H_n)$ be the tangent space to the submanifold $\tilde E\subset\Simil \tilde E$ at the unit.
Then  $\C$ is an ideal in $\lie(\Simil \H_n)$ and $\tilde E$ is a vector subspace of $\lie(\H_n)$ with $[\tilde E,\tilde E]=\C$.
We have $\lie(\H_n)=\tilde E+\C$. Thus,
$$\lie(\Simil \H_n)=(\A^1\oplus\u(n))\zr(\tilde E+\C).$$
For the Lie algebra $\lie(\Simil\tilde E)$ of the Lie group $\Simil \tilde E$ we obtain the decomposition
 $$\lie(\Simil E)=(\A^1\oplus \u(n))\zr\tilde E,$$
where $\A^1=\Real$ is the Lie algebra of the group of real dilations about the origin,
$\tilde E$ is the Lie algebra of the Lie group $\tilde E$.
We see that $A^1$ commutes with $\u(n)$, and $\tilde E$ is an ideal in $\lie(\Simil E)$.

We denote the differential of the projection $\pi:\Simil{\H_n}\to\Simil\tilde E$  also by $\pi$.
Obviously, the linear  map $$\pi:\lie(\Simil{\H_n})\to\lie(\Simil\tilde E)$$ is surjective with the 1-dimensional
kernel $\C$.

\section{The groups $U(1,n+1)_{\Co p_1}$ and  $U(1,n+1)_{<p_1,p_2>}$, their Lie algebras and examples}

Let $S$ be a complex vector
space. Denote by $S_\Real$ the real vector space underlying $S$
and by $J$ the complex structure on $S_\Real$.
The correspondence $S\mapsto (S_\Real,J)$ gives an isomorphism of categories
of complex vector spaces and real vector spaces
with complex structures. For a real vector space $S$  with
a complex structure $J$  we denote by $\tilde S$ the complex vector space given by $(S,J)$.

Let $S$ be a complex vector space. A subspace $S_1\subset S_{\Real}$ is
called {\it complex} if $JS_1=S_1$, where $J$ is the
complex structure on $S_\Real$. A subspace $S_1\subset S_{\Real}$ is
called {\it a  real form} of $S$ if $JS_1\cap S_1=\{0\}$ and
$\dim_\Real S_1=\dim_\Co S$.

Suppose that  $S$ is endowed with a pseudo-Hermitian metric $g$. For $x,y\in S_\Real$
let\\ $\eta(x,y)=\re g(x,y)$. Then $\eta$ is a metric on $S_\Real$
and we have\\ $\eta(Jx,Jy)=\eta(x,y)$ for all
$x,y\in S_\Real$, i.e. $\eta$ is $J$-invariant.
Conversely, for a given real vector space $S$ with a complex structure $J$
and an  $J$-invariant non-degenerate metric $\eta$, let
$g(z,w)=\eta(z,w)+i\eta(z,Jw)$ for all $z,w\in \tilde S$. Then $g$
is a pseudo-Hermitian metric on $\tilde S$. This gives us an isomorphism
of categories of pseudo-Hermitian spaces and real vector spaces endowed
with complex structures and invariant non-degenerate metrics.

A vector subspace $L\subset S$ (resp. $\tilde S$) is called {\it degenerate} if
the restriction of $\eta$ (resp. $g$)  to $L$ is degenerate.
A vector subspace $L\subset S$ (resp. $\tilde S$) is called {\it isotropic} if
the restriction of $\eta$ (resp. $g$)  to $L$ is zero.

Let $S$ be a complex vector space of dimension $n+2$, where $n\geq 0$,
and let $g$ be a pseudo-Hermitian metric on $S$ of signature $(1,n+1)$.
By definition, the pseudo-unitary group $U(1,n+1)$ is the real Lie group
of $g$-invariant automorphisms of $S$, i.e.
$$U(1,n+1)=\{f\in\Aut(S)|g(fz,fw)=g(z,w) \text{ for all } z,w\in S\}.$$
The corresponding Lie algebra consists of $g$-skew symmetric endomorphisms
of $S$, i.e.
$$\u(1,n+1)=\{\xi\in\End(S)|g(\xi z,w)+g(z,\xi w)=0 \text{ for all } z,w\in S\}.$$

Consider the  action of the  group $U(1,n+1)$ on $S_\Real$. Then
$$U(1,n+1)=\{f\in SO(2,2n+2)| Jf=fJ\}$$ and
$$\u(1,n+1)=\{\xi\in\so(2,2n+2)|[\xi,J]=0\}.$$
Here $SO(2,2n+2)$ is the Lie group of $\eta$-orthogonal automorphisms
of $S_\Real$ and $\so(2,2n+2)$ is the corresponding Lie algebra,
which consists of  $\eta$-skew symmetric endomorphisms of $S_\Real$.

Let $\Real^{2,2n+2}$ be a $2n+4$-dimensional real  vector space endowed with a
complex structure $J\in\Aut \Real^{2,2n+2}$ and with a $J$-invariant metric $\eta$ of
signature $(2,2n+2)$.
Let $\Co^{1,n+1}$ be the $n+2$-dimensional complex vector space given by $(\Real^{2,2n+2},J,\eta)$.
Denote by $g$ the  pseudo-Hermitian metric on $\Co^{1,n+1}$ of signature $(1,n+1)$ corresponding to $\eta$.

{\it
We say that a subgroup $G\subset U(1,n+1)\subset SO(n+2,\Co)$ acts  weakly-irreducibly
on  $\Co^{1,n+1}$ if it does not preserve any non-degenerate proper subspace
 of $\Co^{1,n+1}$. We say that a subalgebra  $\g\subset \u(1,n+1)\subset\so(n+2,\Co)$ is  weakly-irreducible
if it does not preserve any non-degenerate proper subspace of $\Co^{1,n+1}$.}


\begin{prop}\label{prop1}
If a subgroup $G\subset U(1,n+1)$ acts weakly-irreducibly on $\Real^{2,2n+2}$,
then $G$ acts weakly-irreducibly on $\Co^{1,n+1}$.
\end{prop}
{\it Proof.} Suppose that  $G\subset U(1,n+1)$ acts weakly-irreducibly on $\Real^{2,2n+2}$
and $G$ preserves a non-degenerate proper subspace $L\subset\Co^{1,n+1}$.
Hence $G$ preserves the subspace $L_\Real\subset \Real^{2,2n+2}$. We claim that
$L_\Real$ is non-degenerate.
Indeed, let $x\in L_\Real\cap (L_\Real)^{\bot_\eta}$.
For any $y\in L$ we have $g(x,y)=\eta(x,y)+i\eta(x,Jy)=0$,
since $JL=L$. Hence, $x\in L\cap L^{\bot_g}$ and $x=0$.
Thus the subspace $L_\Real$ is non-degenerate and we have a contradiction.
$\Box$

The converse to Proposition \ref{prop1} is not true, see Example \ref{exWI2}  below.

For the group $U(1,n+1)$ we have the local  decomposition
$U(1,n+1)=SU(1,n+1)\cdot\mathbb{T}$, where $SU(1,n+1)=U(1,n+1)\cap SL(2+n,\Co)$
and $\mathbb{T}$ is the center of $U(1,n+1)$, which is the 1-dimensional subgroup
generated by the complex structure $J\in U(1,n+1)$.

We fix a basis
$p_1,e_1,...,e_n,q_1$ of $\Co^{1,n+1}$ such that the Gram matrix
of $g$ has the form $${\scriptsize
\left (\begin{array}{ccc}
0 & 0 & 1\\ 0 & E_n & 0 \\ 1 & 0 & 0 \\
\end{array}\right),}$$
 where $E_n$ is the $n$-dimensional identity matrix.
Let $\tilde E\subset \Co^{1,n+1}$ be the vector subspace spanned by $e_1,...,e_n$.
We will consider $\tilde E$ as Hermitian space with the metric
$g|_{\tilde E}.$

Denote by $U(1,n+1)_{\Co p_1}$ the Lie subgroup of $U(1,n+1)$ acting
on $\Co^{1,n+1}$ and  preserving the complex isotropic line $\Co p_1$.
Since $J\in U(1,n+1)_{\Co p_1}$, we have the decomposition
$U(1,n+1)_{\Co p_1}=SU(1,n+1)\cdot\mathbb{T}$, where $SU(1,n+1)_{\Co p_1}=U(1,n+1)_{\Co p_1}\cap SL(2+n,\Co)$.

The Lie algebra $\u(1,n+1)_{\Co p_1}$ of the Lie group $U(1,n+1)_{\Co p_1}$
can be identified with the following matrix algebra
$${\scriptsize \u(1,n+1)_{\Co p_1}=\left\{\left. \left (\begin{array}{ccc}
a & -\overline{z}^t & ic\\ 0 & A & z \\ 0 & 0 & -\overline a \\
\end{array}\right)\right| a\in \Co,\, c\in\Real,\, z\in \tilde E ,\,A \in \u(n)  \right\}.}$$
Here $\u(n)$ is the unitary Lie  algebra of the Hermitian space $\tilde E$.
We identify the above matrix with the quadruple  $(a,A,z,c)$ and
define the following vector subspaces of $\u(1,n+1)_{\Co p_1}:$
$$\A=\{(a,0,0,0)|a\in\Co\},\quad \N=\{(0,0,z,0)|z\in\tilde E\}\quad\text{and}\quad \C=\{(0,0,0,c)|c\in\Real\}.$$
We consider $\u(n)$ as a subalgebra of $\u(1,n+1)_{\Co p_1}$ with the obvious inclusion.
Note that $\C$ is a commutative ideal, which  commutes with
$\u(n)$ and $\N$, and $\A$ is a commutative subalgebra, which commutes
with $\u(n)$. For
$a\in \Co,$ $c\in \Real$,  $z,w\in \tilde E$ and $A \in \u(n)$ we obtain
$$[(a,0,0,0),(0,0,z,c)]=(0,0,\overline{a} z,2\re a \cdot c),\quad
[(0,A,0,0),(0,0,z,0)]=(0,0,Az,0)$$ and
$$[(0,0,z,0),(0,0,w,0)]=(0,0,0,-2\im g(z,w)).$$

We obtain the decomposition  $$\u(1,n+1)_{\Co p_1}=(\A\oplus\u(n))\zr(\N+\C).$$

For the Lie algebra $\su(1,n+1)_{\Co p_1}$ of the Lie group
$SU(1,n+1)_{\Co p_1}$ we have\\
$\su(1,n+1)_{\Co p_1}=\u(1,n+1)_{\Co p_1}\cap\su(1,n+1)=
\{\xi\in\u(1,n+1)_{\Co p_1}|\tr_\Co\xi=0\}=$\\
$\{(a,A,z,c)\in\u(1,n+1)_{\Co p_1}|a-\bar a+\tr_\Co A=0\}$
and $\u(1,n+1)_{\Co p_1}=\su(1,n+1)_{\Co p_1}\oplus \Real J$.

Let $\g\subset\u(1,n+1)\subset\so(n+2,\Co)$ be a weakly-irreducible and not irreducible subalgebra.
Then  $\g$ preserves
a non-degenerate proper  subspace $L\subset\Co^{1,n+1}$. Hence
$\g$ preserves the orthogonal complement $L^{\bot_ g}$ and
the intersection $L\cap L^{\bot_ g}$, which is an isotropic
complex line. Hence $\g$ is conjugated to a weakly-irreducible
subalgebra of $\u(1,n+1)_{\Co p_1}$.

Now we consider the real vector space $\Real^{2,2n+2}$. Let $p_2=Jp_1$,  $f_1=Je_1$,...,
 $f_n=Je_n$ and  $q_2=Jq_1$. Consider the basis $p_1$, $p_2$, $e_1$,...,$e_n$,
 $f_1$,...,$f_n$, $q_1$, $q_2$ of the vector space $\Real^{2,2n+2}$.
With respect to this basis the Gram matrix of the metric
$\eta$  and the complex structure $J$ have the forms
$${\scriptsize \left (\begin{array}{ccccc}
0 &0 &0 & 1 & 0\\
0 &0 &0 & 0& 1\\
0 &0 & E_{2n} & 0 & 0 \\
1 &0 &0 &0 & 0 \\
0 &1 &0 &0 & 0 \\
\end{array}\right)} \text{ \rm{ and }}
{\scriptsize \left (\begin{array}{cccccc}
0&-1 &0 &0 & 0 & 0\\
1& 0 &0 &0 & 0 & 0\\
0 &0 &0& -E_n & 0 & 0 \\
0 &0 &E_n & 0 & 0& 0 \\
0 &0 &0 &0 & 0&-1 \\
0 &0 &0 &0 &1 &0 \\
\end{array}\right),}\text{ \rm{ respectively}}.$$
We consider the vector space $E=\spa_\Real\{e_1,...,e_n,f_1,...,f_n\}\subset \Real^{2,2n+2}$
as an Euclidian  space with the metric $\eta|_E$. Let $E^1=\spar\{e_1,...,e_n\}$ and  $E^2=\spar\{f_1,...,f_n\}$.

We denote by $U(1,n+1)_{<p_1,p_2>}$ the subgroup of $U(1,n+1)$ acting on $\Real^{2,2n+2}$ and
preserving the isotropic $2$-dimensional vector subspace
$\Real p_1\oplus\Real p_2\subset \Real^{2,2n+2}$. The group $U(1,n+1)_{<p_1,p_2>}$ is just
$U(1,n+1)_{\Co p_1}$ acting on $\Real^{2,2n+2}$.
Denote by $SU(1,n+1)_{<p_1,p_2>}$ the group
$SU(1,n+1)_{\Co p_1}$ acting on $\Real^{2,2n+2}$.
We have $U(1,n+1)_{<p_1,p_2>}=SU(1,n+1)_{<p_1,p_2>}\cdot\mathbb{T}$.

The Lie algebra $\u(1,n+1)_{<p_1,p_2>}$ of the Lie group $U(1,n+1)_{<p_1,p_2>}$
can be identified with the  matrix algebra as in Section \ref{w-ir}.
The correspondence between $\u(1,n+1)_{\Co p_1}$
and $\u(1,n+1)_{<p_1,p_2>}$  is given by the identities
$a=a_1+ia_2$, $A=B+iC$ and  $z=z_1+iz_2$.

Now we  consider some  examples (we use denotation of Section \ref{w-ir}).
\begin{ex}\label{exWI1} The subalgebra $\g=\N^1+\C\subset\u(1,n+1)_{<p_1,p_2>}$ is weakly-irreducible.
\end{ex}
{\it Proof.} Suppose that  $\g$ preserves a proper vector subspace $L\subset \Real^{2,2n+2}$.
Let $v=(a_1,a_2,\alpha,\beta,b_1,b_2)\in L$, where we use the decomposition
$$\Real^{2,2n+2}=\Real p_1\oplus\Real p_2\oplus\spar\{e_1,...,e_n\}
\oplus\spar\{f_1,...,f_n\}\oplus \Real q_1\oplus\Real q_2.$$
Applying  an element $(0,0,0,0,z_1,0,0)\in\N^1\subset\g$
with $z_1^tz_1=1$ twice to $v$, we see that $b_1p_1+b_2p_2\in L$.
Applying the element $(0,0,0,0,0,0,1)\in\C\subset\g$ to $v$, we get
$-b_2p_1+b_1p_2\in L$.
Hence if the projection of $L$ to $\Real q_1\oplus\Real q_2$ is non-zero,
then $L$ contains $\Real p_1\oplus\Real p_2$ and the projection of $L^{\bot_\eta}$
to $\Real q_1\oplus\Real q_2$ is zero, moreover, $L^{\bot_\eta}$ is also preserved.
Thus we can assume that the projection of $L$ to $\Real q_1\oplus\Real q_2$ is zero.
Let  $v=(a_1,a_2,\alpha,\beta,0,0)\in L$.
If $\alpha\neq 0$, then applying the element $(0,0,0,0,\alpha,0,0)\in\N^1\subset\g$ to $v$,
we see that $0\neq\alpha^t\alpha p_1 +\alpha^t\beta p_2\in L$, hence  $L$ is degenerate.
The similar statement is for $\beta\neq 0$. If $\alpha=\beta=0$ for all $v\in L$, then $L$
is degenerate. $\Box$
\begin{ex}
The subalgebra $\g=\N^1\subset\u(1,n+1)_{<p_1,p_2>}$ is not weakly-irreducible.
The subalgebra $\g=\N^1\subset\u(1,n+1)_{\Co p_1}$ is weakly-irreducible.
\end{ex}
{\it Proof.}  Suppose that $\g$ preserves a proper subspace $L\subset \Co^{1,n+1}$.
Let $v=(a_1,\alpha,b_1)\in L$ (here we use the decomposition
$\Co^{1,n+1}=\Co p_1\oplus\spac\{e_1,...,e_n\}\oplus\Co q_1$).
Applying an element $(0,0,z_1,0)\in\N^1$
with $\bar z_1^tz_1=1$ twice to $v$, we see that $b_1p_1\in L$. As in Example \ref{exWI1}, we can assume
that the projection of $\g$
to $\Co q_1$ is zero and show that $L$ is degenerate.

The vector subspaces $\spar\{p_1+p_2,e_1+f_1,...,e_n+f_n,q_1+q_2\}\subset \Real^{2,2n+2}$
and\\  $\spar\{p_1-p_2,e_1-f_1,...,e_n-f_n,q_1-q_2\}\subset \Real^{2,2n+2}$
are non-degenerate and preserved by $\g$. $\Box$

\begin{ex}\label{exWI2} The subalgebra $\g=\N^1\oplus\Real J\subset\u(1,n+1)_{<p_1,p_2>}$ is weakly-irreducible.
\end{ex}
{\it Proof.}
Suppose that  $\g$ preserves a proper vector subspace $L\subset \Real^{2,2n+2}$.
Let $v=(a_1,a_2,\alpha,\beta,b_1,b_2)\in L$.
As in Example \ref{exWI1}, we see that $b_1p_1+b_2p_2\in L$.
Applying  the element $J\in\g$ to $b_1p_1+b_2p_2$, we get
$-b_2p_1+b_1p_2\in L$. The end of the proof is as in Example \ref{exWI1}. $\Box$

Let $A^1$, $A^2$, $N^1$, $N^2$ and $C$ be the connected Lie subgroups
of $U(1,n+1)_{\Co p_1}$ corresponding to the subalgebras
$\A^1$, $\A^2$, $\N^1$, $\N^2$ and $\C$ of the Lie algebra
$\u(1,n+1)_{\Co p_1}$. These groups can be identified with the following
groups of matrices

$$\begin{array}{ll}
A^1=\left\{\left. \left (\begin{array}{ccc}
a_1 & 0 & 0\\ 0 & \id & 0\\ 0 & 0 & \frac{1}{a_1}\\
\end{array} \right)\right|\begin{array}{c} a_1\in \Real,\\a_1>0
\end{array} \right \},&
A^2=\left\{\left. \left (\begin{array}{ccc}
e^{ia_2} & 0 & 0\\ 0 & \id & 0\\ 0 & 0 & e^{ia_2}\\
\end{array} \right)\right| a_2\in \Real \right \},\\
N^1=\left\{\left. \left ( \begin{array}{ccc} 1 & - z_1^t & -
\frac{1}{2} z_1^tz_1\\ 0 & \id & z_1\\ 0 & 0& 1 \\
\end{array} \right)\right|z_1\in  E^1 \right \},&
N^2=\left\{\left. \left ( \begin{array}{ccc} 1 &  z_2^t &
\frac{1}{2} z_2^tz_2\\ 0 & \id & z_2\\ 0 & 0& 1 \\
\end{array} \right)\right| z_2\in E^2 \right \},\end{array}$$
$$C=\left\{\left. \left ( \begin{array}{ccc}
1 & 0 & ic\\ 0 & \id & 0\\ 0 & 0 & 1 \\
\end{array} \right)\right|c\in \Real \right \}.$$

We consider the group $U(n)$ as a Lie subgroup of
$U(1,n+1)_{\Co p_1}$ with the inclusion
$$A\in U(n)\mapsto\left (\begin{array}{ccc}1 & 0 & 0\\ 0 & A & 0\\ 0 & 0 & 1 \\\end{array} \right)\in U(1,n+1)_{\Co p_1}.$$

We have the decomposition $U(1,n+1)_{\Co p_1}=(A^1\times A^2 \times U(n))\ZR( N \cdot  C)$.


\section{Action of the group $U(1,n+1)_{\Co p_1}$  on the boundary of the  complex hyperbolic space}
\def\mypar{\thesection~Action of the group $U(1,n+1)_{\Co p_1}$  on $\partial\mathbf{H}^{n+1}_\mathbb{C}$}

As above, let $\Co^{1,n+1}$ be a complex vector space of dimension $n+2$ and $g$ be
a pseudo-Hermitian metric on $\Co^{1,n+1}$ of signature $(1,n+1)$.
A complex line $l\subset \Co^{1,n+1}$ is called {\it negative} if $g(z,z)<0$
for all $z\in l\backslash \{0\}$.
The $n+1$-dimensional {\it complex hyperbolic} space
$\mathbf{H}^{n+1}_\mathbb{C}$ is the
subset of the projective space  $\mathbb{P}\Co^{1,n+1}$ consisting
of all negative lines.

The {\it boundary} $\p\mathbf{H}^{n+1}_\mathbb{C}$
of the complex hyperbolic space
$\mathbf{H}^{n+1}_\mathbb{C}$ is the
subset of the projective space  $\mathbb{P}\Co^{1,n+1}$ consisting of
all  complex isotropic  lines.

We identify $\p\mathbf{H}^{n+1}_\mathbb{C}$
with a $2n+1$-dimensional real sphere in the following way.
Let $p_1,e_1,...,e_n,q_1$ by the basis of $\Co^{1,n+1}$ as above.
We also consider the basis $e_0,e_1,...,e_n,e_{n+1}$,
where  $e_0=\frac{\sqrt{2}}{2}(p_1-q_1)$ and
$e_{n+1}=\frac{\sqrt{2}}{2}(p_1+q_1)$.
With respect to this basis the Gram matrix of $g$ has the form
 $$\left(
\begin{array}{cc}
-1 & 0 \\
0 & E_{n+1}\\
 \end{array} \right).$$

Consider the vector subspace $\tilde{E_1}=\tilde{E}\oplus\mathbb{C} e_{n+1}\subset\tilde E$.
Each isotropic line intersects the affine subspace $e_0+\tilde{E_1}$ at an unique point
and we see that  \begin{multline*}(e_0+\tilde{E_1})\cap\{z\in \Co^{1,n+1}|g(z,z)=0\}=\\
\{z\in \Real^{2,2n+2} | x_0=1,\,y_0=0,\  x_1^2+y_1^2+\cdots +x_{n+1}^2+y_{n+1}^2=1\}\end{multline*}
is a $2n+1$-dimensional unite sphere $S^{2n+1}$.
Here $(x_0+iy_0,...,x_{n+1}+iy_{n+1})$ are the coordinates of a point $z$
with respect to the basis $e_0,e_1,...,e_n,e_{n+1}$.

Let $G\subset U(1,n+1)_{\mathbb{C}p_1}$ be a subgroup.
We identify $\p\mathbf{H}^{n+1}_\mathbb{C}\backslash \{\mathbb{C}p_1\}$
with a Heisenberg space $\mathcal H_n$ and  consider an action of $G$ on
$\mathcal H_n$.
For this let $E_2=\tilde E\oplus i\Real e_{n+1}$ and
$p_0=\Co p_1\cap (e_0+\tilde{E_1})=\sqrt{2}p_1$.  Denote by $s_0$
the stereographic projection
$s_0:S^{2n+1}\backslash\{p_0\}\to E_2$ (here we identify $E_2$ and $e_0+E_2$).
Note that $E_2$ is just a Heisenberg space $\H_n$
with the vertical axis $i\Real e_{n+1}$.
Any $f\in U(1,n+1)_{\Co p_1}$ takes complex isotropic
lines to complex isotropic lines, hence it
acts on  $S^{2n+1}\backslash\{p_0\}=\p\mathbf{H}^{n+1}_\mathbb{C}\backslash \{\mathbb{C}p_1\}$, and
we get the transformation
$\Gamma(f)=s_0\circ f \circ s^{-1}_0$ of $E_2=\H_n$. Moreover, we will see that this
transformation is a Heisenberg similarity transformation.

Now we  consider the basis $p_1,e_1,...,e_n,q_1$ and the  matrices of elements
of $U(1,n+1)_{\Co p_1}$ with respect to this basis.
Computations show that \\ for $a_1\in \Real$, $a_1>0$, the element
$\left (\begin{array}{ccc}
a_1 & 0 & 0\\ 0 & \id & 0\\ 0 & 0 & \frac{1}{a_1} \\
\end{array} \right)$
\begin{tabular}{c}
corresponds to the  real Heisenberg dilation\\
$(z,iu)\mapsto (a_1 z,a_1^2 iu)$,
\end{tabular} where $z\in \tilde E$ and $u\in\Real$;\\
for $a_2\in \Real$ the element
$\left (\begin{array}{ccc}
e^{ia_2} & 0 & 0\\ 0 & \id & 0\\ 0 & 0 & e^{ia_2}\\
\end{array} \right)$
\begin{tabular}{c}
corresponds to the  complex Heisenberg dilation\\
$(z,iu)\mapsto (e^{-ia_2} z,iu)$;
\end{tabular} \\
for $A\in U(n)$ the element
$\left (
\begin{array}{ccc}
1 & 0 & 0\\ 0 & A & 0\\ 0 & 0 & 1 \\
\end{array} \right)$
\begin{tabular}{c}
corresponds the Heisenberg rotation \\
$(z,iu)\mapsto (Az,iu)$;
\end{tabular}\\
for $z_1\in E^1$, $z_2\in E^2$  and $c\in\Real$ the elements\\
$\left ( \begin{array}{ccc}
 1 & -z_1^t &
-\frac{1}{2}z_1^tz_1\\ 0 & \id & z_1\\ 0 & 0& 1 \\
\end{array} \right)$,
$\left ( \begin{array}{ccc}
 1 & z_2^t &
\frac{1}{2}z_2^tz_2\\ 0 & \id & z_2\\ 0 & 0& 1 \\
\end{array} \right)$
 and
$\left (
\begin{array}{ccc}
1 & 0 & ic\\ 0 & \id & 0\\ 0 & 0 & 1 \\
\end{array} \right)$\\
 correspond to the Heisenberg translations\\
$(z,iu)\mapsto (z+z_1,iu+i2\im g(z,z_1)),\quad (z,iu)\mapsto (z+z_2,iu+i2\im g(z,z_2))$\\
and $(z,iu)\mapsto (z,iu+ic)$, respectively.

Thus we have a surjective Lie group
homomorphism $$\Gamma: U(1,n+1)_{\Co p_1}\to \Simil \H_n.$$ The kernel
of $\Gamma$ is the center $\mathbb{T}$ of $U(1,n+1)_{\Co p_1},$
hence the restriction $$\Gamma|_{SU(1,n+1)_{\Co p_1}}:SU(1,n+1)_{\Co p_1}\to \Simil \H_n$$ is a Lie group isomorphism.

We say that an affine subspace $L\subset E$ is {\it complex}
if the corresponding vector subspace is complex,
otherwise we say that $L$ is {\it non-complex}.

Let  $\pi:\Simil{\H_n}\to\Simil\tilde E$ be the obvious projection.
The homomorphism  $\pi$ is surjective and its kernel is 1-dimensional
and  consists of the Heisenberg translations
$(z,iu)\mapsto (z,iu+ic)$.
\begin{theorem}\label{non-compl}
Let $G\subset U(1,n+1)_{\Co p_1}$ be a weakly-irreducible subgroup.
Then
\begin{itemize}
\item[{\rm (1)}] The subgroup $\pi(\Gamma(G))\subset \Simil \tilde E$
does not preserve any proper
complex affine subspace of $E$.
\item[{\rm (2)}] If $\pi(\Gamma(G))\subset \Simil \tilde E$ preserves a proper
non-complex affine subspace $L\subset E$, then the minimal
complex affine subspace of $\tilde E$ containing $L$ is $\tilde E$.
\end{itemize}
\end{theorem}

{\it Proof.} (1)  First we prove that $\pi(\Gamma(G))\subset \Simil \tilde E$
does not preserve any proper
complex  vector  subspace of $E$.
Suppose that $\pi(\Gamma(G))$ preserves a proper complex
vector subspace $L\subset \tilde E$.
Then we have $\pi(\Gamma(G))\subset (\Real^+\times U(L)\times U(L^{\bot_g}))\ZR U$,
where $L^{\bot_g}$ is the orthogonal complement to $L$ in $\tilde E$.
We see that the group $G$ preserves the proper non-degenerate vector subspace
$L^{\bot_g}\subset \Co^{1,n+1}$.

Suppose that $\pi(\Gamma(G))\subset \Simil \tilde E$
preserves a proper
complex affine subspace $L\subset E$. Let $w\in L$
and let $L_0=-w+L$ be the corresponding to $L$ complex vector subspace
of $\tilde E$. Let $f=\left ( \begin{array}{ccc}
 1 &-\overline{w}^t &
-\frac{1}{2}\overline{w}^tw\\ 0 & \id & w\\ 0 & 0& 1 \\
\end{array} \right)\in U(1,n+1)_{\Co p_1}$.
Consider the subgroup $G_1=f^{-1}Gf\subset U(1,n+1)_{\Co p_1}$.
We have $\pi(\Gamma (G_1))=-w\cdot\pi(\Gamma(G))\cdot w$ and
$\pi(\Gamma (G_1))$ preserves the complex vector subspace
$L_0\subset\tilde E$. Hence $G_1$ preserves the complex vector subspace
$L_0^{\bot_g}\subset\Co^{1,n+1},$ which is non-degenerate.
Thus $G$ preserves the proper complex subspace
$f(L_0^{\bot_g})\subset\Co^{1,n+1},$
which is also non-degenerate. This gives us a contradiction.

(2) Suppose that $\pi(\Gamma(G))\subset \Simil \tilde E$
preserves a non-complex affine subspace $L\subset E$.
By the same argument, $G$ preserves the proper complex subspace
$(\spa_{\Co} L_0)^{\bot_g}\subset\Co^{1,n+1},$
which is non-degenerate. Hence, $\spa_{\Co} L_0=\tilde E$.
This proves the theorem. $\Box$

\section{Proof of Theorem \ref{Tw-ir}}
{\bf 1)} First consider the case $n=0$. There are four subalgebras of $\su(1,1)_{<p_1,p_2>}=\A^1\zr\C$:
$$\A^1\zr\C,\qquad \C,\qquad\A^1,\qquad \{(c\gamma,c)\in\A^1\zr\C|c\in\Real\},\text{ where } \gamma\in\Real,\gamma\neq 0.$$
The last subalgebra preserves the non-degenerate proper subspace $\spa\{p_1+p_2+\gamma q_1+\gamma q_2,2p_2-\gamma q_1-\gamma q_2\}\subset\Real^{2,2}$.
The subalgebra $\A^1\subset\su(1,1)_{<p_1,p_2>}$ preserves the non-degenerate proper subspace $\Real p_1\oplus\Real q_1\subset\Real^{2,2}$.

We claim that the subalgebra $\C\subset\su(1,1)_{<p_1,p_2>}$ is weakly-irreducible.
Indeed, suppose that $\C$ preserves a proper subspace $L\subset\Real^{2,2}$. We may assume that $\dim L=1$ or $\dim L=2$
(if $\dim L=3$, then we consider $L^\bot$). Let $\alpha_1 p_1+\alpha_2 p_2+\beta_1 q_1 +\beta_2 q_2\in L$ be a non-zero vector.
Applying the element $(0,1)\in\C$, we get $\beta_1 p_2 -\beta_2 p_1\in L$. If $\dim L=1$, then $\beta_1=\beta_2=0$ and  $L=\Real(\alpha_1 p_1+\alpha_2 p_2)$.
If $\dim L=2$, then $L=\spa\{\alpha_1 p_1+\alpha_2 p_2+\beta_1 q_1 +\beta_2 q_2, \beta_1 p_2 -\beta_2 p_1\}$. In both cases $L$ is degenerate.

Thus the subalgebra $\C\subset\su(1,1)_{<p_1,p_2>}$ is weakly-irreducible. Hence the subalgebra  $\su(1,1)_{<p_1,p_2>}=\A^1\zr\C$
is also weakly-irreducible. Part 1) of Theorem \ref{Tw-ir} is proved.

\vskip0.3cm
{\bf 2)} Suppose that $n\geq 1$.
Theorem \ref{non-compl} shows us that we must find all connected Lie subgroups $F\subset\Simil\tilde E$ that
satisfy the conclusion of Theorem \ref{non-compl}. Then for each subgroup $F$ find all connected subgroups
$G\subset SU(1,n+1)_{<p_1,p_2>}$ with $\pi(\Gamma(G))=F$ and check which of these groups act weakly-irreducibly
on $\Real^{2,2n+2}$. Since all groups are connected, we may  do some steps in terms of their Lie algebras.

{\bf Step 1.} First we describe non-complex vector subspaces $L\subset E$ with $\spac L=\tilde E$.
Let $L\subset E$ be such a subspace and $L_0=L\cap JL$. Let $L_1$ be the orthogonal complement to $L_0$
in $L$. Obviously, $L_0$ is a complex subspace and $L_1\cap JL_1=\{0\}$. Since $\spac L=\tilde E$,
we see that $\spac L_1$ is the orthogonal complement to $L_0$ in $\tilde E$. Thus $L_1$
is a real form of the complex vector space $\spac L_1$. We choose the above  basis $e_1,...,e_n,f_1,...,f_n$ in such a way that
$L_0=\spar\{e_1,...,e_m,f_1,...,f_m\}$ and $L_1=\spar\{e_{m+1},...,e_n\}$, where $m=\dimc L_0=\dimr L-n$ and $0\leq m\leq n$.
We have three cases:

{\bf Case  $m=n$.} $L=\tilde E$ is the whole space;

{\bf Case $m=0$.} $L=\spar\{e_1,...,e_n\}$ is a real form of $\tilde E$;

{\bf Case $0<m<n$.} $L=\spar\{e_1,...,e_m,f_1,...,f_m\}\oplus\spar\{e_{m+1},...,e_n\}$.

{\bf Step 2.} Now
we describe Lie algebras $\f$  of connected  Lie subgroups
$F\subset \Simil \tilde E$ preserving $L$.
Without loss of generality, we may assume that the Lie group $F$ does not  preserve
any proper affine subspace of  $L$, i.e.  $F$ acts irreducibly on $L$, hence
$F$ acts transitively on $L$ (Theorem \ref{Al}). And we can describe all such groups
and their Lie algebras.

{\bf Case $m=n$.} The Lie  algebras corresponding to the transitive similarity transformation groups of type $\A^1$ and $\varphi$ are

$\f^{n,\B,\A^1}=(\A^1\oplus\B)\zr E$, where $\B\subset \u(n)$ is a subalgebra;

$\f^{n,\B,\varphi}=\{x+\varphi(x)| x\in \B\})\zr E,$ where $\varphi:\B\to \A^1$ is a linear map with $\varphi|_{\B'}=0$.

Now consider the  Lie algebras of type  $\psi$. Here we have an orthogonal decomposition $E=W\oplus U$.
Let $W_1=W\cap JW$,  $U_1=U\cap JU$ and let $W_2$ and $U_2$ be the orthogonal complements
to $W_1$ and $U_1$ in $W$ and $U$, respectively. Obviously, $W_1$ and $U_1$ are complex subspaces of $E$,
and $W_2$ and $U_2$ are mutually orthogonal real forms of the complex vector space $\spac W_2=\spac U_2$.
We obtain the orthogonal decomposition of $E$: $E=W_1\oplus U_1\oplus W_2\oplus U_2$. Let $k=\dimc W_1$ and $l=\dimc (W_1\oplus U_1)$.
We choose the basis $e_1,...,e_n,f_1,...,f_n$ in such a way that\\ $W_1=\spar\{e_1,...,e_k,f_1,...,f_k\},$
$U_1=\spar\{e_{k+1},...,e_l,f_{k+1},...,f_l\}$ and\\ $W_2=\spar\{e_{l+1},...,e_n\}$. Then
$U_2=JW_2=\spar\{f_{l+1},...,f_n\}$.

The orthogonal part of a Lie algebra of type $\psi$ is contained in $\so(W)\cap \u(n)$.
Suppose\\ $A=\bigl( \begin{smallmatrix} B&-C\\C&B\end{smallmatrix}\bigr)\in  \so(W)\cap\u(n)$, then
with respect to the decomposition
$$E=E^1_{1,...,k}\oplus   E^1_{k+1,...,l}\oplus E^1_{l+1,...,n}\oplus E^2_{1,...,k}\oplus E^2_{k+1,...,l}\oplus E^2_{l+1,...,n}$$
the matrix $A$ has the form
{\scriptsize$\left(\begin{array}{cccccc}
B_1&0&0&-C_1&0&0\\
0&0&0&0&0&0\\
0&0&0&0&0&0\\
C_1&0&0&B_1&0&0\\
0&0&0&0&0&0\\
0&0&0&0&0&0\\
\end{array}\right)$}. Consequently, $\so(W)\cap \u(n)=\u(W_1)=\u(k).$

Thus the Lie algebras corresponding to  transitively acting groups of type $\psi$ have the form
$$\f^{n,\B,\psi,k,l}=\{x+\psi(x)| x\in \B\}\zr(E^1_{1,...,k}+E^2_{1,...,k}+E^1_{l+1,...,n}),$$
where $0\leq k\leq l\leq n$, $\B\subset\u(k)$ and
 $\psi:\B\to E^1_{k+1,...,l}+E^2_{k+1,...,l}+E^2_{l+1,...,n}$ is a surjective linear map with $\psi|_{\B'}=0$.

{\bf Case $m=0$.} We have $L=L_0=\spar\{e_1,...,e_n\}=E^1$. If a subgroup $F\subset\Simil\tilde E$
preserves $L$, then $F$ is contained in $(A^1\times SO(L)\times SO(L^{\bot_\eta}))\ZR L.$
If $F$ acts transitively on $L$, the we can describe the projection of $F$ to $(A^1\times SO(L))\ZR L$,
but the projection of $F$ to $SO(L)\times SO(L^{\bot_\eta})$ is also contained in
$U(n)$. We have
$$(\so(L)\oplus\so(L^{\bot_\eta}))\cap \u(n)=\left\{\bigl(\begin{smallmatrix} B&0\\0&B\end{smallmatrix}\bigr)| B\in\so(L)\right\}=
\sod(1,...,n).$$
Hence, $(SO(L)\times SO(L^{\bot_\eta}))\cap U(n)=\SOd(1,...,n)$, where  $\SOd(1,...,n)$ is the connected Lie subgroup
of $U(n)$ corresponding to the Lie algebra $\sod(1,...,n)$. Thus the projection of $F$ to
$(A^1\times SO(L))\ZR L$ gives us the full information about $F$.

The Lie  algebras corresponding to the transitive similarity transformation groups of type $\A^1$ and $\varphi$ have the form

$\f^{0,\B,\A^1}=(\A^1\oplus\B)\zr E^1$, where $\B\subset \sod(1,...,n)$ is a subalgebra;

$\f^{0,\B,\varphi}=\{x+\varphi(x)| x\in \B\}\zr E^1,$ where $\varphi:\B\to \A^1$ is a linear map with $\varphi|_{\B'}=0$.

For a  Lie algebra of type $\psi$ we have an orthogonal decomposition $L=W\oplus U$. We choose the vectors
$e_1,...,e_n$ in such a way that $W=\spar\{e_1,...,e_k\}$ and  $U=\spar\{e_{k+1},...,e_n\}$, where $k=\dimr W$.

The Lie algebras corresponding to  transitively acting groups of type $\psi$ have the form
$$\f^{0,\B,\psi,k}=\{x+\psi(x)| x\in \B\})\zr E^1_{1,...,k},$$
where $0<k<n$, $\B\subset\sod(1,...,k)$ and
 $\psi:\B\to E^1_{k+1,...,n}$  is a surjective linear map with $\psi|_{\B'}=0$.

{\bf Case $0<m<n$.}  In this case $L=\spac\{e_1,...,e_m\}\oplus\spar\{e_{m+1},...,e_n\}.$ Hence, $L^{\bot_\eta}=\spar\{f_{m+1},...,f_n\}.$
Suppose that $$A=\bigl( \begin{smallmatrix} B&-C\\C& B\end{smallmatrix}\bigr)\in (\so(L)\oplus\so(L^{\bot_\eta}))\cap\u(n),$$ then
with respect to the decomposition $$E=E^1_{1,...,m}\oplus E^1_{m+1,...,n}\oplus E^2_{1,...,m}\oplus E^2_{m+1,...,n}$$ the element  $A$
has the form{\scriptsize
$\left(\begin{array}{cccc}
B_1&0&-C_1&0\\
0&B_2&0&0\\
C_1&0&B_1&0\\
0&0&0&B_2\\
\end{array}\right)$}.
Consequently, $(\so(L)\oplus\so(L^{\bot_\eta}))\cap\u(n)=\u(m)\oplus\sod(m+1,...,n).$

Thus, as in case m=0, the projection of an $L$-preserving subgroup  $F\subset\Simil\tilde E$  to $\Simil L$ gives us
the full information about $F$.

The Lie  algebras corresponding to the transitive similarity transformation groups of type $\A^1$ and $\varphi$ have the form

$\f^{m,\B,\A^1}=(\A^1\oplus\B)\zr(E^1_{1,...,m}+E^2_{1,...,m}+E^1_{m+1,...,n})$, where
$\B\subset \u(m)\oplus\sod(m+1,...,n)$\\ $\qquad$ is a subalgebra;

$\f^{m,\B,\varphi}=\{x+\varphi(x)| x\in \B\}\zr (E^1_{1,...,m}+E^2_{1,...,m}+E^1_{m+1,...,n}),$\\
$\qquad$ where $\varphi:\B\to \A^1$ is a linear map with $\varphi|_{\B'}=0$.

Now consider the  Lie algebras of  type $\psi$. We have an orthogonal decomposition $L=W\oplus U$.
Let $W_1=W\cap JW$,  $U_1=U\cap JU$, and let $W_2$ and $U_2$ be the orthogonal complements
to $W_1$ and $U_1$ in $W$ and $U$, respectively. We see that $W_1$ and $U_1$ are complex subspaces of
$L_0=L\cap JL=\spac\{e_1,...,e_m\}$. For $W_2$ and $U_2$ we have $W_2\cap JW_2=\{0\}$ and $U_2\cap JU_2=\{0\}$.
Denote by $L_1$ the orthogonal complement to $L_0$ in $L$ and by $L_2$ the orthogonal complement
to $W_1\oplus U_1$ in $L_0$. We get the following orthogonal decompositions:
$L=L_0\oplus L_1=W_1\oplus U_1\oplus W_2\oplus U_2,$ $L_0=W_1\oplus U_1\oplus L_2$ and $L_1\oplus L_2=W_2\oplus U_2$.
Let $k=\dimc W_1$, $l=\dimc(W_1\oplus U_1)$ and $r=\dimr W_2+l$. We see that $0\leq k\leq l\leq m\leq r\leq n$.
We choose the basis $e_1,...,e_n,f_1,...,f_n$   in such a way that
$$W_1=\spac\{e_1,...,e_k\}, U_1=\spac\{e_{k+1},...,e_l\},$$
$$W_2=\spar\{e_{l+1},...,e_m\}\oplus\spar\{e_{m+1},...,e_r\}.$$ Then
$$U_2=\spar\{f_{l+1},...,f_m\}\oplus\spar\{e_{r+1},...,e_n\}.$$

Suppose that $A=\bigl( \begin{smallmatrix} B&-C\\C& B\end{smallmatrix}\bigr)\in (\so(W)\times\so(L^{\bot_\eta}))\cap\u(n)$, then
with respect to the decomposition
$$\begin{array}{rcl}
E&=&E^1_{1,...,k}\oplus E^1_{k+1,...,l}\oplus E^1_{l+1,...,m}\oplus E^1_{m+1,...,r}\oplus E^1_{r+1,...,n}\oplus\\
&&E^2_{1,...,k}\oplus E^2_{k+1,...,l}\oplus E^2_{l+1,...,m}\oplus E^2_{m+1,...,r}\oplus E^2_{r+1,...,n}\end{array}$$
the matrix $A$ has the form{\scriptsize
$\left(\begin{array}{cccccccccc}
B_1&0&0&0&0&-C_1&0&0&0&0\\
0&0&0&0&0&0&0&0&0&0\\
0&0&0&0&0&0&0&0&0&0\\
0&0&0&B_2&0&0&0&0&0&0\\
0&0&0&0&0&0&0&0&0&0\\
C_1&0&0&0&0&B_1&0&0&0&0\\
0&0&0&0&0&0&0&0&0&0\\
0&0&0&0&0&0&0&0&0&0\\
0&0&0&0&0&0&0&0&B_2&0\\
0&0&0&0&0&0&0&0&0&0\\
\end{array}\right)$.}\\
Hence,  $(\so(W)\times\so(L^{\bot_\eta}))\cap\u(n)=\u(k)\oplus\sod(m+1,...,r).$

Thus the Lie algebras corresponding to  transitively acting groups of type $\psi$ have the form
$$\f^{m,\B,\psi,k,l,r}=\{x+\psi(x)| x\in \z(\B)\}\zr(E^1_{1,...,k}+E^2_{1,...,k}+ E^1_{l+1,...,m}+E^1_{m+1,...,r}),$$
where $0<k\leq l\leq m\leq  r\leq n$, $\B\subset\u(k)\oplus\sod(m+1,...,r)$ and\\
 $\psi:\B\to E^1_{k+1,...,l}+E^2_{k+1,...,l}+E^2_{l+1,...,m}+E^1_{r+1,...,n}$ is a surjective linear map with $\psi|_{\B'}=0$.

{\bf Step 3.} Here for each subalgebra $\f\subset \lie(\Simil\tilde E)$ considered above we describe subalgebras $\a$
of the Lie algebra $\lie(\Simil \H_n)$ with $\pi(\a)=\f$.

{\bf Case $m>0$.} Let $\f\subset \lie(\Simil\tilde E)$ be a subalgebra as above with $m>0$. We claim that
if for a subalgebra $\a\subset \lie(\Simil \H_n)$ we have  $\pi(\a)=\f$, then $\a=\f\zr\C$, where $\C=\ker\pi$.
Indeed, the projection $\pi:\lie(\Simil \H_n)\to \lie(\Simil\tilde E)$ is surjective with the kernel $\C$,
consequently, if $\pi(\a)=\f$ and $\a$ does not contain $\C$, then $\a$ has the form $\{x+\zeta(x)|x\in\f\}$ for
some linear map $\zeta:\f\to\C$. Suppose that $\f\neq\f^{m,4,*}$, then we choose $z_1,z_2\in\f\cap E$ with $\im g(z_1,z_2)\neq 0$.
This yields that $[z_1+\zeta(z_1),z_2+\zeta(z_2)]=-2\im g(z_1,z_2)\in\C$. If $\f=\f^{m,4,*}$, then we can choose $z\in\f\cap E$ and
$x\in\z(\B)$ with $\im g(z,\psi(x))\neq 0$, or $x_1,x_2\in\z(\B)$ with $\im g(\psi(x_1),\psi(x_2)\neq 0$.
For each subalgebra $\f\subset \lie(\Simil\tilde E)$ considered above with $m>0$  we define the
Lie subalgebra $$\a^{m,*}=\f^{m,*}\zr\C\subset\lie(\Simil\H_n).$$

{\bf Case m=0.} For each subalgebra $\f\subset \lie(\Simil\tilde E)$ considered above with $m=0$  we have the following two possibilities:

{\bf Subcase 1.} Here we have $\a^{0,*}=\f^{0,*}\zr\C$.

{\bf Subcase 2.} The Lie algebra $\a$ does not contain $\C$. Hence $\a$ has the form $\{x+\zeta(x)|x\in\f\}$ for
some linear map $\zeta:\f\to\C$. For each $\f^{0,*}$ we will find all possible $\zeta$. Since $\a$ is
a Lie algebra, we see that $\zeta$ vanishes on the commutator $\f'$.

{\bf Subcase 2.1.} Consider $\f^{0,\B,\A^1}=(\A^1\oplus\B)\zr E^1$, where $\B\subset \sod(1,...,n)$ is a subalgebra.
Since $(\f^{0,1,\B})'=E^1+\B'$, we obtain $\zeta|_{E^1+\B'}=0$.
Let $a\in \A^1$ and $x\in\z(\B)$. From $[a+\zeta(a),x+\zeta(x)]=a\zeta(x)\in\C$, it follows that $\zeta|_{\z(\B)}=0$.
Thus $\zeta$ can be considered as a linear map $\zeta:\A^1\to\C$.
For any linear map $\zeta:\A^1\to\C$ we consider the Lie algebra
$$\a^{0,\B,\A^1,\zeta}=(\B\oplus\{a+\zeta(a)|a\in \A^1\})\zr E^1.$$

{\bf Subcase 2.2.} Consider $\f^{0,\B,\varphi}$ with $\varphi=0$, i.e.  $\f^{0,\B,0}=\B\zr E^1$, where $\B\subset \sod(1,...,n)$ is a subalgebra.
We have $(\f^{0,2,\B})'=\B'+\spa\{x(u)|x\in\B,u\in E^1\}$.
Choose the vectors $e_1,...,e_n$ so that $E^1_{i_0+1,...,n}$ is the subspace of $E^1$ annihilated
by $\B$ and $E^1_{1,...,i_0}$ is the orthogonal complement to $E^1_{i_0+1,...,n}$ in $E^1$.
The Lie algebra $\B$ is compact, hence $\B$ is totally reducible and $E^1_{1,...,i_0}$ is decomposed into an orthogonal
sum of subspaces, on each of these subspaces $\B$ acts irreducibly. Thus, $\spa\{x(u)|x\in\B,u\in E^1\}=E^1_{1,...,i_0}$,
and $\zeta$ can be considered as a linear map $\zeta:\z(\B)\oplus E^1_{i_0+1,...,n}\to\C$.
For any linear map $\zeta:\B\oplus E^1_{i_0+1,...,n}\to\C$ with $\zeta|_{\B'}=0$ we define the Lie algebra
$$\a^{0,\B,\varphi=0,i_0,\zeta}=\{x+\zeta(x)|x\in\B\oplus E^1_{i_0+1,...,n}\}\zr E^1_{1,...,i_0}.$$

{\bf Subcase 2.3.} Consider $\f^{0,\B,\varphi}=\{x+\varphi(x)| x\in \B\})\zr E^1,$ where $\B\subset \sod(1,...,n)$ and
$\varphi:\B\to \A^1$ is a non-zero linear map with $\varphi|_{\B'}=0$.
As above, we can show that $(\f^{0,\B,\varphi})'=E^1+\B'$. Hence $\zeta$ is a map
$\zeta:\{x+\varphi(x)| x\in \z(\B)\}\to\C$. Denote by the same letter $\zeta$ the linear map
$\zeta:\z(\B)\to\C$ defined by $\zeta(x)=\zeta(x+\varphi(x))$.
Let $x_1,x_2\in \z(\B)$ we have $[x_1+\varphi(x_1)+\zeta(x_1),x_2+\varphi(x_2)+\zeta(x_2)]=\varphi(x_1)\zeta(x_2)-
\varphi(x_2)\zeta(x_1)\in\C$. Hence, $\varphi(x_1)\zeta(x_2)=\varphi(x_2)\zeta(x_1)$. In particular, if $\varphi(x)=0$,
then $\zeta(x)=0$. Hence, $\ker\varphi\subset\ker\zeta$. Conversely, if $\ker\varphi\subset\ker\zeta$, then
$\varphi(x_1)\zeta(x_2)=\varphi(x_2)\zeta(x_1)$ for all $x_1,x_2\in \z(\B)$.

Let $i_0$ be as in Subcase 2.2.
For any linear map $\zeta:\B\oplus E^1_{i_0+1,...,n}\to\C$ such that $\ker\varphi\subset\ker\zeta$ and $\zeta|_{E^1_{i_0+1,...,n}}=0$ we consider the Lie algebra
$$\a^{0,\B,\varphi,i_0,\zeta}=\{x+\varphi(x)+\zeta(x)|x\in\B\oplus E^1_{i_0+1,...,n}\}\zr E^1_{1,...,i_0}.$$
Note that the Lie algebras $\a^{0,\B,\varphi,i_0,\zeta}$  for  $\varphi=0$ and $\varphi\neq 0$ are defined in the same way.

{\bf Subcase 2.4.} As above, for the Lie algebra
$$\f^{0,\B,\psi,k}=\{x+\psi(x)| x\in \B\}\zr E^1_{1,...,k},$$
where $0<k<n$, $\B\subset\sod(1,...,k)$ and  $\psi:\B\to E^1_{k+1,...,n}$
is a surjective linear map with $\psi|_{\B'}=0$, we can prove that the map $\zeta$ vanishes on $\B'+E^1_{1,...,i_0},$
where $i_0$ is as in Subcase 2.2. For any linear map
$\zeta:\B\oplus E^1_{i_0+1,...,k}\to\C$ with $\zeta|_{\B'}=0$ we consider the Lie algebra
$$\f^{0,\B,\psi,,k,i_0,\zeta}=\{x+\psi(x)+\zeta(x)| x\in \B\}\oplus\{u+\zeta(u)| u\in E^1_{i_0+1,...,k}\})\zr E^1_{1,...,i_0}.$$

{\bf Step 4.} For each subalgebra $\a\subset \lie(\Simil\H_n)$ constructed above consider the subalgebra
$\Gamma_0^{-1}(\a)\subset\su(1,n+1)_{<p_1,p_2>}$. Note that we have $\Gamma_0(\N^1_{k,...,l})=E^1_{k,...,l}$,
$\Gamma_0(\N^2_{k,...,l})=E^2_{k,...,l}$, $\Gamma_0(\A^1)=\A^1$ and $\Gamma_0(\C)=\C$, where we consider
$\A^1$ and $\C$  as the  subalgebras of $\su(1,n+1)_{<p_1,p_2>}$ as well as of $\lie(\Simil\H_n)$.
Let $\a^{m,\B}\subset \lie(\Simil\H_n)$ be any subalgebra constructed above with the associated number $0\leq m\leq n$ and the associated subalgebra
$\B\subset\u(m)\oplus\sod(m+1,..,n)=\su(m)\oplus\Real J_m\oplus\sod(m+1,...,n)$. We have $\Gamma|_{\su(n)}=\id_{\su(n)}$, where
$\su(n)$ is considered as the  subalgebras of $\su(1,n+1)_{<p_1,p_2>}$ and of $\lie(\Simil\H_n)$.
Furthermore, since $J_m-\frac{m}{n}J_n\in\su(n)$, we have
$\Gamma_0^{-1}(J_m)=\Gamma_0^{-1}(J_m-\frac{m}{n}J_n)+\frac{m}{n}\Gamma_0^{-1}(J_n)=
J_m-\frac{m}{n}J_n+\frac{m}{n}I_0=J_m-\frac{m}{n}J_n+\frac{m}{n}(0,-\frac{n}{n+2},0,\frac{2}{n+2}E_n,0,0,0)=J_m-\frac{m}{n+2}J.$
Let $\h=\pr_{\u(n)}\Gamma_0^{-1}(\a)=\pr_{\u(n)}\Gamma_0^{-1}(\B)\subset\su(m)\oplus\Real(J_m-\frac{m}{n+2}J_n)$.
Thus, $$\Gamma_0^{-1}(\B)=\{(0,-\frac{1}{2}\tr C,B,C,0,0,0)|\bigl(\begin{smallmatrix}B&-C\\C&B\end{smallmatrix}\bigr)\in\h\}.$$

For the Lie algebra $\a^{m,\h,\psi,k,l,r}$ with $0<m<n$ we have $\B\subset\su(k)\oplus\Real J_k\oplus\sod(m+1,...,r)$,
consequently, $\h\subset\su(k)\oplus\Real (J_k-\frac{k}{n+2}J_n)\oplus\sod(m+1,...,r).$
For the Lie algebra $\a^{n,\h,\psi,k,l}$ we have $\B\subset\su(k)\oplus\Real J_k$,
hence, $\h\subset\su(k)\oplus\Real (J_k-\frac{k}{n+2}J_n).$
For the Lie algebra $\a^{0,\h,\psi,k}$ we have $\B\subset\sod(1,...,k)$,
hence $\h\subset\sod(1,...,k).$ We denote $\Gamma_0^{-1}(\a^{m,\B,*})$ by $\g^{m,\B,*}$.

Thus we obtain a list of subalgebras of $\su(1,n+1)_{<p_1,p_2>}$.
From proposition \ref{prop1} and Theorem \ref{non-compl} it follows that this list contains all weakly-irreducible
subalgebras of $\su(1,n+1)_{<p_1,p_2>}$. Example \ref{exWI2} shows that this list contains also
subalgebras of $\su(1,n+1)_{<p_1,p_2>}$ that are not weakly-irreducible.
Here we verify which of the subalgebras $\Gamma_0^{-1}(\a)\subset\su(1,n+1)_{<p_1,p_2>}$ are weakly-irreducible.

The subalgebras of the form $\g^{m,\h,\A^1}$ and $\g^{m,\h,\varphi}$,  where $0\leq m\leq n$ contain $\N^1+\C$, hence these
subalgebras are weakly-irreducible.

\begin{lem}\label{lem1} The subalgebras of the form $\g^{0,\h,\psi,k}$, $\g^{m,\h,\psi,k,l,r}$ and $\g^{n,\h,\psi,k,l}$ are weakly-irreducible. \end{lem}

{\it Proof.} Let $\g$ be any of these subalgebras and suppose $\g$ preserves a proper vector  subspace $L\subset \Real^{2,2n+2}$.
As in Example \ref{exWI1}, we can show that if $(a_1,a_2,\alpha,\beta,b_1,b_2)\in L$, then we have $-b_2p_1+b_1p_2\in L$
and  $b_1p_1+b_2p_2\in L$. We suppose that $L\subset\Real p_1\oplus \Real p_2\oplus E$.
Let $v=(a_1,a_2,\alpha,\beta,0,0)\in L$. We assume $\alpha\neq 0$. If the projection of $\alpha$ to $W$ is non-trivial,
then we denote this projection by $w$. Applying the element $(0,0,0,0,w,0,0)\in\g$ to $v$, we obtain
$w^t\alpha p_1+w^t\beta p_2\in L$. Consequently  $L$ is degenerate subspace. Suppose that the projections of $\alpha$ and $\beta$ to
$W$ are trivial and the projection of $\alpha$ to $U$ is non-trivial, then there exists an element
$x\in\z(\h)$ such that $\psi(x)=(0,0,0,0,u,0,0),$ where  $u\in U$ is equal to this projection.
The element $x\in\z(\h)$ has the form $x_1+\mu(J_k-\frac{k}{n+2}J_n)$, where
$x_1\in\su(k)\oplus\sod(m+1,...,r)$ (resp. $x_1\in\su(k)$ and $x_1\in\sod(1,...,k)$) for the Lie algebra $\g^{m,\h,\psi,k,l,r}$
(resp. $g^{0,\h,\psi,k}$ and $\g^{n,\h,\psi,k,l}$),  $\mu=0$ or $1$ for the Lie algebras $\g^{m,\h,\psi,k,l,r}$ and $\g^{n,\h,\psi,k,l}$,
and $\mu=0$ for the Lie algebra $\g^{0,\h,\psi,k}$.
Consider the element $X=(0,-\frac{\mu k}{n+2},B_1,C_1+\mu(J_k-\frac{k}{n+2}J_n),u,0,0,0)\in\g,$
where $\bigl(\begin{smallmatrix}B_1&-C_1\\C_1&B_1\end{smallmatrix}\bigr)=x_1.$
Applying the element $X$ to $v$ we get\\
$Xv=(\frac{\mu k}{n+2}a_2+u^t\alpha, -\frac{\mu k}{n+2}a_2+u^t\beta,\frac{\mu k}{n+2}\beta,-\frac{\mu k}{n+2}\alpha,0,0)\in L,\quad\quad\quad
X^2v=$\\
$(\frac{\mu k}{n+2}(-\frac{\mu k}{n+2}a_1+u^t\beta)+\frac{\mu k}{n+2}u^t\beta,
-\frac{\mu k}{n+2}(\frac{\mu k}{n+2}a_2+u^t\alpha)-\frac{\mu k}{n+2}u^t\alpha,
-(\frac{\mu k}{n+2})^2\alpha,-(\frac{\mu k}{n+2})^2\beta,0,0)\in L.$
Hence if $\mu=1$, then $v+(\frac{\mu k}{n+2})^2X^2v=(2\frac{\mu k}{n+2}u^t\beta,2\frac{\mu k}{n+2}u^t\alpha,0,0,0,0)\in L$.
If $\mu=0$, then $Xv=(u^t\alpha,u^t\beta,0,0,0,0)\in L$.

Thus the subspace  $L\subset \Real^{2,2n+2}$ is degenerate  and the lemma is proved. $\Box$

We are left now with the Lie algebras of the form $\g^{0,\h,\A^1,\zeta}$, $\g^{0,\h,\varphi,i_0,\zeta}$ and
$\g^{0,\h,\psi,k,i_0,\zeta}$.

\begin{lem}\label{lt1}
The subalgebra $\g^{0,\h,\A^1,\zeta}\subset\su(1,n+1)_{<p_1,p_2>}$ is not weakly-irreducible.
\end{lem} {\it Proof.} The Lie algebra $\g^{0,\h,\A^1,\zeta}$ preserves the non-degenerate proper vector subspace
$\spar\{p_1+p_2,e_1+f_1,...,e_n+f_n,q_1-\frac{\zeta}{2}p_2+q_2+\frac{\zeta}{2}p_1\}\subset\Real^{2,2n+2}$. $\Box$

\begin{lem}\label{lt2} The subalgebra $\g^{0,\h,\varphi,i_0,\zeta}\subset\su(1,n+1)_{<p_1,p_2>}$ is weakly-irreducible
if and only if $\zeta|_{\z(\h)}\neq 0$ and $\varphi=0$. \end{lem}

{\it Proof.} If $\varphi\neq 0$, then $\g^{0,\h,\varphi,i_0,\zeta}$ preserves
 the non-degenerate proper vector subspace
$\spar\{p_1+p_2,e_1+f_1,...,e_n+f_n,q_1-\frac{\zeta(A)}{2\varphi(A)}p_2+q_2+\frac{\zeta(A)}{2\varphi(A)}p_1\}\subset
\Real^{2,2n+2}$, where $A\in\z(\h)$ is a non-zero element that is orthogonal to $\ker\varphi$.

Now assume that  $\varphi=0$. Suppose that $\zeta|_{\z(\h)}=0$. If $\zeta=0$, then $\g^{0,\h,\varphi,i_0,\zeta}$ preserves
the non-degenerate vector subspaces $\spar\{p_1+p_2,e_1+f_1,...,e_n+f_n,q_1+q_2\}\subset \Real^{2,2n+2}$ and
$\spar\{p_1-p_2,e_1-f_1,...,e_n-f_n,q_1-q_2\}\subset \Real^{2,2n+2}$. If $\zeta\neq 0$, then we choose the vectors
$e_{i_0+1},...,e_{n}$ so that $\zeta|_{\spa\{e_{i_0+1},...,e_{n-1}\}}=0$. Then $\g^{0,2,\h,i_0,\zeta}$ preserves the
non-degenerate vector subspace $\spar\{p_1+p_2,e_1+f_1,...,e_n+f_n,q_1+q_2-2\zeta(e_n)e_n\}\subset \Real^{2,2n+2}$.

Suppose that $\zeta|_{\z(\h)}\neq 0$. Let $x=\bigl(\begin{smallmatrix} B&0\\0&B\end{smallmatrix}\bigr)\in\z(\h)$ with $\zeta(x)\neq 0$.
Since $B\in\so(n)$, we can choose the basis $e_1,...,e_n$ so that $B$ has the form{\scriptsize
$\left(\begin{array}{cccc}
\begin{smallmatrix}0&-\lambda_1\\ \lambda_1&0\end{smallmatrix}& & &0\\
&\ddots& & \\
&  &  \begin{smallmatrix}0&-\lambda_s\\ \lambda_s&0\end{smallmatrix}&\\
0& & & \begin{smallmatrix}0& &\\&\ddots&\\&&0\end{smallmatrix}\end{array}\right),$} where $2s\leq i_0$ and $\lambda_i\neq 0$.

Suppose that $\g^{0,2,\h,i_0,\zeta}$ preserves a proper vector subspace $L\subset \Real^{2,2n+2}$.
Let $v=(a_1,a_2,\alpha,\beta,b_1,b_2)\in L$ and $X=x+\zeta(x)\in\g^{0,\h,\varphi,i_0,\zeta}$.
We have\\ {\scriptsize
$$Xv=\left(\begin{array}{c}-\zeta(x)b_2\\ \zeta(x)b_1\\
-\lambda_{1} \alpha_{2}\\ \lambda_{1} \alpha_{1}\\ \vdots \\-\lambda_{s} \alpha_{s}\\\lambda_{s} \alpha_{s-1}\\
0\\ \vdots\\ 0\\ -\lambda_{1} \beta_{2}\\ \lambda_{1} \beta_{1}\\ \vdots \\-\lambda_{s} \beta_{s}\\
\lambda_{s} \beta_{s-1}\\0\\ \vdots\\ 0\\0\\0\end{array}\right)\in L,\quad
X^3v=-\left(\begin{array}{c}0\\ 0\\
\lambda^2_1 (-\lambda_{1} \alpha_{2})\\ \lambda^2_1(\lambda_{1} \alpha_{1})\\ \vdots \\\lambda^2_s(-\lambda_{s} \alpha_{s})\\ \lambda^2_s(\lambda_{s} \alpha_{s-1})\\
0\\ \vdots\\ 0\\ \lambda^2_1(-\lambda_{1} \beta_{2})\\ \lambda^2_1(\lambda_{1} \beta_{1})\\ \vdots \\\lambda^{2s}_s(-\lambda_{s} \beta_{s})\\
\lambda^2_s(\lambda_{s} \beta_{s-1})\\0\\ \vdots\\ 0\\0\\0\end{array}\right)\in L,
\quad
\cdots,
\quad X^{2s+1}v=-\left(\begin{array}{c}0\\ 0\\
\lambda^{2s}_1 (-\lambda_{1} \alpha_{2})\\ \lambda^{2s}_1(\lambda_{1} \alpha_{1})\\ \vdots \\\lambda^{2s}_s(-\lambda_{s} \alpha_{s})\\
\lambda^{2s}_s(\lambda_{s} \alpha_{s-1})\\
0\\ \vdots\\ 0\\ \lambda^{2s}_1(-\lambda_{1} \beta_{2})\\ \lambda^{2s}_1(\lambda_{1} \beta_{1})\\ \vdots \\\lambda^{2s}_s(-\lambda_{s} \beta_{s})\\
\lambda^{2s}_s(\lambda_{s} \beta_{s-1})\\0\\ \vdots\\ 0\\0\\0\end{array}\right)\in L.$$}

The vector $-\zeta(x)b_2p_1+\zeta(x)b_1p_2$ can by decomposed into a combination of the vectors $Xv$, $X^3v$,...,$X^{2s+1}v\in L$,
hence $-b_2p_1+b_1p_2\in L$.
The end of the proof of the lemma is as in Example \ref{exWI1}. $\Box$

\begin{lem}\label{lt4} The subalgebra $\g^{0,\h,\psi,k,i_0,\zeta}\subset\su(1,n+1)_{<p_1,p_2>}$ is weakly-irreducible
if and only if $\zeta|_{\z(\h)}\neq 0$. \end{lem}

The proof is similar to the proofs of Lemma \ref{lem1} and \ref{lt2}. $\Box$

Now for the Lie algebra $\g^{0,\h,\psi,k,i_0,\zeta}$ we assume that $\zeta|_{\z(\h)}\neq 0$. Obviously, this Lie algebra
has the form $\g^{0,\h,\psi',k',i_0,\zeta'}$, where $\zeta'|_{E^1_{i_0+1,...,k'}}=0$. We will denote this Lie algebra by
$\g^{0,\h,\psi',k,\zeta'}$. Consider the Lie algebra $\g^{0,\h,\varphi,i_0,\zeta}$ such that
$\zeta|_{E^1_{i_0+1,...,n}}\neq 0$, $\zeta|_{\z(\h)}\neq 0$ and $\varphi=0$. Obviously, this is the Lie algebra of the
form $\g^{0,\h,\psi,k,\zeta}$. Thus we may restrict our attention to the Lie algebras $\g^{0,\h,\varphi,i_0,\zeta}$ such
that $\zeta|_{E^1_{i_0+1,...,n}}=0$, $\zeta|_{\z(\h)}\neq 0$ and $\varphi=0$. We denote such Lie algebra by
$\g^{0,\h,\zeta}$.

The theorem is proved. $\Box$

\part{Classification of the holonomy algebras and constructions of metrics}
In this chapter we give the classification of the holonomy algebras of pseudo-K\"ahlerian manifolds of index 2.
For each weakly-irreducible not irreducible holonomy algebra $\g$ we construct a polynomial metric with the holonomy algebra $\g$.
The results are stated in Section \ref{Ber+metr}. In the other two sections we give the proofs of the theorems.

\section{Main results}\label{Ber+metr}


In the following theorem we give the classification of the weakly-irreducible not irreducible
holonomy algebras of pseudo-K\"ahlerian manifolds of index 2.
We use the denotation from Section \ref{w-ir}.

\begin{theorem}\label{hol2n} {\bf 1)} A subalgebra $\g\subset\u(1,1)$ is the weakly-irreducible not irreducible holonomy algebra
of a pseudo-K\"ahlerian manifold of signature $(2,2)$ if and only if $\g$ is conjugated to one of the following
subalgebras of $\u(1,1)_{<p_1,p_2>}$:
\begin{description}

\item[$\hol^1_{n=0}$]$=\u(1,1)_{<p_1,p_2>}$;

\item[$\hol^2_{n=0}$]$=\A^1\oplus\A^2$
{\scriptsize$=\left\{\left. \left (\begin{array}{cccc}
a_1&-a_2  &0 &0\\
a_2&a_1  &0 &0\\
0&0&-a_1&-a_2\\
0&0&a_2&-a_1\\
\end{array}\right)\right|\, a_1,a_2\in \Real \right\}$};

\item[$\hol^{\gamma_1,\gamma_2}_{n=0}$]$=\{(a\gamma_1 ,a\gamma_2 ,0)|a\in\Real\}\zr\C$
{\scriptsize
$=\left\{\left. \left (\begin{array}{cccc}
 a \gamma_1&-a\gamma_2  &0 &-c\\
a\gamma_2&a\gamma_1  &c &0\\
0&0&-a\gamma_1&-a\gamma_2\\
0&0&a\gamma_2&-a\gamma_1\\
\end{array}\right)\right|\, a,c\in \Real \right\}$,}
where $\gamma_1,\gamma_2\in~\Real$.

\end{description}

{\bf 2)} Let $n\geq 1$. Then a subalgebra $\g\subset\u(1,n+1)$ is the weakly-irreducible not irreducible holonomy algebra
of a pseudo-K\"ahlerian manifold of signature $(2,2n+2)$ if and only if $\g$ is conjugated to one of the following
subalgebras of $\u(1,n+1)_{<p_1,p_2>}$:
\begin{description}

\item[$\hol^{m,\un,\A^1,\t\A^2}$]$=(\A^1\oplus \t\A^2\oplus\un)\zr(\N^1+\N^2_{1,...,m}+\C)$
{\scriptsize
$$=\left\{\left. \left
(\begin{array}{cccccccc}
a_1&-a_2 &-z^t_1&-z'^t_1 & -z^t_2&0 &0 &-c\\
a_2&a_1 &z^t_2&0 & -z^t_1&-z'^t_1 &c &0\\
0 &0&B&0&-C&0&z_1&-z_2\\
0&0&0&0&0&-a_2E_{n-m}&z'_1&0\\
0 &0&C&0&B&0&z_2&z_1\\
0&0&0&a_2E_{n-m}&0&0&0&z'_1\\
0&0&0&0&0&0&-a_1&-a_2\\
0&0&0&0&0&0&a_2&-a_1\\
\end{array}\right)
\right|\,
\begin{array}{c}a_1,a_2,c\in \Real,\\
 z_1, z_2\in \Real^m,\\
z_1'\in\Real^{n-m},\\
\bigl(\begin{smallmatrix}B&-C\\C&B\end{smallmatrix}\bigr)\in\un \end{array} \right\},$$}
where  $0\leq m\leq n$ and $\un\subset\u(m)$ is a subalgebra;

\item[$\hol^{m,\un,\A^1,\phi}$]$=\{(a_1,\phi(B,C),B,C,0,0,0)+\phi(B,C)J_{m+1,...,n}|a_1\in\Real,\bigl(\begin{smallmatrix}B&-C\\C&B\end{smallmatrix}\bigr)\in\un\}
\zr(\N^1+\N^2_{m+1,...,n}+\C)$
{\scriptsize
$$=\left\{\left. \left
(\begin{array}{cccccccc}
a_1&-\phi(A) &-z^t_1&-z'^t_1 & -z^t_2&0 &0 &-c\\
\phi(A)&a_1 &z^t_2&0 & -z^t_1&-z'^t_1 &c &0\\
0 &0&B&0&-C&0&z_1&-z_2\\
0&0&0&0&0&-\phi(A)E_{n-m}&z'_1&0\\
0 &0&C&0&B&0&z_2&z_1\\
0&0&0&\phi(A)E_{n-m}&0&0&0&z'_1\\
0&0&0&0&0&0&-a_1&-\phi(A)\\
0&0&0&0&0&0&\phi(A)&-a_1\\
\end{array}\right)
\right|\,
\begin{array}{c}a_1,c\in \Real,\\
 z_1, z_2\in \Real^m,\\
z_1'\in\Real^{n-m},\\ A=\\
\bigl(\begin{smallmatrix}B&-C\\C&B\end{smallmatrix}\bigr)\in\un \end{array} \right\},$$}
where  $0\leq m\leq n$, $\un\subset\u(m)$ is a subalgebra
and $\phi:\un\to\Real$ is a linear map with $\phi|_{\un'}=0$;

\item[$\hol^{m,\un,\varphi,\phi}$]$=\{(\varphi(B,C),\phi(B,C),B,C,0,0,0)+\phi(B,C)J_{m+1,...,n}|\bigl(\begin{smallmatrix}B&-C\\C&B\end{smallmatrix}\bigr)\in\un\}
\zr(\N^1+\N^2_{m+1,...,n}+\C)$
{\scriptsize
$$=\left\{\left. \left
(\begin{array}{cccccccc}
\varphi(A)&-\phi(A) &-z^t_1&-z'^t_1 & -z^t_2&0 &0 &-c\\
\phi(A)&\varphi(A) &z^t_2&0 & -z^t_1&-z'^t_1 &c &0\\
0 &0&B&0&-C&0&z_1&-z_2\\
0&0&0&0&0&-\phi(A)E_{n-m}&z'_1&0\\
0 &0&C&0&B&0&z_2&z_1\\
0&0&0&\phi(A)E_{n-m}&0&0&0&z'_1\\
0&0&0&0&0&0&-\varphi(A)&-\phi(A)\\
0&0&0&0&0&0&\phi(A)&-\varphi(A)\\
\end{array}\right)
\right|\,
\begin{array}{c}c\in \Real,\\
 z_1, z_2\in \Real^m,\\
z_1'\in\Real^{n-m},\\ A=\\
\bigl(\begin{smallmatrix}B&-C\\C&B\end{smallmatrix}\bigr)\in\un \end{array} \right\},$$}
where  $0\leq m\leq n$, $\un\subset\u(m)$ is a subalgebra
and $\varphi,\phi:\un\to\Real$ are  linear maps with $\varphi|_{\un'}=\phi|_{\un'}=0$;

\item[$\hol^{m,\un,\varphi\t\A^2}$]$=\{(\varphi(B,C),a_2,B,C,0,0,0)+a_2J_{m+1,...,n}|a_2\in\Real,\bigl(\begin{smallmatrix}B&-C\\C&B\end{smallmatrix}\bigr)\in\un\}
\zr(\N^1+\N^2_{m+1,...,n}+\C)$
{\scriptsize
$$=\left\{\left. \left
(\begin{array}{cccccccc}
\varphi(A)&-a_2 &-z^t_1&-z'^t_1 & -z^t_2&0 &0 &-c\\
a_2&\varphi(A) &z^t_2&0 & -z^t_1&-z'^t_1 &c &0\\
0 &0&B&0&-C&0&z_1&-z_2\\
0&0&0&0&0&-a_2E_{n-m}&z'_1&0\\
0 &0&C&0&B&0&z_2&z_1\\
0&0&0&a_2E_{n-m}&0&0&0&z'_1\\
0&0&0&0&0&0&-\varphi(A)&-a_2\\
0&0&0&0&0&0&a_2&-\varphi(A)\\
\end{array}\right)
\right|\,
\begin{array}{c}a_2,c\in \Real,\\
 z_1, z_2\in \Real^m,\\
z_1'\in\Real^{n-m},\\
A=\\
\bigl(\begin{smallmatrix}B&-C\\C&B\end{smallmatrix}\bigr)\in\un \end{array} \right\},$$}
where  $0\leq m\leq n$, $\un\subset\u(m)$ is a subalgebra
and $\varphi:\un\to\Real$ is a linear map with $\varphi|_{\un'}=0$;

\item[$\hol^{m,\un,\lambda}$]$=(\{(a_1,\lambda a_1,0,0,0,0,0)+\lambda a_1J_{m+1,...,n}|a_1\in\Real\}\oplus\un)\zr(\N^1+\N^2_{m+1,...,n}+\C)$
{\scriptsize
$$=\left\{\left. \left
(\begin{array}{cccccccc}
a_1&-\lambda a_1 &-z^t_1&-z'^t_1 & -z^t_2&0 &0 &-c\\
\lambda a_1&a_1 &z^t_2&0 & -z^t_1&-z'^t_1 &c &0\\
0 &0&B&0&-C&0&z_1&-z_2\\
0&0&0&0&0&-\lambda a_1E_{n-m}&z'_1&0\\
0 &0&C&0&B&0&z_2&z_1\\
0&0&0&\lambda a_1E_{n-m}&0&0&0&z'_1\\
0&0&0&0&0&0&-a_1&-\lambda a_1\\
0&0&0&0&0&0&\lambda a_1&-a_1\\
\end{array}\right)
\right|\,
\begin{array}{c}a_1,c\in \Real,\\
 z_1, z_2\in \Real^m,\\
z_1'\in\Real^{n-m},\\
\bigl(\begin{smallmatrix}B&-C\\C&B\end{smallmatrix}\bigr)\in\un \end{array} \right\},$$}
where  $0\leq m\leq n$, $\un\subset\u(m)$ is a subalgebra and $\lambda\in\Real$, $\lambda\neq 0$;

\item[$\begin{array}{rl}\hol^{n,\un,\psi,k,l}=&\{(0,0,B,C,\psi_1(B,C),\psi_2(B,C)+\psi_3(B,C),0)|\ \bigl(\begin{smallmatrix}B&-C\\C&B\end{smallmatrix}\bigr)\in\un\}\\
&\zr(\N^1_{1,...,k}+\N^2_{1,...,k}+\N^1_{l+1,...,n}+\C)\end{array}$]
{\scriptsize
$$=\left\{\left. \left (\begin{array}{cccccccccc}
0&0 &-z^t_1&-\psi_1(A)^t&-z'^t_1 & -z^t_2&-\psi_2(A)^t&-\psi_3(A)^t &0 &-c\\
0&0 &z^t_2&\psi_2(A)^t&\psi_3(A)^t & -z^t_1&-\psi_1(A)^t&-z'^t_1 &c &0\\
0 &0&B&0&0&-C&0&0&z_1&-z_2\\
0 &0&0&0&0&0&0&0&\psi_1(A)&-\psi_2(A)\\
0&0&0&0&0&0&0&0&z'_1&-\psi_3(A)\\
0&0&C&0&0&B&0&&z_2&z_1\\
0&0&0&0&0&0&0&0&\psi_2(A)&\psi_1(A)\\
0&0&0&0&0&0&0&0&\psi_3(A)&z'_1\\
0&0&0&0&0&0&0&0&0&0\\
0&0&0&0&0&0&0&0&0&0\\
\end{array}\right)\right|\,
\begin{array}{c}c\in \Real,\\
 z_1, z_2\in \Real^k,\\
 z'_1\in\Real^{n-l},\\
A=\\
 \bigl(\begin{smallmatrix}B&-C\\C&B\end{smallmatrix}\bigr)\\
\in \un\end{array} \right\},$$}
where $0< k\leq l\leq n$,  $\un\subset\u(k)$ is a subalgebra
such that $\dim\z(\un)\geq n+l-2k$,  $\psi:\un\to E^1_{k+1,...,l}\oplus E^2_{k+1,...,l}\oplus E^2_{l+1,...,n}$
is a surjective linear map with $\psi|_{\un'}=0$, $\psi_1=\pr_{E^1_{k+1,...,l}}\circ\psi$, $\psi_2=\pr_{E^2_{k+1,...,l}}\circ\psi$ and
$\psi_3=\pr_{E^2_{l+1,...,n}}\circ\psi$;

\item[$\begin{array}{rl}\hol^{m,\un,\psi,k,l,r}=&
\{(0,0,B,C,\psi_1(B,C)+\psi_4(B,C),\psi_2(B,C)+\psi_3(B,C),0)|\ \bigl(\begin{smallmatrix}B&-C\\C&B\end{smallmatrix}\bigr)\in\un\}\\
&\zr(\N^1_{1,...,k}+\N^2_{1,...,k}+ \N^1_{l+1,...,m}+\N^1_{m+1,...,r}+\C)\end{array}$]
\end{description}
{\scriptsize
$$=\left\{\left. \left (\begin{smallmatrix}
0&0 &-z^t_1&-\psi_1(A)^t&-z'^t_1&-z''^t_1&-\psi_4(A)^t & -z^t_2&-\psi_2(A)^t&-\psi_3(A)^t&0&0 &0 &-c\\
0&0 &z^t_2&\psi_2(A)^t&\psi_3(A)^t&0&0 & -z^t_1&-\psi_1(A)^t&-z'^t_1&-z''^t_1&-\psi_4(A)^t &c &0\\
0&0&B&0&0&0&0&-C&0&0&0&0&z_1&-z_2\\
0&0&0&0&0&0&0&0&0&0&0&\psi_1(A)&-\psi_2(A)\\
0&0&0&0&0&0&0&0&0&0&0&0&z'_1&-\psi_3(A)\\
0&0&0&0&0&0&0&0&0&0&0&0&z''_1&0\\
0&0&0&0&0&0&0&0&0&0&0&0&\psi_4(A)&0\\
0&0&C&0&0&0&0&B&0&0&0&0&z_2&z_1\\
0&0&0&0&0&0&0&0&0&0&0&0&\psi_2(A)&\psi_1(A)\\
0&0&0&0&0&0&0&0&0&0&0&0&\psi_3(A)&z'_1\\
0&0&0&0&0&0&0&0&0&0&0&0&0&z''_1\\
0&0&0&0&0&0&0&0&0&0&0&0&0&\psi_4(A)\\
0&0&0&0&0&0&0&0&0&0&0&0&0&0\\
0&0&0&0&0&0&0&0&0&0&0&0&0&0\\
\end{smallmatrix}\right)\right|\,
\begin{array}{c}c\in \Real,\\
 z_1, z_2\in \Real^k,\\
 z'_1\in\Real^{m-l},\\
z''_1\in\Real^{r-m},\\
\\
A=\\\ \bigl(\begin{smallmatrix}B&-C\\C&B\end{smallmatrix}\bigr)\\\in \un\end{array} \right\},$$}
\begin{description}
\item[] $ $
where $0< k\leq l\leq m\leq r\leq n$, $0<m<n$,  $\un\subset\u(k)$ is a subalgebra
such that $\dim\z(\un)\geq(n+m+l-2k-r)$,  $\psi:\un\to E^1_{k+1,...,l}\oplus E^2_{k+1,...,l}\oplus E^2_{l+1,...,m}\oplus E^1_{r+1,...,n}$
is a surjective linear map with $\psi|_{\un'}=0$, $\psi_1=\pr_{E^1_{k+1,...,l}}\circ\psi$, $\psi_2=\pr_{E^2_{k+1,...,l}}\circ\psi$,
$\psi_3=\pr_{E^2_{l+1,...,m}}\circ\psi$, and $\psi_4=\pr_{E^1_{r+1,...,n}}\circ\psi$.
\end{description}\end{theorem}

\vskip0.2cm We see that to each weakly-irreducible not irreducible holonomy algebras $\hol\subset\u(1,n+1)_{<\Real
p_1,\Real p_2>}$ an integer $0\leq m\leq n$ and a subalgebra $\un=\pr_{\u(m)}\hol\subset\u(m)$ are associated.

Recall that a pseudo-K\"ahlerian manifold is called {\it special pseudo-K\"ahlerian} if its Ricci tensor is zero.
This is equivalent to the inclusion $\hol_x\subset\su(T_xM,g_x,J^M_x)$.

\begin{corol}\label{G6corol1} 1)  A subalgebra $\g\subset\su(1,1)$ is the weakly-irreducible not irreducible holonomy algebra
of a special  pseudo-K\"ahlerian manifold of signature $(2,2)$ if and only if $\g$ is conjugated to  the subalgebra $\C\subset\su(1,1)_{<p_1,p_2>}$
or to $\su(1,1)_{<p_1,p_2>}$.

2) Let $n\geq 1$. Then  a subalgebra $\g\subset\su(1,n+1)$ is the weakly-irreducible not irreducible holonomy algebra
of a special  pseudo-K\"ahlerian manifold of signature $(2,2n+2)$ if and only if $\g$ is conjugated to one of the following
subalgebras of $\su(1,n+1)_{<p_1,p_2>}$:
\begin{description}
\item[] $\hol^{m,\un,\A^1,\phi}$, $\hol^{m,\un,\varphi,\phi}$ with $\phi(B,C)=-\frac{1}{n-m+2}\tr C$;
\item[] $\hol^{n,\un,\psi,k,l}$, $\hol^{m,\un,\psi,k,l,r}$ with $\u\subset\su(k)$.
\end{description}\end{corol}

Now we construct an example of metric with the holonomy algebra $\hol$ for each Lie algebra from Theorem \ref{hol2n}.

Let $0\leq m\leq n$ and $\un\subset\u(m)$ be a subalgebra.
Denote by $L\subset E$  the vector subspace annihilated by $\un$.
We can choose the basis $e_1,...,e_n,f_1,...,f_n$ in such a way that $L=\spa\{e_{n_0+1},...,e_n,f_{n_0+1},...,f_n\}$,
where $\dim L=2(n-n_0)$. Then $\un\subset\u(n_0)$ and $\un$ does not annihilate any proper subspace of $E_{1,...,n_0}$.
Let us consider a basis
  $A_1=\bigl(\begin{smallmatrix}B_1&-C_1\\C_1&B_1\end{smallmatrix}\bigr)$,...,$A_N=\bigl(\begin{smallmatrix}B_N&-C_N\\C_N&B_N\end{smallmatrix}\bigr)$
 of the vector space $\un$ such that $A_1,...,A_{N_1}$ is a basis of the vector space $\un'$ and
$A_{N_1+1},...,A_{N}$ is a basis of the vector space $\z(\un)$ ($N=\dim\un$, $N_1=\dim\un'$).
We denote by $(B^i_{\alpha  j})_{i,j=1}^{n_0}$ and $(C^i_{\alpha  j})_{i,j=1}^{n_0}$  the elements
of the matrices $B_\alpha$ and $C_\alpha$, respectively, where $\alpha=1,...,N$.

Let $x^1,...,x^{2n+4}$ be the standard coordinates  on $\Real^{2n+4}$.
Consider the following metric on $\Real^{2n+4}$: \begin{multline}\label{mkg} g=2dx^1dx^{2n+3}+2dx^2dx^{2n+4}+
\sum^{2n+2}_{i=3}(dx^i)^2+2\sum^{2n+2}_{i=3}u^i dx^{i} dx^{2n+4}\\+f_1\cdot(dx^{2n+3})^2+f_2\cdot(dx^{2n+4})^2
+2f_3dx^{2n+3}dx^{2n+4},\end{multline}
where $u^3$,...,$u^{2n+2}$, $f_1$, $f_2$ and $f_3$ are some functions which  depend on the  holonomy algebra that we wish to obtain.

For the  linear maps $\varphi,\phi:\un\to\Real$ we define the numbers $\varphi_\alpha=\varphi(A_\alpha)$ and
$\phi_\alpha=\phi(A_\alpha)$, $\alpha=N_1+1,...,N$.
For the linear map $\psi_1:\un\to E^1_{k+1,...,l}$ we define the numbers $\psi_{1\alpha i}$ such that
$\psi_1(A_\alpha)=\sum_{i=k+1}^{l} \psi_{1\alpha i}e_i$, $\alpha=N_1+1,...,N$. We define
analogous numbers for the linear maps $\psi_2$, $\psi_3$ and $\psi_4$.

Define the following functions
\begin{align*}
f^0_1=&\sum^{N}_{\alpha=1}\frac{(x^{2n+3})^{\alpha-1}}{(\alpha-1)!}\left(\sum_{i,j=1}^{n_0}\left(B^i_{\alpha  j}x^{i+2} x^{n+j+2}+
\frac{1}{2}C^i_{\alpha  j}x^{i+2} x^{j+2}+\frac{1}{2}C^i_{\alpha  j}x^{n+i+2} x^{n+j+2}\right)\right),\\
f^0_2=&f^0_1+
\sum_{i=1}^{n_0}\left(\left(\sum_{\alpha=1}^N\frac{1}{\alpha!}\sum_{j=1}^{n_0}\left(B^i_{\alpha  j}x^{j+2}-C^i_{\alpha  j}x^{n+j+2}\right)(x^{2n+3})^\alpha\right)\right.^2+\\
&\quad\quad\quad\left.\left(\sum_{\alpha=1}^N\frac{1}{\alpha!}\sum_{j=1}^{n_0}\left(B^i_{\alpha  j}x^{n+j+2}+C^i_{\alpha  j}x^{j+2}\right)(x^{2n+3})^\alpha\right)^2\right),\\
f^0_3=&0.\end{align*}

For the maps $\varphi$ and $\phi$ define the functions
\begin{align*}f^\varphi_1=&-2\sum_{\alpha=N_1+1}^N\frac{1}{\alpha!}\varphi_\alpha x^2(x^{2n+3})^\alpha,\quad
f^\varphi_2=-f^\varphi_1,\quad
f^\varphi_3=2\sum_{\alpha=N_1+1}^N\frac{1}{\alpha!}\varphi_\alpha x^1(x^{2n+3})^\alpha,\\
f^\phi_1=&\sum_{\alpha=N_1+1}^N\phi_\alpha \frac{(x^{2n+3})^{\alpha-1}}{(\alpha-1)!}
\left(-\frac{2}{\alpha} x^1x^{2n+3}+\sum_{i=m+1}^n(x^{i+2})^2\right),
\end{align*}
\begin{align*}
f^\phi_2=&\sum_{\alpha=N_1+1}^N\phi_\alpha \frac{(x^{2n+3})^{\alpha-1}}{(\alpha-1)!}
\left(\frac{2}{\alpha} x^1x^{2n+3}+\sum_{i=m+1}^n(x^{n+i+2})^2\right.\\
&\qquad\qquad\qquad\qquad\qquad\qquad\left.+\frac{(x^{2n+3})^2}{(\alpha+1)\alpha}\sum_{i=m+1}^n((x^{i+2})^2+(x^{n+i+2})^2)\right),\\
f^\phi_3=&\sum_{\alpha=N_1+1}^N\phi_\alpha \frac{(x^{2n+3})^{\alpha-1}}{(\alpha-1)!}
\left(-\frac{2}{\alpha} x^2x^{2n+3}+\sum_{i=m+1}^nx^{i+2}x^{n+i+2}\right)\end{align*}

For $K=N+1,N+2$ we define the functions \begin{align*}
f^{\A^1,K}_1=&-2\frac{1}{K!}x^2(x^{2n+3})^K,\quad f^{\A^1,K}_2=-f^{\A^1,K}_1,\quad f^{\A^1,K}_3=2\frac{1}{K!}x^1(x^{2n+3})^K,\\
f^{\t\A^2,K}_1=&\frac{1}{(K-1)!}(x^{2n+3})^{K-1}\left(-\frac{2}{K}x^1x^{2n+3} +\sum_{i=m+1}^n(x^{i+2})^2\right),\\
f^{\t\A^2,K}_2=&\frac{1}{(K-1)!}(x^{2n+3})^{K-1}\left(\frac{2}{K}x^1x^{2n+3} +\sum_{i=m+1}^n(x^{n+i+2})^2\right.\\
&\qquad\qquad\qquad\qquad\qquad\qquad\left.+
\frac{(x^{2n+3})^2}{(\alpha+1)\alpha}\sum_{i=m+1}^n((x^{i+2})^2+(x^{n+i+2})^2)\right),\\
f^{\t\A^2,K}_3=&\frac{1}{(K-1)!}(x^{2n+3})^{K-1}\left(-\frac{2}{K}x^2x^{2n+3} +\sum_{i=m+1}^nx^{i+2}x^{n+i+2}\right).\end{align*}

For any numbers $1\leq m_1\leq m_2\leq n$ consider the functions
$$\tilde f_{1m_1}^{m_2}=\sum^{m_2}_{i=m_1}(x_{i+2}^2-x_{n+i+2}^2),\quad\tilde f_{2m_1}^{m_2}=-\tilde f_{1m_1}^{m_2},\quad
\tilde f_{3m_1}^{m_2}=2\sum^{m_2}_{i=m_1}x^{i+2}x^{n+i+2}.$$
For any numbers $0\leq m_1\leq m_2\leq n$ and $K\geq N+1$ consider the functions
\begin{align*}\breve f_{1m_1}^{Km_2}=&-\sum^{m_2}_{i=m_1}\frac{2}{(K+i-m_1)!}x^{n+i+2}(x^{2n+3})^{K+i-m_1},\quad \breve f_{2m_1}^{Km_2}=-\breve f_{1m_1}^{Km_2},\\
\breve f_{3m_1}^{Km_2}=&-\sum^{m_2}_{i=m_1}\frac{2}{(K+i-m_1)!}x^{i+2}(x^{2n+3})^{K+i-m_1}.\end{align*}

For the Lie algebra $\hol^{n,\un,\psi,k,l}$ we consider the functions
\begin{align*}
f_1^{n,\psi}=&\sum_{\alpha=N_1+1}^N\frac{2}{\alpha!}\left(\sum_{i=k+3}^{l+2}\psi_{1\alpha i}x^{n+i}-\sum_{i=k+3}^{l+2}\psi_{2\alpha i}x^{i}-
\sum_{i=l+3}^{n+2}\psi_{3\alpha i}x^{i}\right)(x^{2n+3})^\alpha,\\ f_2^{n,\psi}=&-f_1^{n,\psi},\\
f_3^{n,\psi}=&\sum_{\alpha=N_1+1}^N\frac{2}{\alpha!}\left(-\sum_{i=k+3}^{l+2}\psi_{1\alpha i}x^{i}-\sum_{i=k+3}^{l+2}\psi_{2\alpha i}x^{n+i}-
\sum_{i=l+3}^{n+2}\psi_{3\alpha i}x^{n+i}\right)(x^{2n+3})^\alpha.\end{align*}

For the Lie algebra $\hol^{m,\un,\psi,k,l,r}$ we  define the functions
\begin{align*}f_1^{m,\psi}=&\sum_{\alpha=N_1+1}^N\frac{2}{\alpha!}\left(
\sum_{i=k+3}^{l+2}\psi_{1\alpha i}x^{n+i}+\sum_{i=r+3}^{n+2}\psi_{4\alpha i}x^{n+i}\right.\\
&\qquad\qquad\qquad\qquad\left.-\sum_{i=k+3}^{l+2}\psi_{2\alpha i}x^{i}
-\sum_{i=l+3}^{m+2}\psi_{3\alpha i}x^{i}\right)(x^{2n+3})^\alpha,\\ f_2^{m,\psi}=&-f_1^{m,\psi},\\
f_3^{m,\psi}=&\sum_{\alpha=N_1+1}^N\frac{2}{\alpha!}\left(
-\sum_{i=k+3}^{l+2}\psi_{1\alpha i}x^{i}-\sum_{i=r+3}^{n+2}\psi_{4\alpha i}x^{i}\right.\\
&\qquad\qquad\qquad\qquad\left.-\sum_{i=k+3}^{l+2}\psi_{2\alpha i}x^{n+i}-
\sum_{i=l+3}^{m+2}\psi_{3\alpha i}x^{n+i}\right)(x^{2n+3})^\alpha.\end{align*}

Define the functions $u^3$,...,$u^{n_0+2}$, $u^{n+3}$,...,$u^{n+n_0+2}$ as follows:
\begin{align*}
u^i=&\sum_{\alpha=1}^{N}\frac{1}{\alpha!}\left(\sum_{j=1}^{n_0}(B^{i-2}_{\alpha j}x^{j+2}-C^{i-2}_{\alpha j}x^{n+j+2})\right)(x^{2n+3})^\alpha,\\
u^{n+i}=&\sum_{\alpha=1}^{N}\frac{1}{\alpha!}\left(\sum_{j=1}^{n_0}(B^{i-2}_{\alpha j}x^{n+j+2}+C^{i-2}_{\alpha j}x^{j+2})\right)(x^{2n+3})^\alpha,
\end{align*}
where $3\leq i\leq n_0+2$.

For the Lie algebras $\hol^{m,\u,\A^1,\t\A^2}$  and  $\hol^{m,\u,\varphi,\t\A^2}$ we set in addition
\begin{align*}
u^i=&-\frac{1}{(N+2)!}x^{n+i}(x^{2n+3})^{N+2},\\
u^{n+i}=&   \frac{1}{(N+2)!}x^{i}(x^{2n+3})^{N+2},
\end{align*} where $m+3\leq i\leq 2n+2$.

For the Lie algebras $\hol^{m,\u,\A^1,\phi}$  and  $\hol^{m,\u,\varphi,\phi}$ we set
\begin{align*}
u^i=&-\sum_{\alpha=1}^{N}\frac{1}{\alpha!}\phi_\alpha x^{n+i}(x^{2n+3})^{\alpha},\\
u^{n+i}=&\sum_{\alpha=1}^{N}\frac{1}{\alpha!}\phi_\alpha x^{i}(x^{2n+3})^{\alpha},
\end{align*} where $m+3\leq i\leq 2n+2$.

For the Lie algebra $\hol^{m,\u,\lambda}$ we set
\begin{align*}
u^i=&-\frac{1}{(N+1)!}\lambda x^{n+i}(x^{2n+3})^{N+1},\\
u^{n+i}=&   \frac{1}{(N+1)!}\lambda x^{i}(x^{2n+3})^{N+1},
\end{align*} where $m+3\leq i\leq 2n+2$.

We assume that the functions $u^i$ that were not defined now are zero.

If we choose the functions $f_1,$ $f_2$ and $f_3$ such  that
$f_1(0)=f_2(0)=f_3(0)=0$, then for the metric $g$ given by \eqref{mkg} we have
$g_0=\eta$ and we can identify the tangent space to $\Real^{2n+4}$ at $0$ with the vector space
$\Real^{2,2n+2}$.

Note that if  $n=0$, then   $$ g=2dx^1dx^{3}+2dx^2dx^{4}+f_1\cdot(dx^{3})^2+f_2\cdot(dx^{4})^2+2f_3dx^{3}dx^{4}.$$

\begin{theorem}\label{G6metrics}
Let $\hol_0$ be  the holonomy algebra of the metric $g$ at the point $0\in \Real^{2n+4}$.

{\bf 1)} Let $n=0$, then $\hol_0$ depends on the functions $f_1$, $f_2$ and $f_3$ as in Table \ref{tabhol0}.

\begin{tab}\label{tabhol0} Dependence of $\hol_0$ on the functions $f_1$, $f_2$ and $f_3$ for $n=0$
\begin{longtable}{|l|l|}\hline $f_i$, $(i=1,2,3)$ & $\hol_0$\\ \hline
$f_1=-2x^2x^3-x^1(x^3)^2$, $f_2=-f_1$, $f_3=2x^1x^3-x^2(x^3)^2$ &$\hol^{1}_{n=0}$\\ \hline
$f_1=(x^1)^2-(x^2)^2$, $f_2=-f_1$, $f_3=2x^1x^2$ &$\hol^{2}_{n=0}$\\ \hline
$f_1=-2\gamma_1x^2x^3-2\gamma_2x^1x^3$, $f_2=-f_1$, $f_3=2\gamma_1x^1x^3-2\gamma_2x^2x^3$ &$\hol^{\gamma_1,\gamma_2}_{n=0}$\\
& (if $\gamma_1^2+\gamma_2^2\neq 0$)\\ \hline
$f_1=(x^4)^2,$ $f_2=f_3=0$&$\hol^{\gamma_1=0,\gamma_2=0}_{n=0}$\\\hline
\end{longtable}\end{tab}

{\bf 2)} Let $n>0$, then $\hol_0$ depends on the functions $f_1$, $f_2$ and $f_3$ as in Table \ref{tabhol}.

\begin{tab}\label{tabhol} Dependence of $\hol_0$ on the functions $f_1$, $f_2$ and $f_3$ for $n>0$
\begin{longtable}{|l|l|}\hline $f_i$, $(i=1,2,3)$ & $\hol_0$\\ \hline
$f_i=f_i^{\A^1,N+1}+f_i^{\t\A^2,N+2}+f^0_i+\tilde f_{i\,n_0+1}^{m}+\breve f_{i\, m+1}^{N+3\,\,n} $&$\hol^{m,\un,\A^1,\t\A^2}$\\ \hline
$f_i=f_i^{\A^1,N+1}+f_i^{\phi}+f^0_i+\tilde f_{i\,n_0+1}^{m}+\breve f_{i\, m+1}^{N+2\,\,n} $&$\hol^{m,\un,\A^1,\phi}$\\ \hline
$f_i=f_i^{\varphi}+f_i^{\t\A^2,N+1}+f^0_i+\tilde f_{i\,n_0+1}^{m}+\breve f_{i\, m+1}^{N+2\,\,n}$&$\hol^{m,\un,\varphi,\t\A^2}$\\ \hline
$f_i=f_i^{\varphi}+f_i^{\phi}+f^0_i+\tilde f_{i\,n_0+1}^{m}+\breve f_{i\, m+1}^{N+2\,\,n}\breve f_{i\, m+1}^{N+1\,\,n} $&$\hol^{m,\un,\varphi,\phi}$\\ \hline
$f_i=f_i^{\A^1,N+1}+\lambda f_i^{\t\A^2,N+1}+f^0_i+\tilde f_{i\,n_0+1}^{m}+\breve f_{i\, m+1}^{N+2\,\,n}$&$\hol^{m,\un,\lambda}$\\ \hline
$f_i=f^0_i+\tilde f_{i\,n_0+1}^{k}+f_i^{n,\psi}+\breve f_{i\, l+1}^{N+1\,n}  $&$\hol^{n,\un,\psi,k,l}$ (if $\dim\z(\un)\geq n+l-2k$)\\ \hline
$f_i=f^0_i+\tilde f_{i\,n_0+1}^{k}+f_i^{m,\psi}+\breve f_{i\, l+1}^{N+1\,r}  $&$\hol^{m,\un,\psi,k,l,r}$\\ & (if $\dim\z(\un)\geq n+m+l-2k-r$)\\ \hline
\end{longtable}\end{tab}
\end{theorem}


\section{Proof of Theorem \ref{hol2n}}\label{sec3.2.1}

In this section we will prove that the Lie algebras of Theorem \ref{hol2n} exhaust all weakly-irreducible Berger subalgebras of
$\u(1,n+1)_{<p_1,p_2>}$, i.e. all candidates for the holonomy algebras.
The rest of the proof of Theorem \ref{hol2n} will follow from Theorem \ref{G6metrics}.



Now we will describe the spaces of curvature tensors $\R(\g)$ for subalgebras  $\g\subset\u(1,n+1)_{<p_1,p_2>}$.
We will use the following obvious fact. Let $\f_1\subset\f_2\subset \so(r,s)$, then
\begin{equation}\label{svoistvo1} R\in\R(\f_1)\text{ if and only if } R\in\R(\f_2) \text{ and } R(\Real^{r,s}\wedge \Real^{r,s} )\subset
\f_1.\end{equation}

 First we will describe the space $\R(\g^\u)$ for the Lie algebra
{\scriptsize  $$\g^\u=\left\{\left. \left (\begin{array}{ccccc}
a_1&-a_2 &-X^t &0 &-c\\
a_2&a_1 &-Y^t &  c &0\\
0 &0&A&X&Y\\
0&0&0&-a_1&-a_2\\
0&0&0&a_2&-a_1\\
\end{array}\right)\right|\,
\begin{array}{c}
a_1,a_2,c\in \Real,\\ X,Y\in \Real^{2n},\\ A\in\un \end{array}\right\}\subset\so(2,2n+2)_{<p_1,p_2>}.$$}
Here $\un\subset\u(n)$ is a subalgebra and $\so(2,2n+2)_{<p_1,p_2>}$ is the subalgebra of $\so(2,2n+2)$ that preserves
the isotropic plane $\Real p_1\oplus\Real p_2$.

Using the form $\eta$, we identify $\so(2,2n+2)$ with the space $$\Real^{2,2n+2}\wedge \Real^{2,2n+2}=
\spa\{u\wedge v=u\otimes v-v\otimes u|u,v\in\Real^{2,2n+2}\}.$$
The identification is given by the formula $$(u\wedge v)w=\eta(u,w)v-\eta(v,w)u \text{ for all } u,v,w\in\Real^{2,2n+2}.$$
Similarly, we identify $\so(n)$ with $E\wedge E\subset\Real^{2,2n+2}\wedge \Real^{2,2n+2}$.
It is easy to see that the element
{\scriptsize  $$\left (\begin{array}{ccccc}
a_1&-a_2 &-X^t &0 &-c\\
a_2&a_1 &-Y^t &c &0\\
0 &0&A&X&Y\\
0&0&0&-a_1&-a_2\\
0&0&0&a_2&-a_1\\
\end{array}\right)\in\g^\u$$} corresponds to
$$-a_1(\pq)+a_2(\qp)+A+p_1\wedge X+p_2\wedge Y+cp_1 \wedge p_2\in\Real^{2,2n+2}\wedge \Real^{2,2n+2}.$$
Thus we obtain the following decomposition of $\g^\un$:
\begin{equation}\label{decomp}\g^\un=(\Real(\pq)\oplus\Real(\qp)\oplus\un)\zr(p_1\wedge E+p_2\wedge E+\Real p_1\wedge p_2).\end{equation}
By analogy, we have
$$\begin{array}{rl}\u(1,n+1)_{<p_1,p_2>}=&(\Real(\pq)\oplus\Real(\qp)\oplus\u(n))\\&\zr
(\{p_1\wedge x+p_2\wedge Jx|x\in E^1\}+\{p_1\wedge Jx-p_2\wedge x|x\in E^1\}+\Real p_1\wedge p_2).\end{array}$$

The decomposition $\Real^{2,2n+2}=\Real p_1+\Real p_2+E+\Real q_1+\Real q_2$ gives us the decomposition
\begin{multline}\label{decomp3} \Real^{2,2n+2}\wedge\Real^{2,2n+2}=\Real(\pq)+\Real(\qp)+\Real(p_1\wedge q_1-p_2\wedge q_2)+\Real(p_1\wedge q_2+p_2\wedge q_1)\\+
E\wedge E+p_1\wedge E+p_2\wedge E+q_1\wedge E+q_2\wedge E+\Real p_1\wedge p_2+\Real q_1\wedge q_2.\end{multline}

The metric $\eta$ defines the metric $\eta\wedge \eta$ on $\Real^{2,2n+2}\wedge\Real^{2,2n+2}$.
Let $R\in\R(\g^\un)$. It can be proved that \begin{equation}\label{svoistvo3}
\eta\wedge\eta(R(u\wedge v),z\wedge w)=\eta\wedge\eta(R(z\wedge w),u\wedge v)\text{ for all } u,v,z,w\in\Real^{2,2n+2}.\end{equation}
This shows that $R:\Real^{2,2n+2}\wedge\Real^{2,2n+2}\to\g^\un\subset\Real^{2,2n+2}\wedge\Real^{2,2n+2}$ is a symmetric linear map.
Hence $R$ is zero on the orthogonal complement to $\g^\un$ in $\Real^{2,2n+2}\wedge\Real^{2,2n+2}$. In particular,
\begin{equation}\label{b1} R|_{\Real(p_1\wedge q_1-p_2\wedge q_2)+\Real(p_1\wedge q_2+p_2\wedge q_1)+\Real p_1\wedge p_2+p_1\wedge E+p_2\wedge E}=0.\end{equation}
Thus $R$ can be considered as the linear map
$$\begin{array}{rl}R:&\Real(\pq)+\Real(\qp)+E\wedge E+q_1\wedge E+q_2\wedge E+\Real q_1\wedge q_2\\
 &\to\g^\un=(\Real(\pq)\oplus\Real(\qp)\oplus\un)\zr(p_1\wedge E+p_2\wedge E+\Real p_1\wedge p_2).\end{array}$$

Consider the following  set of subsets of $\Real^{2,2n+2}\wedge\Real^{2,2n+2}$,
$$\F=\{\Real(\pq),\Real(\qp),\un,p_1\wedge E,p_2\wedge E,\Real p_1\wedge p_2\}.$$
For any $F\in\F$ we  set $$R_{F}=\pr_{F}\circ R:\Real^{2,2n+2}\wedge\Real^{2,2n+2}\to F,$$
where $\pr_{F}$ is the projection with respect to the decomposition \eqref{decomp}.
Obviously, $$R=\sum_{F\in\F}R_{F}.$$

Note that from \eqref{b1} it follows that $R(p_1\wedge q_1)=R(p_2\wedge q_2)$ and $R(p_1\wedge q_2)=-R(p_2\wedge q_1)$.
Using the Bianchi identity, \eqref{svoistvo3} and \eqref{b1}, it is easy to show that $R$  can be found from Table~3.2.1 on page~\pageref{tabboth},
where on the position $(u\wedge v,\F)$ stays the value $R_{F}(u\wedge v)$.

\addtocounter{tab}{1}

In Table 3.2.1 we have  $x,y\in E$, $\lambda_1,...,\lambda_5\in\Real$, $K_1,K_2,L_1,L_2\in\Hom(\Real,E),R_\un\in\R(\un)$,
$$P_1,P_2\in\P(\un)=\left\{P\in \Hom (E,\un)\left|\,\begin{array}{c}\eta(P(u)v,w)+\eta(P(v)w,u)+\eta(P(w)u,v)=0\\\text{ for all }u,v,w\in E
\end{array}\right\}\right.,$$
$T_1,T_2\in\Hom(E,E)$, $T_1^*=T_1$, $T_2^*=T_2$ and $S\in \Hom(E,E)$ is a linear map such that $S-S^*\in\un$.

It is easy to show that for any elements as above the linear map $R\in\Hom(\Real^{2,2n+2}\wedge\Real^{2,2n+2},\g^\un)$ defined
by  Table 3.2.1 and \eqref{b1} satisfies $R\in\R(\g^\un)$.

For any $0\leq m\leq n$ and $\un\subset\u(m)\oplus\sod(m+1,...,n)$ consider  the  subalgebra
$$\begin{array}{rl}\g^{m,\un}=&(\Real(\pq)\oplus\Real(\qp+J_{m+1,...,n})\oplus\un)\\
&\zr(\{p_1\wedge x+p_2\wedge Jx|x\in E^1\}+\{p_1\wedge Jx-p_2\wedge x|x\in E^1_{1,...,m}\}+\Real p_1\wedge p_2)\\
&\subset\u(1,n+1)_{<p_1,p_2>}.\end{array}$$

For the $\u(n)$-projection of the Lie algebra $\g^{m,\un}$ we have $\pr_{\u(n)}\g^{m,\un}=\un\oplus\Real J_{m+1,...,n}$.
Let us consider a curvature tensor
$R\in\R(\g^{m,\un})$. Since $\g^{m,\un}\subset\g^{\un\oplus\Real J_{m+1,...n}}$, to decompose $R$ we can use \eqref{svoistvo1} and Table  3.2.1.
For $K_1\in\Hom(\Real,E)$ let $K_1^1=\pr_{E^1}\circ K_1$ and $K_1^2=\pr_{E^2}\circ K_1$. Then $K_1=K^1_1+K^2_1$.
For $S\in\Hom(E,E)$ let $$S^{11}=\pr_{E^1}\circ S|_{E^1},\quad S^{12}=\pr_{E^1}\circ S|_{E^2}, \quad S^{21}=\pr_{E^2}\circ S|_{E^1},
\quad S^{22}=\pr_{E^2}\circ S|_{E^2}$$
and extend these linear maps to $E$ mapping the natural complement to zero. We get the decomposition $S=S^{11}+S^{12}+S^{21}+S^{22}$.
For $P_1\in\P(\un)$ let $Q_1=P_1^*\in\Hom(\un,E)$,  $Q_1^1=\pr_{E^1}\circ Q_1$ and $Q_1^2=\pr_{E^2}\circ Q_1$.
Consequently, $Q_1=Q^1_1+Q^2_1$. Consider the analogous  decompositions for the elements $K_2,L_1,L_2\in\Hom(\Real,E)$,
$Q_2=P_2^*\in\Hom(\un,E)$, $T_1,T_2,S^*\in\Hom(E,E)$.
Using the condition $R(\Real^{2,2n+2}\wedge\Real^{2,2n+2})\subset\g^{m,\un}$, we obtain
\begin{align}
\label{b4a} \begin{array}{cccc} K^2_2=JK^1_1, &K^1_2=JK^2_1, &K^1_1(1)\in E^1_{1,...,m},&K^2_1(1)\in E^2_{1,...,m}\end{array}\\
\label{b4b} \begin{array}{cccc} Q^2_2=JQ^1_1, &Q^1_2=JQ^2_1, &Q^2_1(1)\in E^2_{1,...,m}&\end{array}
\end{align}
\begin{multline}
\label{b4c1} T^{11}_1=-JS^{21},\quad T^{21}_1=-JS^{11},\\  S^{21}(E^1_{1,...,m})\subset E^2_{1,...,m},
\quad S^{21}(E^1_{m+1,...,n})\subset E^2_{m+1,...,n}, \end{multline}
\begin{align}
\label{b4c2}\begin{array}{cccc} T^{12}_1=-JS^{22},&T^{22}_1=-JS^{12}, &S^{22}(E)\subset E^2_{1,...,m},& \end{array}\\
\label{b4d1}\begin{array}{cccc} T^{11}_2=JS^{*21},&T^{21}_2=JS^{*11}, &S^{*21}(E)\subset E^2_{1,...,m},& \end{array}\\
\label{b4d2}\begin{array}{cccc} T^{12}_2=JS^{*22},&T^{22}_2=JS^{*12}, &S^{*22}(E)\subset E^2_{1,...,m},& \end{array}\\
\label{b4e} \begin{array}{cccc} L^2_2=JL^1_1, &L^1_2=JL^2_1, &L^2_1(1)\in E^2_{1,...,m}.&\end{array}
\end{align}
From \eqref{b4c1}, \eqref{b4c2}, \eqref{b4d1}, \eqref{b4d2}, and the fact that $T^*_1=T_1$ and $T^*_2=T_2$ it follows that
$$(S^{22})^*=JS^{11}J,\quad (S^{12})^*=JS^{12}J,\quad(S^{21})^*=JS^{21}J,$$
$$S^{11}(E^1)\subset E^1_{1,...,m} ,\quad S^{12}(E^2)\subset E^1_{1,...,m}.\quad S^{21}(E^1_{1,...,m})\subset E^2_{1,...,m}.$$
Since $\pr_{\Real(\qp)}R(q_1\wedge q_2)=\lambda_3(\qp)$, we see that $S^{21}|_{E^1_{m+1,...,n}}=\lambda_3 J|_{E^1_{m+1,...,n}}$.

\begin{figure}[p]\label{tabboth}
\includegraphics[angle=90]{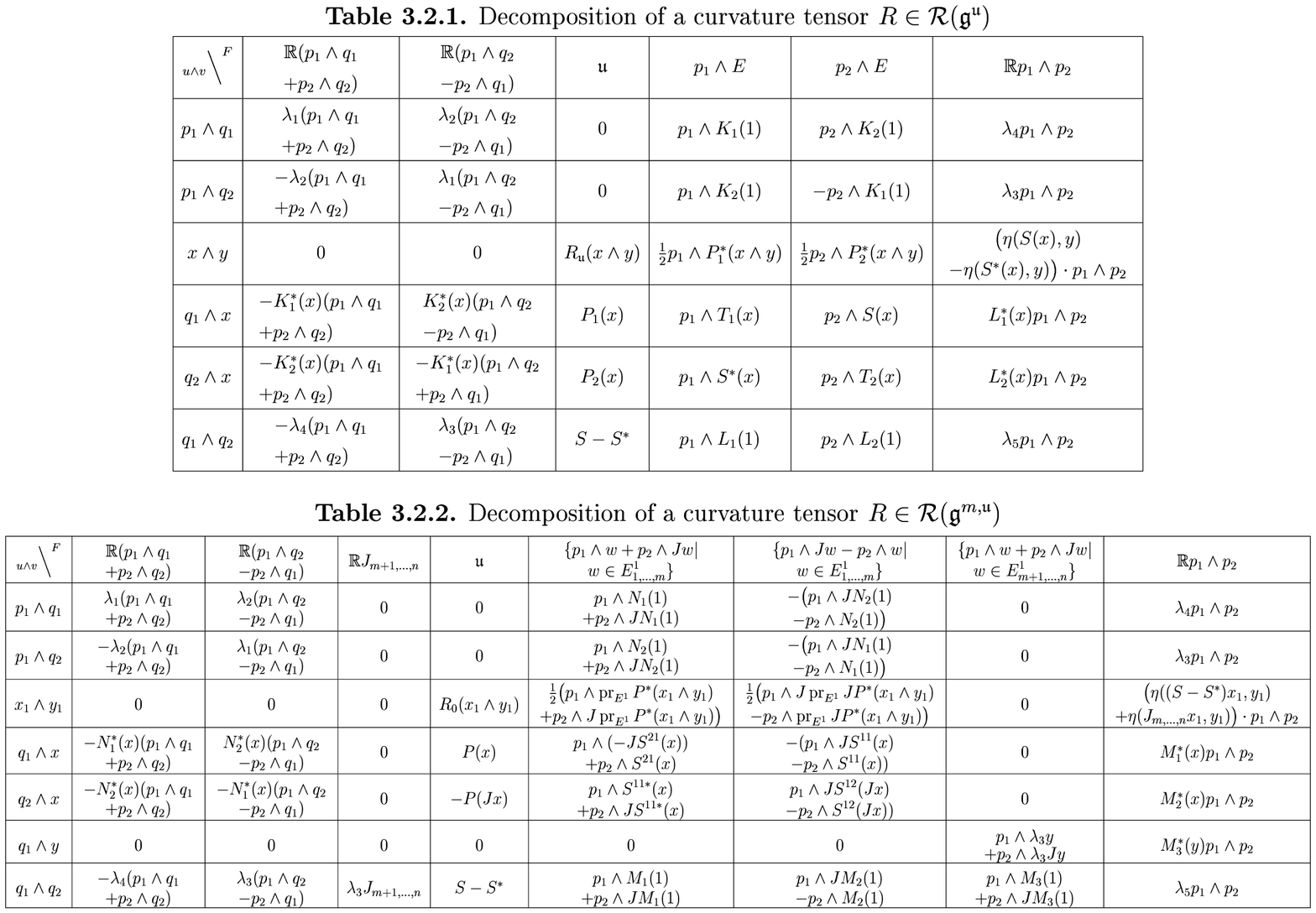}
\end{figure}

Using \eqref{b4b}, we get $Q_2=JQ_1$, i.e. $P_2^*=JP_1^*$. Hence, $P_2=P_1J^*$. Thus, $P_2=-P_1J$.

Since $\pr_{\Real(\qp)}R(E\wedge E)=\{0\}$, from Table 3.2.1 it follows that
$\pr_{\u(n)}R(E\wedge E)\subset\un$.
It is well known that from the inclusion $\un\subset\u(m)\oplus\sod(m+1,...,n)$ it follows that
$\R(\un)=\R(\un\cap\u(m))\oplus\R(\un\cap\sod(m+1,...,n)).$
Moreover, $\R(\sod(m+1,...,n))=\{0\}$. Therefore, $\pr_{\u(n)}R(E\wedge E)\subset\un\cap \u(m)$.
Thus for $R_0=\pr_{\u(n)}\circ R|_{E\wedge E}$ we get $R_0\in\R(\un\cap\u(m))$.

Similarly, since $\pr_{\Real(\qp)}R(q_1\wedge E^1_{1,...,m})=\{0\}$, we see that $\pr_{\u(n)}R(q_1\wedge E)\subset\un$.
In \cite{Gal1} it was proved that
$$\P(\un)=\P(\un\cap\u(m))\oplus\P(\un\cap\sod(m+1,...,n))\quad\text{ and }\quad \P(\sod(m+1,...,n))=\{0\}.$$
Thus,   $P_1\in\P(\un\cap\u(m))$.

From \eqref{svoistvo3} it follows that
\begin{equation}\label{b5a}
R|_{\Real(p_1\wedge q_1-p_2\wedge q_2)+\Real(p_1\wedge q_2+p_2\wedge q_1)+\Real p_1\wedge p_2
+\{q_1\wedge x-q_2\wedge Jx|x\in E^1\}+\{q_1\wedge Jx+q_2\wedge x|x\in E^1_{1,...,m}\}}=0,\end{equation}
\begin{equation}\label{b5b} R|_{
q_1\wedge E^2_{m+1,...,n}+q_2\wedge E^1_{m+1,...,n} +p_1\wedge E+p_2\wedge E}=0.\end{equation}
In particular, $R(q_2\wedge Jx)=R(q_1\wedge x)$ for all $x\in E^1$, and $R(q_1\wedge Jx)=R(q_2\wedge x)$ for all $x\in E^1_{1,...,m}$.

We set the following denotation:
$N_1=K^1_1$, $N_2=JK^2_1$, $M_1=\pr_{E^1_{1,...,m}}\circ L^1_1$, $M_2=-JL^2_1$,   $M_3=\pr_{E^1_{m+1,...,n}}\circ L^1_1$, $P=P_1$,
$S^{21}=\pr_{E^2_{1,...,m}}\circ S|_{E^1_{1,...,m}}$. 
Now the curvature tensor $R\in\R(\g^{m,\un})$ can be found as above from the conditions \eqref{b5a}, \eqref{b5b} and Table 3.2.2.
In this case we assume that
$$\begin{array}{rl}\F=&\{\Real(\pq),\Real(\qp),\Real J_{m+1,...,n},\un,p_1\wedge E,\\
&\{p_1\wedge x+p_2\wedge Jx|x\in E^1_{1,...,m}\},\quad\{p_1\wedge Jx-p_2\wedge x|x\in E^1_{1,...,m}\},\\
&\{p_1\wedge x+p_2\wedge Jx|x\in E^1_{m+1,...,n}\},\Real p_1\wedge p_2\}.\end{array}$$

\addtocounter{tab}{1}

In Table 3.2.2 we have $x_1,y_1\in E$, $x\in E^1_{1,...,m}$, $y\in E^1_{m+1,...,n}$,
$\lambda_1,...,\lambda_5\in\Real$, $N_1,N_2,M_1,M_2\in\Hom(\Real,E^1_{1,...,m}),$ $M_3\in\Hom(\Real,E^1_{m+1,...,n})$,
if $m<n$, then $\lambda_1=\lambda_2=0$ and $N_1=N_2=0$. Furthermore, $R_0\in\R(\un\cap\u(m))$,
$P\in\P(\un\cap\u(m))$,
$S^{11}\in\Hom(E_{1,...,m},E^1_{1,...,m})$, $ S^{21}\in\Hom(E_{1,...,m},E^2_{1,...,m})$,
$S^{12}\in\Hom(E_{1,...,m},E^1_{1,...,m})$,
$S^{11}|_{E^2_{1,...,m}}=0$, $S^{21}|_{E^2_{1,...,m}}=0$, $S^{12}|_{E^1_{1,...,m}}=0$,
$S^{12*}=JS^{12}J$,   $S^{21*}=JS^{21}J$, $S\in\Hom(E_{1,...,m},E_{1,...,m})$, $S= S^{11}+S^{12}+S^{21}+JS^{11*}J$ and $S-S^*\in\un\cap \u(m)$.
Conversely, for  any elements as above the linear map $R\in\Hom(\Real^{2,2n+2}\wedge\Real^{2,2n+2},\g^{m,\un})$ defined
by the Table 3.2.2, \eqref{b5a}  and \eqref{b5b} satisfies $R\in\R(\g^{m,\un})$.

{\bf 1)} Consider the case $n=0$.
\begin{lem} The Lie algebras of Part 1 of Theorem \ref{hol2n} exhaust all weakly-irreducible Berger subalgebras of $\u(1,1)_{<p_1,p_2>}$.\end{lem}
{\it Proof.} Let $R\in\R(\u(1,1)_{<p_1,p_2>})$. As above, $R$ can be found from the conditions $R(p_2\wedge q_2)=R(p_1\wedge q_1)$,
$R(p_2\wedge q_1)=-R(p_1\wedge q_2)$, $R(p_1\wedge p_2)=0$ and the following table:
{\scriptsize
$$\begin{array}{|c|c|c|c|}
\hline
_{u\wedge v}\Big\backslash^{F}&\Real(\pq)&\Real(\qp)&\Real p_1\wedge p_2\\ \hline
p_1\wedge q_1 &\lambda_1(\pq)&\lambda_2(\qp)&\lambda_4 p_1\wedge p_2\\ \hline
p_1\wedge q_2 &-\lambda_2(\pq)&\lambda_1(\qp)&\lambda_3 p_1\wedge p_2\\ \hline
q_1\wedge q_2 &-\lambda_4(\pq)&\lambda_3(\qp)&\lambda_5 p_1\wedge p_2\\ \hline
\end{array}$$}

We will consider all subalgebras of $\u(1,1)_{<p_1,p_2>}$ and check which of these subalgebras are weakly-irreducible Berger subalgebras.
 Let $\g\subset\u(1,1)_{<p_1,p_2>}$ be a subalgebra. We have the following cases:

{\bf Case 1.} $\pr_\C\g=\{0\}$, i.e. $\g\subset\A^1\oplus\A^2$;

{\bf Case 2.} $\C\subset\g$;

{\bf Case 3.} $\C\not\subset\g$ and $\pr_\C\g\neq\{0\}$.

Consider these cases.

{\bf Case 1.} $\pr_\C\g=\{0\}$, i.e. $\g\subset\A^1\oplus\A^2$. We have the following subcases.

{\bf Subcase 1.1.} $\g=\A^1\oplus\A^2$. We claim that $\g$ is a weakly-irreducible Berger subalgebra.
Suppose that $\g$ preserves a non-trivial vector subspace $L\subset\Real^{2,2}$.
Let $\alpha_1 p_1+\alpha_2 p_2+\beta_1 q_1 +\beta_2 q_2\in L$ be a non-zero vector.
Applying the element $(1,0)\in\A^1$, we get $\alpha_1 p_1+\alpha_2 p_2-\beta_1 q_1 -\beta_2 q_2\in L$.
Hence, $\alpha_1 p_1+\alpha_2 p_2\in L$ and $\beta_1 q_1 +\beta_2 q_2\in L$. Applying to these vectors the element $(0,1)\in\A^2$,
we get $\alpha_1 p_2-\alpha_2 p_1\in L$ and $\beta_1 q_2 -\beta_2 q_1\in L$. There are three possibilities:
$L=\Real^{2,2}$, $L=\Real p_1\oplus\Real p_2$ or $L=\Real q_1\oplus\Real q_2$. Hence the subalgebra $\g$ is weakly-irreducible.

Furthermore, $\g$ is spanned be the image of the curvature tensor $R\in\R(\g)$ given by $\lambda_1=1$ and $\lambda_2=\cdots=\lambda_5=0$.

{\bf Subcase 1.2.} $\g=\{(a\gamma_1,a\gamma_2,0)|a\in\Real\}$, where $\gamma_1,\gamma_2\in\Real$.
We claim that this subalgebra is not a  weakly-irreducible Berger subalgebra.
Indeed, if $\gamma_2=0$, then $\g$ preserves the non-degenerate proper subspace $\Real p_1\oplus\Real q_1\subset\Real^{2,2}$.
If $\gamma_1=0$, then $\g$ preserves the non-degenerate proper subspace $\Real (p_1+q_1)\oplus\Real (p_2+q_2)\subset\Real^{2,2}$.
Suppose that $\gamma_1\neq 0$ and $\gamma_2\neq 0$. Let $R\in\R(\g)$.
We have $$R(p_1\wedge q_1)=\lambda_1(\pq)+\lambda_2(\qp)+\lambda_4(p_1\wedge p_2),$$
 $$R(p_1\wedge q_2)=-\lambda_2(\pq)+\lambda_1(\qp)+\lambda_3(p_1\wedge p_2).$$
Hence, $\lambda_3=\lambda_4=0$ and $\frac{\lambda_1}{-\lambda_2}=\frac{\lambda_2}{\lambda_1}$. Therefore, $\lambda_1=\lambda_2=0$.
Moreover, $R(q_1\wedge q_2)=\lambda_5(p_1\wedge p_2)$. Consequently, $\lambda_5=0$. Thus, $R=0$, $\R(\g)=\{0\}$ and $\g$ is not a Berger algebra.

{\bf Case 2.} $\C\subset\g$. We have the following subcases:

{\bf Subcase 2.1.} $\g=(\A^1\oplus\A^2)\zr\C=\u(1,1)_{<p_1,p_2>}$;

{\bf Subcase 2.2.} $\g=\{(a\gamma_1,a\gamma_2,0)|a\in\Real\}\zr C$, where $\gamma_1,\gamma_2\in\Real$.
These subalgebras contain $\C$, hence they are weakly-irreducible (Part 1 of Theorem \ref{Tw-ir}).
These subalgebras are Berger algebras, since any element of these algebras can be obtained as $R(q_1\wedge q_2)$ for some curvature tensor $R$.

{\bf Case 3.} $\C\not\subset\g$ and $\pr_\C\g\neq\{0\}$. Consider the following subcases.

{\bf Subcase 3.1.} $\dim\g=1$, then $\g=\{c\gamma_1,c\gamma_2,c)|c\in\Real\}$, where $\gamma_1,\gamma_2\in\Real$, $\gamma_1\neq 0$ or $\gamma_2\neq 0$.
We claim that $\g$ is not a Berger algebra. Indeed, let $R\in\R(\g)$ by analogy with Subcase 1.2, we have
$\frac{\lambda_1}{-\lambda_2}=\frac{\lambda_2}{\lambda_1}=\frac{\lambda_4}{\lambda_3}$. Hence, $\lambda_1=\lambda_2=0$ and $\lambda_3=\lambda_4=\lambda_5=0$.

{\bf Subcase 3.2.} $\dim\g=2$, then $\pr_{\A^1\oplus\A^2}\g=\A^1\oplus\A^2$ and
$\g=\{(a_1,a_2,\nu(a_1,a_2))|a_1,a_2\in\Real\},$ where $\nu:\Real\oplus\Real\to\Real$ is a non-zero linear map. Let
$\gamma_1,\gamma_2\in\Real$ be  numbers such that $\nu(\gamma_1,\gamma_2)=0$ and  $\nu(-\gamma_2,\gamma_1)=1$. Hence $\g$
has the form $\{(a\gamma_1,a\gamma_2,0)|a\in\Real\}\zr\{-c\gamma_2,c\gamma_2,c)|c\in\Real\}.$ We have
$[(\gamma_1,\gamma_2,0),(-\gamma_2,\gamma_1,1)]=(0,0,2\gamma_1)\in\g$. Therefore, $\gamma_1=0$. Let $\gamma=\gamma_2$.
Thus, $\g=\A^2\oplus\{(c\gamma,0,c)|c\in\Real\}$, $\gamma\neq 0$. This Lie algebra
 is conjugated to the Lie algebra
$\hol_{n=0}^2$. To see this it is enough to choose the new basis
$p_1,p_2,q_1-\frac{1}{2\gamma}p_2,q_2+\frac{1}{2\gamma}p_1$.

The lemma is proved. $\Box$

{\bf 2)} Let $n\geq 1$. We claim that the subalgebras of $\u(1,n+1)_{<p_1,p_2>}$ from the statement of the theorem are weakly-irreducible Berger algebras.
Indeed,  these subalgebras are weakly-irreducible. Table 3.2.2 shows that
all these subalgebras are Berger algebras (any element of each algebra can be obtained as $R(q_1\wedge q_2)$ for proper curvature tensor $R$),
this will follow also from Theorem \ref{G6metrics}.

We must prove that there are no other weakly-irreducible Berger subalgebras of\\ $\u(1,n+1)_{<p_1,p_2>}$.
For this we need all candidates for the weakly-irreducible subalgebras of $\u(1,n+1)_{<p_1,p_2>}$. Recall that in order to classify
weakly-irreducible subalgebras of $\su(1,n+1)_{<p_1,p_2>}$ we considered a Lie algebras homomorphism $\Gamma:\u(1,n+1)_{<p_1,p_2>}\to\lie(\Simil\H_n)$
and its restriction $\Gamma_{\su(1,n+1)_{<p_1,p_2>}}:\su(1,n+1)_{<p_1,p_2>}\to\lie(\Simil\H_n)$ which is an isomorphism.
We proved that if $\g\subset\u(1,n+1)_{<p_1,p_2>}$ is weakly-irreducible, then the subalgebra $\f=\Gamma(\g)\subset\lie(\Simil\H_n)$
satisfies a property. Then we found all  subalgebras $\f\subset\lie(\Simil\H_n)$ that satisfy this property and the isomorphism
$\Gamma_{\su(1,n+1)_{<p_1,p_2>}}$ gave us  the list of candidates for the  weakly-irreducible subalgebras of $\su(1,n+1)_{<p_1,p_2>}$.
Now, since $\ker\Gamma_{\su(1,n+1)_{<p_1,p_2>}}=\Real J$, any weakly-irreducible subalgebra of $\u(1,n+1)_{<p_1,p_2>}$
must have one of the forms $\g$, $\g\oplus\Real J$ or $\g^\xi$, where $\g$ is a  candidate to the weakly-irreducible subalgebras of $\su(1,n+1)_{<p_1,p_2>}$,
$\g^\xi=\{x+\xi(x)J|x\in\g\}$ and $\xi:\g\to \Real$ is a non-zero linear map such that $\g^\xi$ is a Lie algebra. Recall that for $m>0$ all candidates
to the  weakly-irreducible subalgebras of $\su(1,n+1)_{<p_1,p_2>}$ are weakly-irreducible subalgebras.

\begin{lem}\label{lem322n} Let $\g$ be one of the following candidates to the weakly-irreducible
subalgebra of $\su(1,n+1)_{<p_1,p_2>}$:
$$\g^{0,\h,\psi,k},\quad\g^{0,\h,\A^1,\zeta},\quad\g^{0,\h,\varphi,i_0,\zeta},\quad  \g^{0,\h,\psi,k,i_0,\zeta}$$ (these
Lie algebras defined as in the proof of Theorem \ref{Tw-ir}, see Lemmas \ref{lt1}, \ref{lt2} and \ref{lt4}). Then the Lie
algebras of the form $\g$, $\g^\xi$ and $\g\oplus\Real J$ are not Berger algebras.
\end{lem}

{\it Proof.}  Let $\g$ be one of the above Lie algebras and $R\in\R(\g\oplus\Real J)$.
Since  $\h\subset\sod(1,...,n)$, from Table 3.2.2 it follows  that $(\pr_{u(n)}R(\Real^{2,2n+2}\wedge \Real^{2,2n+2}))\cap\h=\{0\}$.
Therefore, if one of the Lie algebras $\g$, $\g^\xi$ and $\g\oplus\Real J$ is a Berger algebra, then $\h=\{0\}$.
This condition holds only for the Lie algebra $\g=\g^{0,\h=\{0\},\A^1,\zeta}$. We claim that $\R(\g\oplus\Real J)=\{0\}$.
Indeed, let $R\in \R(\g\oplus\Real J)$. We decompose $R$ as in Table 3.2.2.  For any $y\in E^1_{1,...,m}$ we have
$R(q_1\wedge y)=p_1\wedge \lambda_3 y+p_2\wedge  \lambda_3 Jy+M_3^*(y)p_1\wedge p_2\in\g$.
Consequently, $\lambda_3=0$ and $M_3=0$. It is easy to show in the same way that all the other components of $R$ are also zero.
Thus the Lie algebras of the form $\g$, $\g^\xi$ and $\g\oplus\Real J$ are not Berger algebras. The lemma is proved. $\Box$

\begin{lem}\label{lem322A} Let $\g\subset\u(1,n+1)_{<p_1,p_2>}$ be a weakly-irreducible  Berger subalgebra with the associated number $m$, $1\leq m<n$.
If $\g$ does not contain the set $\{p_1\wedge w+p_2\wedge Jw|w\in E^1_{m+1,...,n}\}$, then $\pr_{\u(m+1,...,n)}\g=\{0\}$.\end{lem}

{\it Proof.} Consider a curvature tensor $R\in\R(\g)$. We decompose $R$  using Table 3.2.2.
Since $m\geq 1$, we see that $ p_1\wedge p_2\in\g$.
Let $y\in  E^1_{m+1,...,n}$ be a vector such that
$p_1\wedge y+p_2\wedge Jy\not\in\g$, then $R(q_1\wedge y)=p_1\wedge \lambda_3 y+p_2\wedge  \lambda_3 Jy+M_3^*(y)p_1\wedge p_2\in\g$.
Consequently, $\lambda_3=0$. From  Table 3.2.2 it follows that $\pr_{\u(m+1,...,n)}R(\Real^{2,2n+2}\wedge \Real^{2,2n+2})=\{0\}$.
This proves the lemma. $\Box$

\begin{lem}\label{lem322Bn} Let $\g$ be a Lie algebra of the form $\g^{n,\h,\psi,k,l}$ or
$\g^{m,\h,\psi,k,l,r}$. Then
\begin{itemize}
\item[1)] If $\g$ is a Berger algebra, then $\h\subset\su(k)$.

\item[2)] If $\xi:\g\to\Real$ is not zero and $\g^{\xi}$ is a Berger algebra, then there exist elements
 $A=\bigl(\begin{smallmatrix}B&-C\\C&B\end{smallmatrix}\bigr)\in\u(n)$, $z_1,z_2\in\Real^n$ such that
$A+J_{k}-\frac{k}{n+2}J_n\in\z(\h)$, for $x=(0,-\frac{k}{n+2},B,C+J_{k}-\frac{k}{n+2}J_n,z_1,z_2,0)\in\g$ we have
$\xi(x)=\frac{k}{n+2}$,
$\xi$  is zero on the orthogonal complement to  $\Real x$, and the orthogonal complement to $\Real(A+J_{k}-\frac{k}{n+2}J_n)$ in $\u(n)$
is contained in $\u(k)$.
In particular, $x+\xi(x)J=(0,0,B,C+J_{k},z_1,z_2,0)$ and  $\pr_{\u(n)}\g^\xi\subset\u(k)$.

\item[3)] The Lie algebra $\g\oplus\Real J$ is not a Berger algebra.
\end{itemize}
\end{lem}

{\it Proof.} Statements  1) and 3) follow from Lemma \ref{lem322A}.
Let us prove Statement 2).  Suppose that $\g^\xi$ is a Berger algebra for some linear map $\xi:\g\to\Real$.
From Lemma \ref{lem322A}  it follows that  $\pr_{\u(k+1,...,n)}\g^\xi=\{0\}$. This can happen only in the case described in Statement 2).
In this situation the Lie algebra $\g^\xi$ is either  of the form $\hol^{n,\u,\psi_1,k_1,l_1}$ or
$\g^{m,\u,\psi_1,k_1,l_1,r_1}$. The lemma is proved.  $\Box$

Now we have to consider only the Lie algebras  of the form $\g$, $\g^\xi$ and $\g\oplus\Real J$ for
$\g=\g^{m,\h,\A_1}$ and $\g=\g^{m,\h,\varphi}$, where $0\leq m\leq n$ and $\h\subset\su(m)\oplus\Real(J_m-\frac{m}{n+2}J_n)\oplus\sod(m+1,...,n)$.
Lemma \ref{lem322A} yields that if any of these Lie algebras is a Berger algebra, then $\h\subset\su(m)\oplus\Real(J_m-\frac{m}{n+2}J_n)$.

Let us consider the Lie algebra $\g=\g^{m,\h,\A_1}$, where $\h\subset\su(m)\oplus\Real(J_m-\frac{m}{n+2}J_n)$. Obviously, $\g$ is a holonomy algebra of the form
$\hol^{m,\un,\A^1,\phi}$, where $\un=\pr_{\u(m)}\h$ and $\phi:\un\to\Real$ is the map linear map given by
$\phi:\left(\begin{smallmatrix}B&-C\\C&B\end{smallmatrix}\right)\mapsto-\frac{1}{2}\tr C$.
Consider a Lie algebra of the form $\g^\xi$, where $\xi:\g\to\Real$  is a non-zero linear map.
We have $\xi|_{\g'}=\xi|_{\h'\zr(\N^1+\N^{1,...,m}+\C)}=0$, i.e. $\xi$ can be considered as a linear map
$\xi:\A^1\oplus\z(\h)\to\Real$. If $\xi|_{\A^1}\neq 0$ and  $\xi|_{\z(\h)}= 0$, then
$\g^\xi$ is a holonomy algebra of the form $\hol^{m,\un,\varphi,\phi}$ or $\hol^{m,\un,\lambda}$.
Suppose that $\xi|_{\A^1}=0$ and let $A\in\h$ be an element such that $\xi(A)\neq 0$ and $\xi$ is zero
on the orthogonal complement to $\Real A$ in $\h$. Consider the decomposition $A=A_1+a(J_m-\frac{m}{n+2}J_n)$, where
$A_1\in\su(m)$ and $a\in\Real$. If $A_1\neq 0$, then $\g^\xi$ is a holonomy algebra of the form $\hol^{m,\un,\A^1,\phi}$.
Suppose that $A_1=0$. If $\xi(J_m-\frac{m}{n+2}J_n)=-J$, then $\g^\xi$ is a holonomy algebra of the form $\hol^{m,\un,\A^1,\t A^2}$,
otherwise, then $\g^\xi$ is a holonomy algebra of the form $\hol^{m,\un,\A^1,\phi}$.
Obviously, the Lie algebra $\g\oplus\Real J$ is a holonomy algebra  of the form $\hol^{m,\un,\A^1,\t\A_2}$ or $\hol^{m,\un,\A^1,\phi}$.

The case $\g=\g^{m,\h,\varphi}$ can be considered in the same way.

Thus the Lie algebras of the theorem exhaust weakly-irreducible Berger subalgebras of $\u(1,n+1)_{<p_1,p_2>}$.
The proof of the theorem will follow from Theorem \ref{G6metrics}. $\Box$

\section{Proof of Theorem \ref{G6metrics}}

Since the coefficients of each  our metric $g$ are polynomial functions, the Levi-Civita
connection given by $g$ is analytic and the Lie algebra $\hol_0$ is generated by the operators
$$R(X,Y)_0,\nabla R(X,Y;Z_1)_0,\nabla^2 R(X,Y;Z_1;Z_2)_0,...\in\so(T_0\Real^{2n+4},g_0),$$
where  $\nabla^r R(X,Y;Z_1;...;Z_r)=(\nabla_{Z_r}\cdots\nabla_{Z_1}R)(X,Y)$ and
 $X$, $Y$, $Z_1$, $Z_2$,... are vectors at the point $0$.

We consider the case $n>0$. The proof for the case $n=0$ can be obtained by simple computations.

First we consider some general metric and find all covariant derivatives of the curvature tensor of this metric.
Let $1\leq n_0\leq m\leq n$ be integers as in Section \ref{Ber+metr}. We will use the following convention about the ranks of the indices
$$\begin{array}{lll} a,b,c,d=1,...,2n+4,&i,j=3,...,2n+2,&\\
\i,\j=3,...,n_0+2,&\iv,\jv=n_0+3,...,m+2,&\ib,\jb=m+3,...,n+2,\end{array}$$
$$\begin{array}{l}
\ii,\jj=3,...,n_0+2,n+3,...,n+n_0+3,\\\iiv,\jjv=n_0+3,...,m+2,n+n_0+3,...,n+m+2,\\
\iib,\jjb=m+3,...,n+2,n+m+3,...,2n+2.\end{array}$$ We will use the Einstein rule for sums.

We assume that the numbers $B^\i_{\alpha \j}$ and $C^\i_{\alpha \j}$ equal $B^{\i-2}_{\alpha \j-2}$ and $C^{\i-2}_{\alpha \j-2}$, respectively
(here $(B^i_{\alpha j})_{i,j=1}^{n_0}$ and $(C^i_{\alpha j})_{i,j=1}^{n_0}$ are numbers as  in Section \ref{Ber+metr}).
Define the numbers  $A^i_{\alpha j}$  such that $A^\i_{\alpha \j}=B^\i_{\alpha \j}$, $A^{\i+n}_{\alpha \j+n}=B^{\i}_{\alpha \j}$,
$A^{\i+n}_{\alpha \j}=C^{\i}_{\alpha \j}$,  $A^{\i}_{\alpha \j+n}=-C^{\i}_{\alpha \j}$, and $A^i_{\alpha j}=0$ for other $i$ and $j$,
here $\alpha=1,...,N$.

Let $\varphi,\phi:\un\to\Real$ be two linear maps with $\varphi|_{\un'}=\phi|_{\un'}=0$.
Let the numbers $\varphi_\alpha$ and $\phi_\alpha$ ($N_1+1\leq\alpha\leq N$) be as in Section \ref{Ber+metr}.
 If $\varphi=\phi=0$, then we set $N_0=N+2$ and consider some numbers $\varphi_{N+1}$, $\varphi_{N+2}$, $\phi_{N+1}$, $\phi_{N+2}$.
If $\varphi=0$ and $\phi\neq 0$, then we set $N_0=N+1$, $\phi_{N+1}=0$ and consider a number $\varphi_{N+1}$.
If $\varphi\neq 0$ and $\phi= 0$, then we set $N_0=N+1$, $\varphi_{N+1}=0$ and consider a number $\phi_{N+1}$.
If $\varphi\neq 0$ and $\phi\neq 0$, then we set $N_0=N$.
Thus we get some numbers $(\varphi_\alpha)_{\alpha=N_1+1}^{N_0}$ and $(\phi_\alpha)_{\alpha=N_1+1}^{N_0}$.
Consider the following polynomials $$\hat\varphi(x^{2n+3})=\sum_{\alpha=N_1+1}^{N_0}\frac{1}{\alpha!}\varphi_\alpha(x^{2n+3})^{\alpha}\text{ and }
\hat\phi(x^{2n+3})=\sum_{\alpha=N_1+1}^{N_0}\frac{1}{\alpha!}\phi_\alpha(x^{2n+3})^{\alpha}.$$

Consider the metric $g$ given by \eqref{mkg} with the functions
$$f_i=f_i^\varphi+f_i^\phi+\hat f_i(x^\ii,x^{2n+3})+\tilde f_{i\,n_0+1}^m+\breve f_{i\,N_0+1\,\,m+1}^n\qquad (i=1,2,3),$$
where $\tilde f_{i\,n_0+1}^m$ and  $\breve f_{i\,N_0+1\,\,m+1}^n$ are functions as in Section \ref{Ber+metr},
$f_i^\varphi$ and  $f_i^\phi$ are functions defined as in Section \ref{Ber+metr} using $N_0$ instead of $N$,
and $\hat f_i(x^\ii,x^{2n+3})$ are some functions, $\hat f_3(x^\ii,x^{2n+3})=0$.

We assume that $f_1(0)=f_2(0)=f_3(0)=0$, then $g_0=\eta$ and we can identify the tangent space to $\Real^{2n+4}$ at $0$ with the vector space
$\Real^{2,2n+2}$ such that $\frac{\p}{\p x^{1}}|_0=p_1$, $\frac{\p}{\p x^{2}}|_0=p_2$, $\frac{\p}{\p x^{3}}|_0=e_1$,...,$\frac{\p}{\p x^{n+2}}|_0=e_n$,
$\frac{\p}{\p x^{n+3}}|_0=f_1$,...,$\frac{\p}{\p x^{2n+2}}|_0=f_n$, $\frac{\p}{\p x^{2n+3}}|_0=q_1$, $\frac{\p}{\p x^{2n+4}}|_0=q_2$.

For the non-zero Christoffel symbols of the metric $g$ we have
{\scriptsize
\begin{align}
 \Ga^1_{1\,2n+3}&=\frac{1}{2}\frac{\p f_1}{\p x^1},\quad
\Ga^1_{1\,2n+4}=\frac{1}{2}\frac{\p f_3}{\p x^1},\quad
 \Ga^1_{2\,2n+3}=\frac{1}{2}\frac{\p f_1}{\p x^2},\quad
 \Ga^1_{2\,2n+4}=\frac{1}{2}\frac{\p f_3}{\p x^2},\quad
 \Ga^1_{i\,2n+3}=\frac{1}{2}\frac{\p f_1}{\p x^i}, \label{mkga10}\\
\Ga^1_{i\,2n+4}&=\frac{1}{2}\left(-\frac{\p u^i}{\p x^{2n+3}}+\frac{\p f_3}{\p x^i}\right),\quad
\Ga^1_{2n+3\,2n+3}=\frac{1}{2}\left(\frac{\p f_1}{\p x^{2n+3}}
       +f_3\frac{\p f_1}{\p x^{2}}+f_1\frac{\p f_1}{\p x^1}\right),\label{mkga20}\\
\Ga^1_{2n+3\,2n+4}&=\frac{1}{2}\left(f_3\frac{\p f_3}{\p x^{2}}+f_1\frac{\p f_3}{\p x^1}\right),\quad
\Ga^1_{2n+4\,2n+4}=\frac{1}{2}\left(-\frac{\p f_2}{\p x^{2n+3}}
       +f_3\frac{\p f_2}{\p x^{2}}+f_1\frac{\p f_2}{\p x^1}\right),\label{mkga30}\\
\Ga^2_{1\,2n+3}&=\frac{1}{2}\frac{\p f_3}{\p x^1},\quad
\Ga^2_{1\,2n+4}=\frac{1}{2}\frac{\p f_2}{\p x^1},\quad
 \Ga^2_{2\,2n+3}=\frac{1}{2}\frac{\p f_3}{\p x^2},\quad
 \Ga^2_{2\,2n+4}=\frac{1}{2}\frac{\p f_2}{\p x^2},\quad
  \Ga^2_{i\,j}=\frac{1}{2}\left(\frac{\p u^i}{\p x^j}+\frac{\p u^j}{\p x^i}\right),\quad \label{mkga40}\\
\Ga^2_{i\,2n+3}&=\frac{1}{2}\left(\frac{\p u^i}{\p x^{2n+3}}+\frac{\p f_3}{\p x^i}\right),\quad
\Ga^1_{i\,2n+4}=\frac{1}{2}\left(\sum_{j=3}^{2n+2}u^j\left(\frac{\p u^i}{\p x^j}-\frac{\p u^j}{\p x^i}\right)
   +\frac{\p f_2}{\p x^i}\right),\label{mkga50}\\
\Ga^2_{2n+3\,2n+3}&=\frac{1}{2}\left(2\frac{\p f_3}{\p x^{2n+3}}+ \sum_{i=3}^{2n+2}u^i\frac{\p f_1}{\p
x^i}+\left(f_2-\sum_{i=3}^{2n+2}(u^i)^2\right)\frac{\p f_1}{\p x^2}+ f_3\frac{\p f_1}{\p x^1}\right),\label{mkga55}
\end{align}
\begin{align}
\Ga^2_{2n+3\,2n+4}&=\frac{1}{2}\left(\sum_{i=3}^{2n+2}u^i\left(-\frac{\p u^i}{\p x^{2n+3}}+\frac{\p f_3}{\p x^{i}}\right)+
\frac{\p f_2}{\p x^{2n+3}}+ \left(f_2-\sum_{i=3}^{2n+2}(u^i)^2\right)\frac{\p f_3}{\p x^2}+f_3\frac{\p f_3}{\p
x^1}\right),\label{mkga60}\\ \Ga^2_{2n+4\,2n+4}&=\frac{1}{2}\left(\sum_{i=3}^{2n+2}u^i\frac{\p f_2}{\p x^{i}}+
\left(f_2-\sum_{i=3}^{2n+2}(u^i)^2\right)\frac{\p f_2}{\p x^2}+f_3\frac{\p f_2}{\p x^{1}}\right),\quad
 \Ga^i_{j\,2n+4}=\frac{1}{2}\left(\frac{\p u^i}{\p x^j}-\frac{\p u^j}{\p x^i}\right), \label{mkga80}\\
 \Ga^i_{2n+3\,2n+3}&=\frac{1}{2}\left(-\frac{\p f_1}{\p x^i}+u^i\frac{\p f_1}{\p x^2}\right),\quad
\Ga^i_{2n+3\,2n+4}=\frac{1}{2}\left(\frac{\p u^i}{\p x^{2n+3}}-\frac{\p f_3}{\p x^i}+u^i\frac{\p f_3}{\p x^2}\right),\label{mkga90}\\
 \Ga^i_{2n+4\,2n+4}&=\frac{1}{2}\left(-\frac{\p f_2}{\p x^i}+u^i\frac{\p f_2}{\p x^2}\right),\quad
 \Ga^{2n+3}_{2n+3\,2n+3}=-\frac{1}{2}\frac{\p f_1}{\p x^1},\quad
 \Ga^{2n+3}_{2n+3\,2n+4}=-\frac{1}{2}\frac{\p f_3}{\p x^1}, \label{mkga100}\\
 \Ga^{2n+3}_{2n+4\,2n+4}&=-\frac{1}{2}\frac{\p f_2}{\p x^1},\quad
 \Ga^{2n+4}_{2n+3\,2n+3}=-\frac{1}{2}\frac{\p f_1}{\p x^2},\quad
 \Ga^{2n+4}_{2n+3\,2n+4}=-\frac{1}{2}\frac{\p f_3}{\p x^2},\quad
 \Ga^{2n+4}_{2n+4\,2n+4}=-\frac{1}{2}\frac{\p f_2}{\p x^2}.\label{mkga110}
\end{align}}

Note that if $a\not\in\{1,2\}$, then  $\Ga^a_{1b}=\Ga^a_{2b}=0$. This means that the holonomy algebra $\hol_0$ of  the metric $g$ at the point $0$
preserves the vector subspace $\Real p_1\oplus \Real p_2\subset \Real^{2,2n+2}=T_0\Real^{2n+4}$, hence $\hol_0$ is contained in $\so(2,2n+2)_{<p_1,p_2>}$, where
{\scriptsize
$$\so(2,2n+2)_{<p_1,p_2>}=\left\{\left.\left (\begin{array}{ccccc}
a_{11}&a_{12} &-X^t &0 &-c\\
a_{21}&a_{22} &-Y^t &  c &0\\
0 &0&A&X&Y\\
0&0&0&-a_{11}&-a_{21}\\
0&0&0&-a_{12}&-a_{22}\\\end{array}\right)\right|\,\begin{array}{c}
\left(\begin{smallmatrix}a_{11}&a_{12}\\a_{21}&a_{22}\end{smallmatrix}\right)\in\gl(2),\\c\in \Real,\\ X,Y\in \Real^{2n},\\ A\in\so(2n) \end{array}\right\}.$$}
In particular, it is enough to compute the following components of the covariant derivatives of the curvature tensor:
$R^a_{bcd;a_1;a_2;...}$, $R^a_{icd;a_1;a_2;...}$, $R^i_{jcd;a_1;a_2;...}$, where $a,b\in\{1,2\}$ and  $3\leq i,j\leq 2n+2$.
We have
{\scriptsize
\begin{multline}
R^1_{1\,2n+3\,2n+4}=R^2_{2\,2n+3\,2n+4}=\sum_{\alpha=N_1+1}^{N_0}\frac{1}{(\alpha-1)!}\varphi_\alpha(x^{2n+3})^{\alpha-1},\\
 R^1_{1ab}=R^2_{2ab}=0      \text{ if } \{a\}\cup\{b\}\neq\{2n+3,2n+4\},\label{mkR10}\end{multline}
\begin{multline}
R^2_{1\,2n+3\,2n+4}=-R^1_{2\,2n+3\,2n+4}=\sum_{\alpha=N_1+1}^{N_0}\frac{1}{(\alpha-1)!}\phi_\alpha(x^{2n+3})^{\alpha-1},\\
R^2_{1ab}=R^1_{2ab}=0       \text{ if } \{a\}\cup\{b\}\neq\{2n+3,2n+4\},\label{mkR20}
\end{multline}
\begin{align}
R^1_{ij\,2n+3}&=\frac{1}{2}\frac{\p^2 f_1}{\p x^i\p x^j}-\frac{1}{4}\left(\frac{\p u^i}{\p x^j}+\frac{\p u^j}{\p x^i}\right)\frac{\p f_1}{\p x^2},\\
R^1_{ij\,2n+4}&=-\frac{1}{2}\frac{\p^2 u^i}{\p x^j\p x^{2n+3}}+\frac{1}{2}\frac{\p^2 f_3}{\p x^i\p x^j}
-\frac{1}{4}\left(\frac{\p u^i}{\p x^j}+\frac{\p u^j}{\p x^i}\right)\frac{\p f_3}{\p x^2},\label{mkR30}\\
R^2_{ij\,2n+3}&=-\frac{1}{2}\frac{\p^2 u^j}{\p x^i\p x^{2n+3}}+\frac{1}{2}\frac{\p^2 f_3}{\p x^i\p x^j}
-\frac{1}{4}\left(\frac{\p u^i}{\p x^j}+\frac{\p u^j}{\p x^i}\right)\frac{\p f_3}{\p x^2},\label{mkR40}
\end{align}
\begin{multline}R^2_{ij\,2n+4}=\frac{1}{4}\sum_{i_1=3}^{2n+2}\left(\frac{\p u^{i_1}}{\p x^i}-\frac{\p u^{i}}{\p x^{i_1}}\right)
\left(\frac{\p u^{j}}{\p x^{i_1}}-\frac{\p u^{i_1}}{\p x^{j}}\right)+\frac{1}{2}\sum_{i_1=3}^{2n+2}
u^{i_1}\left(\frac{\p^2 u^{i}}{\p x^{i_1}\p x^j}-\frac{\p^2 u^{i_1}}{\p x^{i}\p x^j}\right)\\
+\frac{1}{2}\frac{\p^2 f_2}{\p x^i\p x^j}-\frac{1}{4}\left(\frac{\p u^i}{\p x^j}+\frac{\p u^j}{\p x^i}\right)\frac{\p f_2}{\p x^2},\label{mkR50}
\end{multline}
\begin{multline}
R^1_{i\,2n+3\,2n+4}=\frac{1}{4}\left(-2\frac{\p^2 u^j}{(\p x^{2n+3})^2}+2\frac{\p^2 f_3}{\p x^i\p x^{2n+3}}-
\sum_{j=3}^{2n+2}\left(\frac{\p u^i}{\p x^j}-\frac{\p u^j}{\p x^i}\right)\frac{\p f_1}{\p x^j}\right.\\
+\left(\sum_{j=3}^{2n+2}u^j\left(\frac{\p u^i}{\p x^j}-\frac{\p u^j}{\p x^i}\right)+\frac{\p f_2}{\p x^i}\right)\frac{\p f_1}{\p x^2}\\
\left.-\frac{\p u^i}{\p x^{2n+3}}\frac{\p f_3}{\p x^2}-\frac{\p u^i}{\p x^{2n+3}}\frac{\p f_1}{\p x^1}-\frac{\p f_1}{\p x^i}\frac{\p f_3}{\p x^1}
-\frac{\p f_3}{\p x^i}\frac{\p f_3}{\p x^2}+\frac{\p f_3}{\p x^i}\frac{\p f_1}{\p x^1}\right),
\label{mkR60}\end{multline}
\begin{multline}
R^2_{i\,2n+3\,2n+4}=\frac{1}{4}\left(2\sum_{j=3}^{2n+2}u^j\left(\frac{\p^2 u^i}{\p x^j\p x^{2n+3}}-\frac{\p^2 u^j}{\p x^i\p x^{2n+3}}\right)
+\sum_{j=3}^{2n+2}
\left(\frac{\p u^j}{\p x^{2n+3}}-\frac{\p f_3}{\p x^j} \right)
\left(\frac{\p u^i}{\p x^j}-\frac{\p u^j}{\p x^i}\right)\right.\\
+\frac{\p f_2}{\p x^i}\frac{\p f_3}{\p x^2}
+\frac{\p^2 f_2}{\p x^i\p x^{2n+3}}
\left.-\frac{\p u^i}{\p x^{2n+3}}\frac{\p f_2}{\p x^2}-\frac{\p u^i}{\p x^{2n+3}}\frac{\p f_3}{\p x^1}-\frac{\p f_1}{\p x^{i}}\frac{\p f_2}{\p x^1}
-\frac{\p f_3}{\p x^i}\frac{\p f_2}{\p x^2}-\frac{\p f_3}{\p x^i}\frac{\p f_3}{\p x^1}\right),
\label{mkR65}\end{multline}
\begin{align}
R^i_{j\,2n+3\,2n+4}=&\sum_{\alpha=1}^{N}\frac{1}{(\alpha-1)!}A^i_{\alpha j}(x^{2n+3})^{\alpha-1},\label{mkR70}\\
R^\ib_{\jjb\,2n+3\,2n+4}=&-R^\jjb_{\ib\,2n+3\,2n+4}=\delta^{\ib+n}_{\jjb}\sum_{\alpha=N_1+1}^{N_0}\frac{1}{(\alpha-1)!}
\phi_\alpha(x^{2n+3})^{\alpha-1},
\quad R^\ib_{\jb\,2n+3\,2n+4}=R^\iib_{\jjb\,2n+3\,2n+4}=0,\label{n47a}\\
R^i_{jab}=&0   \text{ if } \{a\}\cup\{b\}\neq\{2n+3,2n+4\},\label{n47b}\\
R^i_{j\,2n+3\,2n+4}=&0   \text{ if } i,j\not\in\{3,...,n_0+2,n+3,...,n+n_0+2\} \text{ or } i,j\not\in\{m+3,...,n+2,n+m+3,...,2n+2\}.\label{n47c}
\end{align}}

Specifically,
\begin{align}
R^1_{\iv\,\jv\,2n+3}&=R^2_{n+\iv\,\jv\,2n+3}=\delta_{\iv\jv},\quad R^1_{n+\iv\,\jv\,2n+3}=-R^2_{\iv\,\jv\,2n+3}=0,\label{mkR80}\\
R^1_{\iv\,\jv\,2n+4}&=R^2_{n+\iv\,\jv\,2n+4}=0,\quad R^1_{n+\iv\,\jv\,2n+4}=-R^2_{\iv\,\jv\,2n+4}=-\delta_{\iv\jv},\label{mkR90}\\
R^1_{\iiv ab}&=R^2_{\iiv ab}=0, \text{ if } \{a\}\cup\{b\}\not\subset \{\jjv,2n+3\}\cup \{\jjv,2n+4\} \text{ for some }\jjv, \label{mkR100}
\end{align}
\begin{align}
R^1_{\ib\,2n+3\,2n+4}&=R^2_{n+\ib\,2n+3\,2n+4}=\sum_{\jb=m+1}^{n}\frac{1}{(N_0+\jb-m-1)!}(x^{2n+2})^{N_0+\jb-m-1},\label{mkR105}\\
R^1_{n+\ib\,2n+3\,2n+4}&=R^2_{\ib\,2n+3\,2n+4}=0,\quad R^1_{\iib ab}=R^2_{\iib ab}=0 \text{ if } \{a\}\cup\{b\}\neq\{2n+3,2n+4\}.\label{mkR110}
\end{align}
We also wish to have \begin{equation}
R^1_{\i ab}=R^2_{n+\i\, ab}\text{ and }R^1_{n+\i\, ab}=-R^2_{\i\, ab}.\end{equation}
The computations shows that these equalities hold if we choose
\begin{equation} \hat f_i=f^0_i\quad (i=1,2,3),\end{equation}
where the functions $f^0_i$ are as in Section \ref{Ber+metr}.
In particular,
\begin{equation} R^1_{\ii\,\jj\,2n+4}=\sum_{\alpha=1}^{N}\frac{1}{(\alpha-1)!}A^\ii_{\alpha \jj}(x^{2n+3})^{\alpha-1}. \label{mkR120}\end{equation}
Thus,
\begin{equation} R^1_{iab}=R^2_{n+i\,ab}\text{ and } R^1_{n+i\,ab}=-R^2_{i\,ab}, \text{ where } 3\leq i\leq n+2 \label{mkR130}\end{equation}

To compute  the covariant derivatives of the curvature tensor we will need the following Christoffel symbols
\begin{align} \Ga^a_{1b}&=\Ga^a_{2b}=0 \text{ if } a\not\in\{1,2\},\quad \Ga^{2n+3}_{ia}=\Ga^{2n+4}_{ia}=0,\quad \Ga^i_{jj_1}=\Ga^i_{j \,2n+3}=0, \label{mkga200}\\
\Ga^\ii_{\jj\,2n+3}&=\sum_{\alpha=1}^N\frac{1}{\alpha!}A^\ii_{\alpha \jj}(x^{2n+3})^\alpha,\label{mkga210}
\quad \Ga^\ib_{\jjb\,2n+3}=-\Ga^\jjb_{\ib\,2n+3}=\delta^{\ib+n}_\jjb\hat\phi,\end{align}
\begin{multline}\Ga^{2n+3}_{2n+3\,2n+3}=-\Ga^{2n+3}_{2n+4\,2n+4}=\Ga^{2n+4}_{2n+3\,2n+4}=\hat\phi,\\
-\Ga^{2n+3}_{2n+3\,2n+4}=\Ga^{2n+4}_{2n+3\,2n+3}=-\Ga^{2n+4}_{2n+4\,2n+4}=\hat\varphi.\label{mkga230}\end{multline}

\begin{lem}\label{mklem1}
We have
\begin{itemize}
\item[{\rm 1)}] $R^\ii_{\jj ab;a_1;...;a_r}=\sum_{t\in T_{ab a_1...a_r}}z_tA^\ii_{t\jj}$, where $T_{ab a_1...a_r}$ is a finite set of indices,
$z_t$ are functions and $A_t\in\un$;

\item[{\rm 2)}] If $\varphi\neq 0$, then
$R^1_{1ab;a_1;\cdots;a_r}=R^2_{2ab;a_1;\cdots;a_r}=\sum_{t\in T_{ab a_1...a_r}}z_t\varphi(A_t)$;

If $\varphi=0$, then $R^1_{1ab;a_1;\cdots;a_r}=R^2_{2ab;a_1;\cdots;a_r}$ and this equals $0$ for $0\leq r\leq N-1$;

\item[{\rm 3)}] If $\phi\neq 0$, then
$R^2_{1ab;a_1;\cdots;a_r}=-R^1_{2ab;a_1;\cdots;a_r}=\sum_{t\in T_{ab a_1...a_r}}z_t\phi(A_t)$,
$\,\,\,\,R^\ib_{\jjb ab;a_1;\cdots;a_r}=-R^\jjb_{\ib ab;a_1;\cdots;a_r}=\delta^{\ib+n}_\jjb\sum_{t\in T_{ab a_1...a_r}}z_t\phi(A_t)$;

If $\phi=0$, then $R^2_{1ab;a_1;\cdots;a_r}=-R^1_{2ab;a_1;\cdots;a_r}$, $R^\ib_{\jjb ab;a_1;\cdots;a_r}=-R^\jjb_{\ib ab;a_1;\cdots;a_r}$
and these components equals $0$ for $0\leq r\leq N-1$;

\item[{\rm 4)}] $R^1_{iab;a_1;...;a_r}=R^2_{n+i\,ab;a_1;...;a_r}$ and $R^1_{n+i\,ab;a_1;...;a_r}=-R^2_{i\,ab;a_1;...;a_r}$, where $3\leq i\leq n+2$;

\item[{\rm 5)}] $R^1_{n+\ib\,ab;a_1;...;a_r}=R^2_{\ib ab;a_1;...;a_r}=0$;

\item[{\rm 6)}] $R^i_{jab;a_1;...;a_r}=0   \text{ if } i,j\not\in\{3,...,n_0+2,n+3,...,n+n_0+2\} \text{ or } i,j\not\in\{m+3,...,n+2,n+m+3,...,2n+2\}$.

\end{itemize}\end{lem}

{\it Proof.} The lemma can be easily proved using the induction and equalities \eqref{mkR10}, \eqref{mkR20}, (\ref{mkR70}--\ref{n47c}), \eqref{mkR110} and \eqref{mkR130}.

For example, let us prove  Part 1) and Part 2) of the lemma  for $a_r=2n+4$. Let $r\geq 1$
Suppose that the lemma is true for all $s<r$. Let $a_r=2n+4$. Suppose that $\varphi\neq 0$.

We have\\ $R^i_{jbc;a_1;...;a_{r-1};2n+4}=$
$$\begin{array}{l}\frac{\p R^i_{jbc;a_1;...;a_{r-1}}}{\p x^{2n+4}}+\Ga^i_{a\,2n+4} R^a_{jbc;a_1;...;a_{r-1}}-
\Ga^a_{j\,2n+4} R^i_{abc;a_1;...;a_{r-1}}-\Ga^a_{b\,2n+4} R^i_{jac;a_1;...;a_{r-1}}\\-\Ga^a_{c\,2n+4} R^i_{jba;a_1;...;a_{r-1}}
-\Ga^a_{a_1\,2n+4} R^i_{jbc;a;a_2...;a_{r-1}}-...-\Ga^a_{a_{r-1}\,2n+4} R^i_{jbc;a_1;...;a_{r-2};a}.\end{array}$$
Since  $\hol_0\subset\so(2,2n+2)_{<p_1,p_2>}$, we have
$R^{2n+3}_{jbc;a_1;...;a_{r-1}}=R^{2n+4}_{jbc;a_1;...;a_{r-1}}=R^i_{1bc;a_1;...;a_{r-1}}=R^i_{2bc;a_1;...;a_{r-1}}=0$.
Using this, \eqref{mkga200}, \eqref{mkga210} and the induction hypotheses, we get\\
$R^i_{jbc;a_1;...;a_{r-1};2n+4}=$
$$\begin{array}{l}\sum_{t\in T_{bc a_1...a_{r-1}}}\frac{\p z_t}{\p x^{2n+4}}A^i_{tj}+
\sum_{\alpha=1}^{N}\frac{1}{\alpha!}\sum_{t\in T_{bc a_1...a_r}}(x^{2n+3})^\alpha z_t[A_\alpha,A_t]_j^i\\
-\sum_{a=1}^{2n+4}\sum_{t\in T_{ac;a_1;...;a_{r-1}}}\Ga^a_{b\,2n+4}z_t A_{tj}^i-\sum_{a=1}^{2n+4}\sum_{t\in T_{ba;a_1;...;a_{r-1}}}\Ga^a_{c\,2n+4}z_t A_{tj}^i\\
-\sum_{a=1}^{2n+4}\sum_{t\in T_{bc;a;a_2;...;a_{r-1}}}\Ga^a_{a_1\,2n+4}z_t A_{tj}^i-...
-\sum_{a=1}^{2n+4}\sum_{t\in T_{bc;a_1;...;a_{r-2};a}}\Ga^a_{a_{r-1}\,2n+4}z_t A_{tj}^i.\end{array}$$

We also have\\$R^1_{1bc;a_1;...;a_{r-1};2n+4}=$
 $$\begin{array}{l}\frac{\p R^1_{1bc;a_1;...;a_{r-1}}}{\p x^{2n+4}}+\Ga^1_{a\,2n+4} R^a_{1bc;a_1;...;a_{r-1}}-
\Ga^a_{1\,2n+4} R^1_{abc;a_1;...;a_{r-1}}-\Ga^a_{b\,2n+4} R^1_{1ac;a_1;...;a_{r-1}}\\-\Ga^a_{c\,2n+4} R^1_{1ba;a_1;...;a_{r-1}}
-\Ga^a_{a_1\,2n+4} R^1_{1bc;a;a_2...;a_{r-1}}-...-\Ga^a_{a_{r-1}\,2n+4} R^1_{1bc;a_1;...;a_{r-2};a}.\end{array}$$
Since  $\hol_0\so(2,2n+2)_{<p_1,p_2>}$, we see that
$R^{a}_{1bc;a_1;...;a_{r-1}}=\Ga^a_{1\,2n+4}=0$ if $a\not\in\{1,2\}$. Hence,
$$\begin{array}{l}\Ga^1_{a\,2n+4} R^a_{1bc;a_1;...;a_{r-1}}- \Ga^a_{1\,2n+4} R^1_{abc;a_1;...;a_{r-1}}\\
=\Ga^1_{2\,2n+4} R^2_{1bc;a_1;...;a_{r-1}}-\Ga^2_{1\,2n+4} R^1_{2bc;a_1;...;a_{r-1}}\\
=(\Ga^1_{2\,2n+4}+\Ga^2_{1\,2n+4})R^2_{1bc;a_1;...;a_{r-1}}=0,\end{array}$$ where we used  Statement 3) for $r-1$
and the fact that $\Ga^1_{2\,2n+4}+\Ga^2_{1\,2n+4}=0$.
Thus,\\ $R^1_{1bc;a_1;...;a_{r-1};2n+4}=$
$$\begin{array}{l}\sum_{t\in T_{bc a_1...a_{r-1}}}\frac{\p z_t}{\p x^{2n+4}}\varphi(A_t)\\-
\sum_{a=1}^{2n+4}\sum_{t\in T_{ac;a_1;...;a_{r-1}}}\Ga^a_{b\,2n+4}z_t \varphi(A_t)
-\sum_{a=1}^{2n+4}\sum_{t\in T_{ba;a_1;...;a_{r-1}}}\Ga^a_{c\,2n+4}z_t \varphi(A_t)\\-
\sum_{a=1}^{2n+4}\sum_{t\in T_{bc;a;a_2;...;a_{r-1}}}\Ga^a_{a_1\,2n+4}z_t \varphi(A_t)-...
-\sum_{a=1}^{2n+4}\sum_{t\in T_{bc;a_1;...;a_{r-2};a}}\Ga^a_{a_{r-1}\,2n+4}z_t \varphi(A_t).\end{array}$$
In the same way we can compute $R^2_{2bc;a_1;...;a_{r-1};2n+4}$. Now the statement follows from the induction hypotheses and the fact that $\varphi|_{\un'}=0$.
$\Box$

\begin{lem}\label{mklem2}
If one of the numbers $a$, $b$, $a_1$,...,$a_r$ belongs to the set $\{1,...,2n+2\}$, then $$R^i_{jab;a_1;...;a_r}=0.$$
\end{lem}
{\it Proof.} To prove the lemma it is enough to prove  the following 3 statements
\begin{itemize}
\item[{\rm 1)}] If $1\leq c\leq 2n+2$, then $\frac{\p R^i_{jab;a_1;...;a_{r-1}}}{\p x^c}=0$;

\item[{\rm 2)}] If $1\leq c\leq 2n+2$, then $R^i_{jab;a_1;...;a_{r-1};c}=0$;

\item[{\rm 3)}] If for fixed $a$, $b$, $a_1$,...,$a_{r-1}$ (and for  all $1\leq i,j\leq j$)
we have $R^i_{jab;a_1;...;a_{r-1}}=0$, then $R^i_{jab;a_1;...;a_{r-1};a_r}=0$.
These statements can be proved  using the induction,  \eqref{mkR70} and  \eqref{mkga200}. $\Box$
\end{itemize}

\begin{lem}\label{mklem3} For all $0\leq r\leq N-1$ we have
$$R^\ii_{\jj\,2n+3\,2n+4;2n+3;...;2n+3(r \text{ times})}(0)=A^\ii_{r+1\,\jj}+\sum_{\alpha=1}^{r}\mu_{r\alpha} A^\ii_{\alpha \jj},$$ where $\mu_{r\alpha}$ are some numbers.
\end{lem}
{\it Proof.} We will prove the following three statements
 \begin{itemize}
\item[{\rm 1)}] If $0\leq r\leq N_1-1$, then $$R^\ii_{\jj\,2n+3\,2n+4;2n+3;...;2n+3(r \text{ times})}=\sum_{\alpha=r+1}^N\frac{1}{(\alpha-r-1)!}
A^\ii_{\alpha \jj}(x^{2n+3})^{\alpha-r-1}+
\sum_{\beta=1}^{N}H_\beta(x^{2n+3})A_{\beta \jj}^\ii,$$ where $H_\beta(x^{2n+3})$ are polynomials of $x^{2n+3}$ such that
$H_\beta(0)=H'_\beta(0)=\cdots=H_\beta^{(N_1-r)}(0)=0$;

\item[{\rm 2)}] For all $0\leq r\leq N-1$ we have $$R^\ii_{\jj\, 2n+3\,2n+4;a_1;...;a_r}=\sum_{\alpha=1}^{r+1} G_\alpha(x^{2n+3})A^\ii_{\alpha \jj}+
\sum_{\alpha=r+2}^{N} F_\alpha(x^{2n+3})A^\ii_{\alpha \jj},$$ where $G_\alpha(x^{2n+3})$ and $F_\alpha(x^{2n+3})$  are polynomials of $x^{2n+3}$ such that
$F_{r+2}(0)=0$, $F_{r+3}(0)=F'_{r+3}(0)=0$,...,$F_N(0)=F'_N(0)=\cdots=F_N^{(N-r-2)}(0)=0$;

\item[{\rm 3)}] If $N_1\leq r\leq N-1$, then\\ $R^\ii_{\jj\,2n+3\,2n+4;2n+3;...;2n+3(r \text{ times})}=$
$$\sum_{\alpha=r+1}^N\frac{1}{(\alpha-r-1)!}A^\ii_{\alpha \jj}(x^{2n+3})^{\alpha-r-1}+\sum_{\beta=1}^{N_1}W_\beta(x^{2n+3})A_{\beta \jj}^\ii+
\sum_{\gamma=N_1+1}^{N}Q_\gamma(x^{2n+3})A_{\gamma \jj}^\ii,$$ where $W_\beta(x^{2n+3})$ and $Q_\gamma(x^{2n+3})$ polynomials of $x^{2n+3}$ such that
$Q_{r+1}(0)=0$, $Q_{r+2}(0)=Q'_{r+2}(0)=0$,...,$Q_N(0)=Q'_N(0)=\cdots=Q_N^{(N-r-1)}(0)=0$.
\end{itemize}

Let us prove Statement 1). For $r=0$ the statement follows from \eqref{mkR70}. Let $r>0$ and suppose that  Statement 1)  holds for all $s<r$.
We have\\
$R^\ii_{\jj\,2n+3\,2n+4;2n+3;...;2n+3(r \text{ times})}=$
$$\begin{array}{l}\frac{\p R^\ii_{\jj\,2n+3\,2n+4;2n+3;...;2n+3(r-1 \text{ times})}}{\p x^{2n+3}}+
\Ga^\ii_{a\,2n+3}R^a_{\jj\,2n+3\,2n+4;2n+3;...;2n+3(r-1 \text{ times})}\\
-\Ga^a_{\jj\,2n+3}R^\ii_{a\,2n+3\,2n+4;2n+3;...;2n+3(r-1 \text{ times})}\\-\Ga^a_{2n+3\,2n+3}R^\ii_{\jj\,a\,2n+4;2n+3;...;2n+3(r-1 \text{ times})}
-\Ga^a_{2n+4\,2n+3}R^\ii_{\jj\,2n+3\,a;2n+3;...;2n+3(r-1 \text{ times})}\\
-\Ga^a_{2n+3\,2n+3}R^\ii_{\jj\,2n+3\,2n+4;a;2n+3;...;2n+3(r-2 \text{ times})}
\\-\cdots-\Ga^a_{2n+3\,2n+3}R^\ii_{\jj\,2n+3\,2n+4;2n+3;...;2n+3(r-2 \text{ times});a}.\end{array}$$
Using Lemma \ref{mklem2}, \eqref{mkga200} and \eqref{mkga230} we get
\begin{multline} R^\ii_{\jj\,2n+3\,2n+4;2n+3;...;2n+3(r \text{ times})}=\\
\begin{array}{l}\frac{\p R^\ii_{\jj\,2n+3\,2n+4;2n+3;...;2n+3(r-1 \text{ times})}}{\p x^{2n+3}}
-(r+1)\hat\phi(x^{2n+3})R^\ii_{\jj\,2n+3\,2n+4;2n+3;...;2n+3(r-1 \text{ times})}\\
-\hat\varphi(x^{2n+3})(R^\ii_{\jj\,2n+3\,2n+4;2n+4;2n+3;...;2n+3(r-2 \text{ times})}\\
+\cdots+R^\ii_{\jj\,2n+3\,2n+4;2n+3;...;2n+3(r-2 \text{ times});2n+4}).\end{array}
\label{mk400}\end{multline}
Statement 1) follows from the induction hypotheses and the fact that $\hat\varphi(0)=\hat\varphi'(0)=\cdots=\hat\varphi^{(N_1-1)}(0)=0$.

Statement 2) can by proved by analogy (by Lemma \ref{mklem2} we may  consider only $a_r=2n+3$ and $2n+4$).
Statement 3) can be proved using \eqref{mk400} and Statement 2).

From Statement 1) it follows that for  $0\leq r\leq N_1-1$ we have
$$R^\ii_{\jj\,2n+3\,2n+4;2n+3;...;2n+3(r \text{ times})}(0)=A^\ii_{r+1\,\jj}.$$ The end of the proof of the lemma follows from this and Statement 3). $\Box$

From Statement 1) of Lemma \ref{mklem1} and Lemma \ref{mklem3} it follows that $\pr_{\so(2m)}\hol_0=\un$.

\begin{lem}\label{mklem4} For all $0\leq r\leq N-1$ we have
$$R^1_{\ii\jj\,2n+4;2n+3;...;2n+3(r \text{ times})}(0)=A^\ii_{r+1\,\jj}+\sum_{\alpha=1}^{r}\nu_{r\alpha} A^\ii_{\alpha \jj},$$ where $\nu_{r\alpha}$ are some numbers.
\end{lem}
{\it Proof.} The lemma can be proved by analogy to the proof of Lemma \ref{mklem3} using \eqref{mkR120}. $\Box$

Using \eqref{mkR105}, we  get
\begin{multline}\label{mk410}
R^1_{\ib\,2n+3\,2n+4;2n+3;...;2n+3(r \text{ times})}(0)=R^2_{n+\ib\,2n+3\,2n+4;2n+3;...;2n+3(r \text{ times})}(0)\\
=\delta_{\ib\,r-N_0+m+3} \text{ for all } N_0\leq r\leq N_0+n-m-1.
\end{multline}

Consider the Lie algebra $\hol^{m,\un,\A^1,\t\A^2}$. We take  $\varphi=\phi=0$, $N_0=N+2$, $\varphi_{N+1}=1$, $\varphi_{N+2}=0$, $\phi_{N+1}=0$ and $\phi_{N+2}=1$.
Then the above metric coincides with the metric from Table \ref{tabhol}   for the Lie algebra $\hol^{m,\un,\A^1,\t\A^2}$. From Lemma \ref{mklem1} it follows that
$\hol_0\subset\hol^{m,\un,\A^1,\t\A^2}$. From Lemma  \ref{mklem3} it follows that $\un\subset\hol_0$.
From \eqref{mkR80}, \eqref{mkR90}, \eqref{mk410}, Lemma \ref{mklem4} and the fact that $\un$ does not annihilate
any proper subspace of $E_{1,...,n_0}$ it follows that $\N^1+\N^2\subset\hol_0$.
It can be shown that $$R^1_{1\,2n+3\,2n+4;2n+3;...;2n+3(N+1 \text{ times})}(0)=R^2_{2\,2n+3\,2n+4;2n+3;...;2n+3(N+1 \text{ times})}(0)=1,$$ $$
R^2_{1\,2n+3\,2n+4;2n+3;...;2n+3(N+1 \text{ times})}(0)=R^1_{2\,2n+3\,2n+4;2n+3;...;2n+3(N+1 \text{ times})}(0)=0,$$
$$R^1_{1\,2n+3\,2n+4;2n+3;...;2n+3(N+2 \text{ times})}(0)=R^2_{2\,2n+3\,2n+4;2n+3;...;2n+3(N+2 \text{ times})}(0)=0,$$
$$R^2_{1\,2n+3\,2n+4;2n+3;...;2n+3(N+2 \text{ times})}(0)=-R^1_{2\,2n+3\,2n+4;2n+3;...;2n+3(N+2 \text{ times})}(0)=1,$$
$$R^\ib_{\jjb\,2n+3\,2n+4;2n+3;...;2n+3(N+2 \text{ times})}(0)=-R^\jjb_{\ib\,2n+3\,2n+4;2n+3;...;2n+3(N+2 \text{ times})}(0)=\delta^{\ib+n}_\jjb.$$

Hence, $\A^1+\t\A^2\subset\hol_0$. The inclusion $\C\subset\hol_0$ is obvious. Thus, $\hol_0=\hol^{n,\un,\A^1,\t\A^2}$.

The Lie algebras  $\hol^{m,\un,\A^1,\phi}$, $\hol^{m,\un,\varphi,\t\A^2}$, $\hol^{m,\un,\varphi,\phi}$ and $\hol^{m,\un,\lambda}$,
can be considered in the same way.

Now we are left with the Lie algebras $\hol^{n,\un,\psi,k,l}$ and $\hol^{m,\un,\psi,k,l,r}$. For them we can use the following lemma.
Let ${\bar i}=n_0+3,...,n+2,n+n_0+3,...,2n+2$.

\begin{lem}\label{mklem5} Consider the following metrics on  $\Real^{2n+4}$
$$\begin{array}{rl} g=&2dx^1dx^{2n+3}+2dx^2dx^{2n+4}+
\sum^{2n+2}_{i=3}(dx^i)^2+2\sum^{2n+2}_{i=3}u^i(x^\ii,x^{2n+3})dx^{i} dx^{2n+4}\\&
+(\hat f_1(x^\ii,x^{2n+3})+\bar f_1(x^{\bar i},x^{2n+3}))(dx^{2n+3})^2+(\hat f_2(x^\ii,x^{2n+3})+\bar f_2(x^{\bar i},x^{2n+3})) (dx^{2n+4})^2\\
&+2(\hat f_3(x^\ii,x^{2n+3})+\bar f_3(x^{\bar i},x^{2n+3})) dx^{2n+3}dx^{2n+4},\\
g_1=&2dx^1dx^{2n+3}+2dx^2dx^{2n+4}+\sum^{2n+2}_{i=3}(dx^i)^2+2\sum^{2n+2}_{i=3}u^i(x^\ii,x^{2n+3})dx^{i} dx^{2n+4}\\&
+\hat f_1(x^\ii,x^{2n+3})(dx^{2n+3})^2+\hat f_2(x^\ii,x^{2n+3})(dx^{2n+4})^2+2\hat f_3(x^\ii,x^{2n+3})dx^{2n+3}dx^{2n+4},\\
g_2=&2dx^1dx^{2n+3}+2dx^2dx^{2n+4}+\sum^{2n+2}_{i=3}(dx^i)^2\\&
+\bar f_1(x^{\bar i},x^{2n+3})(dx^{2n+3})^2+\bar f_2(x^{\bar i},x^{2n+3}) (dx^{2n+4})^2+2\bar f_3(x^{\bar i},x^{2n+3}) dx^{2n+3}dx^{2n+4}.\end{array}$$
Let $R$, $\hat R$ and $\bar R$ be the corresponding curvature tensors, then $R=\hat R+\bar R$.
\end{lem}

{\it Proof.} Using  equalities \eqref{mkga10} -- \eqref{mkga110}, it is easy to see that for the corresponding Christoffel symbols
we have $\Ga^a_{bc}=\hat\Ga^a_{bc}+\bar\Ga^a_{bc}$, $\hat\Ga^a_{bc}\bar\Ga^{a_1}_{ab_1}=0$ and $\bar\Ga^a_{bc}\hat\Ga^{a_1}_{ab_1}=0$. The proof of the lemma
follows from the formula for the curvature tensor. $\Box$

Consider the Lie algebra $\hol^{n,\un,\psi,k,l}$. We have $f_i=f^0_i+\tilde f_{i\,n_0+1}^{k}+f_i^{n,\psi}+\breve f_{i\,
l+1}^{N+1\,n}$ ($i=1,2,3$). From \eqref{mkga10} and \eqref{mkga40} it follows that the holonomy algebra $\hol_0$
annihilats the vectors $p_1$ and $p_2$.

By Lemma \ref{mklem5}, it is enough to compute the covariant derivatives of the curvature tensor $\hat R$ of the metric $g_1$ with
$f_i=f^0_i+\tilde f_{i\,n_0+1}^{k}$ ($i=1,2,3$)
and the covariant derivatives of the curvature tensor $\bar R$ of the metric $g_2$ with
$f_i=f_i^{n,\psi}+\breve f_{i\, l+1}^{N+1\,n} $ ($i=1,2,3$).

Consider the curvature tensor $\hat R$. Set $\varphi=\phi=0$ and $N_0=N$, then we can use the above computations.

As in Lemma \ref{mklem3}, we can show that for all $0\leq r\leq N-1$ it holds
$$\hat R^i_{j\,2n+3\,2n+4;2n+3;...;2n+3(r \text{ times})}(0)=A^i_{r+1\,j},$$ and $\hat R^i_{j\,2n+3\,2n+4;2n+3;...;2n+3(r \text{ times})}(0)=0$ for $r\geq N$.

From \eqref{mkga230} and proof of Lemma \ref{mklem1} it follows that
$$\hat R^i_{jbc;a_1;...;a_{r-1};2n+4}=\sum_{\alpha=1}^{N}\frac{1}{\alpha!}(x^{2n+3})^\alpha [A_\alpha,\hat R(b,c;a_1;...;a_{r-1})]_j^i.$$
Hence if at least one of the numbers $a_1,...,a_r$ equals $2n+4$, then  $\hat R(b,c;a_1;...;a_{r})\in\un'$.
We can use also   Lemma \ref{mklem2}, \eqref{mkR80} and \eqref{mkR90}.

For the curvature tensor $\bar R$ we have $\bar R^i_{jbc;a_1;...;a_r}=0.$ Let $3\leq i\leq n+2$, then

\begin{multline*}\bar R^1_{i\,2n+3\,2n+4;2n+3;...;2n+3( r\text{ times})}(0)=
\bar R^2_{i+n\,2n+3\,2n+4;2n+3;...;2n+3( r\text{ times})}(0)\\=\left\{\begin{array}{ll}
-\psi_{1\,r+1\,i},&\text{if } k+3\leq i\leq l+2\text{ and } N_1\leq r\leq N-1,\\
1,&\text{if } l+3\leq i\leq n+2\text{ and }  r= N+i-l-3,\\
0,&\text{ else,}\end{array}\right.\end{multline*}
\begin{multline*}\bar R^1_{i+n\,2n+3\,2n+4;2n+3;...;2n+3( r\text{ times})}(0)=
-\bar R^2_{i\,2n+3\,2n+4;2n+3;...;2n+3( r\text{ times})}(0)\\
=\left\{\begin{array}{ll}
\psi_{2\,r+1\,i},&\text{if } k+3\leq i\leq l+2\text{ and } N_1\leq r\leq N-1,\\
\psi_{3\,r+1\,i},&\text{if } l+3\leq i\leq n+2\text{ and } N_1\leq r\leq N-1,\\
0,&\text{ else.}\end{array}\right.\end{multline*}
It can be also proved  that $\bar R^1_{jbc;a_1;...;a_r}=\bar R^2_{jbc;a_1;...;a_r}=0$, if $\{b\}\cup\{c\}\neq\{2n+3,2n+4\}$ or
$\{a_1\}\cup\cdots\cup\{a_r\}\neq\{2n+3\}$.

Now it is easy to see that $\hol_0=\hol^{n,\un,\psi,k,l}$. The Lie algebra $\hol^{m,\un,\psi,k,l,r}$ can be considered in the same way.

The theorem is proved. $\Box$.


\part{Applications and examples}

In this chapter we consider some examples and applications.
In Section \ref{secLie} we give examples of  real  4-dimensional Lie groups with left-invariant \pK metrics and find their holonomy algebras.
In Section \ref{secsym} we use our classification of holonomy algebras
to give a new proof for the classification of simply connected \pK symmetric spaces of index 2 with weakly-irreducible not
irreducible holonomy algebras.
Finally in Section \ref{secSasaki} we consider time-like cones over Lorentzian Sasaki manifolds.
These cones are \pK manifolds of index 2. We describe the local Wu decomposition of the cone in terms of the initial Lorentzian Sasaki manifold
and we find all possible weakly-irreducible not irreducible holonomy algebras of the cones.

\section{Examples of 4-dimensional Lie groups with left-invariant pseudo-K\"ahlerian metrics}\label{secLie}

Let $G$ be a Lie group endowed with a left-invariant  metric  $g$
and let $\g$ be the Lie algebra of $G$. We will consider $\g$ as the Lie algebra of
left-invariant vector fields on $G$ and as the tangent space at the identity  $e\in G$.
Let $X,Y,Z\in\g$. Since $X,Y,Z$ and  $g$ are left-invariant, from the Koszul formulae it follows that the Levi-Civita connection on
$(G,g)$ is given by
\begin{equation}\label{KoszulLie} 2g(\n_XY,Z)=g([X,Y],Z)+g([Z,X],Y)+g(X,[Z,Y]),
\end{equation}
where $X,Y,Z\in\g$. In particular, we see that  the vector field $\n_XY$ is also left-invariant.
Hence $\n_X$ can be considered as the linear operator $\n_X:\g\to\g$. Obviously, $\n_X\in\so(\g,g)$.
For the curvature tensor $R$ of $(G,g)$ at the point $e\in G$ we have
\begin{equation}\label{RLie} R(X,Y)=[\n_X,\n_Y]-\n_{[X,Y]},\end{equation} where $X,Y\in\g$.
From Theorem \ref{analythol} and the expressions for the covariant derivatives of the curvature tensor
 it follows that  the holonomy algebra $\hol_e$ at the point $e\in G$ is given by
\begin{equation}\label{holLie} \hol_e=\m_0+[\m_1,\m_0]+[\m_1,[\m_1,\m_0]]+\cdots,   \end{equation}
where $$\m_0=\spa\{R(X,Y)|X,Y\in\g\}\quad\text{and}\quad \m_1=\spa\{\n_X|X\in\g\}.$$

Now we consider 4-dimensional Lie algebras with the basis $p_1,p_2,q_1,q_2$ and with the  metric that has  the Gram  matrix
${\scriptsize\left(\begin{array}{cccc}0&0&1&0\\0&0&0&1\\1&0&0&0\\0&1&0&0\end{array}\right)}$ with respect to this basis.
Define the following Lie algebras by giving their non-zero brackets:
\begin{itemize}
\item[$\g_1$:]
$[p_1,q_1]=p_1+q_2,\quad [p_1,q_2]=-p_2-q_1,\quad [p_2,q_1]=p_2+q_1,\quad [p_2,q_2]=p_1+q_2;$

\item[$\g_2$:] $[p_1,q_2]=p_1,\quad [p_2,q_1]=-p_1,\quad  [q_1,q_2]=p_1+q_1.$

\end{itemize}

\begin{ex}\label{eLie1} The holonomy algebras of the Levi-Civita connections on the simply connected
Lie groups corresponding to the Lie algebras $\g_1$, $\g_2$ and $\g_3$ are given in the following table:

\centerline{ \begin{tabular}{|c|c|}\hline$\g$&$\hol_e$\\\hline $\g_1$&$\hol^2_{n=0}$\\\hline
$\g_2$&$\hol^{\gamma_1=0,\gamma_2=1}_{n=0}$\\\hline
\end{tabular}}
\end{ex}

{\it Proof.} Consider the Lie algebra $\g_1$. From \eqref{KoszulLie} it follows that
$$\n_{p_1}=\n_{q_2}= {\scriptsize\left(\begin{array}{cccc}0&-1&0&0\\1&0&0&0\\0&0&0&-1\\0&0&1&0\end{array}\right)},\quad
\n_{p_2}=\n_{q_1}={\scriptsize\left(\begin{array}{cccc}-1&0&0&0\\0&-1&0&0\\0&0&1&0\\0&0&0&1\end{array}\right)}.$$
From this and \eqref{RLie} it follows that $R(p_1,p_2)=R(q_1,q_2)=0$,
$$R(p_1,q_1)=R(p_2,q_2)={\scriptsize\left(\begin{array}{cccc}0&2&0&0\\-2&0&0&0\\0&0&0&2\\0&0&-2&0\end{array}\right)},\quad
R(p_1,q_2)=R(q_1,p_2)={\scriptsize\left(\begin{array}{cccc}2&0&0&0\\0&2&0&0\\0&0&-2&0\\0&0&0&-2\end{array}\right)}.$$
The statement for the Lie algebra $\g_1$ follows from \eqref{holLie}.

Similarly, for $\g_2$ we have
$$\n_{p_1}= {\scriptsize\left(\begin{array}{cccc}0&0&0&1\\0&0&-1&0\\0&0&0&0\\0&0&0&0\end{array}\right)},\quad
\n_{q_1}={\scriptsize\left(\begin{array}{cccc}0&1&0&1\\-1&0&-1&0\\0&0&0&1\\0&0&-1&0\end{array}\right)},\quad
\n_{p_2}=\n_{q_2}=0$$
and $R(p_1,p_2)=R(p_1,q_1)=R(p_2,q_2)=0$,
$$R(p_1,q_2)=-R(p_2,q_1)={\scriptsize\left(\begin{array}{cccc}0&0&0&-1\\0&0&1&0\\0&0&0&0\\0&0&0&0\end{array}\right)},\quad
R(q_1,q_2)={\scriptsize\left(\begin{array}{cccc}0&-1&0&-2\\1&0&2&0\\0&0&0&-1\\0&0&1&0\end{array}\right)}.$$

For $\g_3$ we obtain

$$\n_{p_1}= {\scriptsize\left(\begin{array}{cccc}-\gamma&-\gamma&0&1\\\gamma&-\gamma&-1&0\\0&0&\gamma&-\gamma\\0&0&\gamma&\gamma\end{array}\right)},\quad
\n_{p_2}= {\scriptsize\left(\begin{array}{cccc}-\gamma&\gamma&0&1\\-\gamma&-\gamma&-1&0\\0&0&\gamma&\gamma\\0&0&-\gamma&\gamma\end{array}\right)},$$
$$\n_{q_1}= {\scriptsize\left(\begin{array}{cccc}\gamma&\gamma+1&0&-1\\-\gamma-1&\gamma&1&0\\0&0&-\gamma&\gamma+1\\0&0&-\gamma-1&-\gamma\end{array}\right)},\quad
\n_{q_2}= {\scriptsize\left(\begin{array}{cccc}-\gamma&\gamma+1&0&1\\-\gamma-1&-\gamma&-1&0\\0&0&\gamma&\gamma+1\\0&0&-\gamma-1&\gamma\end{array}\right)}$$
and $R(p_1,p_2)=0$,
$$R(p_1,q_1)=R(p_2,q_2)=\gamma R(q_1,q_2)=2\gamma(1+2\gamma) {\scriptsize\left(\begin{array}{cccc}0&-1&0&0\\1&0&0&0\\0&0&0&-1\\0&0&1&0\end{array}\right)},$$
$$R(p_1,q_2)=-R(p_2,q_1)=2(1+2\gamma)  {\scriptsize\left(\begin{array}{cccc}\gamma&0&0&-1\\0&\gamma&1&0\\0&0&-\gamma&0\\0&0&0&-\gamma\end{array}\right)}.$$

The theorem is proved. $\Box$


\section{Pseudo-K\"ahlerian symmetric spaces of index 2}\label{secsym}

In this section we use our classification of holonomy algebras to give a new proof for the classification of simply
connected \pK symmetric spaces of index 2 with weakly-irreducible not irreducible holonomy algebras. This classification
follows from results of M.~Berger \cite{Ber57} and of  I.~Kath and M.~Olbrich  \cite{KathOldrich}. The facts about
pseudo-Riemannian symmetric spaces that we use here can be found in the books \cite{Be,Helgason,K-N}.

Recall that a pseudo-Riemannian manifold $(M,g)$ is called {\it a symmetric space} if the geodesic symmetry $s_p$ with
respect to any point $p\in M$ is a globally defined isometry. If a pseudo-K\"ahlerian manifold is a symmetric space, than
it is called {\it a pseudo-K\"ahlerian symmetric spaces}. A pseudo-Riemannian symmetric space is geodesically  complete
and its curvature tensor is parallel, i.e.  $\nabla R=0$. Any pseudo-Riemannian symmetric space admits a connected
transitive Lie group of isometries $G$ and can be identified with the factor $G/K$, where $K\subset G$ is the stabilizer
of a point $o\in M$. Moreover, there exists an involutive automorphism $\sigma:G\to G$ such that $S^0\subset K\subset S$,
where $S$ is the set of the fixed points of the automorphism $\sigma$ and $S^0$ is the connected identity component of the
set $S$. For the group $G$ can be chosen the group of transvections  of $(M,g)$, i.e. the group $\{s_p\circ s_q|p,q\in
M\}$. Conversely,  any such triple $(G,K,\sigma)$ under the condition that the homogeneous space $G/K$ admits a
$G$-invariant pseudo-Riemannian metric $g$ defines the pseudo-Riemannian symmetric space $(G/H,g)$. A symmetric space
$(M,g)$ is called {\it semi-simple} if the group of transvections of $(M,g)$ is semi-simple.

Let  $\g$ be a Lie algebra with an involutive automorphism $\sigma:\g\to\g$. Denote by $\g_\pm$ the eigenspaces of $\sigma$ corresponding
to the eigenvalues $\pm 1$.
We get the decomposition of $\g$ into the direct sum of the vector subspaces $$\g=\g_++\g_-.$$
For the Lie brackets of the Lie algebra $\g$ we have
$$[\g_+,\g_+]\subset\g_+,\quad[\g_+,\g_-]\subset\g_-,\quad [\g_-,\g_-]\subset\g_+.$$
We see that $\g_+\subset\g$ is a subalgebra and $\ad_{\g_+}|_{\g_-}:\g_+\to\gl(\g_-)$ is a representation.

Any simply  connected pseudo-Riemannian symmetric space $(M,g)$ of signature $(r,s)$ is uniquely defined by a triple
$(\g,\sigma,\eta)$ ({\it a symmetric triple}), where
$\g$ is a Lie algebra, $\sigma$ is an involutive automorphism of $\g$ and
$\eta$ is a non-degenerate $\ad_{\g_+}$-invariant  symmetric bilinear form of signature $(r,s)$ on the vector space  $\g_-$.

The vector space $\g_-$ can be identified with the tangent space to $M$ at a fixed point $o\in M$.
Then the form $\eta$ is identified with the form $g_o$. The curvature tensor $R_o$ of the manifold $(M,g)$ is given by
$$R_o(x,y)z=-[[x,y],z],$$
where $x,y,z\in\g_-=T_oM$. Since any symmetric space is an analytic manifold and $\nabla R=0$,
for the holonomy algebra of the manifold $(M,g)$ at the point $o\in M$ we have
$$\hol=\spa\{R_o(T_oM,T_oM)\}=[\g_-,\g_-].$$

Let $\h\subset\so(r,s)$ be a subalgebra. The tensorial extension of the  representation $\h\hookrightarrow\so(r,s)$
defines a representation of $\h$ in the space of curvature tensors $\R(\h)$ as follows
$$A\cdot R= R_A,\quad R_A(x,y)z= [R(x\wedge y),A]+R(Ax\wedge y)+R(x\wedge Ay),$$
where $A\in\h$, $R\in\R(\h)$ and $x,y,z\in \Real^{r,s}$.
Consider the following space
$$\R_0(\h)=\{R\in\R(\h)|R_A=0 \text{ for all } A\in\h\}.$$
Let $(M,g)$ be a symmetric space and $R$ its curvature tensor. Since $\nabla R=0$, we have $R_o\in\R_0(\hol)$.
Conversely, let $\h\subset\so(r,s)$ be a subalgebra and $R\in\R_0(\h)$.
Consider the Lie algebra $\g=\h+\Real^{r,s}$ with the Lie brackets
$$[x,y]=-R(x,y),\quad [A,x]=Ax,\quad [A,B]=A\circ B-B\circ A,$$
where $x,y\in\Real^{r,s}$ and $A,B\in\h$.
Define the involutive  automorphism $\sigma$ of $\g$ by $\sigma(A+x)=A-x,$ where $A\in \h$, $x\in \Real^{r,s}$.
We have $\g_+=\h$ and $\g_-=\Real^{r,s}$.
Obviously, the pseudo-Euclidian metric $\eta$ on $\Real^{r,s}$ is $\ad_{g_+}$-invariant.
We get the symmetric triple $(\g,\sigma,\eta)$ that defines us a simply connected  pseudo-Riemannian symmetric space of signature $(r,s)$.

For any symmetric triple $(\g,\sigma,\eta)$, consider the new symmetric triple $(\g_1=[\g_-,\g_-]+\g_-,\sigma|_{\g_1},\eta)$.
The simply connected symmetric spaces
corresponding to the both triples have at the origins the same metrics and curvature tensors,
hence these symmetric spaces are isometric.
The holonomy algebras of these symmetric spaces coincide with $g_{1+}=[\g_-,\g_-]$.

Thus any pseudo-Riemannian symmetric space $(M,g)$ of signature $(r,s)$ defines a pair $(\hol,R)$ ({\it a symmetric pare}),
where $\hol\subset\so(r,s)$ is a subalgebra, $R\in\R_0(\hol)$ and we have $\hol=R(\Real^{r,s}\wedge\Real^{r,s})$.
Conversely, any such pair defines a simply connected pseudo-Riemannian symmetric space.

{\it An isomorphism }of two symmetric triples $f:(\g_1,\sigma_1,\eta_1)\to (\g_2,\sigma_2,\eta_2)$ is a Lie algebra
isomorphism $f:\g_1\to\g_2$ such that $f\circ\sigma_1=\sigma_2\circ f$ and the induced map $f:\g_{1-}\to\g_{2-}$
is an isometry of the pseudo-Euclidean spaces $(\g_{1-},\eta_1)$ and $(\g_{2-},\eta_2)$.
Isomorphic symmetric triples define isometric pseudo-Riemannian symmetric spaces.
We say that two symmetric pairs are {\it isomorphic} if they define isomorphic symmetric triples.
For a positive real number $c\in\Real$, the symmetric triples $(\g,\sigma,\eta)$ and $(\g,\sigma,c\eta)$
define diffeomorphic simply connected symmetric spaces and the metrics of these spaces differ by the factor $c$.

Let $(\hol,R)$ be a non-trivial symmetric  pair (i.e. $R\neq 0$) and $(\g=\hol+\Real^{r,s},\sigma,\eta)$
be the corresponding symmetric triple. For any positive real number $c\in \Real$
consider the symmetric pair $(\hol,cR)$ and the  symmetric triple
$(\g_c=\hol+\Real^{r,s},\sigma,c\eta).$ Obviously, the map $f:\g\to\g_c$ given by
$f(A+x)=A+\frac{1}{\sqrt{c}}x,$ where $A\in\hol,$ $x\in\Real^{r,s}$, is an isomorphism
of the symmetric triples $(\g,\sigma,\eta)$ and $(\g_c,\sigma,c\eta)$.
Note that the  symmetric pairs $(\hol,R)$  and $(\hol,-R)$ can define non-isomorphic
symmetric spaces, e.g. in the Riemannian case if $\hol$ is irreducible.

From the Wu theorem it follows that any simply connected pseudo-Riemannian  symmetric space is a product of a pseudo-Euclidean space and
of  simply connected pseudo-Riemannian  symmetric spaces with weakly-irreducible holonomy algebras.
Thus the classification problem  of simply connected pseudo-Riemannian  symmetric spaces is reduced to the case of spaces with weakly-irreducible
holonomy algebras. The classification of simply connected pseudo-Riemannian  symmetric spaces with irreducible holonomy
algebras was obtained  by M.~Berger in \cite{Ber57}. Thus we are left with the case of symmetric  spaces with weakly-irreducible not irreducible
holonomy algebras.
From the above it follows that to solve the classification problem for a signature $(r,s)$ we need a classification
of weakly-irreducible not irreducible holonomy algebras $\hol\subset\so(r,s)$ for each of  which there exists an
$R\in\R_0(\hol)$ such that $R(\Real^{r,s}\wedge \Real^{r,s})=\hol$.
Then each line $L\subset \R_0(\hol)$ such that for any non-zero  $R\in L$
we have $R(\Real^{r,s}\wedge \Real^{r,s})=\hol$ will give us two symmetric pairs $(\hol,R)$  and $(\hol,-R)$ (where $R\in L$ is fixed).
Checking which of these pairs define non-isometric symmetric triples, we will get a classification of
simply connected pseudo-Riemannian  symmetric spaces of signature $(r,s)$.

Now consider simply connected pseudo-K\"ahlerian  symmetric spaces of index 2.
Recall that we have a fixed  complex structure $J$ on $\Real^{2,2n+2}$ and a  basis $p_1$, $p_2$, $e_1$,...,$e_n$,  $f_1$,...,$f_n$,
$q_1$, $q_2$ of  $\Real^{2,2n+2}$
such that the Gram matrix of the pseudo-Euclidian  metric $\eta$  and the matrix of $J$  have the forms\\
$${\scriptsize \left (\begin{array}{ccccc}
0 &0 &0 & 1 & 0\\
0 &0 &0 & 0& 1\\
0 &0 & E_{2n} & 0 & 0 \\
1 &0 &0 &0 & 0 \\
0 &1 &0 &0 & 0 \\
\end{array}\right)}\quad\text{ and }\quad{\scriptsize \left (\begin{array}{cccccc}
0&-1 &0 &0 & 0 & 0\\
1& 0 &0 &0 & 0 & 0\\
0 &0 &0& -E_n & 0 & 0 \\
0 &0 &E_n & 0 & 0& 0 \\
0 &0 &0 &0 & 0&-1 \\
0 &0 &0 &0 &1 &0 \\
\end{array}\right),}\text{ respectively.}$$
We denote by $\u(1,n+1)_{<p_1,p_2>}$ the subalgebra of $\u(1,n+1)$ that preserves the $J$-invariant 2-dimensional
isotropic  subspace  $\Real p_1\oplus\Real p_2\subset\Real^{2,2n+2}$. This Lie algebra
can be identified with the following matrix algebra
$${\scriptsize  \u(1,n+1)_{<p_1,p_2>}=\left\{\left. \left (\begin{array}{cccccc}
a_1&-a_2 &-z^t_1 & -z^t_2 &0 &-c\\
a_2&a_1 &z^t_2 & -z^t_1 &c &0\\
0 &0&B&-C&z_1&-z_2\\
0 &0&C&B&z_2&z_1\\
0&0&0&0&-a_1&-a_2\\
0&0&0&0&a_2&-a_1\\
\end{array}\right)\right|\,
\begin{array}{c}
a_1,a_2,c\in \Real,\\ z_1, z_2\in \Real^n,\\
\bigl(\begin{smallmatrix}B&-C\\C&B\end{smallmatrix}\bigr)\in \u(n) \end{array} \right\}.}$$
Using the form $\eta$, we identify the Lie algebra $\so(2,2n+2)$ with the space $\Real^{2,2n+2}\wedge \Real^{2,2n+2}$ of bivectors.
Then the element of $\u(1,n+1)_{<p_1,p_2>}$ given by the above matrix is identified with the bivector
$$-a_1(\pq)+a_2(\qp)+A+p_1\wedge z_1+p_2\wedge Jz_1+
p_1\wedge Jz_2-p_2\wedge z_2+ cp_1 \wedge p_2,$$
where $A\in\Real^{2n}\wedge \Real^{2n} $ corresponds to $\bigl(\begin{smallmatrix}B&-C\\C&B\end{smallmatrix}\bigr)\in\u(n)$.

Now we can formulate the main theorem of this section.

\begin{theorem}\label{sym} If  $(\hol,R)$ is a symmetric pair associated with a simply connected
\pK symmetric space of signature $(2,2n+2)$ $(n\geq 0)$
with weakly-irreducible not irreducible holonomy algebra, then $(\hol,R)$ is isomorphic exactly to one of the
symmetric pairs of the following list (we define curvature tensors giving their non-zero values):

\begin{description}

\item[1.] n=0: \begin{description}
\item[(a)] $\hol_1=\hol^{\gamma_1=0,\gamma_2=0}_{n=0}=\Real p_1\wedge p_2=
\Real{\scriptsize \left(\begin{array}{cccc}0&0&0&-1\\0&0&1&0\\0&0&0&0\\0&0&0&0\end{array}\right)},$\\
$R_{\lambda_5=1}(q_1\wedge q_2)=p_1\wedge p_2;$
\item[(b)] $(\hol_1,-R_{\lambda_5=1})$;

\item[(c)] $\hol_{1c}=\hol^{2}_{n=0}=\Real (p_1\wedge q_1+p_2\wedge q_2)\oplus\Real (p_1\wedge q_2-p_2\wedge q_1)$\\ $=
{\scriptsize \left.\left\{
\left(\begin{array}{cccc}a_1&-a_2&0&0\\a_2&a_1&0&0\\0&0&-a_1&-a_2\\0&0&a_2&-a_1\end{array}\right) \right|
a_1,a_2\in\Real\right\} } ,$\\ $R_{\lambda_1=-\frac{1}{2}}(p_1\wedge q_1)=R_{\lambda_1=-\frac{1}{2}}(p_2\wedge
q_2)=-\frac{1}{2}(p_1\wedge q_1+p_2\wedge q_2),$\\ $R_{\lambda_1=-\frac{1}{2}}(p_1\wedge
q_2)=-R_{\lambda_1=-\frac{1}{2}}(p_2\wedge q_1)=-\frac{1}{2}(p_1\wedge q_2-p_2\wedge q_1);$

\item[(d)] $(\hol_{1c},-2R_{\lambda_1=-\frac{1}{2}})$;

\item[(e)] $\hol_{1c}$,\\
$R_{\lambda_2=1}(p_1\wedge q_1)=R_{\lambda_2=1}(p_2\wedge q_2)=p_1\wedge q_2-p_2\wedge q_1,$\\ $R_{\lambda_2=1}(p_1\wedge
q_2)=-R_{\lambda_2=1}(p_2\wedge q_1)=-(p_1\wedge q_1+p_2\wedge q_2);$

\end{description}

\item[2.] n=1: $\hol_2=\hol^{m=n-1=0,\{0\},\varphi=0,\phi=0}
=\Real(p_1\wedge e_1+p_2\wedge f_1)\oplus\Real p_1\wedge p_2$ $$={\scriptsize\left.\left\{ \left(\begin{array}{cccccc}
0&0&-z_1^t&0&0&-c\\ 0&0&0&-z_1^t&c&0\\ 0&0&0&0&z_1&0\\ 0&0&0&0&0&z_1\\ 0&0&0&0&0&0\\
0&0&0&0&0&0\end{array}\right)\right|z_1,c\in\Real\right\}},$$ $R_{M_3(1)=e_1}(q_1\wedge q_2)=p_1\wedge e_1+p_2\wedge
f_1,$\\ $R_{M_3(1)=e_1}(q_1\wedge e_1)=R_{M_3(1)=e_1}(q_2\wedge f_1)=p_1\wedge p_2$;

\item[3.]$n\geq 0$:  $\hol_3=\hol^{m,\{\Real J_m\},\varphi=0,\phi=2}$\\
$=\Real(2J-J_m)\zr(\spar\{p_1\wedge e_i+p_2\wedge f_i,\,p_1\wedge f_j-p_2\wedge e_j|\,1\leq i\leq n,\,m+1\leq j\leq n \}+\Real p_1\wedge p_2)$
$${\scriptsize
=\left\{\left. \left (\begin{array}{cccccccc}
0&-2a &-z^t_1&-z'^t_1 & -z^t_2&0 &0 &-c\\
2a&0 &z^t_2&0 & -z^t_1&-z'^t_1 &c &0\\
0 &0&0&0&-aE_m&0&z_1&-z_2\\
0&0&0&0&0&-2aE_{n-m}&z'_1&0\\
0&0&aE_m&0&0&0&z_2&z_1\\
0&0&0   &2aE_{n-m}&0&0&0&z'_1\\
0&0&0&0&0&0&0&-2a\\
0&0&0&0&0&0&2a&0\\
\end{array}\right)\right|\,
\begin{array}{c}a,c\in \Real,\\
 z_1, z_2\in \Real^m,\\
 z'_1\in\Real^{n-m},
\end{array} \right\},}$$
where $0\leq m\leq n$,\\
$R=R_1+R_{\lambda_5}$, $R_1=R_{\lambda_3=-2}+R_{S=-\frac{1}{2}J_m}$,  \\
$R_1(p_1\wedge q_2)=-R_1(p_2\wedge q_1)=-2 p_1\wedge p_2,$\\
$R_1(e_i\wedge f_i)=-p_1\wedge p_2$, $1\leq i\leq m$,\\
$R_1(e_j\wedge f_j)=-2p_1\wedge p_2$, $m+1\leq j\leq n$,\\
$R_1(q_1\wedge e_i)=R_1(q_2\wedge f_i)=-\frac{1}{2}(p_1\wedge e_i+p_2\wedge f_i)$, $1\leq i\leq m$,\\
$R_1(q_2\wedge e_i)=-R_1(q_1\wedge f_i)=\frac{1}{2}(p_1\wedge f_i-p_2\wedge e_i)$, $1\leq i\leq m$,\\
$R_1(q_1\wedge e_j)=R_1(q_2\wedge f_j)=-2(p_1\wedge e_j+p_2\wedge f_j)$, $m+1\leq j \leq n$,\\
$R_1(q_1\wedge q_2)=J_m-2J$,\\
$R_{\lambda_5}=\lambda_5 p_1\wedge p_2$, where $\lambda_5\in\Real$;

\item[4.] $n\geq 0$:  $(\hol_3,-R_1-R_{\lambda_5})$.
\end{description}
 \end{theorem}

{\bf Remark.} The simply connected symmetric spaces corresponding to the symmetric triples (c), (d) and (e) are
semi-simple and they correspond to the spaces $SL(2,\Co)/{\Co^*}$, $SO(3,\Co)/{\Co^*}$ and $Sp(2,\Co)/{\Co^*}$ from
\cite{Ber57}. The simply connected symmetric spaces corresponding to the other symmetric triples exhaust all simply
connected  \pK symmetric spaces of index 2 that are not semi-simple (they were classified in \cite{KathOldrich}).

{\bf Proof of Theorem \ref{sym}.} As the first step   we will find all symmetric pairs $(\hol,R)$, where
$\hol\subset\u(1,n+1)_{<p_1,p_2>}$ is a  weakly-irreducible not irreducible holonomy algebra.
On the second step we will check which of the obtained symmetric pairs define isometric
simply connected symmetric spaces.

{\bf Step 1.} Fix a holonomy algebra $\hol\subset\u(1,n+1)_{<p_1,p_2>}$ and a  curvature tensor
$R\in\R_0(\hol)$ such that $R(\Real^{2,2n+2}\wedge\Real^{2,2n+2})=\hol$.
Let $0\leq m\leq n$ be the associated number and $\u\subset\u(m)$ the unitary part of $\hol$.
Using (\ref{svoistvo1}) and Table 3.2.2 we can decompose $R$ into several components,
we will use the corresponding denotation. For any element $A$ as in Table 3.2.2, we denote by
$R_A$ the curvature tensor obtained by  Table 3.2.2 under the condition that the other elements
of Table 3.2.2  are zero.

It is easy to see that for the Lie algebra $\hol^2_{n=0}$ we have $\R_0(\hol^2_{n=0})=\R(\hol^2_{n=0})=\Real
R_{\lambda_1=1}\oplus \Real R_{\lambda_2=1}$. Let $\hol$ be any other holonomy algebra, then $\hol$ containes the ideal
$\Real p_1\wedge p_2$.

\begin{lem}\label{symL1} For the components of the curvature tensor $R$ we have

$$\begin{array}{lll} \lambda_1=\lambda_2=\lambda_4=0,& N_1=N_2=0,& R_0=0,\\
P=0,&  S^{11}=0,& S^{12}=\frac{\lambda_3}{4}J|_{E^2_{1,...,m}},\\
S^{21}=\frac{\lambda_3}{4}J|_{E^1_{1,...,m}},&S-S^*=\frac{\lambda_3}{2}J_m.\end{array}$$
Furthermore,  if $\lambda_3\neq 0$, then $M_1=M_2=0$.
\end{lem}

{\it Proof.}
Using the fact that $R\in\R_0(\hol)$, \eqref{b5a} and  Table 3.2.2, we get
$$\begin{array}{rl}0&=R_{p_1\wedge p_2}(q_1\wedge q_1)=[R(q_1\wedge q_2),p_1\wedge p_2]+R(p_2\wedge q_2)+R(q_1\wedge (-p_1))\\
&=[R(q_1\wedge q_2),p_1\wedge p_2]+2R(p_1\wedge q_1)\\
&=2\lambda_4p_1\wedge p_2+2(\lambda_1\pq+\lambda_2\qp+p_1\wedge N_1(1)+p_2\wedge JN_1(1)\\
&\,\,-(p_1\wedge JN_2(1)-p_2\wedge N_2(1))+
\lambda_4 p_1\wedge p_2.\end{array}$$ This implies that $\lambda_1=\lambda_2=\lambda_4=0$ and  $N_1=N_2=0$.

Consequently if $n=0$, then for the Lie algebra $\hol$ we have two possibilities: either $\hol=\hol_1$ or
$\hol=\hol_{n=0}^{\gamma_1=0,\gamma_2=1}=\Real(\qp)\oplus\Real p_1\wedge p_2.$ It is easy to check that $\R_0(\hol_1)=\Real R_{\lambda_5=1}$
and $\R_0(\hol_{n=0}^{\gamma_1=0,\gamma_2=1})=\Real R_{\lambda_3=-2}\oplus\Real R_{\lambda_5=1}$.

Suppose that $n\geq 1$. Let $m_1=k$ if $\hol$ is a Lie algebra of type $\hol^{n,\u,\psi,k,l}$ or $\hol^{m,\u,\psi,k,l,r}$, and
$m_1=m$ if $\hol$ is a Lie algebras of another type.

Let $w\in E^1_{1,...,m_1}$. From Theorem \ref{hol2n} it follows that
$p_1\wedge w+p_2\wedge Jw\in \hol$. We have
$$\begin{array}{rl}0&=\pr_\u(R_{p_1\wedge w+p_2\wedge Jw}(q_1\wedge q_2))\\
&=\pr_\u([R(q_1\wedge q_2),p_1\wedge w+p_2\wedge Jw]+R(w\wedge q_2)+R(q_1\wedge Jw))\\
&=-2\pr_u(R(q_2\wedge w))=2P(Jw).\end{array}$$ Hence, $P|_{E^2_{1,...,m}}=0$. Similarly,  the condition
$\pr_\u(R_{p_1\wedge Jw-p_2\wedge w}(q_1\wedge q_2))=0$ implies that
$P|_{E^1_{1,...,m}}=0$. Thus, $P=0$.

For $x\in E^1_{1,...,m_1}+E^2_{1,...,m_1}$ we have
$$\begin{array}{rl}0&=\pr_\u(R_{p_1\wedge w+p_2\wedge Jw}(q_1\wedge x))\\
&=\pr_\u([R(q_1\wedge x),p_1\wedge w+p_2\wedge Jw]+R(w\wedge x)+R(q_1\wedge (-\eta(w,x)p_1+\eta(Jw,x)p_2)))\\
&=R_0(w,x).\end{array}$$ Similarly, the condition $\pr_\u(R_{p_1\wedge Jw-p_2\wedge w}(q_1\wedge x))=0$
yields that $R_0(Jw,x)=0$. Thus, $R_0=0$.


For $w\in E^1_{1,...,m_1}$  we have
$$\begin{array}{rl}0&=\pr_{\N^1_{1,...,m}+\N^2_{1,...,m}}(R_{{p_1\wedge w+p_2\wedge Jw}}(q_1\wedge q_2))\\
&=\pr_{\N^1_{1,...,m}+\N^2_{1,...,m}}([R(q_1\wedge q_2),p_1\wedge w+p_2\wedge Jw]+R(w\wedge q_2)+R(q_1\wedge Jw))\\
&=[\lambda_3(\pq)+S-S^*,p_1\wedge w+p_2\wedge Jw]-2\pr_{\N^1_{1,...,m}+\N^2_{1,...,m}}(R(q_2\wedge w))\\
&=\lambda_3(p_2\wedge w-p_1\wedge Jw)+p_1\wedge ((S^{11}-S^{*11})w+(S^{21}-S^{*21})w)\\
&\quad +p_2\wedge ((S^{12}-S^{*12})Jw+(S^{22}-S^{*22})Jw)\\
&\quad -2(p_1\wedge S^{11*}(w)+p_2\wedge JS^{11*}(w)+p_1\wedge JS^{12}(Jw)-p_2\wedge S^{12}(Jw)).\end{array}$$
Considering the terms from $p_1\wedge E^1$, $p_1\wedge E^2$, $p_2\wedge E^1$, $p_2\wedge E^2$
and using the formulas
$$S^{*21}=S^{12*}=JS^{12}J,\quad S^{*12}=S^{21*}=JS^{21}J,\quad S^{22}=JS^{11*}J,\quad  S^{*22}=S^{22*}=JS^{11}J$$
we  get
\begin{equation}\label{symeq1}S^{11}=3S^{*11}\quad \text{ and }\quad\lambda_3 Jw-S^{21}w+3JS^{12}Jw.\end{equation}
Likewise, from the equality $\pr_{\N^1_{1,...,m_1}+\N^2_{1,...,m_1}}(R_{{p_1\wedge Jw-p_2\wedge w}}(q_1\wedge q_2))=0$
it follows that
\begin{equation}\label{symeq2} 3S^{11}=S^{*11}\quad\text{ and }\quad\lambda_3 w+S^{12}Jw+3JS^{21}w=0.\end{equation}

Suppose that the Lie algebra $\hol$ is of type $\hol^{n,\u,\psi,k,l}$ or $\hol^{m,\u,\psi,k,l,r}$. Then from Table 3.2.2
it follows that $\lambda_3=0$. This, \eqref{symeq1} and \eqref{symeq2} yield that $S=0$. Therefore, $\u=\{0\}$
and the Lie algebra $\hol$ can not be neither of type $\hol^{n,\u,\psi,k,l}$ nor of type $\hol^{m,\u,\psi,k,l,r}$.

From \eqref{symeq1} and \eqref{symeq2} it follows that
$$S^{11}=0,\quad S^{12}=\frac{\lambda_3}{4}J|_{E^2_{1,...,m}},\quad S^{21}=\frac{\lambda_3}{4}J|_{E^1_{1,...,m}}.$$
Thus,
$$\begin{array}{rl}S-S^*&=S^{12}+S^{21}-S^{*12}+S^{*21}\\
&=\frac{\lambda_3}{4}J|_{E^2_{1,...,m}}+\frac{\lambda_3}{4}J|_{E^1_{1,...,m}}+
+\frac{\lambda_3}{4}J|_{E^1_{1,...,m}}+\frac{\lambda_3}{4}J|_{E^2_{1,...,m}}=\frac{\lambda_3}{2}J_m.\end{array}$$

Suppose that $\lambda_3\neq 0$. From the proof of the lemma and Theorem \ref{hol2n} it follows that
$\hol=\hol_3$, where $0\leq m\leq n$.
Hence, $2J-J_m\in\hol$. For $w\in E^1_{1,...,m}$ we have
$$\begin{array}{rl}0=&R_{2J-J_m}(q_1\wedge w)=[R(q_1\wedge w),2J-J_m]+R(2q_2\wedge x)+R(q_1\wedge J_mx)\\
=& [R(q_1\wedge w) ,2J-J_m]+2R(q_2\wedge x)-R(q_2\wedge x)=[R(q_1\wedge w),2J-J_m]+R(q_1\wedge x)\\
=&[p_1\wedge \frac{1}{2}x+p_2\wedge \frac{1}{2}Jx+M_1^*(x)p_1\wedge p_2,2J-J_m]-p_1\wedge \frac{1}{2}Jx+p_2\wedge \frac{1}{2}x+M_2^*(x)p_1\wedge p_2\\
=&M_2^*(x)p_1\wedge p_2.\end{array}$$
Therefore, $M_2=0$. Similarly, $M_1=0$. The lemma is proved. $\Box$

As all symmetric pairs for $n=0$ were found in Lemma \ref{symL1}, suppose that $n\geq 1$.
If $\lambda_3=0$, then $S=0$. From Lemma \ref{symL1}, Table 3.2.2 and the fact that $R(\Real^{2,2n+2}\wedge\Real^{2,2n+2})=\hol$
it follows that $\hol=\hol_2$. It is easy to check that
$\R_0(\hol_2)=\Real R_{M_3(2)=e_1}\oplus\Real R_{\lambda_5=1}$.

Suppose that  $\lambda_3\neq 0$. From Lemma \ref{symL1} and Theorem \ref{hol2n} it follows that
$\hol=\hol_3$. It is easily shown that
$\R_0(\hol_3)=\Real R_{1}\oplus\{R_{M_3}\}\oplus\Real R_{\lambda_5=1}$.
Note that if $n=0$, then the Lie algebra $\hol_3$ coincides with the Lie algebra $\hol_{n=0}^{\gamma_1=0,\gamma_2=1}$.

{\bf Step 2.} Now to complete the proof of the theorem, we must prove the following statements:
\begin{itemize}

\item[1)] the symmetric pairs $(\hol_1,R_{\lambda_5=1})$ and
$(\hol_1,-R_{\lambda_5=1})$ are not isomorphic;

\item[2)] the symmetric pairs $(\hol_2,R_{M_3(1)=e_1})$, $(\hol_2,R_{M_3(1)=e_1}+R_{\lambda_5})$ and
$(\hol_2,-R_{M_3(1)=e_1}+R_{\lambda_5})$ are isomorphic for any $\lambda_5\in\Real$;

\item[3)] the symmetric pairs $(\hol_3,R_1+R_{\lambda_5})$ and
$(\hol_3,R_1+R_{\lambda_5}+R_{M_3})$ are isomorphic for any  $\lambda_5\in\Real$ and  $M_3\in\Hom(\Real,E^1_{m+1,...,n})$;

\item[4)] the symmetric pairs $(\hol_3,R_1+R_{\lambda_5})$ and
$(\hol_3,R_1+R_{\lambda_5'})$ are not isomorphic if $\lambda_5\neq \lambda_5'$;

\item[5)] the symmetric pairs $(\hol_3,-R_1-R_{\lambda_5})$ and
$(\hol_3,-R_1-R_{\lambda_5'})$ are not isomorphic if $\lambda_5\neq \lambda_5'$;

\item[6)] the symmetric pairs $(\hol_3,R_1+R_{\lambda_5})$ and
$(\hol_3,-R_1+R_{\lambda_5'})$ are not isomorphic for any $\lambda_5,\lambda_5'\in\Real$;

\item[7)] if $\lambda_1\neq 0$, then the symmetric pairs $(\hol_{1c},R_{\lambda_1}+R_{\lambda_2})$ and
$(\hol_{1c},R_{\lambda_1})$ are  isomorphic for any $\lambda_2\in\Real$;

\item[8)] the symmetric pairs $(\hol_{1c},R_{\lambda_2=1})$ and $(\hol_{1c},-R_{\lambda_2=1})$
 are  isomorphic;

\item[9)] if $\lambda_1>0$, then any two of the symmetric pairs $(\hol_{1c},R_{\lambda_1})$,
$(\hol_{1c},-R_{\lambda_1})$ and $(\hol_{1c},R_{\lambda_2=1})$ are not isomorphic.

\end{itemize}

Let us prove these statements.

1) It is easy to check that the sectional curvature of the curvature tensor $R_{\lambda_5=1}\in\R(\hol_1)$ is non-zero and non-negative.
In converse, the sectional curvature of the curvature tensor $-R_{\lambda_5=1}\in\R(\hol_1)$ is non-zero and non-positive.
For instance, $$\eta(R_{\lambda_5=1}((p_1+q_1)\wedge(p_2+q_2))(p_1+q_1),(p_2+q_2))=1$$ and
$$\eta(-R_{\lambda_5=1}((p_1+q_1)\wedge(p_2+q_2))(p_1+q_1),(p_2+q_2))=-1.$$

2) To show that  the symmetric pairs
$(\hol_2,R_{M_3(1)=e_1})$ and $(\hol_2,R_{M_3(1)=e_1}+R_{\lambda_5})$ are isomorphic
consider the corresponding symmetric triples $(\g_1=\hol_2+\Real^{2,4},\sigma_1,\eta)$ and
$(\g_2=\hol_2+\Real^{2,4},\sigma_2,\eta)$ and consider the Lie algebra isomorphism $f:\g_1\to\g_2$ given by
$$\begin{array}{llll}f|_{\hol}=\id_{\hol},& f(p_1)=p_1,& f(p_2)=p_2,& f(e_1)=\lambda_5p_2+e_1,\\
f(f_1)=-\lambda_5 p_1+f_1,& f(q_1)=q_1+\lambda_5 f_1,&   f(q_2)=q_2-\lambda_5 e_1.&\end{array}$$
Similarly we can show that the symmetric pairs
$(\hol_2,-R_{M_3(1)=e_1})$ and $(\hol_2,-R_{M_3(1)=e_1}+R_{\lambda_5})$ are isomorphic.
Let $(\g_3,\sigma_3,\eta)$ be the symmetric triple corresponding to the symmetric pair $(\hol_2,-R_{M_3(1)=e_1})$.
The isomorphism $f:\g_1\to\g_3$ given by
$$\begin{array}{llll}f(p_1\wedge e_1+p_2\wedge f_2)=-(p_1\wedge e_1+p_2\wedge f_2),& f(p_1)=p_1,& f(p_2)=p_2,& f(e_1)=-e_1,\\
R(p_1\wedge p_2)=p_1\wedge p_2,& f(f_1)=-f_1,& f(q_1)=q_1,&   f(q_2)=q_2\end{array}$$
implies that the symmetric pairs $(\hol_2,R_{M_3(1)=e_1})$ and $(\hol_2,-R_{M_3(1)=e_1})$
are isomorphic.

3) Let  $(\g_1,\sigma_1,\eta)$ and $(\g_2,\sigma_2,\eta)$ be the symmetric triples corresponding to the symmetric pairs
   $(\hol_3,R_1+R_{\lambda_5}+R_{M_3})$ and
$(\hol_3,R_1+R_{\lambda_5})$, respectively. Choosing in the proper way the vectors $e_{m+1},...,e_n$ of the basis,
we can assume that $M_3(1)=\alpha e_{m+1}$, where $\alpha\in\Real$. Consider the Lie algebra isomorphism $f:\g_1\to\g_2$ given by
$$\begin{array}{llll}f|_{\hol}=\id_{\hol},& f(p_1)=p_1,& f(p_2)=p_2,& f(e_1)=\frac{\alpha}{2}p_2+e_1,\\ f(f_1)=-\frac{\alpha}{2}p_1+f_1,&
f(q_1)=q_1+\frac{\alpha}{2}f_1,&   f(q_2)=q_2-\frac{\alpha}{2}e_1.&\end{array}$$
The isomorphism $f$ implies that the symmetric pairs are isomorphic.

4),5) the proof is by the same argument as in the proof of Statement 1).

6) To prove that the symmetric pairs $(\hol_3,R_1+R_{\lambda_5})$ and $(\hol_3,-R_1+R_{\lambda_5'})$ are not isomorphic
for any $\lambda_5,\lambda_5'\in\Real$ let us compute the corresponding Ricci curvatures. Since $R_{\lambda_5}$ takes
values in $\su(1.n+1)$, we have $\Ric(R_1+R_{\lambda_5})=\Ric(R_1)$. It is known that for the Ricci curvature of a
pseudo-K\"ahlerian manifold holds $\Ric(X,Y)=-\tr JR(X,JY)$. Therefore,
$$\Ric(R_1)(q_1,q_1)=\Ric(R_1)(q_2,q_2)=-\tr(J(J_m-2J))=m-2n-4<0.$$ Moreover, $\Ric(R_1)$ is zero on the other basis
vectors, since $R_1$  on these vectors takes values in $\su(1.n+1)$. On the other hand,
$$\Ric(-R_1)(q_1,q_1)=\Ric(-R_1)(q_2,q_2)=2n-m+4>0$$ and $\Ric(-R_1)$ is zero on the other basis vectors. This proves
Statement 6).

7) Suppose that $\lambda_1\neq 0$. Let $(\g_1,\sigma_1,\eta)$ and $(\g_2,\sigma_2,\eta)$ be the symmetric triples
corresponding to the symmetric pairs
    $(\hol_{1c},R_{\lambda_1}+R_{\lambda_2})$ and
$(\hol_{1c},R_{\lambda_1})$. To prove Statement 7) consider the Lie algebra isomorphism $f:\g_1\to\g_2$ given by
$$\begin{array}{lll}f|_{\hol}=\id_{\hol},& f(p_1)=p_1,& f(p_2)=\frac{\lambda_2}{\lambda_1}p_1+p_2,\\
 f(q_1)=q_1-\frac{\lambda_2}{\lambda_1}q_2,&   f(q_2)=q_2.&\end{array}$$

8) Let $(\g_1,\sigma_1,\eta)$ and $(\g_2,\sigma_2,\eta)$ be the symmetric triples corresponding to the symmetric pairs
    $(\hol_{1c},R_{\lambda_2=1})$ and $(\hol_{1c},-R_{\lambda_2=1})$. To prove the statement  consider the Lie algebra isomorphism $f:\g_1\to\g_2$ given by
$$\begin{array}{lll}f|_{\hol}=\id_{\hol},& f(p_1)=p_2,& f(p_2)=p_1,\\
 f(q_1)=q_2,&   f(q_2)=q_1.&\end{array}$$

9) To prove the statement note the following: the only non-vanishing values of $\Ric(R_{\lambda_1})$ are
$$\Ric(R_{\lambda_1})(p_1,q_1)=\Ric(R_{\lambda_1})(p_2,q_2)=2\lambda_1>0;$$ the only non-vanishing values of
$\Ric(-R_{\lambda_1})$ are $$\Ric(-R_{\lambda_1})(p_1,q_1)=\Ric(-R_{\lambda_1})(p_2,q_2)=-2\lambda_1<0;$$
$$\Ric(R_{\lambda_2=1})(p_2,q_1)=2>0,\quad \Ric(R_{\lambda_2=1})(p_1,q_2)=-2<0.$$

The theorem is proved. $\Box$

Recall that  a pseudo-Riemannian manifold $(M,g)$ is called {\it locally symmetric} if the geodesic symmetry with respect to any point $p\in M$
is a locally defined isometry. This condition is equivalent to $\nabla R=0$. For the holonomy algebra of a locally symmetric
pseudo-Riemannian manifold $(M,g)$  at a point $p\in M$ we have $\hol_p=\spa\{R_p(T_pM,T_pM)\}$. The symmetric pair $(\hol_p,R_p)$ defines
a simply connected pseudo-Riemannian symmetric space with the same holonomy algebra. We get the following corollary.

\begin{corol} The  weakly-irreducible not irreducible  holonomy  algebras of locally symmetric
pseudo-K\"ahlerian manifolds of signature $(2,2n+2)$ are exhausted by $\hol^2_{n=0}$,
$\hol^{\gamma_1=0,\gamma_2=0}_{n=0}$, $\hol^{m=n-1=0,\{0\},\varphi=0,\phi=0}$ and $\hol^{m,\{\Real
J_m\},\varphi=0,\phi=2}$.
\end{corol}

\section{Holonomy of time-like cones over  Lorentzian Sasaki manifolds}\label{secSasaki}

\subsection{Definitions and results}

Let $(M,g)$ be a pseudo-Riemannian manifold and $c=\pm $. The pseudo-Riemannian manifold $$(\t M^c=\Real^+\times M,\t
g^c=cdr^2+r^2g)$$ is called {\it the space-like cone over} $(M,g)$ if $c=+$ and it is called {\it the time-like cone over}
$(M,g)$ if $c=-$. One of the reason why the pseudo-Riemannian cones are of interest is the following. In \cite{Baer}
Ch.~B\"ar proved that it real  Killing spinors on a pseudo-Riemannian spin manifold $(M,g)$ correspond to parallel spinors
on $(\t M^+,\t g^+)$. Similarly,  imaginary  Killing spinors on a pseudo-Riemannian spin manifold $(M,g)$ correspond to
parallel spinors on $(\t M^-,\t g^-)$, see \cite{Boh03}. In \cite{Gallot79} S.~Gallot proved that {\it the space-like cone
over a complete Riemannian manifold is indecomposable or flat}. This statement does not hold for pseudo-Riemannian
manifolds and it also does not hold for non-complete Riemannian manifolds, see examples in \cite{Cone}.

We will denote by  $\p_r$  the radial vector field on  $(\t M^\pm,\t g^\pm)$.

An odd-dimensional Lorentzian manifold is called {\it a Sasaki manifold} if its time-like cone $(\t M^-,\t g^-)$ is a pseudo-K\"ahlerian manifold.
A Lorentzian manifold $(M,g)$ is a Sasaki manifold if and only if
there exists a vector field $\xi$ on $(M,g)$ such that
\begin{itemize} \item[1.] $\xi$ is a Killing vector field with $g(\xi,\xi)=-1$.
\item[2.] The map $J=-\na\xi:TM\to TM$ satisfies  $$J^2X=-X-g(X,\xi)\xi \text{ and } (\na_X J)(Y)=-g(X,Y)\xi+g(Y,\xi)X,$$
where $X,Y$ are vector fields on $M$, see \cite{Helga,HelgaFelipe}.
\end{itemize}
Such vector field $\xi$ on $M$ is called {\it a Sasaki}  field.
The Sasaki field $\xi$ and the pseudo-K\"ahlerian structure $\t J$ on $(\t M^-,\t g^-)$ are related by
$$\t J X=\t\na_X\xi\quad \text{ and }\quad \xi=\t J(r\p_r),$$
where $X$ is a vector field on $\t M$. A Lorentzian Sasaki manifold with given Sasaki field $\xi$ is denoted by $(M,g;\xi)$.
A Lorentzian Sasaki manifold $(M,g)$ is an  Einstein manifold if and only if the time-like cone over it is Ricci-flat.

Lorentzian Sasaki manifolds are of interest in many respects. For instance,  they admit twistor spinors \cite{Helga,HelgaFelipe,Kat99,Boh03}.
For us they provide examples of manifolds with the holonomy algebras contained in $\u(1,n+1)$ and $\su(1,n+1)$.

An odd-dimensional Riemannian manifold $(M,g)$ is called a Sasaki manifold if its space-like cone $(\t M^+,\t g^+)$ is a K\"ahlerian manifold.
This is equivalent to existence of a Killing vector field $\xi$ on $(M,g)$ with $g(\xi,\xi)=1$ that satisfies conditions similar to the above ones.

The standard example for regular Lorentzian Sasaki manifolds are $S^1$-bundles over Riemannian K\"ahler Einstein spaces of negative scalar curvature,
see  \cite{Helga,HelgaFelipe}. The next example gives a way to construct   Lorentzian Sasaki manifolds.

\begin{ex} \label{ex3} Let $(M_1,g_1;\xi_1)$ be a Lorentzian Sasaki manifold and $(M_2,g_2;\xi_2)$ a Riemannian Sasaki manifold.
Then the manifold
$$(M=\bR^+\times M_1\times M_2,g=ds^2+\ch^2(s)g_1+\sh^2(s)g_2;\xi_1+\xi_2)$$
is a Lorentzian Sasaki manifold and the time-like cone over $(M,g)$ is a product of truncated  cones over   $(M_1,g_1)$ and $(M_2,g_2)$,
$$(\tilde M^-= ((1,+\infty) \times M_1)\times((0,1) \times M_2), \t g^-=(-dr_1^2 + r_1^2 g_1)+(dr_2^2 + r_2^2 g_2)).$$
\end{ex}

{\it Proof.} Consider the  product of the truncated cones
$$(\tilde M^-_1\times\t M^+_1= ((1,+\infty) \times M_1)\times((0,1) \times M_2), \t g^-_1+\t g_2^+=(-dr_1^2 + r_1^2 g_1)+(dr_2^2 + r_2^2 g_2)).$$
Consider the functions $$r=\sqrt{r_1^2-r_2^2}\in\bR^+,\quad s=\arcth\left(\frac{r_2}{r_1}\right)\in\bR^+ .$$
The functions  $r$ and $s$ give the diffeomorphism $(1,+\infty)\times (0,1)\to \bR^+\times \bR^+$.
For   $\t M^-_1\times \t M_2^+$ we get
$$\t M_1^-\times \t M_2^+=\bR^+\times \bR^+\times M_1\times M_2$$
and $$\t g^-_1+\t g^+_2=-dr^2+r^2 (ds^2+\ch^2(s)g_1+\sh^2(s)g_2).$$
Obviously, the manifold $(\t M^-,\t g^-)$ is pseudo-K\"ahlerian, i.e. the manifold $(M,g)$ is a Lorentzian Sasaki manifold.
For its Sasaki field we obtain
$$\xi=\t J(r\p_r)=\t J(r\ch(s)\p_{r_1}+r\sh(s)\p_{r_2})=\t J(r_1\p_{r_1}+r_2\p_{r_2})=\xi_1+\xi_2.$$
$\Box$

The following theorem is converse to Example \ref{ex3} and it gives a local description in the case of decomposable
cone over a Lorentzian Sasaki manifold.

\begin{theorem} \label{Sasaki1}
Let $(M,g)$ be a Lorentzian Sasaki manifold  and\\ $(\t{M}^-=\bR^+\times M,\t{g}^-=-dr^2+r^2g)$  the cone over $(M,g)$.
Suppose that the holonomy algebra of $(\t{M}^-,\t{g}^-)$ preserves a non-degenerate proper subspace.
Then  there exists a dense open submanifold $U_1\subset M$  that is locally isometric to a manifold of the form
$$W=(a,b)\times N_1\times N_2,\quad (a,b)\subset\bR^+$$
with the metric $$g_W=ds^2+ \ch^2(s) g_1+\sh^2(s)g_2,$$
where $(N_1,g_1)$ is a Lorentzian Sasaki manifold and $(N_2,g_2)$ is a Riemannian Sasaki manifold.

Moreover, any point  $(r,x)\in\bR^+\times U_1\subset\t{M}$  has a neighborhood of the form
$$((a_1,b_1)\times N_1)\times ((a_2,b_2)\times N_2),\quad (a_1,b_1),(a_2,b_2)\subset\bR^+$$
with the metric $$(-dt_1^2+t_1^2g_1)+(dt_2^2+t_2^2g_2).$$
\end{theorem}

From now on we study Lorentzian Sasaki manifold such that the holonomy algebra of its time-like cone preserves a degenerate subspace.

\begin{theorem}\label{Sasaki2} Let $(M,g)$ be a Lorentzian Sasaki manifold  and\\ $(\t{M}^-=\bR^+\times M,\t{g}^-=-dr^2+r^2g)$  the time-like cone over $(M,g)$
with the pseudo-K\"ahlerian structure $\t J$.
If the holonomy algebra of $(\t{M}^-,\t{g}^-)$ preserves a 2-dimensional isotropic $\t J$-invariant  subspace,
then it annihilates this subspace, i.e. in this situation there exist locally on $\t M^-$ two isotropic parallel vector fields $p_1$ and $p_2=\t J p_1$.
\end{theorem}

\begin{theorem}\label{Sasaki3} Let $(M,g)$ be a $2n+3$-dimensional
 Lorentzian Sasaki manifold  and\\ $(\t{M}^-=\bR^+\times M,\t{g}^-=-dr^2+r^2g)$  the time-like cone over $(M,g)$
with the pseudo-K\"ahlerian structure $\t J$.
Suppose that the holonomy algebra $\hol(\t M^-)$  of $(\t{M}^-,\t{g}^-)$ annihilates a 2-dimensional isotropic $\t J$-invariant  subspace,
i.e. there exist locally on $\t M^-$ two isotropic parallel vector fields $p_1$ and $p_2=\t J p_1$.
Then

I. There is a dense open subset $U_1\subset M$ such that each point $z_0\in U_1$ has an open neighbourhood $U_{z_0}\subset U_1$ of the form
 $$U_{z_0}=(a,b)\times N,\quad a\in\bR\cup\{-\infty\},\quad b\in\bR\cup\{+\infty\},\quad a<b$$
and  the metric $g|_{U_{z_0}}$ is given by
$$g|_{U_{z_0}}=ds^2+e^{-2s}g_N,$$ where $(N,g_N)$ is a Lorentzian manifold with a parallel isotropic vector field.

II.  There exist coordinates $x,y,\hat x,\hat y,x_1,...,x_{2n}$ on $\Real^+\times U_{z_0}\subset\t M^-$ such that $\hat x,\hat y,x_1,...,x_{2n}$ are
coordinates on $N$ and the metric $\t g^-|_{\Real^+\times U_{z_0}}$ has the form
$$\t g^-|_{\Real^+\times U_{z_0}}=2dxdy+y^2g_N=2dxdy+y^2(2d\hat xd\hat y+g_1),$$
where $g_1$ is a family  of K\"ahlerian metrics on the integral manifolds corresponding
to the coordinates $x_1,...,x_{2n}$ depending on $\hat y$. In these coordinates, $p_1=\p_x$ and $p_2=\hat y \p_x-\frac{1}{y}\p_{\hat x}$.

III. Denote by $E$ the distribution on $U_{z_0}$ generated by $\p_{x_1},...,\p_{x_{2n}}$. There exists a vector field $\hat X\in E$ that does not depend on $x,y,\hat x$ and  such that
the pseudo-K\"ahlerian structure $\t J$ is given by
$$\t J p_1=p_2,\quad \t J p_2=-p_1,$$
$$\t J Y=\lam(Y)p_1+\lam(J_E Y)p_2+ J_E Y,$$
where $Y\in E$, $\lam(Y)=-\frac{y}{1+\hat y^2}g_1(Y,\hat y J_E\hat X+\hat X)$ and
\beq J_E Y=\pr_E\na^N_Y\hat X-\hat yY-(1+\hat y^2)\pr_E\na^N_Y\p_{\hat y}\label{Jusl1}\eeq
is a $\hat y$-family  of K\"ahlerian structures on the integral manifolds of the distribution $E$,
$$\begin{array}{rl}\t J\p_y=&\frac{\yy}{2(1+\yy^2)}g_1(\hat X,\hat X)p_1-\frac{1}{2(1+\yy^2)}g_1(\hat X,\hat X)p_2+\frac{1}{y}\hat X-\yy\p_y-\frac{1+\yy^2}{y}\pyy,\\
\t J\p_{\hat y}=&-\frac{y}{2(1+\yy^2)}g_1(\hat X,\hat X)p_1+\frac{1}{1+\yy^2}(J_E\hat X-\hat y\hat X)+y\p_y+\hat y\pyy. \end{array}.$$
The following conditions are satisfied
\beq\label{Jusl2} J_E\pr_E\na^N_{\p_{\hat y}}Y=\pr_E\na^N_{\p_{\hat y}}J_EY\text{ for all }Y\in E,\eeq
\beq\label{Jusl3} \pr_E\na^N_{\p_{\hat y}}\hat X=\frac{1}{1+\hat y^2}(\hat y\hat X-J_E\hat X).\eeq
The Sasaki field of the Sasaki structure of $(M,g)$ is given by
$$\xi=-\yy\p_s+\left(\frac{e^{2s}}{2}+\frac{g_1(\hat X,\hat X)}{2(1+\yy^2)}\right)\p_{\hat x}+\hat X-(1+\yy^2)\pyy.$$

IV.  If the holonomy algebra $\hol$ of $(\Real^+\times U_{z_0},\t g^-|_{\Real^+\times U_{z_0}})$ is weakly-irreducible,
then the holonomy algebra $\hol(N)$ of $(N,g_N)$ is weakly-irreducible and
\begin{itemize}\item[1)] if $\hol(N)$ is of type $\g^{2,\u},$ $\u\subset\u(n)$ (see Section \ref{sec1.3}), then $\hol$ is
of type $\hol^{n,\u,\varphi=0,\phi=0}$;
\item[2)] if $\hol(N)$ is of type $\g^{4,\u,k,\psi},$  then $\hol$ is
of type $\hol^{n,\u,\psi,k,l}$ for some $l$, $k\leq l\leq n$.\end{itemize}

V. Conversely, for any Lorentzian manifold $(N,g_N=2d\hat xd\hat y+g_1)$, where $g_1$ is as above with  a given vector field $\hat X\in TN$
and a family of K\"ahler structures $J_E$ satisfying $\p_{\hat x}\hat X=0$, $g_N(\hat X,\p_{\hat x})=g_N(\hat X,\p_{\hat y})=0$,
 (\ref{Jusl1}), (\ref{Jusl2}) and (\ref{Jusl3}), the manifold
$((a,b)\times N,ds^2+e^{-2s}g_N)$ is a Lorentzian Sasaki manifold with the Sasaki field $\xi$ as above and
the holonomy algebra of the  time-like cone $(\t M^-,\t g^-)$ is as in Statement IV.
\end{theorem}

\begin{ex}\label{ExS2} If in the situation of Theorem \ref{Sasaki3}, the metric $g_1$ has the form $$g_1=f(\hat y)g_0,$$
where $f(\hat y)$ is a positive function and  $g_0$ is a K\"ahlerian metric depending on the coordinates $x_1,...,x_{2n}$, then
$f=const$.
In this  case the holonomy algebra of  $(\Real^+\times U_{z_0},\t g^-|_{\Real^+\times U_{z_0}})$ is not weakly-irreducible.
\end{ex}

\subsection{Proof of Theorem \ref{Sasaki1}}

Using the Koszul formulae, it is easy to prove the following lemma.
\blem\label{lemSasaki1} Let $(M,g)$ be a pseudo-Riemannian manifold, $c=\pm$ and  $(\t M^c,\t g^c)$ a cone over $(M,g)$.
The Levi-Civita connection $\tilde \nabla^c$ of the metric $\tilde g^c$ is given by
$$ \begin{array}{lcllcl} \tilde \nabla^c_{\p_r} \p_r&=&0,&  \tilde \nabla^c_{\p_r}X&=&\p_rX+\frac{1}{r}X,\\
\tilde \nabla^c_X\p_r&=&\frac{1}{r} X,& \tilde \nabla^c_XY&=&\nabla_XY-cr g(X,Y)\p_r,
\end{array}
$$ where $\nabla$ is the  Levi-Civita connection of  the metric $g$ and $X,Y\in TM$ are
considered as vector fields on $M$ depending on the parameter $r$.

The curvature tensors $\tilde R^c$ and $R$ of the connections $\tilde \nabla^c$ and $\nabla$ are related by
$$\begin{array}{l}
\t R^c(\cdot,\cdot)\p_r=\t R^c(\p_r,\cdot)\cdot=0,\\
\tilde R^c(X,Y)Z=R(X,Y)Z-cR_1(X,Y)Z, \end{array}$$
where $R_1(X,Y)Z=g(Y,Z)X-g(X,Z)Y$, $X,Y,Z\in TM$.
\elem

We will need also the following  more general lemma. \blem\label{lemSasaki1A} Let $(M,g)$ be a pseudo-Riemannian manifold,
$c=\pm$ and  $(M^{c,f}=\Real\times M,\t g^{c,f}=cdr^2+f^2(r)g)$ a warped product manifold over $(M,g)$. The Levi-Civita
connection $\tilde \nabla^{c,f}$ of the metric $\tilde g^{c,f}$ is given by $$ \begin{array}{lcllcl}  \nabla^{c,f}_{\p_r}
\p_r&=&0,&   \nabla^{c,f}_{\p_r}X&=&\p_rX+\frac{f'(r)}{2f(r)}X,\\ \nabla^{c,f}_X\p_r&=&\frac{f'(r)}{2f(r)}X,&
\nabla^{c,f}_XY&=&\nabla_XY-c\frac{f'(r)}{2} g(X,Y)\p_r,
\end{array}
$$ where $\nabla$ is the Levi-Civita connection of  the metric $g$ and $X,Y\in TM$ are
considered as vector fields on $M$ depending on the parameter $r$.
\elem

In the proof of the theorem we will write $\t M$, $\t g$ and $\t \na$ instead of $\t M^-$, $\t g^-$ and $\t \na^-$, respectively.

Suppose that the  holonomy algebra $\hol_x$ of $(\t{M},\t{g})$ at a point $x\in\t{M}$ is not weakly-irreducible, that is $T_x\t{M}$
is a sum $T_x\t{M}=(V_1)_x\oplus (V_2)_x$ of two non-degenerate $\hol_x$-invariant orthogonal subspaces.
They define locally two parallel non-degenerate distributions $V_1$ and $V_2$. For simplicity we assume that $V_1$ and $V_2$ are defined
globally (this happens, for example, if $M$ is simply connected).
Denote by  $X_1$ and $X_2$ the projections of the vector field $\p_r$ to the distributions $V_1$ and $V_2$, respectively. We have
\beq\label{e1} \p_r=X_1+X_2.\eeq
We  decompose the vectors $X_1$ and $X_2$ with respect to the decomposition $T\t{M}=T\bR\oplus TM$,
\beq\label{e2} X_1=\alpha\p_r+X,\quad  X_2=(1-\alpha)\p_r-X,\eeq
where $\alpha$ is a function on $\t{M}$ and $X$ is a vector field on $\t{M}$ tangent to $M$. Since $\t g(\p_r,\p_r)=-1$ and $\t g(X_1,X_2)=0$,
we have \beq\label{e3}\t{g}(X,X)=\alpha^2-\alpha,\quad \t{g}(X_1,X_1)=-\alpha, \quad \t{g}(X_2,X_2)=\alpha-1.\eeq

\blem \label{lemA3} The open subset $U=\{x|\alpha(x)\neq 0,1 \}\subset\t{M}$ is dense.\elem
{\it Proof.}  Suppose that $\alpha=1$ on an open subspace  $V\subset\t{M}$. We claim that  $\p_r\in V_1$ on $V$.
Indeed, on $V$  we have $$X_1=\p_r+X,\,\,\, X_2=-X\,\,\,\,\text{ and   }\,\,\,\,\, \t{g}(X,X)=0.$$
We show that $X=0$.
Let $Y_2\in V_2$. We obtain the decomposition $Y_2=\lam\p_r+Y,$ where $\lam$ is a function on $\t{M}$ and $Y\in TM$.
Then  $$\t{\n}_{Y_2}X_1=\t{\n}_{\lam\p_r+Y}(\p_r+X)=\frac{1}{r}Y+\t{\n}_{Y_2}X.$$
Note that $Y=Y_2-\lam X_1+\lam X$. Hence, $$\t{\n}_{Y_2}X_1=\frac{1}{r}(Y_2-\lam X_1+\lam X)+\t{\n}_{Y_2}X.$$
Since $X,\t{\n}_{Y_2}X,Y_2\in V_2$ and $\t{\n}_{Y_2}X_1\in V_1$, we see that $$\frac{1}{r}(Y_2+\lam X)+\t{\n}_{Y_2}X=0.$$
From $\t{g}(X,X)=0$, it follows that $\t{g}(\t{\n}_{Y_2}X,X)=0$. Thus we get $\t{g} (Y_2,X)=0$ for all $Y_2\in V_2$.
Since $V_2$ in non-degenerate, we conclude that $X=0$. Hence, $\p_r\in V_1$.

Thus for any $Y_2\in  V_2$ it holds $\t{\n}_{Y_2}\p_r=\frac{1}{r}Y_2$.
Since the distribution  $V_1$ is parallel and $\p_r\in V_1$, we see that $Y_2=0$. Therefore, $V_2=0$.
$\Box$

Further on, we will consider the dense open submanifold $U\subset\t{M}$. By the  definition of $U$,
the vector fields $X_1$, $X_2$ and $X$  are nowhere isotropic on $U$.
For $i=1,2$ let $E_i\subset V_i$ be the subdistribution of $V_i$ orthogonal to $X_i$. Denote by $L$ the distribution of lines on $U$ generated
by the vector field $X$.
We obtain on $U$ the orthogonal decomposition $$T\t{M}=T\bR\oplus L\oplus E_1\oplus E_2.$$

\blem \label{lemA4} Let $Y_1\in  E_1$ and $Y_2\in  E_2$, then
\begin{itemize}
\item[1.] $Y_1\alpha=Y_2\alpha=\p_r\alpha=0.$
\item[2.] $\t{\n}_{Y_1}X=\frac{1-\alpha}{r}Y_1$, $\t{\n}_{Y_2}X=-\frac{\alpha}{r}Y_2$.
\item[3.] $\t{\n}_{\p r}X=\p_r X+\frac{1}{r}X=0$.
\item[4.] $\t{\n}_{X}X=\left(\frac{(1-\alpha)^2}{r}-X\alpha \right)X_1+\left(\frac{\alpha^2}{r}-X\alpha \right)X_2$.
\end{itemize}\elem

{\it Proof.} Using Lemma \ref{lemSasaki1}, we get $$\begin{array}{l}
\t{\n}_{Y_1}X_1=(Y_1\alpha)\p_r+\frac{\alpha}{r}Y_1+\n_{Y_1}X,\\
\t{\n}_{Y_1}X_2=-(Y_1\alpha)\p_r+\frac{1-\alpha}{r}Y_1-\n_{Y_1}X.\end{array}$$ Since $Y_1\in  E_1\subset  V_1$ and  the
distributions $V_1$, $V_2$ are parallel, we see that $$(Y_1\alpha)\p_r+\frac{\alpha}{r}Y_1+\n_{Y_1}X\in V_1\cap
V_2=\{0\}.$$ This yields that $Y_1\alpha=0$ and $\n_{Y_1}X=\frac{1-\alpha}{r}Y_1$. The other claims can be proved
similarly. $\Box$

Since $\p_r\alpha=0$, the function $\alpha$ is a function on $M$.
Note that  $U=\bR^+\times U_1$, where $$U_1=\{x\in M|\alpha(x)\neq 0,1 \}\subset M.$$

Claim 3 of Lemma \ref{lemA4} shows that $X= \frac{1}{r}\t X$, where $\t X$ is a vector field on the manifold $M$.
Hence the distributions $L$ and $E=E_1\oplus E_2$ do not depend on $r$ and  can be considered
as distributions on $M$. Claim 2 of Lemma \ref{lemA4} shows that the distributions $E_1$ and $E_2$ also
do not depend on $r$. We get on $U_1$  the orthogonal decompositions $$TM=L\oplus E,\quad E=E_1\oplus E_2.$$

\blem \label{lemA5} The function $\alpha$ satisfies  the following differential equation on $U_1$:
     $$\t X\alpha=2(\alpha-\alpha^2).$$\elem
{\it Proof.} Since $\t g(X,X)=\alpha^2-\alpha$, we see that $$2\t g(\t \na_{X}X,X)=2\alpha X\alpha-X\alpha.$$
Using this and  Claim 4 of Lemma \ref{lemA4}, we get
$$(2\alpha-1)\left(X\alpha-\frac{2}{r}(\alpha-\alpha^2)\right)=0.$$
If $\alpha=\frac{1}{2}$, then $\t{\n}_{X}X=\frac{1}{4r}\p_r$ and $\n_{\t X}\t X=\frac{r}{2}\p_r$. The last equality is impossible. $\Box$

From Lemma  \ref{lemA5} it follows that if $t$ is a coordinate on $M$ corresponding to the vector field $\t X$, then
$$\alpha(t)=\frac{e^{2t}}{e^{2t}+c},$$ where $c$ is a constant.

From Lemmas \ref{lemA4} and \ref{lemA5} it follows that
$$\t{\n}_{X}X=-\frac{\alpha-\alpha^2}{r}\p_r+\frac{1-2\alpha}{r}X.$$
On the subset $U_1\subset M$ we get the following
\beq\label{e4}  \n_{\t X}\t X=(1-2\alpha)\t X,\quad \n_{Y_1}\t X=(1-\alpha)Y_1,\quad \n_{Y_2}\t X=-\alpha Y_2,\eeq
where $Y_1\in E_1$ and $Y_2\in E_2$.

\blem \label{lemA6} The distributions $E_1,E_2,E=E_1\oplus E_2\subset TM$ defined on $U_1\subset M$ are involutive.
Let $x\in U_1$ and $M_x\subset U_1$ be the maximal connected integral submanifold of the distribution $E$, then
the distributions $E_1|_{M_x},E_2|_{M_x}\subset TM_x=E|_{M_x}$ are parallel.
\elem  {\it Proof.}
Let $Y_1,Z_1\in E_1$. Since $g(Y_1,\t X)=g(Z_1,\t X)=0$, we have
$$\begin{array}{rl}g([Y_1,Z_1],\t X)=&g(\n_{Y_1}Z_1-\n_{Z_1}Y_1,\t X)=g(Y_1,\n_{Z_1}\t X)-g(Z_1,\n_{Y_1}\t X)\\
=&g(Y_1,(1-\alpha)Z_1)-g(Z_1,(1-\alpha)Y_1)=0.\end{array}$$
Let $Y_2\in E_2$. Since $E_1\subset V_1$, $E_2\subset V_2$ and the distributions $V_1, V_2\subset T\t{M}$ are parallel,
we see that
$$\begin{array}{rl}g([Y_1,Z_1],Y_2)=&g(\n_{Y_1}Z_1-\n_{Z_1}Y_1,Y_2)\\
=&\frac{1}{r^2}\t{g}(\t{\n}_{Y_1}Z_1+r\t{g}(Y_1,Z_1)\p_r-\t\n_{Z_1}Y_1-r\t{g}(Z_1,Y_1)\p_r,Y_2)\\
=&\frac{1}{r^2} (\t{g}(\t{\n}_{Z_1}Y_1,Y_2)-\t{g}(\t{\n}_{Y_1}Z_1,Y_2))=0.\end{array}$$
This proves that the distribution $E_1\subset TM$ is involutive. The proof of the other statements is analogous. $\Box$

Note that in Example \ref{ex3} we have the following
$$V_1=T\t{M_1}^-,\quad V_2=T\t{M_2}^+,\quad E_1=T M_1,\quad E_2=T M_2,$$
$$X_1=\ch(s)\p_{r_1},\quad X_2=\sh(s)\p_{r_2 },\quad \alpha=\ch^2(s)>1,\quad X=-\frac{1}{r}\sh(s)\ch(s)\p_s.$$

Since the manifold $(\t M,\t g)$ is pseudo-K\"ahlerian, the distributions $V_1$ and $V_2$ are $\t J$-invariant.
Without loss of generality, we can assume that
the restriction of the metric $\t g$ to the distribution $V_1$ has index 2 and  the restriction of the metric $\t g$ to the distribution $V_2$
is positively definite. From (\ref{e2}) it follows that $\alpha\geq 1$.

Let $x\in U_1$. Then $\alpha>1$ on some open subset $W\subset U_1$ containing the point $x$.
From Lemma \ref{lemA5} it follows that the gradient of the function $\alpha$ is nowhere zero on $W$.
From Lemma \ref{lemA4} it follows that on $W$ the gradient of the function $\alpha$ is proportional to the vector field
$\t X|_{W}$. Hence we can assume that  $W$ is diffeomorphic to the product $(a,b)\times N$, where $(a,b)\subset\bR$ and $N$ is diffeomorphic to
the level sets of the function $\alpha$. Not also that the level sets of the function $\alpha$ are
integral submanifolds of the involutive distribution $E$. Denote such a manifold passing through a point $y\in W$ by $W_y$.
Since $\t X$ is orthogonal to $E$ and  $Z(g(\t X,\t X))=0$ for all $Z\in E$, the metric $g|_W$ must have the following form
$$g|_W=ds^2+g_N,$$ where $g_N$ is a family  of pseudo-Riemannian metrics on the  submanifolds $W_y$ depending on the parameter $s$.
We  suppose that $\p_s=\frac{\t X}{\sqrt{g(\t X,\t X)}}$.

From Lemma \ref{lemA6} and the Wu theorem
it follows that each manifold $W_y$ is locally a product of two pseudo-Riemannian manifolds.
For $Y_1,Z_1\in E_1$ and $Y_2,Z_2\in E_2$ in virtue Lemma of \ref{lemA4}  we have
$$(L_{\t X}g)(Y_1,Z_1)=2(1-\alpha)g(Y_1,Z_1),\quad (L_{\t X}g)(Y_1,Y_2)=0,\quad (L_{\t X}g)(Y_2,Z_2)=2\alpha g(Y_2,Z_2).$$
This means that the  one-parameter group of local diffeomorphisms of $W$ generated by the vector field $\t X$
preserves the Wu decomposition of the manifolds $W_y$. Hence  the manifold $N$
can be locally decomposed into the direct product of two manifolds $N_1$ and $N_2$ that are diffeomorphic to integral
manifolds of the distributions $E_1$ and $E_2$.
We can choose $W$ in such a way that $N=N_1\times N_2$.
Therefore, $$g_N=h_1+h_2,$$ where $h_1$ and $h_2$ are metrics on the integral manifolds of
the distributions $E_1$ and $E_2$ depending on $s$.

From Lemma  \ref{lemA5} it follows that the function $\alpha$ satisfies the following differential equation
$$\p_s\alpha=-2\sqrt{\alpha^2-\alpha}.$$
Hence, $$\alpha=\frac{(e^{-2s}+c_1)^2}{4e^{-2s}c_1}=\ch^2(s+\ln\sqrt{c_1}),$$ where $c_1$ is a constant.
We can assume that $c_1=1$.  We see that $s$ is defined on an interval $(a,b)\subset\Real^+$.

Let $Y_1,Z_1\in E_1$ be vector fields on $W$ such that $[Y_1,\p_s]=[Z_1,\p_s]=0$. From (\ref{e4}) it follows that
$\n_{Y_1}\p_s=-\frac{\sqrt{\alpha^2-\alpha}}{\alpha}Y_1$. From the Koszul formula it follows that
$2g(\n_{Y_1}\p_s,Z_1)=\p_sg(Y_1,Z_1)$. Thus,  $$-2\th(s)g(Y_1,Z_1)=\p_sg(Y_1,Z_1).$$ This means that $$h_1=\ch^2(s)g_1,$$
where $g_1$ does not depend on $s$. Similarly, $$h_2=\sh^2(s)g_2,$$ where $g_2$ does not depend on $s$.

For the cone over $W$ we get
$$\bR^+\times W=\bR^+\times (a,b)\times N_1\times N_2$$
and $$\t{g}|_{\bR^+\times W}=-dr^2+r^2(ds^2+\ch^2(s)g_1+\sh^2(s)g_2).$$

Consider the functions $t_1=r\ch(s)$, $t_2=r\sh(s)$.
They define a diffeomorphism from $\bR^+\times (a,b)$ onto a subset $V\subset \bR^+\times\bR^+$.

Let $(r,y)\in\bR^+\times W\subset\t{M}$, then there exist a subset
$(a_1,b_1)\times (a_2,b_2)\subset V,$ where $(a_1,b_1),(a_2,b_2)\subset\bR^+$ and $r\in(a_1,b_1)$.

On the subset $$((a_1,b_1)\times N_1)\times ((a_2,b_2)\times N_2)\subset\bR^+\times W$$
the metric $\t{g}$ has the form
$$(-dt_1^2+t_1^2g_1)+(dt_2^2+t_2^2g_2).$$

Since $(\t M, \t g)$ is a pseudo-K\"ahlerian manifold of index 2, we see that the manifold $((a_1,b_1)\times N_1,-dt_1^2+t_1^2g_1)$ is pseudo-K\"ahlerian of index 2 and
the manifold $((a_2,b_2)\times N_2,dt_2^2+t_2^2g_2)$ is K\"ahlerian.

 Theorem \ref{Sasaki1} is proved.
$\Box$


\subsection{Proof of Theorem \ref{Sasaki2}}
In the proof we will write $\t M$, $\t g$ and $\t \na$ instead of $\t M^-$, $\t g^-$ and $\t \na^-$, respectively.
Let $z\in \t M$.  Suppose that the holonomy algebra $\hol_z$  of $(\t{M},\t{g})$ preserves a 2-dimensional isotropic $\t J_z$-invariant  subspace $L_z\subset T_z\t M$.
Then in a neighbourhood $U_z\subset \t M$ of the point $z$ there exists a parallel isotropic $\t J$-invariant distribution $L$ of rang 2.
Let $p_1$ be a nowhere vanishing vector field from the distribution $L$ (we can assume that $U_z$ is enough small). Then the vector fields $p_1$ and $p_2=\t Jp_1$ generate
the distribution $L$.
We have the decompositions $$p_1=\alpha_1\p_r+Z_1,\quad p_2=\alpha_2\p_r+Z_2,$$
where $\alpha_1$, $\alpha_2$ are functions on $U_z$ and $Z_1,Z_2\in TM|_{U_z}$.
Consider the open subset $$U'=\{x\in U_z|\alpha_1(x)\neq 0\text{ or } \alpha_2(x)\neq 0\}.$$
We claim that the subset $U'$ is dense in $U_z$.
Indeed, suppose that $\alpha_1=\alpha_2=0$ on an open subset
$V\subset U_z$, then $p_1=Z_1$ and $p_2=Z_2$ on $V$.
Consequently, $$\t{\n}_Yp_2=\t{\n}_YZ_2=\n_YZ_2+rg(Y,Z_2)\p_r$$ for any $Y\in TM|_{V}$.
Since $L$ is parallel, there exist  two 1-forms $\lam_1$ and $\lam_2$ on $V$ such that
$$\t{\n}_Yp_2=\lam_1(Y) p_1+\lam_2(Y)p_2.$$
We get $\n_YZ_2+rg(Y,Z_2)\p_r=\lam_1(Y) Z_1+\lam_2(Y)Z_2$. This shows that $g(Y,Z_2)=0$ for all $Y\in TM|_{V}$, i.e. $p_2=0$ and  we get a contradiction.

Let $x\in U'$. Then $\alpha_1(x)\neq 0$  or $\alpha_2(x)\neq 0$. Assume that $\alpha_1(x)\neq 0$.
Since the holonomy algebra $\hol_x$ preserves $L_x$, for any $X,Y\in T_x\t M$, there exist numbers $C_1(X,Y),...,C_4(X,Y)$.
Such that $$\t R_x(X,Y)p_{1x}=C_1(X,Y)p_{1x}+C_2(X,Y)p_{2x} \text{ and } \t R_x(X,Y)p_{2x}=C_3(X,Y)p_{1x}+C_4(X,Y)p_{2x}.$$
Since $\t R_x(X,Y)\in\u(T_x\t M,\t g)$, we see that $C_1(X,Y)=C_4(X,Y)$ and $C_2(X,Y)=-C_3(X,Y)$.
From Lemma \ref{lemSasaki1} it follows that $$\t R_x(X,Y)p_{1x}= R_x(X,Y)Z_{1x} \text{ and } \t R_x(X,Y)p_{2x}= R_x(X,Y)Z_{2x}.$$
Hence,       $$\begin{array}{l} C_1(X,Y)(\alpha_1(x)\p_r(x)+Z_{1x})-C_3(X,Y)(\alpha_2(x)\p_r(x)+Z_{2x})=R_x(X,Y)Z_{1x},\\
C_3(X,Y)(\alpha_1(x)\p_r(x)+Z_{1x})+C_1(X,Y)(\alpha_2(x)\p_r(x)+Z_{2x})=R_x(X,Y)Z_{2x}.\end{array}$$
Consequently,
$$C_1(X,Y)\alpha_1(x)-C_3(X,Y)\alpha_2(x)=0\text{ and } C_3(X,Y)\alpha_1(x)+C_1(X,Y)\alpha_2(x)=0.$$
Thus, $(C_1^2(X,Y)+C_3^2(X,Y))\alpha_1(x)=0$ and $C_1(X,Y)=C_3(X,Y)=0$.
We see that $\t R_x(X,Y)$ annihilates the vector space $L_x$. Hence, $\t R(\cdot,\cdot)L=0$ on $U_z$.

Consider any curve $\gamma:[a,b]\to \t M$ such that $\gamma(a)=z$ and denote by
$\tau_\gamma:T_z\t M\to T_{\gamma(b)}\t M$ the parallel displacement along $\gamma$.
Since the points $z$ and $\gamma(b)$ belong to a connected simply connected open subset of $\t M$, we see that
the parallel distribution  $L$ is defined at the point $\gamma(b)$. Hence
for any $X,Y\in T_{\gamma(b)}\t M$ we have
$$\t R(X,Y)\tau_\gamma(p_{1z})=0\text{ and } \t R(X,Y)\tau_\gamma(p_{2z})=0.$$
From this and the Ambrose-Singer theorem it follows that $\hol_z$ annihilates the vector subspace  $L_z\subset T_z\t M$.  $\Box$

\subsection{Proof of Theorem \ref{Sasaki3}}

In the proof we will write $\t M$, $\t g$ and $\t \na$ instead of $\t M^-$, $\t g^-$ and $\t \na^-$, respectively.

I.  Suppose that the holonomy algebra of $(\t M,\t g)$ preserves a 2-dimensional isotropic $\t J$-invariant subspace.
Then for each point $z\in\t M$ there exists an open neighbourhood $V_z$ of $z$ and two parallel isotropic vector fields
$p_1$ and $p_2=\t J p_1$ defined on $V_z$.  Consider the decompositions
$$p_1=\alpha_1\p_r+Z_1,\quad p_2=\alpha_2\p_r+Z_2,$$ where $\alpha_1,$ $\alpha_2$ are functions on $V_z$ and $Z_1,Z_2\in TM|_{V_z}$.
Note that \beq\label{eL1}\t{g}(Z_1,Z_1)=\alpha_1^2,\quad \t{g}(Z_2,Z_2)=\alpha_2^2,\quad \t{g}(Z_1,Z_2)=\alpha_1\alpha_2,\eeq and $Z_1,$ $Z_2$ are nowhere vanishing.
We claim that the open subset $$U_z=\{x\in V_z|\alpha_1(x)\neq 0\text{ and }\alpha_2(x)\neq 0\}\subset V_z$$ is dense in $V_z$.
Indeed, suppose that $\alpha_1(x)=0$ or $\alpha_2(x)=0$ for each point $x$ of an open subset $V\subset V_z$.
If for some $y\in V$ we have $\alpha_1(y)\neq 0$, then $\alpha_2=0$ on some open subset $V_1\subset V$ containing the point $y$ (and vice verse).
Hence, we can assume that $\alpha_1=0$ on some open subset $V_1\subset V$.
On $V_1$ we have  $p_1=Z_1$. Therefore,
$$0=\t{\n}_Yp_1=\t{\n}_YZ_1=\n_YZ_1+rg(Y,Z_1)\p_r$$ for all $Y\in TM|_{V_1}\subset T\t M|_{V_1}$.
This yields that  $g(Y,Z_1)=0$ for all $Y\in TM|_{V_1}$ and we get a contradiction.

\blem \label{lemL2} Let $Y\in TM|_{U_z}\subset T\t M|_{U_z}$. For $i=1,2$  we have
\begin{itemize}
\item[1.] $\p_r\alpha=0,$ $Y\alpha_i=-rg(Y,Z_i)$.
\item[2.] $\t{\n}_{Y}Z_i=-\frac{\alpha}{r}Y$.
\item[3.] $ \t{\n}_{\p r} Z_i=\p_r Z_i+\frac{1}{r}Z_i=0$, i.e. $Z_i=\frac{1}{r}\t Z_i$, where $\t Z_i$ is a locally defined vector field on $M$.
\item[4.] $\t Z_i\alpha_i=-\alpha_i^2$.
\end{itemize}\elem
{\it Proof.} Claims 1-3 follow from the fact that $\t{\n} p_1=\t\n p_2=0$. Claim 4 follows from (\ref{eL1}) and  Claim 1 of the lemma. $\Box$

Let $z=(r_0,z_0),$ where $r_0\in\Real^+$ and $z_0\in M$.
From Lemma \ref{lemL2} it follows that $\alpha_1$ and $\alpha_2$ can be considered as locally defined
functions on $M$ and $U_z=\Real^+\times U_{z_0}$, where $U_{z_0}\subset M$ is a neighbourhood of the point $z_0$.

From Lemma \ref{lemL2} it follows  that on $U_{z_0}$ the gradient of the function $\alpha_1$ is equal to
the vector field $-\t Z_1$.
Hence the manifold $U_{z_0}$ is diffeomorphic to the product $(a,b)\times N$, where $N$ is a manifold diffeomorphic to
the level sets of the function $\alpha_1$. Note also that the level sets of the function $\alpha_1$ are orthogonal to the vector field $\t Z_1$.
Consequently  the metric $g$ must have the following form
$$g=ds^2+g_1,$$ where $g_1$ is an $s$-family of Lorentzian metrics on the level sets of the function $\alpha_1$,
and $$\p_s=\frac{\t Z}{\alpha_1}.$$

From Lemma  \ref{lemL2} it follows that the function $\alpha$ satisfies the following differential equation
$$\p_s\alpha_1=-\alpha_1.$$
Hence, $$\alpha_1(s)=c_1e^{-s},$$ where $c_1\in \bR$ is a constant.
Changing $s$, we can assume that  $c_1=\pm 1$. Both cases are similar and we suppose that $c_1=-1$, i.e. $\alpha_1=-e^{-s}$
Note that $(a,b)=-\ln(\inf_{U_{z_0}}(-\alpha_1),\sup_{U_{z_0}}(-\alpha_1))$.

Let $Y_1,Y_2\in TM$ be vector fields orthogonal to $\p_s$ and such that $[Y_1,\p_s]=[Y_2,\p_s]=0$. From Lemma \ref{lemL2}
it follows that $\n_{Y_1}\p_s=-Y_1$. From the Koszul formula it follows that  $2g(\n_{Y_1}\p_s,Y_2)=\p_sg(Y_1,Y_2)$. Thus,
\beq -2g_1(Y_1,Y_2)=\p_sg_1(Y_1,Y_2).\label{fS1}\eeq This means that $$g_1=e^{-2s}g_N,$$ where the metric  $g_N$ does not
depend on $s$.

Thus we get the decompositions $$U_{z}=\bR^+\times (a,b) \times N$$
and $$\t{g}=-dr^2+r^2(ds^2+e^{-2s}g_N).$$

Consider the diffeomorphism
$$\bR^+\times\bR^-\to \bR^+\times\bR$$
given by $$(x,y)\mapsto\left(\sqrt{-2xy},\ln\sqrt{-\frac{2x}{y}}\right).$$
The inverse diffeomorphism has the form
$$(r,s)\mapsto\left(\frac{r}{2}e^s,-re^{-s}\right).$$

We have
$$\begin{array}{l} \p_x=e^{-s}\p_r+\frac{e^{-s}}{r}\p_s,\\
\p_y=-\frac{e^{s}}{2}\p_r+\frac{e^{s}}{2r}\p_s.\end{array}$$
On $U_{z}$  the metric $\t g$ has the form
$$\t g=2dxdy+y^2g_N.$$
Note that \beq\label{y=ra} y=r\alpha_1.\eeq

II. Let us prove the following lemma.

\blem\label{S31} The vector field $\hat p_2=r(\alpha_2 p_1-\alpha_1 p_2)$ is a vector field on $(N,g_N)$ and it is parallel
on $(N,g_N)$.\elem
{\it Proof.} We have $$\hat p_2=r(\alpha_2 p_1-\alpha_1 p_2)=r(\alpha_2 Z_1-\alpha_1 Z_2)=\alpha_2 \t Z_1-\alpha_1 \t Z_2.$$
From (\ref{eL1}) it follows that $\t g(\hat p_2,\p_s)=0$. Obviously, $\t g(\hat p_2,\p_r)=0$.
We claim that $\p_r\pp=\p_s\pp=0$. Indeed, from Lemma \ref{lemL2} it follows that $\p_r\pp=\p_r(\alpha_2 \t Z_1-\alpha_1 \t Z_2)=0$.
Furthermore, using Lemma \ref{lemL2}, we get
$$\begin{array}{rl}\t\n_{\p_s}\pp=&\t\n_{\p_s}(r(\alpha_2 p_1-\alpha_1 p_2))=r(\p_s\alpha_2)p_1-r(\p_s\alpha_1)p_2\\
=&\frac{r}{\alpha_1}(\t Z_1\alpha_2)p_1-\frac{r}{\alpha_1}(\t Z_1\alpha_1)p_2=\frac{r}{\alpha_1}(-\alpha_1\alpha_2)p_1-\frac{r}{\alpha_1}(-\alpha_1^2)p_2=-\pp.\end{array}$$
On the other hand, from Lemmas \ref{lemL2} and \ref{lemSasaki1A} it follows that
$$\t\n_{\p_s}\pp=\n_{\p_s}\pp=\n^N_{\p_s}\pp=\p_s\pp-\pp.$$
Thus, $\p_s\pp=0$.

Let $Y\in TN|_{U_{z_0}}$
Using Lemmas \ref{lemSasaki1A}, \ref{lemL2} and the fact that $g(Y,Z_1)=0$, we get
$$\begin{array}{rl}\n^N_Y\pp=&\n_Y\pp-e^{-2s}g_N(Y,\pp)\p_s=\t\n_Y\pp-rg(Y,\pp)\p_r- g(Y,\pp)\p_s\\
=&\t\n_Yr(\alpha_2 p_1-\alpha_1p_2)-rg(Y,\pp)\p_r- g(Y,\pp)\p_s\\
=&r(Y\alpha_2)p_1-r(Y\alpha_1)p_2-rg(Y,\pp)\p_r- g(Y,\pp)\p_s\\
=&-r^2g(Y,Z_2)(e^{-s}\p_r+\frac{e^{-s}}{r}\p_s)-rg(Y,r(\alpha_2Z_1-\alpha_1Z_2))\p_r-g(Y,r(\alpha_2Z_1-\alpha_1Z_2))\p_s\\
=&-r^2g(Y,Z_2)(\alpha_1\p_r+\frac{\alpha_1}{r}\p_s)+\alpha_1r^2g(Y,Z_2)\p_r+\alpha_1rg(Y,Z_2)\p_s=0.\end{array}$$
The lemma is proved. $\Box$

In \cite{Sch74} R.~Schimming proved that in this situation there exist local coordinates
$\hat x, \hat y,x_1,...,x_{2n}$ on $(N,g_N)$  such that $\p_{\hat x}=\pp$ and the metric $g_N$ has the form
$$g_N=2d\hat xd\hat y+g_1,$$ where $g_1$ is a $y$-family  of Riemannian metrics
on the integral manifolds corresponding to the coordinates $x_1,...,x_{2n}$. We can choose the open set $U_{z_0}$
in such a way that the above coordinates are global on $N$.
We see that
$$\t g|_{U_z}=\t g|_{\Real^+\times U_{z_0}}=2dxdy+y^2g_N=2dxdy+y^2(2d\hat xd\hat y+g_1).$$

We will need the following lemma, which can be proved using the Koszul formulae.
\blem\label{lemSasaki1B} Let $(N,g_N)$ be a pseudo-Riemannian manifold, $(\t M,\t g)$ the pseudo-Riemannian manifold of the form
$(\t M=\Real\times\Real^+\times N,\t g=2dxdy+y^2g_N)$.
Then the non-trivial covariant derivatives $\tilde \nabla$ of the metric $\tilde g$ are the following
$$ \begin{array}{l}  \tilde \nabla_{\p_y}X=\p_yX+\frac{1}{y}X,\\
\tilde \nabla_X\p_y=\frac{1}{y} X,\\
\tilde \nabla_XY=\nabla^N_XY-y g(X,Y)\p_x,
\end{array} $$ where $\nabla^N$ is the Levi-Civita connection of  the metric $g_N$ and $X,Y\in TN$ are
considered as vector fields on $N$ depending on the parameter $y$.  \elem

Using Lemma \ref{lemL2},  it is easy to prove the following lemma.
\begin{lem}\label{dejst} We have $$\begin{array}{llll}\grad_M\alpha_1=-\t Z_1,& \grad_M\alpha_2=-\t Z_2,& \py\alpha_1=\frac{1}{2r},&  \py\alpha_2=\frac{\alpha_2}{2r\alpha_1},\\
&\grad_M\frac{\alpha_2}{\alpha_1}=\frac{1}{\alpha_1^2}\pp,& \pyy\alpha_1=0,\end{array}$$
in particular, $\pyy\frac{\alpha_2}{\alpha_1}=1$ and we may assume that $\yy=\frac{\alpha_2}{\alpha_1}$.
\end{lem}

III. Now we will find the necessary and sufficient conditions for the metric $\t g|_{U_z}$ to admit a pseudo-K\"ahlerian structure
$\t J$ such that $p_2=\t Jp_1$. Suppose that $\t g|_{U_z}$  admits such structure $\t J$. Denote by $L$ the distribution
generated by the vector fields $p_1$, $p_2$ and by $E$ the distribution generated by $\p_{x_1},...,\p_{x_{2n}}$.
Obviously, $L^\bot=L\oplus E$. Since $L$ is $\t J$-invariant and parallel, we see that $L\oplus E$ is also  $\t J$-invariant
and parallel. Hence, there exist differential 1-forms $\lam$ and $\lam_1$ such that for any $Y\in E$ we have
$$\t JY=\lam(Y)p_1+\lam_1(Y)p_2+J_E Y,$$ where $J_E:E\to E$, $J_EY=\pr_E\t JY$.
Since $\t J$ is $\t g$-invariant, i.e. $\t g(\t J\cdot,\t J\cdot)=\t g(\cdot,\cdot)$, we see that
$J_E$ is $g_1$-invariant. Since $\t J^2=-1$, we have
$$-Y=\lam(Y)p_2-\lam_1(Y)p_1+\lam(J_EY)p_1+\lam_1(J_EY)p_2+J_E^2 Y.$$
Hence, $\lam_1(Y)=\lam(J_E Y)$, $J_E^2=-1$
\beq\label{tJ}\t JY=\lam(Y)p_1+\lam(J_EY)p_2+J_E Y.\eeq

Suppose that $$\t J\p_y=a_1p_1+b_1p_2+X+c_1\p_y+d_2\pyy,$$
where $a_1,b_1,c_1,d_1\in\Real$ and $X\in E$.
Using the  conditions $\t g(\J\py,\py)=0$, $\t g(p_1,\py)=1$, $\t g(p_2,\py)=\frac{\alpha_2}{\alpha_1}$ and the fact that $\J$ is $\t g$-invariant, we get
$$a_1=-b_1\frac{\alpha_2}{\alpha_1},\quad c_1=- \frac{\alpha_2}{\alpha_1},\quad d_1=-\frac{\alpha_1^2+\alpha_2^2}{r\alpha_1^3}.$$
Similarly, from the conditions $\t g(\J\pyy,\pyy)=0$, $\t g(p_1,\py)=0$ and $\t g(p_2,\py)=-r\alpha_1$ it follows that
$$\J\pyy=a_2p_1+X_2+r\alpha_1\py+\frac{\alpha_2}{\alpha_1}\pyy,$$ where $a_2\in\Real$ and $X_2\in E$.

Since $\J^2\py=-\py$, we see that
\beq\label{Ss1}b_1=\frac{\alpha_1}{2\alpha_2}\lam(J_EX),\quad a_2=\frac{r\alpha^3}{\alpha_1^2+\alpha_2^2}\left(\lam(X)+\frac{\alpha_2^2-\alpha_1^2}{\alpha_1^2}b_1\right),\quad
X_2=\frac{r\alpha^3}{\alpha_1^2+\alpha_2^2}\left(J_EX-\frac{\alpha_2}{\alpha_1}X\right).\eeq
The condition $\J^2\pyy=-\pyy$ does not give us anything new.

For the  Levi-Civita connection $\n^N$ on $(N,g_N)$ let $\n^E=\pr_E\n^N|_E$ be the induced connection on the distribution $E$ and let $h:TN\times  E\to\Real$
be the bilinear map such that for any $Z\in TN$ and $Y\in E$ holds $\pr_{\Real \pp}\n^N_ZY=h(Z,Y)\pp.$ Thus, for any $Z\in TN$ and $Y\in E$ we have
\beq\label{n*}\n^N_ZY=h(Z,Y)\pp+\n^E_ZY.\eeq

Let $Y\in E$, then
$$\begin{array}{rl}
0=&\t g(Y,\py)=\t g(\J Y,\J\py)\\
=&\t g(\lam(Y)p_1+\lam(J_EY)p_2+J_EY,-b_1\frac{\alpha_2}{\alpha_1}p_1+b_1p_2+X-\frac{\alpha_2}{\alpha_1}\py-\frac{\alpha_1^2+\alpha_2^2}{r\alpha_1^3}\pyy)\\
=&\t g(J_EY,X)-\frac{\alpha_2}{\alpha_1}\lam(Y)+\frac{\alpha_1^2+\alpha^2_2}{\alpha^2_1}\lam(J_EY)-\frac{\alpha^2_2}{\alpha^2_1}\lam(J_EY)\\
=&r^2\alpha_1^2g_1(J_EY,X)-\frac{\alpha_2}{\alpha_1}\lam(Y)+\lam(J_E Y).\end{array}$$
Substituting $J_EY$ for $Y$, we get
$$0=-r^2\alpha_1^2g_1(Y,X)-\frac{\alpha_2}{\alpha_1}\lam(J_EY)-\lam(Y).$$
From the last two equalities it follows that
\beq\lam(Y)=-\frac{\alpha_1^3r^2}{\alpha_1^2+\alpha_2^2}g_1(Y,\alpha_2J_EX+\alpha_1X).\label{Ss9A}\eeq
Let $\hat X=\frac{1}{y}X$. Using (\ref{y=ra}) and Lemma \ref{dejst}, we get
\beq\lam(Y)=-\frac{y}{1+\yy^2}g_1(Y,\yy J_E\hat X+\hat X).\label{Ss9A1}\eeq
Not also that
\begin{align}& \alpha_1\lam(Y)+\alpha_2\lam(J_EY)=-r^2\alpha_1^3g_1(Y,X),\label{Ss9D}\\
&\alpha_2\lam(Y)-\alpha_1\lam(J_EY)=r^2\alpha_1^3g_1(JY,X). \label{Ss9D1}\end{align}

The condition $\t g(\J Y,\J\pyy)=0$ does not give us  anything new.

From (\ref{Ss9A}), (\ref{y=ra}) and Lemma \ref{dejst} it follows that
\beq b_1=-\frac{\alpha_1^4r^2}{2(\alpha_1^2+\alpha_2^2)}g_1(X,X)=-\frac{y^2}{2(1+\yy^2)}g_1(X,X)=-\frac{1}{2(1+\yy^2)}g_1(\hat X,\hat X),\eeq
\beq a_2=-\frac{\alpha_1^5r^3}{2(\alpha_1^2+\alpha_2^2)}g_1(X,X)=-\frac{y^3}{2(1+\yy^2)}g_1(X,X)=-\frac{y}{2(1+\yy^2)}g_1(\hat X,\hat X),\eeq

Now we check the condition $\t\n\t J=0$, i.e. $\t\n_{Y_1}\J Y_2=\J\t\n_{Y_1}Y_2$ for all $Y_1,Y_2\in T\t M$.

Let $Y\in E$, then
$$\J\t\n_Y\py=\J\left(\frac{1}{y}Y\right)=\frac{1}{y}(\lam(Y)p_1+\lam(J_EY)p_2+J_EY).$$
On the other hand,
$$\begin{array}{rl}\t\n_Y\J\py=&\t\n_Y(-b_1\frac{\alpha_2}{\alpha_1}p_1+b_1p_2+X-\frac{\alpha_2}{\alpha_1}\py-\frac{\alpha_1^2+\alpha_2^2}{r\alpha_1^3}\pyy)\\
=&-\frac{\alpha_2}{\alpha_1}Yb_1p_1+Yb_1p_2+\t\n_YX-\frac{\alpha_2}{\alpha_1}\n_Y\py-\frac{\alpha_1^2+\alpha_2^2}{r\alpha_1^3}\t\n_Y\pyy\\
=&-\frac{\alpha_2}{\alpha_1}Yb_1p_1+Yb_1p_2+\t\n_YX-\frac{\alpha_2}{y\alpha_1}Y-\frac{\alpha_1^2+\alpha_2^2}{r\alpha_1^3}\n^N_Y\pyy.\end{array}$$
Note that $$\begin{array}{rl}\t\n_YX=&\n^N_YX-yg_1(Y,X)p_1=\n^E_YX+h(Y,X)\pp-yg_1(Y,X)p_1\\=&\n^E_YX+h(Y,X)r(\alpha_2p_1-\alpha_1p_2)-yg_1(Y,X)p_1.\end{array}$$
Using the Koszul formulae, it is easy to check that $\n^N_Y\pyy\in E$.
Hence the condition  $\J\t\n_Y\py=\t\n_Y\J\py$ is equivalent to the following three conditions
\begin{align}\label{Ss29}& -\frac{1}{y}\lam(Y)-\frac{\alpha_2}{\alpha_1}Yb_1+h(Y,X)r\alpha_2-yg_1(Y,X)=0,\\
\label{Ss30}& -\frac{1}{y}\lam(J_EY)+Yb_1-h(Y,X)r\alpha_1=0,\\
\label{Ss28}&\frac{1}{y}J_EY=\n^E_YX-\frac{\alpha_2}{y\alpha_1}Y-\frac{\alpha_1^2+\alpha_2^2}{r\alpha_1^3}\n^N_Y\py.\end{align}
Adding (\ref{Ss29}) multiplied by $\alpha_1$ and (\ref{Ss30}) multiplied by $\alpha_2$, we get
$$-\frac{1}{y}(\alpha_1\lam(Y)+\alpha_2\lam(J_EY))-\alpha_1yg_1(Y,X)=0.$$
This equality follows from (\ref{Ss9D}) and (\ref{y=ra}).
Adding (\ref{Ss29}) multiplied by $\alpha_2$ and (\ref{Ss30}) multiplied by $-\alpha_1$, we get
$$-\frac{1}{y}(\alpha_2\lam(Y)-\alpha_1\lam(J_EY))-\frac{\alpha_1^2+\alpha_2^2}{\alpha_1}Yb_1+h(Y,X)r(\alpha_1^2+\alpha_2^2)-y\alpha_2g_1(Y,X)=0.$$
Using the Koszul formulae, it is easy to check that $h(Y,X)=-g_1(\n^N_Y\pyy,X)$. Using (\ref{Ss9D1}) and (\ref{y=ra}),  we see that the above equality
follows from (\ref{Ss28}).

Thus the condition $\J\t\n_Y\py=\t\n_Y\J\py$ is equivalent to (\ref{Ss28}). From (\ref{y=ra}) and Lemma \ref{dejst} it follows that
(\ref{Ss28}) is equivalent to
\beq J_EY=\n^E_Y\hat X-\yy Y-(1+\yy^2)\n^N_Y\pyy.\label{Ss28A}\eeq

Let $Y_1,Y_2\in E$, then
$$\begin{array}{rl}\J\t\n_{Y_1}Y_2&=\J(\n^N_{Y_1}Y_2-yg_N(Y_1,Y_2)p_1)=\J(\n^E_{Y_1}Y_2+h(Y_1,Y_2)\pp-yg_N(Y_1,Y_2)p_1)\\
&=\lam(\n^E_{Y_1}Y_2)p_1+\lam(J_E\n^E_{Y_1}Y_2)p_2+J_E\n^E_{Y_1}Y_2+h(Y_1,Y_2)r(\alpha_2p_2+r\alpha_1p_1)-yg_1(Y_1,Y_2)p_2.\end{array}$$
On the other hand,
$$\begin{array}{rl}\t\n_{Y_1}\J Y_2&=\t\n_{Y_1}(\lam(Y_2)p_1+\lam(J_EY_2)p_2+J_EY_2)\\
&=Y_2\lam(Y_2)p_1+Y_2\lam(J_EY_2)p_2+\t\n_{Y_1}J_EY_2\\
&=Y_2\lam(Y_2)p_1+Y_2\lam(J_EY_2)p_2+\n^N_{Y_1}J_EY_2-yg_1(Y_1,J_EY_2)p_1\\
&=Y_2\lam(Y_2)p_1+Y_2\lam(J_EY_2)p_2+\n^E_{Y_1}J_EY_2+h(Y_1,J_EY_2)r(\alpha_2p_1-\alpha_1p_2)-yg_1(Y_1,J_EY_2)p_1.\end{array}$$
 Thus the condition $\J\t\n_{Y_1}Y_2=\t\n_{Y_1}\J Y_2$ is equivalent to the following three conditions
\begin{align}&\lam(\n^E_{Y_1}Y_2)+r\alpha_1h(Y_1,Y_2)=Y_1\lam(Y_2)-r\alpha_1g_1(Y_1,J_EY_2)+h(Y_1,J_EY_2)r\alpha_2,\label{Ss201}\\
&\lam(J_E\n^E_{Y_1}Y_2)+r\alpha_2h(Y_1,Y_2)-r\alpha_1g(Y_1,Y_2)=Y_1\lam(J_EY_2)-h(Y_1,J_EY_2)r\alpha_1,\label{Ss202}\\
&J_E\n^E_{Y_1}Y_2=\n^E_{Y_1}J_EY_2.\label{Ss203}\end{align}
Since the  Levi-Civita connections on the integral manifolds of the distribution  $E$ coincide with the restrictions of $\n^E$, we see that
the integral manifolds of the distribution  $E$ are K\"ahlerian.

Adding (\ref{Ss201}) multiplied by $\alpha_1$ and (\ref{Ss202}) multiplied by $\alpha_2$ and using (\ref{Ss9D}), we get
\begin{multline}\label{Ss201A} -r^2\alpha_1^3g_1(\n^E_{Y_1}Y_2,X)+r(\alpha_1^2+\alpha_2^2)h(Y_1,Y_2)-r\alpha_1\alpha_2g_1(Y_1,Y_2)=\\
Y_1(-r^2\alpha_1^3g_1(Y_2,X))-r\alpha_1^2g_1(Y_1,J_EY_2).\end{multline}
Adding (\ref{Ss201}) multiplied by $\alpha_2$ and (\ref{Ss202}) multiplied by $-\alpha_1$ and using (\ref{Ss9D1}), we get
\begin{multline}\label{Ss201B} r^2\alpha_1^3g_1(J_E\n^E_{Y_1}Y_2,X)+r\alpha_1^2g_1(Y_1,Y_2)=\\
Y_1(r^2\alpha_1^3g_1(J_EY_2,X))-r\alpha_1\alpha_2g_1(Y_1,J_EY_2)+r(\alpha_1^2+\alpha_2^2)h(Y_1,J_EY_2).\end{multline}
It is easy to check that (\ref{Ss201A}) follows from (\ref{Ss28}). Substituting   $J_E Y_2$ for $Y_2$ in (\ref{Ss201A}) and using
(\ref{Ss203}), we see that (\ref{Ss201A}) and (\ref{Ss201B}) are equivalent.
Thus if (\ref{Ss28}) is true, then  the condition $\J\t\n_{Y_1}Y_2=\t\n_{Y_1}\J Y_2$ is equivalent to (\ref{Ss203}).

Similarly we can show that the condition $\J\t\n_Y\pyy=\t\n_Y\J\pyy$ is equivalent to (\ref{Ss28}).

Let $Y\in E$. From the Koszul formulae it follows that  $h(Y,\pyy)=0$. Hence,
$$\begin{array}{rl}\t\n_{\pyy}\J Y=&\t\n_{\pyy}(\lam(Y)p_1+\lam(J_E Y)p_2+J_EY)\\
=&\pyy\lam(Y)p_1+\pyy\lam(J_E Y)p_2+\n^E_{\pyy}J_EY.\end{array}$$
Furthermore, $$\begin{array}{rl}\J\t\n_{\pyy} Y=\J\n^E_{\pyy}Y=\lam(\n^E_{\pyy}Y)p_1+\lam(J_E \n^E_{\pyy}Y)p_2+J_E\n^E_{\pyy}Y.\end{array}$$
The condition $\J\t\n_{\pyy}Y=\t\n_{\pyy}\J Y$ is equivalent to the following three conditions
\begin{align}\label{Ss351}&\pyy\lam(Y)=\lam(\n^E_{\pyy}Y),\\
\label{Ss352}&\pyy\lam(\J_EY)=\lam(\J_E\n^E_{\pyy}Y),\\
\label{Ss353}&\J_E\n^E_{\pyy}Y=\n^E_{\pyy}J_EY.\end{align}
Adding (\ref{Ss351}) multiplied by $\alpha_1$ and (\ref{Ss352}) multiplied by $\alpha_2$ and using (\ref{Ss9D}), we get
$$\alpha_1\pyy\lam(Y)+\alpha_2\pyy\lam(J_EY)=-r^2\alpha_1^3g_1(\n^E_{\pyy}Y,X).$$
Hence,
$$\alpha_1\pyy\lam(Y)+\pyy(\alpha_2\lam(J_EY))-\lam(J_EY)\pyy \alpha_2=-r^2\alpha_1^3g_1(\n^E_{\pyy}Y,X).$$
Using Lemma \ref{dejst} and (\ref{Ss9D}), we get
$$\pyy(-r^2\alpha_1^3g_1(Y,X))-\alpha_1\lam(J_EY)=-r^2\alpha_1^3g_1(\n^E_{\pyy}Y,X).$$
Since $\nabla^Ng_1=0$, we have
$$-r^2\alpha_1^3(g_1(\n^E_{\pyy}Y,X)+g_1(Y,\n^E_{\pyy}X))-\alpha_1\lam(J_EY)=-r^2\alpha_1^3g_1(\n^E_{\pyy}Y,X).$$
Hence, $$-r^2\alpha_1^3g_1(Y,\n^E_{\pyy}X)-\alpha_1\lam(J_EY)=0.$$
Using (\ref{Ss9A}), we get
$$g_1(Y,\n^E_{\pyy}X)=\frac{\alpha_1}{\alpha_1^2+\alpha_2^2}(g_1(Y,\alpha_2X-\alpha_1J_EX)).$$
Since $g_1$ is non-degenerate, we see that
\beq\label{Ss35}\n^E_{\pyy}X=\frac{\alpha_1}{\alpha_1^2+\alpha_2^2}(\alpha_2X-\alpha_1J_EX)=\frac{1}{1+\yy^2}(\yy X-J_EX).\eeq
Adding (\ref{Ss351}) multiplied by $\alpha_2$ and (\ref{Ss352}) multiplied by $-\alpha_1$ and using (\ref{Ss9D}), we obtain the same condition.

Let $Y\in E$,  $\py Y=0$, then $$\begin{array}{rl}\t\n_{\py}\J Y=&\t\n_{\py}(\lam(Y)p_1+\lam(J_E Y)p_2+J_EY)\\
=&\py\lam(Y)p_1+\py\lam(J_E Y)p_2+\py(J_EY)+\frac{1}{y}J_EY.\end{array}$$
Furthermore, $$\begin{array}{rl}\J\t\n_{\py} Y=\J(\frac{1}{y}Y) =\frac{1}{y}(\lam(Y)p_1+\lam(J_E Y)p_2+J_EY).\end{array}$$
The condition $\J\t\n_{\py}Y=\t\n_{\py}\J Y$ is equivalent to the following three conditions
\begin{align}\label{Ss361}&\py\lam(Y)=\frac{1}{y}\lam(Y),\\
\label{Ss362}&\py\lam(\J_EY)=\frac{1}{y}\lam(J_EY),\\
\label{Ss363}&\py(J_E Y)=0.\end{align}
Condition (\ref{Ss363}) shows that $J_E$ does not depend on $y$.
Adding (\ref{Ss361}) multiplied by $\alpha_1$ and (\ref{Ss362}) multiplied by $\alpha_2$ and using (\ref{Ss9D}), we get
$$\alpha_1\py\lam(Y)+\alpha_2\py\lam(J_EY)=-\frac{1}{y}r^2\alpha_1^3g_1(Y,X).$$
Hence, $$\py(\alpha_1\lam(Y)+\alpha_2\lam(J_EY))-\py(\alpha_1)\lam(Y)-\py(\alpha_2)\lam(J_EY)=-r\alpha_1^2g_1(Y,X).$$
Using Lemma \ref{dejst} and (\ref{Ss9D}), we get
$$\py(-r^2\alpha_1^3g_1(Y,X))-\frac{1}{2r}\lam(Y)-\frac{\alpha_2}{2r\alpha_1}\lam(J_EY)=-r\alpha_1^2g_1(Y,X).$$
Hence, $$-\frac{5}{2}r\alpha_1^2g_1(Y,X)-r^2\alpha_1^3\py g_1(Y,X)+\frac{1}{2r\alpha_1}(\alpha_1\lam(Y)+\alpha_2\lam(J_EY))=-r\alpha_1^2g_1(Y,X).$$
Consequently, $$\py g_1(Y,X)=-\frac{1}{y}g_1(Y,X).$$
Since $\t\n \t g=0$, we have
$$\begin{array}{rl}\py g_1(Y,X)=&\py (\frac{1}{y^2}\t g(Y,X))=-\frac{2}{y^3}\t g(Y,X)+\frac{1}{y^2}\py\t g(Y,X)\\
=&-\frac{2}{y^3}\t g(Y,X)+\frac{1}{y^2}\t g(\t\n_{\py}Y,X)+\frac{1}{y^2}\t g(Y,\t\n_{\py}X)\\
=&-\frac{2}{y^3}\t g(Y,X)+\frac{1}{y^2}\t g(\frac{1}{y}Y,X)+\frac{1}{y^2}\t g(Y,\py X+\frac{1}{y}X)\\
=&g_1(Y,\py X).\end{array}$$
Therefore, $$g_1(Y,\py X)=-\frac{1}{y}g_1(Y,X).$$
Thus, $$\py X=-\frac{1}{y}X.$$
Since $X=\frac{1}{y}\hat X$, we see that $\py\hat X=0$, i.e. the vector field  $\hat X$ does not depend on $y$.
Similarly we can show that $\hat X$ does not depend on $x$ and $\hat x$.

For the Sasaki field we have
$$\begin{array}{rl}\xi=&\J(r\p_r)=r\J\left(\frac{e^s}{2}\p_x-e^{-s}\p_y\right)=r\J\left(-\frac{1}{2\alpha_1}p_1+\alpha_1\p_y\right)=-\frac{r}{2\alpha_1}p_2+y\J\p_y\\
=&-\frac{r}{2\alpha_1}p_2 +\frac{y\yy}{2(1+\yy^2)}g_1(\hat X,\hat X)p_1-\frac{y}{2(1+\yy^2)}g_1(\hat X,\hat X)p_2+\hat X-y\yy\p_y-(1+\yy^2)\pyy\\
=&-\left(\frac{r}{2\alpha_1}+\frac{yg_1(\hat X,\hat X)}{2(1+\yy^2)}\right)\left(\yy p_1-\frac{1}{y}\pp\right)
+\frac{y\yy}{2(1+\yy^2)}g_1(\hat X,\hat X)p_1+\hat X-y\yy\p_y-(1+\yy^2)\pyy\\
=&\frac{1}{y}\left(\frac{r}{2\alpha_1}+\frac{yg_1(\hat X,\hat X)}{2(1+\yy^2)}\right)\pp -\frac{r\yy}{2\alpha_1}p_1+\hat X-y\yy\p_y-(1+\yy^2)\pyy\\
=&\frac{1}{y}\left(\frac{r}{2\alpha_1}+\frac{yg_1(\hat X,\hat X)}{2(1+\yy^2)}\right)\p_{\hat x} -\frac{r\yy}{2\alpha_1}\left(\alpha_1\p_r+\frac{\alpha_1}{r}\p_s\right)+\hat X
-y\yy\left(-\frac{1}{2\alpha_1}\p_r+\frac{1}{2r\alpha_1}\p_s \right)-(1+\yy^2)\pyy\\
=&-\yy\p_s+\left(\frac{1}{2\alpha_1^2}+\frac{g_1(\hat X,\hat X)}{2(1+\yy^2)}\right)\p_{\hat x}+\hat X-(1+\yy^2)\pyy.\end{array}$$

IV. Let $z_0=(\s,w)\in U_{z_0}=(a,b)\times N$.
Let $\hol_{z}$ be the holonomy algebra of the manifold $(U_{z},\t g|_{U_{z}})$ at the point $z\in U_z$ and
$\hol_w(N)$ be the holonomy algebra of the manifold $(N,g_N)$ at the point $w\in N$.
Denote by $\g$ the subalgebra of $\u(T_zU_z,\t g_z)$ that annihilates the vectors $p_{1z}$ and $p_{2z}$. Obviously, $\hol_z\subset\g$.
Consider a basis $p_{1z},p_{2z},\p_{x_1}|_z,...,\p_{x_{2n}}|_z,\hat q_2,\hat q_1$ of $T_zU_z$, where $\hat q_1$ and $\hat q_2$ are any
vectors orthogonal to $E_z$ and such that $\t g(p_{1z},\hat q_1)=\t g(p_{2z},\hat q_2)=1$ and other values of $\t g(\hat q_1,\cdot)$ and $\t g(\hat q_2,\cdot)$ are zero.
We claim that using this basis we can identify the Lie algebra $\g$ with the following matrix algebra
\beq\g=\left\{\left.\left( {\scriptsize\begin{array}{ccccc}
0&0 &X_1^t &c &0\\
0&0 &X_2^t &  0 &-c\\
0 &0&A&-X_2&-X_1\\
0&0&0&0&0\\
0&0&0&0&0\\
\end{array}}\right)\right|\,
\begin{array}{c}
c\in\Real\\
X_2\in E_z,\\
X_1=-J_EX_2+\frac{1}{y(z)(1+\yy(z)^2)}A(-\yy(z) X_z+J_EX_z),\\
A\in\u(E_z,\t g_z) \end{array}\right\}.\label{Sa1}\eeq
Indeed,  any element  $B\in \so(T_zU_z,\t g_z)$ such that $Bp_{1z}=Bp_{2z}=0$ has in the above basis the matrix  as above with
any $X_1\in E_z$ and $A\in\so(E_z,\t g_z)$.
Let $Y\in E_z$, then
$$\begin{array}{rl}\J B\J Y=&\J B(\lam(Y)p_{1z}+\lam(J_{E}Y)p_{2z}+J_{E}Y)=\J(\t g(X_1,J_{E}Y)p_{1z}+\t g(X_2,J_{E}Y)p_{2z}+AJ_EY)\\
=&\t g(X_1,J_{E}Y)p_{2z}-\t g(X_2,J_{E}Y)p_{1z}+\lam(AJ_EY)p_{1z}+\lam(J_EAJ_EY)p_{2z}+J_EAJ_EY,\\
BY=&\t g(X_1,Y)p_{1z}+\t g(X_2,Y)p_{2z}+AY.\end{array}$$
From the condition $B\in\u(T_zU_z,\t g_z)$, i.e. $\J B\J Y=-BY$, and (\ref{tJ}) it follows that $A\in\u(E_z,\t g_z)$ and
$$\lam(AY)=\t g(X_2-J_EX_1,Y).$$
Using (\ref{Ss9D1}), we get
$$-\frac{y(z)}{1+\yy(z)^2}g_1(AY,\yy(z) J_EX+X)=y(z)^2g_1(X_2-J_EX_1,Y).$$
Hence, $$\frac{y(z)}{1+\yy(z)^2}g_1(Y,A(\yy(z)J_EX+X))=y(z)^2g_1(X_2-J_EX_1,Y).$$
Thus, $$X_2=J_EX_1+\frac{1}{y(z)(1+\yy(z)^2)}A(\yy(z)J_EX+X)$$   and $$ X_1=-J_EX_2+\frac{1}{y(z)(1+\yy(z)^2)}A(-\yy(z) X_z+J_EX_z).$$

We claim that in the basis $p_{1z},\hat p_{2z},\p_{x_1}|_z,...,\p_{x_{2n}}|_z,\frac{1}{y^2(z)}\pyy|_z,\py|_z$ we have
\beq\g=\left\{\left.\left( {\scriptsize\begin{array}{ccccc}
0&0 &X_1^t+\yy(z)X_2^t &c &0\\
0&0 &-\frac{1}{y(z)}X_2^t &  0 &-c\\
0 &0&A&\frac{1}{y(z)}X_2&-X_1-\yy(z)X_2\\
0&0&0&0&0\\
0&0&0&0&0\\
\end{array}}\right)\right|\, {\scriptsize
\begin{array}{c}
c\in\Real\\
X_2\in E_z,\\
X_1=-J_EX_2+\frac{1}{y(z)(1+\yy(z)^2)}A(-\yy(z) X_z+J_EX_z),\\
A\in\u(E_z,\t g_z) \end{array}}\right\}.\label{Sa2}\eeq
Indeed, for $Y\in E_z$ we have
$$\begin{array}{rl}B(Y)=&\t g(X_1,Y)p_{1z}+\t g(X_2,Y)p_{2z}+AY=\t g(X_1,Y)p_{1z}+\t g(X_2,Y)\left(\yy(z)p_{1z}-\frac{1}{y(z)}\hat p_{2z}\right)+AY\\
=&\t g(X_1+\yy(z)X_2,Y)p_{1z}-\frac{1}{y(z)}\t g(X_2,Y)\hat p_{2z}+AY.\end{array}$$

Consider the vector subspace $E'=\{Y'=Y-\lam(J_EY)p_{1z}|Y\in E_z\}\subset T_zU_z$.
Let $Y\in E_z$, then
$$\begin{array}{rl}\J(Y-\lam(J_EY)p_{1z})=&\lam(Y)p_{1z}+\lam(J_EY)p_{2z}+J_EY-\lam(J_EY)p_{2z}\\
=&J_EY-\lam(J_EJ_EY)p_{1z},\end{array}$$
i.e. $E'$ is $\t J$-invariant.
Let $p_{1z},p_{2z},e_1,...,e_{2n}, q_2,q_1$ be a basis of $T_zU_z$, where $e_1,...,e_{2n}$ is a basis of $E'$ and $q_1$, $q_2$ are any
vectors orthogonal to $E'$ and such that $\t g(p_{1z},q_1)=\t g(p_{2z},q_2)=1$ and other values of $\t g(q_1,\cdot)$ and $\t g(q_2,\cdot)$ are zero.
We claim that using this basis we can identify the Lie algebra $\g$ with the following matrix algebra
\beq\g=\left\{\left.\left( {\scriptsize\begin{array}{ccccc}
0&0 &-J_E X_2^t &c &0\\
0&0 &X_2^t &  0 &-c\\
0 &0&A&-X_2&J_EX_2\\
0&0&0&0&0\\
0&0&0&0&0\\
\end{array}}\right)\right|\,
\begin{array}{c}
c\in\Real,\\
X_2\in E_z,\\
A\in\u(E_z,\t g_z) \end{array}\right\}.\label{Sa2A}\eeq
Indeed, for $Y\in E_z$ we have\\
$B(Y-\lam(J_EY)p_{1z})$
$$\begin{array}{crl}
&=&\t g(X_1,Y)p_{1z}+\t g(X_2,Y)p_{2z}+AY\\
&=&\t g(X_1,Y)p_{1z}+\t g(X_2,Y)p_{2z}+AY-\lam(J_EAY)p_{1z}+\lam(J_EAY)p_{1z}\\
&=&\left(\t g(X_1,Y)-\frac{y(z)}{1+\yy^2(z)}g_1(J_EAY,\yy(z)J_EX+X)\right)p_{1z}+\t g(X_2,Y)p_{2z}+AY-\lam(J_EAY)p_{1z}\\
&=&\t g(-J_EX_2,Y)p_{1z}+\t g(X_2,Y)p_{2z}+AY-\lam(J_EAY)p_{1z}.\end{array}$$

In \cite{Leistner05a} T.~Leistner proved that the projection of $\hol_z$ to
$\so(T_w,g_{Nw})\subset\so(T_z,\t g_z)$ coincides with $\hol_w(N)$.
From this,  (\ref{Sa2}) and (\ref{Sa2A}) it follows that
\beq\hol_z=\left\{\left.\left( {\scriptsize\begin{array}{ccccc}
0&0 &-J_E X_2^t &ca &0\\
0&0 &X_2^t &  0 &-ca\\
0 &0&A&-X_2&J_EX_2\\
0&0&0&0&0\\
0&0&0&0&0\\
\end{array}}\right)\right|\, \left({\scriptsize\begin{array}{ccc}0&X_2^t&0\\0&A&-X_2\\0&0&0\end{array}}\right)\in\hol_w(N),\,\, c\in\Real \right\},\label{Sa4}\eeq
where $a=0$ or $a=1$.

\begin{lem}\label{holz}
Suppose that the holonomy algebra $\hol_z$ is weakly-irreducible. Then $\hol_w(N)$ is also weakly-irreducible.\end{lem}
{\it Proof.} Suppose that the holonomy algebra $\hol_w(N)$ is not weakly-irreducible. From the Wu theorem it follows that
the manifold $(N,g)$ is locally a product of a Lorentzian manifold and of a Riemannian manifold. Changing $U_z$, we can assume
that this decomposition is global, i.e. $(N=N_1\times N_2,g_N=g_{N_1}+g_{N_2})$, where $(N_1,g_{N_1})$ is  a Lorentzian manifold and $(N_2,g_{N_2})$
is a Riemannian manifold. In particular, $\pp\in TN_1\subset TN$. Since $\pp$ is parallel, it does not change in the directions of $TN_2$,
i.e. $\pp$ is a vector field on $N_1$.
Let $\dim N_1=n_1+2$.
Applying the theorem from \cite{Sch74} to the manifold $(N_1,g_{N_1})$, we  see that
there exist   coordinates $\hat x, \hat y,x_1,...,x_{n_1}$ on $(N_1,g_{N_1})$ such that $\p_{\hat x}=\pp$ and the metric $g_{N_1}$
has the form
$$g_{N_1}=2d\hat x d\hat y+g_2,$$
where  $g_2$ is a $\hat y$-family  of Riemannian metrics
on the integral manifolds corresponding to the coordinates $x_1,...,x_{n_1}$.
We can assume that these coordinates are global and that there exist global coordinates
$x_{n_1+1},...,x_{2n}$ on $N_2$. We get $$g_N=g_{N_1}+g_{N_2}=2d\hat x d\hat y+g_2+g_{N_2}.$$
Let $E_1\subset T_{N_1}\subset TN$ be the distribution generated by $\p_{x_1},...,\p_{x_{n_1}}$.
In our previous notation, $g_1=g_2+g_{N_2}$ and  $E=E_1\oplus TN_2$.
Since the integral manifolds of the distribution $E$ are K\"ahlerian, we see that the integral manifolds of the distribution $E_2$  are K\"ahlerian
and the manifold $(N_2,g_2)$ is K\"ahlerian. In particular, $n_1$ is even.
Consider the following basis of $E_z,$ $e_1=\p_{x_1}|_w,...,e_{2n}=\p_{x_{2n}}|_w.$
For the holonomy algebra $\hol_w(N)$ we get
\beq\hol_w(N)=\left\{\left.\left( {\scriptsize\begin{array}{cccc}
0&-J_E X_2^t &0&0 \\
0&A_1&0&-X_2\\
0&0&A_2&0\\
0&0&0&0\\
\end{array}}\right)\right|\, \left({\scriptsize\begin{array}{ccc}0&X_2^t&0\\0&A_1&-X_2\\0&0&0\end{array}}\right)\in\hol(N_1),\,\, A_2\in\hol(N_2) \right\}.\label{Sa5}\eeq
From this, (\ref{Sa4}) and the fact that $X_2,J_EX_2\in E_1$ it follows that the holonomy algebra
$\hol_z$ preserves the vector subspace $T_zN_2\subset T_zU_z$, i.e. it is not weakly-irreducible. $\Box$

From (\ref{Sa4}) it follows that  if $\hol_w(N)$ is of type $\g^{2,\u},$  then $\hol_z$ is
of type $\hol^{n,\u,\varphi=0,\phi=0}$;     if $\hol_w(N)$ is of type $\g^{4,\u,k,\psi},$  then $\hol_z$ is
of type $\hol^{n,\u,\psi,k,l}$ for some $l$, $k\leq l\leq n$. This proves Part IV of the theorem.
The proof of Part V follows from the proof of Part III.

The theorem is proved. $\Box$

\subsection{Proof of the statement of Example \ref{ExS2}}

Suppose that the metric $g_1$ has the form $g_1=f(\yy)g_0$, where $g_0$ does not depend on $\yy$.
Let $Y\in E,$ $\pyy Y=0$. From (\ref{Ss28A}) and Lemma \ref{lemSasaki1A} it follows that
\beq J_EY=\n^E_Y\hat X+FY,\label{SE21}\eeq
where $F=-\yy-(1+\yy^2)\frac{f'(\yy)}{2f(\yy)}$.
Applying $J_E$, we get
\beq -Y=\n^E_YJ_E\hat X+FJ_EY,\label{SE22}\eeq
Adding (\ref{SE21}) multiplied by $F$ and (\ref{SE22}), we get.
$$-(1+F^2)Y=\n^E_Y(F\hat X+J_E\hat X).$$
Hence, $$-Y=\n^0_Y\left(\frac{1}{1+F^2}(F\hat X+J_E\hat X)\right),$$
where $\n^0$ is the Levi-Civita connection of the metric $g_0$.
Since $\pyy Y=0$, we see that \beq\pyy\left(\frac{1}{1+F^2}(F\hat X+J_E\hat X)\right)=0.\label{SE3}\eeq
From (\ref{Ss35}) and Lemma \ref{lemSasaki1A} it follows that
$$\pyy \hat X=-\frac{f'(\yy)}{2f(\yy)}\hat X+\frac{\yy}{1+\yy^2}\hat X-\frac{1}{1+\yy^2}J_E\hat X.$$
Similarly, $$\pyy J_E\hat X=-\frac{f'(\yy)}{2f(\yy)}J_E\hat X+\frac{\yy}{1+\yy^2}J_E \hat X+\frac{1}{1+\yy^2}\hat X.$$
Hence the equality (\ref{SE3}) ie equivalent to the following two equations
\begin{align}
& -\frac{2F^2F'}{1+F^2}+F'-\frac{f'}{2f}+\frac{F\yy}{1+\yy^2}+\frac{1}{1+\yy^2}=0,\label{SE41}\\
& -\frac{2FF'}{1+F^2}-\frac{F}{1+\yy^2}-\frac{f'}{2f}+\frac{\yy}{1+\yy^2}=0.\label{SE42}
\end{align}
From (\ref{SE42}) and the definition of $F$  it follows that \beq\label{SE5} \frac{FF'}{1+F^2}=\frac{\yy}{1+\yy^2}.\eeq
Hence, \beq\label{SE6} 1+F^2=c(1+\yy^2),\eeq where $c$ is a constant.
From (\ref{SE41}) and (\ref{SE5}) it follows that $F'=-c$. Hence, $F=-c\yy+c_1$, where $c_1$ is a constant.
From this and $(\ref{SE6})$ it follows that $c=1$ and $c_1=0$. Thus, $f'(\yy)=0$ and $f(\yy)=const$.
Obviously, the holonomy algebra of $(N,g_N)$ is not weakly-irreducible. Thus $\hol$ is also not weakly-irreducible.  $\Box$

\bibliographystyle{unsrt}

\begin{thebibliography}{90}

\bibitem{Al} D.~V.~Alekseevsky, {\it Riemannian manifolds  with exceptional
holonomy groups}, Funksional Anal. i Prilozhen. 2 (2), 1-10, 1968.

\bibitem{Al2} D.~V.~Alekseevsky, {\it Homogeneous Riemannian manifolds
of negative curvature}, Mat. sb. N1, 93-117, 1975.

\bibitem{A-V-S} D.~V.~Alekseevsky, E.B. Vinberg, A.S. Solodovnikov,
{\it Geometry of spaces of constant curvature},
Itogi nauki i tehn. VINITI. Sovr. pr. mat. Fund. napr. 29, 5-146, 1988.

\bibitem{Cone} D.~V.~Alekseevsky, V.~Cort\'es, A.~S.~Galaev, T.~Leistner, {\it Cones over pseudo-Riemannian manifolds and their holonomy}, in preparation.

\bibitem{Am-Si} W.~Ambrose, I.~M.~Singer {\it A theorem on
holonomy}, Trans. Amer. Math. Soc. 79(1953),  428-443.

\bibitem{As} V.~V.~Astrahantsev, {\it On the holonomy groups of
4-dimensional pseudo-Riemannian manifolds,} Matematicheskie
zametki, 9, N1(1971),  59-66.

\bibitem{HelgaInes} H.~Baum, I.~Kath, {\it Parallel spinors and holonomy groups on pseudo-Riemannian spin manifolds,}
SFB 288 Preprint No 276, 1997.

\bibitem{Helga} H.~Baum, {\it Conformal Killing spinors and special geometric structures in Lorentzian geometry -- a survey}, Preprint, 2001.

\bibitem{HelgaFelipe} H.~Baum, F.~Leitner,
{\it The twistor equation in Lorentzian spin geometry,}
Math. Z. 247, No.4, 795-812 (2004).

\bibitem{Baer} Ch.~B\"ar, {\it Real Killing spinors and holonomy,} Commun. Math. Phys., 154(3), 509-521, 1993.

\bibitem{B-I} L.~Berard~Bergery, A.~Ikemakhen,
{\it On the Holonomy of Lorentzian Manifolds}, Proceeding of
symposia in pure math., volume 54, 27 -- 40, 1993.


\bibitem{B-Inn}  L.~Berard~Bergery, A.~Ikemakhen, {\it
 Sur l'holonomie des vari\'et\'es
pseudo-riemanniennes de signature (n,n)},
Bull. Soc. Math. France. 1997. T.125, f1, 93--114.



\bibitem{Be} A.~L.~Besse, {\it Einstein manifolds}, Springer-Verlag,
Berlin-Heidelberg-New York, 1987.

\bibitem{Ber} M.~Berger, {\it Sur les groupers d'holonomie des
vari\'et\'es \`aconnexion affine et des
vari\'et\'es riemanniennes}, Bull. Soc. Math. France 83 (1955), 279-330.


\bibitem{Ber57} M.~Berger,  {\it  Les espace sym\'etriques
non compacts}, Ann. Sci. \'ecole Norm. Sup. 1957. V.74.
P.85--177.

\bibitem{Boh03} Ch.~Bohle. {\it Killing spinors on Lorentzian manifolds.}  J. Geometry and Physics 45 (2003), 285-308.

\bibitem{Bo-Li} A.~Borel, A.~Lichnerowicz {\it Groupes d'holonomie des
vari\'et\'es riemanniennes},
C. R. Acad. Sci. Paris. 1952. V.234. P.279--300.


\bibitem{Boubel} Ch.~Boubel, {\it Sur l'holonomie des vari\'et\'es
pseudo-riemanniennes}, PhD thesis, Universit\'e Henri Poincar\'e, Nancy, 2000.

\bibitem{Boubel1} Ch.~Boubel, {\it On the holonomy of Lorentzian metrics}, Pr\'epublication de l'ENS Lyon no. 323 (2004).


\bibitem{Bo-Ze}  Ch.~Boubel, A.~Zeghib, {\it  Dynamics of some
Lie subgroups of $O(n,1)$, applications}
Pr\'epublication de l'ENS Lyon. 2003. N 315.


\bibitem{Bryant} R.~Bryant, {\it Metrics with exceptional holonomy}, Ann. of Math. (2) 126 (1987), 525-576.

\bibitem{Bryant1} R.~Bryant, {\it Recent advances in the theory of
holonomy}, S\'eminaire Bourbaki 51 $\grave{e}$me
ann\'ee. 1998-99. $n^o$ 861.

\bibitem{Br3} R.~Bryant, {\it  Classical, exceptional, and exotic holonomies}: a status report, in Actes de la Table Ronde de
G\'eom\'etrie Diff\'erentielle en l'Honneur de Marcel ee e Berger, So c. Math. France, 1996, 93-166.


\bibitem{Car1} E.~Cartan, {\it Les groupes projectifs qui laissent invariante aucune multiplicit\'e plane}, Bull. Soc. Math. Fr. 41, 53-96.(1913) ou Qeuvres
compl\'etes t. I, vol.1, P.355-398.

\bibitem{Car2} E.~Cartan, {\it Les groupes projectifs continus r\'eele qui ne  laissent invariante aucune multiplicit\'e plane},
J. Math. pures et applicqu\'ees 10, 149-186 (1914)
ou Qeuvres compl\'etes t. I, vol.1, P.493-530.


\bibitem{Car5} E.~Cartan, {\it Les groupes
r\'eels simples finis et continus,}  Ann. Scient. Ecol.
Norm. Sup. 1914. V.31. P.263--355 ou Oeuvres compl$\grave{e}$tes
T.III.  659-746 et 799--824.

\bibitem{Car6} E.~Cartan, {\it  Sur une classe remarquable
d'espaces de Riemann}, Bull. Soc. math. France 1926. V.54. P.214--264,
1927. V.55. P.114-134 ou Oeuvres compl$\grave{e}$tes
T.I, V.2. P.587--659.

\bibitem{Car8} E.~Cartan, {\it Les groupes d' holonomie des
espaces g\'en\'eralis\'es}, Acta. Math. 1926. V.
48. P.1--42 ou Oeuvres compl$\grave{e}$tes T.III. V.2. P.997--1038.

\bibitem{CMS} Q.-S.~Chi,  S.~Merkulov,  L.~Schwachhofer {\it On the existence of infinite series of exotic holonomies},
Inventiones Math. 126 (1996), 391--411.


\bibitem{Rha} G.~DeRham, {\it
Sur la r\'eductibilit\'ed'um espace de Riemann},  Comm.
Math. Helv. 1952. V.26. P.328--344.


\bibitem{Disc-Ol} A.~J.~Di~Scala, C.~lmos, {\it The
geometry of homogeneous submanifolds of hyperbolic space},
Math. Z. 237, 199 -- 209 (2001).

\bibitem{Gal1} A.~S.~Galaev, {\it The spaces of curvature tensors for
holonomy algebras of Lorentzian manifolds}, Diff. Geom. and its Applications 22, 1--18 (2005).

\bibitem{Gal2} A.~S.~Galaev, {\it Isometry groups of Lobachevskian spaces, similarity transformation
groups of Euclidean spaces and  Lorentzian holonomy groups}, arXiv:math.DG/0404426, 2004.

\bibitem{G4} A.~S.~Galaev, {\it Remark on  holonomy groups of pseudo-Riemannian manifolds of signature (2,n+2)}, Arxiv:math.DG/0406397, 2004.

\bibitem{Gal5} A.~S.~Galaev, {\it Metrics that realize all Lorentzian holonomy algebras},
International Journal of Geometric Methods in Modern Physics, Vol. 3,  Nos. 5$\&$6, p. 1025--1045 (2006).


\bibitem{Gallot79} S.~Gallot. \'Equations  diff\'erentielles caract\'eristiques de la sph\`ere. (French) Ann. Sci. \'Ecole Norm. Sup.
(4) 12 (1979), no. 2, 235--267.

\bibitem{Gold} W.~M.~Goldman, {\it Complex hyperbolic geometry}, Clarendon Press, Oxford, 1999.

\bibitem{H-O}  J.~Hano, H.~Ozeki, {\it  On the holonomy group of
linear connections},  Nagoya Math. J. 1956. V.10. P.97--100.

\bibitem{Helgason} S.~Helgason, {\it Differential geometry, Lie groups and symmetric spaces}, Academic Press, 1978.

\bibitem{Ik} A.~Ikemakhen, {\it Examples of indecomposable non-irreducible
Lorentzian manifolds}, Ann. Sci. Math. Qu\'ebec
20(1996), no. 1, 53-66.

\bibitem{Ik22} A.~Ikemakhen, {\it Sur l'holonomie des vari\'et\'es
pseudo-riemanniennes de signature (2,2+n)}, Publ. Mat. 43, no. 1, 55--84, 1999.

\bibitem{Jo} D.~Joyce, {\it Compact manifolds with special holonomy},
Oxford University Press, 2000.

\bibitem{Kar} F.~I.~Karpelevich, {\it Transitive surfaces for semisimple subgroups
of isometry groups of symmetric spaces}, Dokl Ak. Nauk SSSR, XCIII, N3, 1953.

\bibitem{Kat99} I.~Kath, {\it Killing spinors on pseudo-Riemannian manifolds}, 1999, Habilitationsschrift, Humboldt-Universit\"at zu Berlin.

\bibitem{KathOldrich} I.~Kath, M.~Olbrich, {\it On the structure of pseudo-Riemannian symmetric spaces}, arXiv:math.DG/0408249, 2004.

\bibitem{K-N} S.~Kobayashi, K. Nomizu, {\it Foundations of differential geometry,} volumes 1,2, Interscience Wiley, New York, 1963, 1967.

\bibitem{Le1} T.~Leistner, {\it Berger algebras, weak-Berger algebras and
Lorentzian holonomy}, sfb 288-preprint no. 567, 2002.

\bibitem{Le2} T.~Leistner, {\it Towards a classification of Lorentzian holonomy
groups}, arXiv: math.DG/0305139, 2003.

\bibitem{Le3} T.~Leistner, {\it Towards a classification of Lorentzian holonomy
groups. Part II: Semisimple, non-simple weak-Berger algebras}, arXiv:math.DG/0309274, 2003.

\bibitem{Le4} T.~Leistner, {\it
Holonomy and parallel spinors in Lorentzian geometry,}
PhD thesis, Humboldt-Universit\"at zu Berlin, 2003.

\bibitem{Leistner05a} T.~Leistner, {\it Conformal holonomy of C-spaces, Ricci-flat, and Lorentzian manifolds.} arXiv:math.DG/0501239.

\bibitem{MS99} S.~Merkulov, L.~Schwachh\"ofer {\it
Classification of irreducible holonomies of torsion-free affine
connections},  Ann. Math. 1999. V.150. P.77--149.

\bibitem{LS5} K.~Sfetsos, D.~Zoakos, {\it
Supersymmetry and
Lorentzian holonomy in various dimensions,} Preprint, 2004.

\bibitem{Sim} J.~Simons, {\it  On the transitivity of holonomy
systems}, Annals of Math. September 1962. V.76(2). P.213--234.

\bibitem{Sch74} R.~Schimming, {\it Riemannsche R\"aume mit ebenfrontiger und mit ebener Symmetrie,}
Mathematische Nachrichten, 59, S.128-162, 1974.

\bibitem{Sch99} L.~Schwachh\"ofer L, {\it  On the classification of
holonomy representations,} Habilitationsschrift,
Mathematisches Institut der Universit\"at Leipzig, 1999.


\bibitem{V-O} E.~B.~Vinberg, A.~L.~Onishchik, {\it Seminar on Lie groups and
algebraic groups}, Moscow, URSS, 1995.

\bibitem{Walker} A.~G.~Walker, {\it On parallel fields of
partially null vector spaces},  Quart. J. of Math. September 1949.
V.20. P.135--145.

\bibitem{Wang} M.~Y.~Wang, {\it  Parallel spinors and parallel
forms,}  Ann. Global Anal. Geom. 1989. V.7(1). P.59--68.


\bibitem{Wu} H.~Wu, {\it Holonomy groups of indefinite metrics}, Pacific
journal of math., 20, 351 -- 382, 1967.

\end{thebibliography}

\end{document}